\documentclass{pspum-l} \usepackage{amsmath, amsthm, amscd, amsfonts,amssymb} 
\numberwithin{equation}{section} 
\newcommand{\teq}{\arabic{section}.\arabic{equation}}
\newcommand{\teql}{\Alph{section}.\arabic{equation}}

\newcommand{\sqr}[2]{{\vcenter{\vbox{\hrule height.#2pt\hbox{\vrule width.#2pt
height#1pt \kern#1pt\vrule width.#2pt}\hrule height.#2pt}}}}

\newcounter{eqcount}

\renewcommand{\labelenumi}{{{\rm (\teq \alph{enumi})}}} 
\newenvironment{edesc}{\refstepcounter{equation}\begin{enumerate}}%
{\end{enumerate}}
\newenvironment{triv}{\refstepcounter{equation}\begin{list}%
{{\hbox{\rm(\teq)\ }}} \item }{\end{list}}
\newenvironment{trivl}{\refstepcounter{equation}\begin{list}%
{{\hbox{\rm(\teql)\ }}} \item }{\end{list}}

\newcommand{\ring}[1]{{\mathbb #1}}
\newcommand\bZ{{\ring{Z}}}

\newcommand\bC{{\ring{C}}} \newcommand\bR{{\ring{R}}}
\newcommand\bF{{\ring{F}}} \newcommand\bQ{{\ring{Q}}}
\newcommand\bH{{\ring{H}}}
\newcommand{\csp}[1]{{\mathbb #1}}
\newcommand{\tsp}[1]{{\mathcal #1}}

\newcommand{\esp}[1]{{\mathcal #1}}

\newcommand{\bM}{\csp{M}}
\newcommand{\prP}{\csp{P}}
\newcommand{\afA}{\csp{A}}
\newcommand{\sA}{{\esp{A}}} 
 
\newcommand{\sO}{{\tsp{O}}} 
\newcommand{\sQ}{\tsp{Q}}
\newcommand{\sP}{{\tsp {P}}} 
\newcommand{\sL}{{\tsp {L}}} 
\newcommand{\sT}{{\tsp {T}}} \newcommand{\sH}{{\tsp {H}}}
 \newcommand{\sY}{{\tsp {Y}}}
\newcommand{\sM}{{\tsp {M}}} 
 \newcommand{\sE}{{\tsp {E}}}

\newcommand{\bT}{{\csp {T}}} 

\newcommand{\eql}[2]{{\rm (\ref{#1}\ref{#2})}} 

\newcommand{\vect}[1]{{\pmb #1}} 
\newcommand{\ba}{{\vect{a}}} \newcommand{\bg}{\vect{g}}
 \newcommand{\be}{{\vect{e}}}
\newcommand{\bp}{{\vect{p}}} \newcommand{\bx}{{\vect{x}}}
 \newcommand{\bw}{{\vect{w}}}
 \newcommand{\bz}{{\vect{z}}}
\newcommand{\bm}{{\vect{m}}}  
\newcommand{\bh}{{\vect{h}}}  \newcommand{\bu}{{\vect{u}}}
\newcommand{\row}[2]{{#1_1,\ldots,#1_{#2}}}

\newcommand{\smatrix}[4]{{\big(\begin{array}{cc}
\!\lower2pt\hbox{$\scriptstyle#1$} &\lower2pt\hbox{$\scriptstyle#2$}\!
\\\! \raise2pt\hbox{$\scriptstyle#3$} &\raise2pt\hbox{$\scriptstyle#4$}
\!\end{array}\big)}}
\newcommand{\scases}[5]{{#1 =\left\{\begin{array}{ll} #2 &
\mbox{for $#3$}\\ #4 & \mbox{for $#5$}.\end{array}\right.}}

\newcommand{\texto}[1]{{\textr{#1}}}
\newcommand{\GL}{\texto{GL}} \newcommand{\SL}{\texto{SL}}
 \newcommand{\ind}{\texto{ind}}
\newcommand{\PSL}{\texto{PSL}} \newcommand{\PGL}{\texto{PGL}}
 
\newcommand{\Hom}{\texto{Hom}} \renewcommand{\ni}{\texto{Ni}}
\newcommand{\Spec}{\texto{Spec}} \newcommand{\Pic}{\texto{Pic}}
\newcommand{\textr}[1]{{\text{\rm #1}}}
\newcommand{\tr}{\textr{tr}} \newcommand{\ord}{\textr{ord}}
\newcommand{\abs}{\textr{abs}}  \newcommand{\cyc}{\textr{cyc}}
 
\newcommand{\pC}{{\textr{C}}} \newcommand{\inn}{\textr{in}}
 \newcommand{\Aut}{\textr{Aut}}
\newcommand{\pr}{\textr{pr}}

\newcommand{\norm}{{\triangleleft\,}}
\newcommand{\RET}{{\text{\rm Riemann's Existence Theorem}}}
\newcommand{\BCL}{{\text{\rm Branch Cycle Lemma}}}
\newcommand{\IGP}{{\text{\rm Inverse Galois Problem}}}

\newcommand{\alg}{\texto{alg}}

\newcommand{\rd}{\texto{rd}}
\newcommand{\ari}{\texto{ar}}

\newcommand{\tG}[1]{{}_{#1}\tilde G}

\newcommand{\bbQ}{{\bar{\ring{Q}}}}
\newcommand{\GAP}{{\bf GAP}}

\newcommand{\sph}{{\vphantom 1}}

\newcommand{\textb}[1]{{\text{\bf #1}}}
\newcommand{\bfC}{{\textb{C}}}
\newcommand{\longmapright}[2]{\smash{\mathop{\hbox to
#2pt{\rightarrowfill}}\limits^{#1}}}
\newcommand{\longmapleft}[2]{\smash{\mathop{\hbox to
#2pt{\leftarrowfill}}\limits^{#1}}}
\newcommand{\mapdown}[1]{\Big\downarrow\rlap{$\vcenter{\hbox{$\scriptstyle#1$}}
$}} \newcommand{\lmapdown}[1]{\llap{$\vcenter{\hbox{$\scriptstyle{#1}$}}
$}\Big\downarrow}
 
\newcommand{\mapright}[1]{\smash{\mathop{\longrightarrow}\limits^{#1}}}
\newcommand{\mapleft}[1]{\smash{\mathop{\longleftarrow}\limits^{#1}}}

\newcommand{\np}{{+}}   \newcommand{\nm}{{-}}
\newcommand\sem{\setminus}

\newcommand{\lrang}[1]{{\langle #1\rangle}}
\newcommand{\blrang}[1]{{\big< #1\big>}}

\newcommand{\eqdef}{\stackrel{\text{\rm def}}{=}}
\newcommand{\bsl}{\backslash}

\newfont{\sevenrm}{cmr7}
\newfont{\bsevenrm}{cmbx7}
\newfont{\mathseven}{cmsy7}
\newfont{\bigmath}{cmsy10 scaled 1200}
\newfont{\fiverm}{cmr5}
\newfont{\bfiverm}{cmbx5}
\newfont{\hel}{cmbx10 scaled 1200}
\newfont{\eu}{eufb10}
\newfont{\sseu}{eufm5}
\newfont{\seu}{eufm7}
\newfont{\Cal}{cmmib10}
\newfont{\sCal}{cmmib7}

\theoremstyle{plain}
\newtheorem{thm}{Theorem}[section] 
\newtheorem{lem}[thm]{Lemma}
\newtheorem{princ}[thm]{Principle}

\newtheorem{prop}[thm]{Proposition}
\newtheorem{cor}[thm]{Corollary}

\newtheorem{res}[thm]{Result}
\theoremstyle{definition}
\newtheorem{defn}[thm]{Definition}
\newtheorem{exmp}[thm]{Example}

\newtheorem{quest}[thm]{Question}
\newtheorem{prob}[thm]{Problem}

\theoremstyle{remark}
\newtheorem{rem}[thm]{Remark}

\newcommand{\xs}{\times^s\!}
\newcommand{\wsp}{{$\,$---$\,$}} 

\newcommand{\comm}[1]{{}}
 
\let\phi=\varphi

\newcommand{\cm}{{\eu C\kern-5pt{\raise2.1pt\hbox{$.$}}\,\,}}  
\newcommand{\C}{\texto{C}} \newcommand{\Spin}{\texto{Spin}} 
\newcommand{\hfl}[2]{\smash{\mathop{\hbox to 12mm{\leftarrowfill}} 
\limits^{\scriptstyle#1}_{\scriptstyle#2}}} 
 
\newcommand{\one}{{\pmb 1}} 
\newcommand{\Ind}{{\text{\rm Ind}}} 
\newcommand{\sh}{{\text{\bf sh}}}
\newcommand{\comp}{{\text{\bf cm}}} 
\newcommand{\wid}{{\text{\bf wd}}} 
\newcommand{\mpr}{{\text{\bf mp}}} 
\renewcommand{\GAP}{{\cite{GAP}}}
\newcommand{\Col}[2]{{\{#1\}_{#2=0}^\infty}}
\newcommand{\rest}{\text{\rm res}^*}
\newcommand{\hk}{{\hat \kappa}}
\newcommand{\Cusp}{{\texto{Cusp}}}
\newcommand{\Sp}{{\text{\rm Sp}}}
\newcommand{\psigma}{{\pmb \sigma}}
\renewcommand{\tr}{\textr{t}}

\begin{document} 
\font\eightrm=cmr8  \font\eightit=cmsl8 \let\it=\sl  
 
\title[Hurwitz monodromy and Modular Towers]{Hurwitz monodromy, spin separation 
\\ and higher levels of a Modular Tower} 
 
\date{\today} 
\newcommand\rk{{\text{\rm rk}}} 
\newcommand{\TG}[2]{{{}_{#1}^{#2}\tilde{G}}} 
\newcommand{\vsmatrix}[4]{{\left(\smallmatrix #1 & #2 \\ #3 & #4  
\\ \endsmallmatrix \right)}} 
\newcommand{\tbg}[1]{{{}^{#1}\tilde\bg}} 
\newcommand\Inn{{\text{\rm Inn}}} 

\author[P.~Bailey]{Paul Bailey}
 
\author[M.~Fried]{Michael D.~Fried} 
 \address{UC Irvine, Irvine, CA 92697, USA}
\email{pbailey@math.uci.edu} 
\address{UC Irvine, Irvine, CA 92697, USA}
\email{mfried@math.uci.edu}


\begin{abstract} Each finite $p$-perfect group $G$ ($p$ a prime) has a 
universal central $p$-extension coming from the $p$ part of its  
{\sl Schur multiplier\/}. Serre gave a
Stiefel-Whitney class approach to analyzing spin covers of   alternating groups 
($p=2$)
aimed at  geometric covering space problems that included their
regular realization for the Inverse Galois Problem. 

A special case of a general result is that any finite simple group with a 
nontrivial $p$ part to its Schur
multiplier has an infinite string of perfect centerless group
covers exhibiting nontrivial Schur multipliers for the
prime $p$. Sequences of moduli spaces of curves attached to
$G$ and $p$, called {\sl Modular Towers}, capture the
geometry of these many appearances of Schur multipliers
in degeneration phenomena of {\sl Harbater-Mumford cover
representatives}. These are modular curve tower generalizations. So, they 
inspire conjectures akin to Serre's open image theorem,   
including that at suitably high levels we expect no rational
points. 

Guided by two papers of Serre's, these
cases reveal common appearance of spin
structures producing $\theta$-nulls on these moduli
spaces. The results immediately apply
to all the expected Inverse Galois topics. This 
includes systematic exposure of 
moduli spaces having points where the
field of moduli is a field of definition and
other points where it is not. \end{abstract}

\thanks{Thm.~\ref{presH4}, from a 1987 preprint of M.~Fried, 
includes observations of J.~Thompson. D.~Semmen corrected a faulty first version 
of
Prop.~\ref{fratExt} and made the computations of Lem.~\ref{hatP-P}. Support for 
Fried
came from NSF
\#DMS-9970676 and  a senior research Alexander von Humboldt award.} 
\maketitle 

\setcounter{tocdepth}{2} 
\tableofcontents
\listoftables

\newcommand{\x}{{\tilde x}}  \newcommand{\A}{{{}_2\tilde 
A}} \newcommand{\qA}[2]{{{}_2^{#1}\tilde A_{#2}}} 
\newcommand{\F}{{{}_2\tilde F}} 

\section{Overview, related fundamental groups and preliminaries}  Here is one 
corollary of
a result from the early 1990s. There is an exact sequence  (\cite{FrVAnnals}, 
discussed in \cite[Cor.~6.14]{MM}
and
\cite[Cor.~10.30]{VB}): 
\begin{equation} \label{startGK} 1\to \tilde F_\omega\to G_\bQ\to 
\prod_{n=2}^\infty S_n\to 1.
\end{equation}
The group on the left is the profree group on a countable number of
generators. The group on the right is the direct product of the
symmetric groups, one copy for each integer: \eqref{startGK}  catches the 
absolute Galois group $G_\bQ$ of
$\bQ$  between two known groups. Suppose a subfield $K$
of
$\bar \bQ$ has $G_K$ a projective (profinite) group. Further, \cite{FrVAnnals} 
conjectures: Then $G_K$ is
pro-free (on  countably many generators) if and only if $K$ is {\sl 
Hilbertian\/}.  The conjecture 
generalizes
\eqref{startGK} and  Shafarevich's conjecture that the cyclotomic closure of
$\bQ$ has a pro-free absolute Galois group. 

So, \eqref{startGK} is a
positive statement about the Inverse Galois Problem.  Still, it works by 
maneuvering around the
nonprojectiveness of $G_\bQ$ as a projective profinite group. Modular Towers 
captures, in 
moduli space properties,  implications of this nonprojective  nature of $G_\bQ$. 
Its stems from constructing 
from each finite group $G$, a projective profinite group $\tilde G$ whose 
quotients naturally entangle $G$ with many
classical spaces: modular curves (\S\ref{modCurveIllus}) and spaces of Prym 
varieties (\S\ref{monhalfcan} and
\S\ref{modtowPullback}), to mention just two special cases. This
$\tilde G$ started its arithmetic geometry life on  a different brand of problem 
\cite[Chap.~21]{FrJ}. 

\cite{FrMT} gave definitions, motivations and applications surrounding
Modular Towers starting from the special case of modular curve towers. This 
paper continues that, with a
direct attack on properties of higher Modular Tower levels. 
Our models for applications include the {\sl Open Image
Theorem\/} on modular curve towers in 
\cite{SeAbell-adic} (details in \S\ref{outMon}). In Modular Tower's language 
this is the (four)
{\sl dihedral group involution  realizations\/} case (\cite{Fr-Se2} and 
\cite{Fr-Se1} have {\sl very\/}  elementary
explanations). Modular Towers comes from profinite group theory attached to any
finite group, a prime dividing its order and a choice of conjugacy classes in 
the group. This moduli approach uses 
a different type of group theory \wsp modular representations \wsp than the 
homogeneous space approach typical of modular curves or Siegel upper half
space. Still, these spaces appear here with a more elementary look and a new set 
of applications. 

It shows especially when we compare the monodromy groups of the tower levels 
with those for
modular curves. A key property used in \cite{SeAbell-adic}:  Levels of a modular 
curve tower attached
to the prime $p$ ($p\ge 5$) have a projective sequence of Frattini covers as 
monodromy groups.  Results similar to this should
hold  also for Modular Towers (Ques.~\ref{pFrattiniTower}). That is part of our 
suggestion that there is an Open
Image Theorem for Modular Towers.  We first briefly
summarize results; then set up preliminaries to explain them
in detail. The name {\sl Harbater-Mumford representatives\/} appears so often, 
it could well have been part of the title. Its precise
definition is as a type of element in a Nielsen class. Geometrically its meaning 
for Hurwitz spaces is an especially detectable
degeneration of curves in a family. The name derives from two papers: 
\cite{HaMC} and
\cite{MuDC}. \S\ref{genHMreps} recognizes there are two uses made of H-M reps. 
It suggests generalizations of these based on $p'$
properties. 

The tough question for Modular Towers is not if it is useful or connected well 
to the rest of mathematics. 
Rather, it is if its problems are sufficiently tractable for progress. By 
proving our Main Conjecture in convincing
cases, with analytic detail  often applied to modular curves, we assure rapid 
progress is possible. The Main
Conjecture  shows regular realizations of significant finite groups as Galois 
groups will not appear
serendipitously. Indeed, the Hurwitz space approach gives meaning to considering 
{\sl where\/} to look for them. Therefore, we have
put a  summary statement for those who ask:
\begin{quote} So, where are those regular
realizations (\S\ref{galSummary})? \end{quote} 

\subsection{Quick summary of results} 
Four branch point Hurwitz families produce quotients of the
upper half  plane by finite index subgroups of
$\PSL_2(\bZ)$. The Hurwitz monodromy group 
$H_4$ a group $\bar M_4$, identified with $\PSL_2(\bZ)$, as a natural  quotient. 
The
structure of a {\sl Modular Tower\/} gives sequences of 
$j$-line covers formed by natural moduli problems. These are immensely more 
numerous than
the special case of modular curve towers though they have both similarities to 
them and
exotic new properties generally.   

There is a {\sl Modular Tower\/} for any finite group $G$, prime $p$
dividing
$|G|$ and collection 
$\bfC$ of ($r\ge 3$)  $p'$ conjugacy classes. Then $G$ has a universal {\sl 
$p$-Frattini\/} profinite group cover, the profinite limit  of a sequence 
$\{G_k\}_{k=0}^\infty$ of $p$-extensions of $G=G_0$. Each {\sl level\/} $k\ge 0$ 
has 
a moduli space $\sH_k$ in the (reduced) Modular Tower. When $r=4$ all components 
of 
the $\sH_k\,$s are  $j$-line covers.    Then, {\sl reduced inner\/}
Hurwitz spaces generalize the  modular curve sequence 
$\{Y_1(p^{k+1})\}_{k=0}^\infty$. 
Analyzing 
the action of $H_4$ (and $\bar M_4$) 
on Nielsen classes attached to each Modular Tower level gives their detailed 
structure. 
If there is $k\ge 0$ with all  $\sH_k$ components of genus at least 2, then
for each number field $K$ there is 
$k_{G_0,\bfC,K}$ with $\sH_k(K)=\emptyset$ for $k\ge
k_{G_0,\bfC,K}$  (Thm.~\ref{thm-rbound}).

Main Conjecture (\cite[Main 
Conj.~0.1]{FrKMTIG}, Prob.~\ref{MPMT}, \S\ref{diophMT}): For $r\ge 4$, $G$ 
centerless and
$p$-perfect (\S\ref{pperfSec}), and an explicit
$k'_{G_0,\bfC}=k_0'$, all 
$\sH_k$ components have {\sl general type\/} if $k\ge k'_{G_0,\bfC}$. So, if 
$r=4$,  
the genus of $\sH_k$ components, $k\ge k'_0$, exceeds 1.  \S\ref{usebcl} reminds 
how 
finding regular realizations of the collection  $\{G_k\}_{k=0}^\infty$
forces this problem on the levels of a Modular Tower. The methods of the paper 
pretty much restrict 
to the case $r=4$. We use $p=2$ and $G=A_5$ to find  considerable about 
excluding $r=4$ for
general
$(G,p)$. \S\ref{indSteps} states the unsolved problems for $r=4$ framed by the 
genus computations
of \S\ref{ellipFix}. The effect of increasing
components in going to higher levels of a Modular Tower justifies concentration 
on analyzing the component structures at level 1 in our main examples.

For
$G_0=A_5$, 
$p=2$ and  $\bfC_{3^4}$, four 3-cycle conjugacy classes, $k'_0=1$
works. Level
$k=0$ has one (genus 0) component.  Level $k=1$ is the moduli of Galois covers 
with
group
${}_2^1\tilde  A_5=G_1$, a maximal Frattini extension of
$A_5$ with exponent 2 kernel, and order 3 branching. A $\bQ$  regular
realization of
$({}_2^1\tilde A_5,\bfC_{3^4})$ produces a $\bQ$ point on level 1.  
Level 1 components ($\sH_1^+$ and and $\sH_1^-$, of genuses 12  and 9) 
correspond to values
of a  {\sl spin\/}  cover invariant (as in \cite[Part III]{FrMT} and  Serre
\cite{SeTheta}). Pairings on Modular Tower cusps appear in
the symmetric
$\sh$-{\sl incidence matrix}. Each component corresponds to a mapping class 
orbit on {\sl
reduced Nielsen  classes\/}, whose notation is $\ni({}_2^1\tilde 
A_5,\bfC_{3^4})$
($\ni(G_k\bfC)$ for the level $k$ Nielsen class generally). 

Level 1 has one component of
real points containing all {\sl Harbater-Mumford\/} (H-M) and {\sl near H-M\/}  
reps.~(Thm.~\ref{oneOrbitkappa}). The  finitely many ($\PGL_2(\bC)=\PSL_2(\bC)$ 
equivalence)
$({}_2^1\tilde A_5,\bfC_{3^4})$ regular $\bQ$ realizations  are on the genus 12 
component. 
Our computation for real points emphasizes that cusps also interpret
as elements of a Nielsen class. This combines geometry with group computations. 
H-M
cusps are real points on the reduced moduli space representing total degeneracy 
of the
moduli of curve covers upon approach to the cusp.  The near H-M cusps, also 
real, are a
Modular Tower phenomenon. The corresponding cusp at one level below in the 
tower, is an
H-M rep. Above level 0, H-M and near H-M cusps are the only real cusps. This 
holds for any
Modular Tower when the prime
$p$ is 2. 

Given $\alpha_0\in H^2(G_0,\bZ/p)$ there is a universal sequence of nontrivial 
elements 
$\{\alpha_k\in H^2(G_k,\bZ/p)\}_{k=0}^\infty$. Properties are in  
Prop.~\ref{RkGk}. This is
a special case of
\cite[Schur multipliers results 3.3]{FrKMTIG}. An embedding  of a group $G$ in 
the
alternating group $A_n$ ($n\ge 4$) is {\sl spin separating\/} if the lift $\hat 
G$ to
the spin cover
$\Spin_n=\hat A_n\to A_n$ is nonsplit over
$G$ (see Def.~\ref{spinSepRepDef}). Such a central extension represents an
element $\alpha \in H^2(G,\bZ/2)$. Given an a priori $\alpha$, we also
say the embedding of $G$ in $A_n$ is $\alpha$-spin separating. When
$G_0=G$, this produces a sequence $\{\alpha_k\in H^2(G_k,\bZ/2)\}_{k=0}^\infty$ 
(as above
for $p=2$).   So, each  
$\alpha_k$ represents a central extension $T_k'\to G_k$ with $\bZ/2$ as kernel. 
Finding 
spin separating representations for these extensions has applications that 
include
identifying projective systems of components on certain Modular Tower levels. 
These
conjecturally traces to the evenness and oddness of associated $\theta$-nulls on
these components (\S\ref{hcan}). 

Prop.~\ref{expInv} produces spin separating representations of
$G_1$ for which $\hat G_1=T_1'$ is $\alpha_1$. Let $M$ be the kernel of
$G_1\to A_5$ and use that 
$G_1$ factors through  the spin cover
$\SL_2(\bZ/5)\to A_5$. A special case of a general operating principle applies 
next. With
$\bg\in \ni(A_5,\bfC_{3^4})$ an H-M rep.~at level 0, let $\bar M_{\bg}$ be the 
subgroup of
$\bar M_4$ stabilizing $\bg$. This  
allows us (a significant special case of a general situation) to regard 
$\bar M_{\bg}$ acting on
$(a,b)\in M^2$ with  both $a,b$ over $-1\in \SL_2(\bZ/5)$ and lying  
in a certain linear subspace of $M\times M$. There are two
orbits: 
$a$ and $b$ in the same $A_5$ orbit (action on $M$); and $a$ and $b$ in 
different $A_5$ orbits (Prop.~\ref{expInv}). 
This gives information on all levels of the same tower, and toward the Main 
Conjecture when $r=4$
for all groups $G_0=G$ (\S\ref{indSteps}). 

Finite covers 
of the $r$-punctured Riemann sphere $\prP^1_z=\bC\,\cup\{\infty\} $ have an 
attached  group $G$ and unordered set of conjugacy classes $\bfC$ of $G$. The
braid group  approach to the \IGP\ combines moduli space geometry with finite
group representations to  find  $(G,\bfC)$ regular realizations as rational
points on Hurwitz spaces. \cite{ Fr-Schconf}  includes an update on applying
this to classical problems without obvious  connection to the \IGP.

The history of the \IGP\ shows it is
hard to find $\bQ$  realizations of most finite groups $G$. Such a realization 
is a
quotient of $\pi_1(U_\bz,z_0)$: $U_\bz$ is $\prP^1_z$ with $r$ distinct points
$\bz$ removed.  To systematically investigate  the diophantine 
difficulties we take large  (maximal {\sl Frattini}) quotients
$\tG p$ of
$\pi_1(U_{\bz,z_0})^\alg$ covering $G$ instead of one finite group $G$. {\sl 
Modular 
Towers\/}
are moduli space sequences encoding 
quotients of $\pi_1(U_\bz,z_0)$ isomorphic to $\tG p$ as $\bz$ varies. This 
generalizes  classical modular curve
sequences, and their traditional  questions, through an \IGP\ formulation.
The modular curve case starts with a dihedral  group $D_p=\bZ_p\xs\{ \pm 1\}
$, and a natural projective limit of the dihedral groups $\{  D_{p^{ k+1} } \} 
_{ k=0}
^\infty$. In the generalization,  the maximal
$p$-Frattini (universal 
$p$-Frattini; \S \ref{frattGps}) cover $\tG p$ of any finite group  $G$ replaces 
$\bZ_p\xs\{ \pm
1\}
$. 

Levels of a Modular Tower correspond
to characteristic quotients $\{G_k\}_{ k=0}^\infty$ of $\tG p$, 
replacing  $\{  D_{p^{ k+1} } \}
_{ k=0} ^\infty$. Finding $\bfC$ regular realizations of the groups 
$\{G_k\}_{k=0}^\infty$
is to finding rational points on levels of a Modular Tower as finding regular
involution realizations of $\{D_{p^{k+1}}\}_{k=0}^\infty$ is to finding rational
points on the modular curves $\{X_1(p^{k+1})\}$ (\S\ref{HurView}). Though the 
Modular Tower
definition supports many general properties of modular curves, when $r>4$
the moduli space levels are algebraic varieties of dimension $r-3>1$.  
By conjecture (under the
$p$-perfect hypothesis) we expect a structured disappearance of rational points 
at
high levels in generalization of modular curves. This holds by
proof in many cases beyond those of modular curves from this paper. 

Modular Towers joins number theory  approaches to $G_\bQ$ actions on one hand to 
the \IGP\ braid group approach. Profinite systems of points on Modular Towers 
reveal
relations in the Lie group  of $G_\bQ$ attached to a
$p$-adic representation. Yet, they benefit from the explicit tools of  finite
group theory (modular representations).

Our concentration is on Modular Towers for groups distinct from dihedral 
groups. Many have geometric properties not appearing for towers of modular 
curves.
Especially they can have several connected components. The Hurwitz viewpoint 
applied to towers of modular curves shows how some points of general Modular 
Towers
work.  Examples:
\S\ref{HurView}  on cusps of the modular curves 
$X_0(p^{ k+1} )$ and $X_1(p^{ k+1} )$; and \S\ref{fineModDil}  comparing the 
moduli problems from isogenies of elliptic curves with the Hurwitz space 
interpretation of 
these spaces. Each Modular Tower level maps naturally to a space of abelian 
varieties,
given by the jacobians (with dimension inductively computed from 
\S\ref{jawareMT}; each far from
simple) of the curve
covers parametrized by the points of the moduli space. For example, level $1$ of 
our
$(A_5,\bfC_{3^4},p=2)$ Modular Tower maps to a one-dimensional (nonconstant) 
space of Jacobians
having dimension $1+2^7$. 
\S\ref{jawareMT} considers how to use this to measure how far is a Modular Tower 
from being a tower of modular
curves.  

\subsection{Classical paths on $U_\bz$}  \label{classGens}  Let $\bz=(\row 
z r)$ be distinct points (punctures on $\prP^1_z$) and denote $\prP^1_z\setminus 
\bz$ by
$U_\bz$.
\RET\ produces 
$(G,\bfC)$ covers. Classical 
generators come from paths, $P_i=\gamma_i\delta_i\gamma_i^{-1} $, $i=1,\dots,r$, 
satisfying these conditions:  
\begin{edesc}  \label{classGen}  \item $\delta_i$ clockwise bounds a disc 
$\Delta_i$ 
around $z_i$, $\Delta_i\setminus\{z_i\} \subset U_\bz$; \item $\gamma_i$ starts 
at $z_0$ 
and ends at some point on $\delta_i$; and \item \label{classGenc}  excluding 
their 
beginning and end points, $P_i$ and $P_j$ never meet if $i\ne j$.  \end{edesc}  
Let $\row 
{\bar g}  r$ be the respective homotopy classes of $\row P r$. Use 
\eql{classGen}{classGenc}  to attach an order of clockwise emanation to $\row P 
r$.  Assume
it is in the  order of their subscript numbering. Then: \begin{triv}  
\label{prod-one} 
$\bar g_1\bar  g_2\dots \bar g_r=1$: The product-one condition. \end{triv} A 
surjective
homomorphism
$\psi:\pi_1(U_\bz,z_0)\to G$ produces a cover with monodromy group $G$.
More precise data stipulates that a system of classical generators   
map to 
$\bfC$. Such a $\psi$ is a { \sl geometric\/}  $(G,\bfC)$ representation.  

If $\alpha$ is an automorphism of $\pi_1(U_\bz,z_0)$, it sends generators to new 
generators, changing $\psi$ to $\psi\circ\alpha$. An inner automorphism of 
$\pi_1(U_\bz,z_0)$ produces a cover equivalent to the old cover. We use moduli 
of  covers, so equivalence two homomorphisms if they differ by an inner 
automorphism. Only automorphisms of $\pi_1(U_\bz,z_0)$ from the Hurwitz 
monodromy 
group $H_r$  send classical generators to classical generators (\eqref{hwgp}: 
possibly changing the 
intrinsic order of the paths). Such automorphisms arise by deforming  
$(\bz,z_0)$. They permute the conjugacy classes of the $\bar g_i\,$s in 
$\pi_1(U_\bz,z_0)$. Given $\psi$, $\bfC$ is a well defined $H_r$ invariant.

Suppose $z_0$ and the unordered set $\bz$ have definition field $\bQ$.  
Lem.~\ref{bcl}, the \BCL, gives a {\sl necessary\/}  ({\sl
rational union\/}) condition for an  arithmetic $(G,\bfC)$ (over $\bQ$) 
representation; $\psi$ factors through a representation of the 
arithmetic fundamental  group (\S \ref{IGMot}).  Such a cover is a 
$(G,\bfC)$ (regular) realization (over $\bQ$).  

For every finite group $G$ and prime $p$ 
dividing $|G|$, there is a universal $p$-Frattini  
(\S \ref{frattGps}) cover ${}_p\tilde G$ of $G$. Generators of $G$ lift 
automatically to generators
of any finite quotient of ${}_p\tilde G$ factoring through $G$. This and other 
properties give all
finite quotients of ${}_p\tilde G$ a touching resemblance to $G$. A {\sl Modular 
Tower\/} 
(\S\ref{defMT}) is the system of moduli spaces encoding all {\sl $p$-Frattini
$(G,\bfC)$\/}  covers of $\prP^1_z$  (\S\ref{usebcl}).  Reducing by 
$\PSL_2(\bC)$ 
action gives a {\sl reduced\/}  Modular Tower. Let 
$Y_1(p^k)$ be the modular curve $X_1(p^k)$ without its cusps. With $G$ a 
dihedral 
group $D_p$ ($p$ odd) and $\bfC$ four repetitions of the involution conjugacy 
class,  the
(reduced) Modular Tower is the curve sequence  $$ \cdots \to Y_1(p^{k+1}) 
\to Y_1(p^k) \to \cdots \to Y^1(p)\to \prP^1_j\setminus \{\infty\}.$$ Our main 
results
compute properties of an
$A_5$ Modular  Tower attached to $\bfC_{3^4} $, four repetitions of elements of 
order 3 
as
associated  conjugacy classes $\bfC_{3^4} $. \S\ref{arithmetic}  extends the
general  discussion of Modular Towers beyond \cite{FrMT}  and \cite{FrKMTIG}.

\subsection{Results of this paper} When $r=4$, reduced Modular Tower levels are 
curves 
covering the classical $j$-line.  As they are upper half plane quotients by 
finite index 
subgroups of $\SL_2(\bZ)$, this case shows how reduced Modular Towers generalize 
properties of modular curves. The groups $H_4$ (Hurwitz monodromy group) and 
$\bar
M_4$ (reduced mapping class group) from Thm.~\ref{presH4} apply to 
make reduced Modular Tower computations such as these. \begin{edesc} 
\item Branch cycle descriptions and {\sl cusp\/} data for these $j$-line covers 
(\S\ref{j-covers}). \item A Klein 4-group $\sQ''$ test for reduced Hurwitz 
spaces being
{\sl b-fine\/} or fine moduli spaces  (Prop.~\ref{redHurFM}). 
\end{edesc} The group $\sQ''$ is a quotient of a normal
quaternion subgroup of $H_4$. We use the term {\sl b(irational)-fine\/} to mean 
that
excluding the locus over $j=0$ or 1 (or $\infty$), the reduced space is a fine 
moduli
space. The geometric description of these groups in Prop.~\ref{isotopyM4} 
presents
arithmetic conclusions when the action of
$\sQ''$ on Nielsen classes is trivial. These  are convenient for applications to 
the Inverse
Galois Problem (see
\S\ref{applyA5} and Rem.~\ref{course-finemoduli}).  

The universal $p$-Frattini cover of a finite group $G$ lies behind the 
construction of Modular
Towers. Any triple $(G,p,\bfC)$ with $p$ a prime dividing $|G|$ and $\bfC$ some 
$p'$ conjugacy
classes in $G$ produces a tower of moduli spaces.  \S\ref{startA5MT} is brief
about the relevant group theory, relying on results from \cite{FrKMTIG}. Still, 
\S\ref{frattiniPre} reminds of the geometry of inaccessibility that makes 
Frattini covers
intriguing. The moduli space phenomena in the
case $G=A_5$ and $p=2$ interprets modular representation theory in the geometry 
of the
particular Modular Tower. 

\subsubsection{2-Frattini covers of $A_5$} 
\S\ref{startA5MT}  describes the  characteristic quotients
of the universal 2-Frattini cover of
$A_5$.  It  produces the sequence of groups $\{
G_k\}_{k=0}^\infty$ defining the levels of the  Modular
Tower of this paper where $G_0=A_5$. 
\S\ref{realPtEx}   uses  $j$-line branch cycles (\S\ref{j-covers}) to draw 
diophantine
conclusions about  regular realizations. The main  information comes from 
analyzing
components of real points  at level 0 and 1 of the $(A_5,\bfC_{ 3^4})$ Modular 
Tower.

\subsubsection{$\bR$-cover and $\bR$-Brauer points} \label{cover-brauer} If 
$\phi:X\to \prP^1_z$ is
a (4 branch point) cover over a field
$K$ in a given Nielsen class $\ni$ (\S\ref{setupNC}), then it produces a
$K$ point on the corresponding {\sl reduced Hurwitz space\/} $\sH^\rd$ 
(Def.~\ref{redHurSpace}).
Call this a {\sl $K$-cover point}. The branch locus $\bz$ of a $K$-cover point 
gives a $K$
point on $\prP^1_j\setminus
\{\infty\}$. Detailed analysis of the case $K=\bR$ is especially helpful when 
$p=2$. There
are three kinds of
$K$ points on $\sH^\rd$: 
\begin{edesc}  \item $K$-cover points (with $K$ understood, call these cover
points); 
\item {\sl $K$-Brauer\/} points (not $K$-cover points) given by $K$ covers
$\phi:X\to Y$ with
$Y$ a conic in
$\prP^2$ over
$K$ (\S\ref{redHtoH}); and \item $K$ points that are neither $K$-cover or $K$-
Brauer points. 
\end{edesc} 

When $\sH^\rd$ is a fine (resp.~b-fine) moduli space (\S\ref{fineModSect}), 
there are only
$K$-cover and $K$-Brauer points (resp.~excluding points over $j=0$ or 1).  
\S\ref{2A5Cusps} interprets how cusps attach to orbits of an element 
$\gamma_\infty$
in the Hurwitz monodromy group quotient we call $\bar M_4$. It also explains  H-M and
near H-M representatives (of a Nielsen class).  These correspond to cusps at the
end of components of $\bR$ points on levels $k\ge1$ of our main Modular Tower. 
These
(eight) $\bR$ points are exactly the 
$\bR$-cover points at level $k=1$ of the $(A_5,\bfC_{3^4})$ Modular Tower lying 
above any $j_0\in
(1,\infty)\subset \prP^1_j\setminus \{\infty\}$. 

There are also eight $\bR$-Brauer points, {\sl complements\/} of H-M and near H-M
reps.~(\S\ref{cuspNotat}). This gives the full set of
$\bR$ points on  the genus 12 component of 
$\sH(G_1,\bfC_{3^4})^{\inn,\rd}=\sH_1$ we
label
$\sH^+_1$.  Prop.~\ref{oneOrbitkappa} shows there is only one component of $\bR$ 
points on $\sH_1$.
So, $\sH^+_1$ contains all $\bQ$ points at level 1. 

The $k$th level ($k\ge
1)$,
$\sH_k^\rd$, $k\ge 1$, of this reduced Modular Tower has a component 
$\sH^{\rd,+}_k$ with
both H-M reps. and near H-M reps.~(Prop.~\ref{nearHMbraid}). There is a 
nontrivial central
Frattini extension
$T_k'\to G_k$ (Prop.~\ref{RkGk}).  Near H-M reps.~also correspond to $T_k'$ 
covers whose
field of moduli over
$\bR$ is
$\bR$, though $\bR$ is not a field of definition. There are also $\bR$-Brauer 
points
on $\sH^{\rd,+}_k$ (Prop.~\ref{HMnearHMLevel}). The production of an infinite 
number of
significantly different examples of these types of points from one moduli 
problem is a
general application of Modular Towers.

Real points on 
$\sH(A_5,\bfC_{3^4})^{\abs,\rd}$ (the absolute reduced Hurwitz space; 
equivalence
classes of degree 5 covers) and the inner reduced space  
$\sH(A_5\bfC_{3^4})^{\inn,\rd}$ are all
cover  points (Lem.~\ref{A5absR} and Lem.~\ref{A5innR}). This is a special case 
of the
general Prop.~\ref{isotopyM4} when the group $\sQ''$ acts trivially on Nielsen 
classes.  

\subsubsection{Analysis of level 1 and the Main Conjecture}  \S\ref{compA5MT}  
lists
geometric properties of level 1 of the  reduced $(A_5,\bfC_{ 3^4})$
Modular Tower ($p=2$). 

Level 0 is an irreducible  genus 0
curve with infinitely many $\bQ$ points \cite{Fr3-4BP}. A
dense set of $\bQ$  points produces a $\bQ$ regular
realization of $(A_5,\bfC_{3^4})$; each gives a $\bQ$ 
regular realization of the $A_5$ {\sl spin\/} cover $\Spin_5$ which identifies 
with $\SL_2(\bZ/5)$
(Rem.~\ref{plifting}). The level 1 group 
${}_2^1\tilde  A_5=G_1$ covers $A_5=G_0$, factoring through the spin
cover.  Thm.~\ref{twoOrbits}  shows level 1 has two absolutely irreducible $\bQ$
components  of  genuses 12 and 9. These correspond to two $H_4$
orbits,
$O_1^+$ and 
$O_1^-$,
in  the Nielsen class $\ni({}_2^1\tilde A_5,\bfC_{3^4})^\inn$ 
(\S\ref{defNielClass}).

Further, all $\bR$ (and therefore $\bQ$) points  lie only in the genus 12
component, which has just one connected component of real points
(Prop.~\ref{oneOrbitkappa}).  This  starts with showing the 16 {\sl H-M\/}  and 
16
{\sl near H-M\/}  representatives in the (inner --- not reduced; see
\eqref{nearHM} and \S\ref{cuspsH1})
Nielsen class fall in one
$H_4$ orbit.   Together these account for the eight longest (length-20) {\sl 
cusp
widths\/}  on the  corresponding {\sl reduced\/} Hurwitz space component.

Computations in 
\S\ref{2A5Cusps}--\S\ref{compA5MT}  show how everything organizes around these 
cusps lying at  the end of real components on $\sH_1$. This separates two
components of 
$\sH_1^\rd$ according to the cusps lying on each. The computer program \GAP\ 
helped 
guide our results (\S\ref{gaplist}), which include \GAP-less arguments 
accounting 
for all phenomena of this case. This reveals the group theoretic obstructions to 
proving the 
main conjecture (Prob.~\ref{MPMT}; for $r=4$).  We still need Falting's Theorem
to conclude there are only finitely many $(G_1,\bfC_{3^4})$ realizations over 
$\bQ$ (or over any
number field). 

Many Modular Towers should resemble this one, especially in the critical 
analysis of
H-M and near H-M reps. Detailed forms, however, of the Main Conjecture are
unlikely to be easy as we see from sketches of other examples in 
Ex.~\ref{A4C34-wstory} and
Ex.~\ref{A5C54-wstory}.  H-M and near H-M reps.~in these examples work for us, 
though they
differ in detail from our main example.  
 
\subsubsection{Reduced equivalence classes}
When $r=4$, analyzing the $\bar M_4=\PSL_2(\bZ)$ action on reduced Nielsen 
classes has a short 
summary (Prop.~\ref{j-Line}). Suppose $\bg=(g_1,g_2,g_3,g_4)\in G^4$ (where 
$G=G_k$ for some $k$) represents an element in a reduced Nielsen class. Then 
$\bg$ 
satisfies the product-one condition: $g_1g_2g_3g_4=1$. We have the following 
actions. 
\begin{edesc} \item The {\sl
shift\/}, 
$\gamma_1=\sh$  of \S\ref{shiftOp}, maps $(g_1,g_2,g_3,g_4)$ to 
$(g_2,g_3,g_4,g_1)$. 
\item The {\sl middle twist\/}  $q_2=\gamma_\infty$ maps $\bg$ to 
$(g_1,g_2g_3g_2^{-1} 
,g_2,g_4)$. \end{edesc}
Orbits of $\gamma_\infty$ correspond to cusps of a $j$-line cover. Each
$\lrang{\sh,\gamma_\infty} =\bar M_4$ (a mapping class group quotient) orbit
corresponds to a cover of the $j$-line. \S\ref{shincidence}  introduces  the  
(symmetric) 
{\sl $\sh$-incidence matrix\/}  summarizing a {\sl pairing\/}  using $\sh$ on
$\gamma_\infty$ orbits.  This matrix shows the $\bar M_4$ orbit structure of the
genus 12 component of $\sH_1^\rd$ (\S\ref{shinc1}). This  orbit contains all H-M
(or near H-M) orbits and their shifts  (Def.~\ref{hmorbit}).  

This illustrates analyzing the $\bar M_4$ and $\gamma_\infty$ orbits on reduced 
Nielsen classes, including arithmetic properties of cusps. Since 
components are moduli spaces this means analyzing degeneration of objects in the 
moduli 
space on approach (over $\bR$ or $\bQ_p$; \cite{FrLInv}) to the cusps. It 
generalizes 
analysis of elliptic curve degeneration in approaching a cusp of a modular 
curve. These are 
moduli spaces of curves. Here as in \S\ref{mochp}, Modular Towers
involve nilpotent fundamental groups that detect degeneration of curves, not 
just their Jacobians. \S\ref{genrShInc} shows the idea
of $\sh$-incidence works for general $r$. It is a powerful tool for simplifying the computation
of braid orbits on Nielsen classes, though we haven't had time to explore its geometry
interpretations for general $r$. 

\subsubsection{\IGP\ and obstructed components}  \label{IGPobstruct}
Up to w-equivalence \eql{covEq}{covEqw} there are but finitely many four branch 
point $\bQ$
realizations of $({}_2^1\tilde  A_5,\bfC_{3^4})$, a diophantine result
(Thm.~\ref{twoOrbits}).   Should any exist, they correspond to rational points 
on the genus 12 
(level 1) component  of the Modular Tower. From Lem.~\ref{4HM4NHM},  $\sQ''$ 
acts faithfully at
level 1. So, $\bQ$ points on
$\sH_1^\rd$ correspond  either to $\bQ$ regular realizations of the 
representation cover for
the Schur multiplier for
${}_2^1\tilde A_5$ ($\bQ$-cover points) or to $\bQ$-Brauer points 
(Lem.~\ref{redCocycle}). Which
is a diophantine subtlety.

Only the genus 12 orbit $O_1^+$ is unobstructed for the {\sl big\/}  $H_r$ 
invariant $\nu(O_1^+)$ of \cite[Def.~3.11]{FrMT}  (\S\ref{schMultInv}). The 
component 
corresponding to $O_1^-$ is { \sl obstructed\/}  \cite[\S0.C]{FrKMTIG}. The 
general meaning of
obstructed \wsp applied to a component (at level, say, $k$) \wsp is that it has 
nothing above it at
level
$k+1$.  This applies either to the geometric moduli space component or to the 
$H_r$ orbit
on a Nielsen class associated to it. Prop.~\ref{compSerre}  shows
$\nu(O_i)$,
$i=1,2$, distinguishes the components. That  is, $\nu(O_1^-)$ doesn't contain 1 
and the next levels
of the Modular Tower have {\sl  nothing\/}  over it.   The $\nu$ invariant 
separates most  Modular
Tower  components in this paper. A more general idea of {\sl lifting 
invariant\/} appears in
\S\ref{NielSep}. The more general idea is easier to understand, though not so
efficient for actually separating components. 

\subsubsection{Clifford algebras and Serre's invariant}  \cite{SeTheta}  
interpreted 
a case of this, for alternating groups. The ingredients included a Clifford 
algebra, and 
half-canonical classes ($\theta$-characteristics). Serre's invariant is the 
special case
of \S\ref{NielSep} using pullbacks of subgroups of alternating groups to their 
spin
covers.  It appears nontrivially at level  0 of each 
$A_n$  Modular Tower with $p=2$ and 3-cycle conjugacy classes  (Table 
\ref{An3rConst},
Prop.~\ref{3355obst},
\cite[Ex.~III.12]{FrMT}). 

The phrase {\sl spin separating representation\/}
refers to presenting a group inside the alternating group so it has a nonsplit
central extension through pull back to the spin cover 
(Def.~\ref{spinSepRepDef}). 
\S\ref{cliffAlgLift} uses a  Clifford algebra to produce  spin separating 
representations
of  $G_1={}_2^1\tilde A_5$. 
\S\ref{hcan} outlines the potential $\theta$-nulls (giving automorphic functions 
on
the Modular Tower levels) coming from spin separating representations 
(details in the expansion of \cite{FrLInv}). 

This explains part of this paper's
most mysterious  phenomenon: The two $\bar M_4$ orbits on $\ni({}_2^1\tilde
A_5,\bfC_{3^4})^{\inn,\rd}$  have the same image groups, and degrees. This is
despite their being different $\bar M_4$ representations distinguished just by 
the 
number and length of their $\gamma_\infty$  orbits. Appearance of the Clifford 
algebra dominates details of
our study of this  Modular Tower level 1 through component separation. 

We haven't yet shown a similar
{\sl separated components\/} phenomenon happens beyond level 1 of our main 
Modular Tower (see \S
\ref{shinc1}). Still, \S\ref{nofineMod} relates the subtly different geometry of 
the (real) cusps
associated with H-M reps.~versus near H-M reps.~to the phenomenon of fields of 
moduli versus field
of definitions attached to realizing {\sl Frattini central extensions}. This 
Schur
multiplier result (Prop.~\ref{nearHMbraid}) happens at all levels of the Modular 
Tower (certain
to generalize to most Modular Towers). This is a weaker form of the geometry of 
spin
separation and  an example of how Frattini central extensions (Schur 
multipliers) of perfect groups
affect the geometry of all levels of any Modular Tower. 

The two level 1 representations of
$\bar M_4$, with their relation to spin separating representations has 
ingredients like those in 
\cite{GKS}. So, \S\ref{GKS} uses this analogy to find  a precise
measure of the difference between the two representations. 

\subsection{General results and divergences with modular curves}   Let $K$ be a 
field.  A {\sl point\/}  on
a Modular Tower is a projective  system of points 
$\{\bp_k\}_{k+0}^\infty$: $\bp_{k+1}$  on the $k+{1}$st level maps to $\bp_k$ at 
level $k$.  
Similarly,  define 
points on reduced Modular Towers. If all $\bp_k\,$s have definition field 
$K$, then  this sequence 
defines a $K$ point. Suppose $K$ is a number field  (more generally, $K\le \bar 
\bQ$ has 
an infinite number of places with finite residue class field). Then,
Thm.~\ref{thm-rbound}   implies a Modular Tower (of {\sl inner\/}  Hurwitz 
spaces) has
no $K$ points.

\begin{prob}[Main Problem on Modular Towers] \label{MPMT}  Assume $(G,\bfC,p)$ 
is 
data for a 
Modular Tower with $G$ centerless and  $p$-perfect (Def.~\ref{pperfect}).
For
$k$ large,  show all 
components of a Modular Tower at level $k$ have general type.  Let $K$ be a 
number 
field. Find 
explicit large $k$ so the $k$th level of a Modular Tower of inner Hurwitz spaces 
contains 
no $K$ points (\cite{FrKMTIG}; see  \S\ref{usebcl}).\end{prob}  

\subsubsection{Fine moduli} \label{fineModSect}
The hypotheses on $(G,\bfC,p)$ say that $G$ has no $\bZ/p$  
quotient. The centerless condition  
ensures a Modular Tower of inner Hurwitz spaces consists of fine moduli spaces 
(Prop.~\ref{fineMod}). Our Main Example $(A_5,\bfC_{3^4} ,p=2)$ has all inner 
levels fine moduli 
spaces. Its reduced spaces are also, if $k\ge 1$ (Lem.~\ref{4HM4NHM}). Level 0, 
however, of this (inner and reduced) Modular Tower is not a fine or b-fine 
moduli space 
(\S\ref{znezero}).   

\subsubsection{Modular curve-like aspects}  Suppose $G=P\xs H$ where $H$ is a 
$p'$-group 
(\S\ref{notation}) acting through automorphisms on a $p$-group $P$. Then, $\tG 
p$ is $\tilde P\xs
H$  with $\tilde P$ 
the pro-free pro-$p$ group on the minimal number of generators for $P$ 
(Rem.~\ref{extH}  or 
\cite[Chap.~21]{ FrJ}). Refer to this as {\sl $p$-split data\/}  for a Modular 
Tower. 
\begin{exmp} [$H$ 
acting on a lattice $L$] \label{HL}  Assume a finite group $H$ acts on $L=\bZ^n$ 
(possibly the integers of a number field). For each prime $p$ not dividing
$|H|$,  form $H_p=L/pL\xs 
H$. Let $\tilde F_{p,n}$ be the pro-free pro-$p$ group of rank $n$ and let 
$\bfC$ be conjugacy
classes in
$H$. Form the Modular Tower for 
$(H_p,\bfC,p)$, based on the characteristic $p$-Frattini quotients of 
${}_p\tilde  H_p=\tilde F_{p,n}\xs H$ (directly from those of $\tilde F_{p,n}$).  

This has a modular curve-like property: You can vary $p$. Modular 
curves are 
the case $n=1$, $H=\bZ/2$ with $\bfC$ four repetitions of the nontrivial element 
of $H$. Whenever  
$n=1$, the analysis is similar to modular curves. Compare the case $n=1$, 
$H=\bZ/3$, $\bfC$ two
repetitions  each of the conjugacy classes of $\pm 1\in \bZ/3$ and $p\equiv 1 
\bmod 3$ with
the  harder case where $n=2$, $p=2$ in Ex.~\ref{A4C34-wstory}.  
\end{exmp} 

\begin{exmp}[$G$ any simple group] \label{SG}  Characteristic 
$\tG p$ quotients  $\{G_k\}_{k=0}^\infty$ are 
perfect and centerless (Prop.~\ref{fineMod}). For each $k\ge 1$,
the natural map $G_{k+1}\to G_k$ maps through the exponent $p$ part of the 
universal central 
extension $\hat G_k$ of $G$.  \end{exmp}

Modular Towers for Ex.~\ref{SG}  depart from towers of moduli spaces in the 
literature.
This paper emphasizes new group and geometry phenomena amid 
comfortable similarities to modular curves even for the many Modular Towers 
covering the 
classical $j$-line. This paper aims to show, when $r=4$, the genus of 
reduced Modular Tower components grows as the levels rise. 

Modular curve-like aspects in Ex.~\ref{SG} guide finding information 
about cusp 
behavior. The  phenomenon of {\sl obstructed components\/}, appearing in our 
main
example, is the most non-modular  curve-like aspect of general  Modular 
Towers. 
Suppose in Ex.~\ref{HL}, for a given $p$, $H$ acts irreducibly on $L/pL$. Then, 
the 
Nielsen class 
$\ni(H_p,\bfC)$ is nonempty (Lemma \ref{HLcont}).  Deciding in 
Ex.~\ref{HL}  if there are 
obstructed  components ($H_r$ orbits) at high levels above a  level 0 $H_r$ 
orbit seems 
difficult. 

When $r=4$,  this reflects  on the genuses of the corresponding $j$-line covers.    
Even
when 
$n=2$, $p=2$, 
$H=\bZ/3$, and $\bfC$ consists of several copies of the two nontrivial conjugacy 
classes of 
$\bZ/3$, obstructed components appear at all level $k$ for $r$ suitably large 
(dependent on
$k$). This special case of \cite[Conway-Parker App.]{FrVMS} uses the nontrivial
Schur multiplier appearing at {\sl each\/} level (Prop.~\ref{RkGk}  generalizing 
\cite[\S4]{FrKMTIG}). The case
$G_0=(\bZ/2\times\bZ/2)\xs 
\bZ/3=A_4\ (\le A_5)$ 
has elements like our main example, because  their  universal $2$-Frattini 
covers are close 
(\S\ref{pfrattini}). Yet, for $\bfC$ consisting of $r=4$ ($r$ not large)  3-
cycle conjugacy classes,
the early levels of the Modular Towers for $A_4$ and $A_5$ ($p=2$) are not at 
all alike
(Ex.~\ref{A4C34-wstory}).  

These calculations collect the key points for considering a full proof of the 
Main
Modular Tower Conjecture when $r=4$. Using this paper, \cite[Main  
Theorem]{FrLInv}  will show
what to expect of all  Modular Tower 
levels from $A_n$ and $r$ ($\ge n-1$) 3-cycles as conjugacy classes. 

\subsubsection{Immediate extensions when $r=4$}  A full result on $A_5$, $p=2$ 
and $r=4$
would  contribute to the 
following questions on the \IGP. \begin{edesc}  \label{compA5p2r4}  \item 
\label{compA5p2r4a}  Does any 
other reduced $A_5$ Modular Tower have $\bQ$ points at level 1? \item 
\label{compA5p2r4b}  Is there any 
set of four conjugacy classes $\bfC$ for ${}_2^1\tilde A_5$ where the reduced 
inner 
Hurwitz space contains 
infinitely many (any?) $\bQ$ points? \end{edesc}  An example best illustrates 
the serious
points about Question 
\eql{compA5p2r4}{compA5p2r4a}. Let $\C_5^+$ be the conjugacy class of 
$(1\,2\,3\,4\,5)$, $\C_5^-$ the conjugacy class of $(1\,3\,5\,2\,4)$ (with 
$\C_3$ the
conjugacy class of  3-cycles). 

From the \BCL, the only positive possibilities for \eql{compA5p2r4}{compA5p2r4b} 
are Modular Towers
for $A_5$ with 
$\bfC=(\C_3,\C_3,\C_5^+,\ 
C_5^-)$ or $(\C_5^+,\C_5^+,\C_5^-,\C_5^-)$. For both, the level 0 Modular Tower 
has two 
components, with exactly one obstructed (nothing above it at level 1). The 
respective 
unobstructed $H_4$ 
orbits contain the following representatives: \begin{equation}  \label{othbfC}  
\begin{array} {rl} \text{ a)\ } 
&((2\,5\,4),(4\,5\,3),(1\,2\,3\,4\,5),(2\,1\,5\,3\,4)) \\ \text{ b)\ }&
((1\,2\,3\,4\,5),(5\,4\,3\,2\,1),(2\,1\,4\,3\,5), (5\,3\,4\,1\,2)).\end{array} 
\end{equation} 
Prop.~\ref{3355obst} explains this and why only the (\ref{othbfC}b) 
Modular Tower has real points at level 1. Rem.~\ref{notinftyA5} explains the one 
mystery
left on answering \eql{compA5p2r4}{compA5p2r4b}, about genus 1 curves at level 1 
of the
reduced Modular Tower for $(\C_5^+,\C_5^+,\C_5^-,\C_5^-)$. 

Modular Towers use $p'$ conjugacy classes. If we allow conjugacy classes where 
$p$ 
divides the elements' 
orders, then Hurwitz spaces attached to $G_k$ and these classes must have $r$ 
large if the 
components are 
to have definition field $\bQ$ (Thm.~\ref{thm-rbound}). There appear, however, 
to be 
nontrivial cases of 
$\bQ$ components for $r=4$ and ${}_2^1\tilde A_5$ in 
\eql{compA5p2r4}{compA5p2r4b} using
conjugacy  classes of order 4. Completing a yes answer to the following question 
requires
completing  these cases.  
\begin{prob} Are there only finitely many ${} _2^1\tilde A_5$ realizations (over 
$\bQ$; 
up to 
w-equivalence) with at most four branch points?  \end{prob} For, however, the 
group 
${}_2^2\tilde A_5$, if 
$r=4$, the finitely many possible realizations can only fall on {\sl level 2\/}  
of the Modular 
Tower for 
$(A_5,\bfC_{3^4})$ or for (\ref{othbfC}b).


\subsubsection{Extending results to $r\ge 5$} \label{rge5extens}
\cite{Me} produced $\bQ$ covers in the Nielsen class $\ni(A_n,\bfC_{3^{ n-1}})$ 
when 
$n\ge 5$ is odd (see also \cite[\S9.3]{Se-GT}).  \cite[Thm.~1]{FrLInv} shows 
there is 
exactly one component here. So, the results combine to show $\bQ$ points are 
dense in 
$\sH(A_n,\bfC_{3^{ n-1}})^\inn$: There are many $(A_n,\bfC_{3^{ n-1}})$ 
realizations 
over $\bQ$ that produce $\Spin_n$ realizations.

\begin{prob} \label{Anprob} Does  Thm.~\ref{twoOrbits} generalize to say 
$\sH({}_2^1\tilde A_n,\bfC_{3^{n-1}})^{\inn,\rd}(\bQ)$ is finite for $n\ge 5$ 
odd? 
\end{prob} For $n$ even this follows from \cite[Thm.~1]{FrLInv} showing
$\sH({}_2^1\tilde A_n,\bfC_{3^{n-1}})^{\inn,\rd}$ is empty. Falting's Theorem 
\cite{FaMordCon} works when $n=5$ because the moduli spaces are 
curves. There is no replacement yet for Falting's Theorem for $n\ge 6$.

Here is a special case of the problem that would produce positive results for 
the Inverse Galois
Problem. 

\begin{prob} \label{AnprobPos}  
Is $\sH({}_2^k\tilde A_n,\bfC_{3^{r}})^{\inn,\rd}(\bQ)\ne \emptyset$
for $r$ large (dependent on  $n,k$). \end{prob}

Ex.~\ref{full3cycleList} describes the complete set of components for 
$\sH({}_2^1\tilde
A_n,\bfC_{3^{r}})^{\inn,\rd}$ (see Table \ref{shincni0}). The procedure of
\cite[Thm.~1]{FrLInv} uses the many embeddings of $A_k$ in $A_n$ for $n\ge k$. 
Any such embedding 
extends to an embedding of the level $k$ characteristic 2-Frattini cover of 
$A_k$ into the level
$k$ characteristic 2-Frattini cover of $A_n$ (special case of 
Lem.~\ref{fratIsom}). So,
Thm.~\ref{schMultInv}  describing the two  components of
$\sH({}_2^1\tilde A_5,\bfC_{3^{4}})^{\inn,\rd}$ immediately gives information on 
components
of $\sH({}_2^1\tilde A_n,\bfC_{3^{r}})^{\inn,\rd}$ for all $r\ge n-1$ and $n\ge 
5$. 
\cite[Conway-Parker App.]{FrVMS} proves for $n,k$ fixed and $r$ large each 
element of the Schur
multiplier of ${}_2^k\tilde A_n$ determines a component. This is a special
case of a general result, though here we are precise about the components
for all values of $r$. Combining this with 
\S\ref{indSteps} hints at describing components for higher levels of Modular 
Towers for {\sl all\/}
alternating groups. 

There are many
$j$-line covers. It is significant to find simple invariants distinguishing 
Modular Towers
from general quotients of the upper half plane by a subgroup of
$\PSL_2(\bZ)$. As this paper and its consequents shows, they have far more in 
common with
modular curves than does a general such quotient. The structure of a symmetric
integral $\sh$-incidence matrix makes such a distinction, effectively
capturing complicated data about Modular Tower levels. Example: It displays each 
connected component  (Lem.~\ref{shincBlocks}).  The cusp pairing for $r=4$ 
extends to the case $r\ge 5$ (\S\ref{genrShInc}). 

\subsubsection{Significance of cusps}
Thm.~\ref{twoOrbits}  results from information about the level 1 cusps of this 
Modular 
Tower. Most interesting are cusps attached to H-M representatives. \cite[\S 
III.F]{FrMT} explains
the moduli manifestation of total degeneration around such a cusp.  The 
compactification is
through equivalence classes of {\sl specialization sequences\/}. Especially 
short specialization sequences detect H-M reps.~from
total degeneration. This is a very
different compactification from that of the  {\sl stable compactification
theorem\/} cover versions, like the log structures \cite{WeFM} uses.  
Specialization sequences 
require no extension of the base field. So they are compatible with absolute 
Galois group
actions. H-M reps.~give the type of moduli degeneration  useful for a  {\sl 
tangential basepoint\/} 
(language of Ihara-Matsumoto
\cite{IharMat}  and Nakamura 
\cite{tanBasePts}). These  also use compactifications not quite like stable 
compactification (see
\cite[p.~163, Comments C.5]{FrMT}). Likely for $r=4$ there is not a huge 
difference between these
compactifications, though that is unlikely for $r\ge 5$. 

For Modular Towers over
$\bQ$ with $p=2$ and
$r=4$ (or higher) consider  pull back to the configuration space 
$\prP^4\setminus D_4$ of branch
points of these  covers. Then, a real point on $\sH^\inn_k$  with $k\ge 1$ 
corresponds to a
cover with  complex conjugate pairs of branch points (Thm.~\ref{thm-rbound}). 
So, one must recast
relations on the image  of $G_\bQ$ coming from the Grothendieck-Teichm\"uller 
group. More
appropriate for this situation than  four real points is complex conjugate pairs 
of branch points
(compare with 
\cite[p.~107]{IharIntCong}). For other problems one would want configurations 
with a complex
conjugate pair and two real points of branch points (as in \S\ref{covjvalues}). 
This adapts for 
$G_\bQ$ acting on quotients  of the profinite completion of $H_4$ corresponding 
to projective
limits of the monodromy groups of monodromy group cover from Modular Towers with
$p=2$.  

Here we connect to \cite{BeFr} at level 0. An example result gives exactly {\sl 
three\/} degree five $(A_5,\bfC_{3^4})$ covers of $\prP^1_z$ (up to w-
equivalence 
\eql{covEq}{covEqw}), over $\bQ$ with branch points in $\bQ$. There are three 
non-cusp 
points in the order 12 group generated by the cusps on the {\sl absolute\/} 
reduced space 
attached to $(A_5,\bfC_{3^4})$.  These three covers correspond to those three 
non-cusp 
points. Mazur's explicit bound on torsion points shows, among all 
$(A_5,\bfC_{3^4})$ 
realizations  (over $\bQ$),  exactly three w-equivalence classes can have 
rational branch 
points (\cite{BeFr}).  A result of Serre  (Prop.~\ref{serLift}; \cite{SeLiftAn}) 
allows 
computing cusp widths in the $({}_2^1A_5,\bfC_{3^4})$ case. This includes the 
length-20 cusps attached to H-M and near H-M representatives. Computing the 
genus of 
curves at all levels of a Modular Tower depends on such formulas.

\subsubsection{Characterizing divergences with modular curves} We 
show many phenomena for Modular Tower levels of simple groups not appearing in 
modular curves. There is a precise dividing line between two types of Modular 
Towers, with modular curves a model for
one type.
\cite{Fr-Se1} told the story of how modular curves are essentially the moduli 
for dihedral group realizations with four
involution conjugacy classes. The absolute (resp.~inner) Hurwitz spaces 
correspond to the curves $X_0(p^{k+1})$
(resp.~$X_1(p^{k+1})$). The case of general $r\ge 4$ is about Hurwitz spaces 
associated to modular
curve like covers of the moduli of hyperelliptic curves of genus $g=\frac{r-
2}2$. A generalization of this situation
would include $G_0=bZ/p^t\xs A$ with $A$ a(n abelian) subgroup of $(\bZ/p^t)^*$ 
(acting naturally on $G_0$) and a
collection of conjugacy classes
$\bfC$ whose elements generate $A$.  For, however, inner Hurwitz
spaces to be over $\bQ$, $\bfC$ must be a rational union of conjugacy classes 
(trivial case of the branch cycle argument).
This puts a lower bound on $r$. For example, if $A$ is cyclic and order $d$, 
then $\bfC$ must contain a minimum of
$\phi(d)$ conjugacy classes.  That is why one rarely sees these spaces in the 
classical context, though they share with
modular curves the attribute of having 1-dimensional $p$-Frattini modules. 

For a $p$-perfect centerless group, the characteristic $p$-Frattini $G_0$ module 
$M_0$ 
has dimension 1 if and only if $\one_{G_k}$ never appears in the Loewy layers of 
$M_k$ for {\sl any\/} $k$
(\S\ref{onesAppearNGP} and \S\ref{genFrattini}). Appearance of those 
$\one_{G_k}\,$s could mean 
obstructed components, uncertain location of cusps in components and related 
moduli interpretations of spaces 
for covers whose field of moduli is not a field of definition of representing 
points of the levels of a Modular Tower.
These useful geometric phenomena, present Modular Towers as a new tool for
investigating still untouched  mysteries. That is our major theme. Further, 
given
$k$, if the conjugacy classes in
$\bfC$ repeat often enough (as a function of $k$),
\cite[App.]{FrVMS} implies these complications will occur at level $k$ of the 
Modular Tower. A big mystery is
whether they occur at infinitely many levels of {\sl any\/} Modular Tower for a 
$p$-perfect centerless group.

\subsubsection{The Open Image Theorem and $j$-awareness}
\S~\ref{openIm}  
formulates how to extend part of Serre's {\sl Open Image Theorem\/}  
\cite{SeElliptic}  (for just the prime $p$)
to  Modular 
Towers. This extension would be a  tool for finding a precise lower bound on $k$ 
for 
which  higher levels of a Modular Tower have  no rational points. Ex.~\ref{HL} 
is appropriate for generalizations of Mazur-Merel. 

\begin{prob} Consider $\sH(L/p^kL\xs H,\bfC)^{\inn,\rd} =\sH_{L/p^kL\xs H,\bfC} 
$ 
with $p$ running over primes not dividing $|H|$ and integers $k\ge 0$. Find 
$m(H,\bfC)$  so  
$\sH_{L/p^kL\xs H,\bfC} $ has general type if $p+k$ exceeds $m(H,\bfC)$. If 
$r=4$, and 
$d\ge 1$, find 
explicit $m(H,\bfC,d)'$ with $\sH_{L/p^kL\xs H,\bfC} (K)$ empty for $k+p> 
m(G,\bfC,d)'$ and 
$[K:\bQ]\le d$. \end{prob}  Mazur-Merel is the case $H=\bZ/2=\{\pm 1\} $ and 
$L=\bZ$ 
($-1$ acting by 
multiplication), $r=4$ and classes all the nontrivial element in $H$.

We continue with the case $r=4$. What is the analog of Serre's result for 
nonintegral $j_0\in
\prP^1_j\setminus
\{\infty\}$.  We interpret that to say, for $j_0\in \bar \bQ$ not an algebraic 
integer, the action
of $G_{\bQ(j_0)}$ on  projective systems of points on 
$\{X_1(p^{k+1})\}_{k=0}^\infty$ lying over $j_0$ has an open
orbit. The ingredients in our generalizing statement (still a conjecture) is a 
Modular Tower meaning to $j_0$ being
suitably ($\ell$-adically) close to $\infty$ and the Ihara-Matsumoto-Wewers 
approach to using tangential base points.
One topic we couldn't resist was what we call $j$-awareness. It would be trivial 
to generalize the Open Image Theorem if each
appropriate Modular Tower was close to being a tower of modular curves. That 
would  mean there is some cover $Y\to \prP^1_j$ so
the pullback of a tower of modular curves dominates the levels of the given 
Modular Tower. Checking this
out on the main examples of this paper (with $p=2$) repeatedly called for 
analyzing Prym varieties   (Ex.~\ref{SnAnbcycles}
and 
\S\ref{jawareMT}).

\subsubsection{Developing a $p$-adic theory} \label{mochp} We compare our
approach to $\bR$ points with that of \cite[p.~978]{MochpUnif}. Recall the
complex conjugate of a complex manifold. Suppose 
$\{ (U_\alpha,
\phi_\alpha)\}_{\alpha\in I}$ is an atlas for a complex manifold $X$. For
simplicity, assume $X$ is a 1-dimensional complex manifold. Create a new
manifold by composing each $\phi_\alpha: U_\alpha\to \bC$ with complex
conjugation. Call the resulting map $\phi^*_\alpha$. 

\begin{lem} The atlas $\{ (U_\alpha,\phi^*_\alpha)\}_{\alpha\in I}$ is a new 
complex manifold
structure $X^*$ on the set $X$. \end{lem}

\begin{proof} Denote by $\bar z_\alpha(x)$ the value of $\phi^*_\alpha(x)$. 
Compare the
transition functions $\phi^*_\beta\circ (\phi^*_\alpha)^{-1}(\bar z_\alpha)=\bar 
z_\beta$
with the function $\phi_\beta\circ \phi_\alpha^{-1}(z_\alpha)$ as $\bar 
z_\alpha$
varies over the complex conjugate of $z_\alpha$ running over 
$\phi_\alpha(U_\alpha\cap U_\beta)$. The effect of the former is this: $$\bar
z_\alpha\mapsto z_\alpha\mapsto
\phi_\alpha^{-1}(z_\alpha)
\mapsto
\phi_\beta(\phi_\alpha^{-1}(z_\alpha))=z_\beta \mapsto \bar z_\beta.$$ Suppose 
$f$ is a
local expression of the transition function $\phi_\beta\circ \phi_\alpha^{-
1}(z_\alpha)$
as a power series (about the origin) in $z_\alpha$. Then, the power series 
expressions
for
$\phi^*_\beta\circ (\phi^*_\alpha)^{-1}(\bar z_\alpha)$ comes by applying 
complex
conjugation to the coefficients of $f$, and so the resulting function is 
analytic in $\bar z_\alpha$. 
\end{proof}

\begin{defn} A complex manifold $X$ has an {\sl $\bR$ structure\/} if $X^*$ is
analytically isomorphic to $X$. 
\end{defn}

Apply this to a cover $\phi: X\to \prP^1_z$ using the complex structure from the
$z$-sphere, as in Prop.~\ref{compTest}. You can also use the complex structure
from uniformization by the upper half plane $\bH$ as in  
\cite[p.~978]{MochpUnif}. Keep in mind: $X^*$ is the same set as
$X$. For $x\in X$, there is a natural map
$x\mapsto x^*\in X^*$ where we regard $x^*$ as the complex conjugate of $x$,  
though in
this formulation $x^*$ is the same point on the set $X$. Also, the lower half 
plane $\bH^*$
naturally uniformizes $X^*$. An
$\bR$ structure induces
$\psi:X^*\to X$, which in turn induces $\tilde \psi: \bH^*\to \bH$. Of course, 
$\psi$
induces $\psi^*: X\to X^*$. Mochizuki makes simplifying assumptions: There is 
$x^*$ 
with $\phi(x^*)=x$ and $\psi^*\circ \psi$ the identity map. 

Take  $\sE(X,x)^{\alg}$ to be all meromorphic algebraic
functions in a neighborhood of $x$ that extend analytically  along each path in 
$X$
(\S\ref{groth}).  Then,
$\sE(X,x)^{\alg}/\bR(z)$ is Galois with profinite group $\pi_1(X,x)^{\alg,\bR}$ 
forming a split
sequence $\pi_1(X,x)^{\alg}\le \pi_1(X,x)^{\alg,\bR}\to G(\bC/\bR)$ (as in 
Prop.~\ref{splitseq}).
The virtues of this \cite{gunInd} inspired approach are these. 

\begin{edesc} \item A splitting of $G(\bC/\bR)$ given by $\tilde \psi\circ C$, 
an anti-holomorphic
involution  $\smatrix a b c d\in GL_2^-$ (determinant -1) acting on $\tau\in 
\bH$ by
$\tau\mapsto \frac{a\tau+b}{c\tau+d}$. 
\item So,  $\pi_1(X,x)^{\alg,\bR}$ acts through the matrix group $\GL_2(\bR)$ 
extending
$\rho:\pi_1(X,x)^{\alg}\to \PSL_2(\bR)$. 
\item Extend $\rho$ by $\PSL_2(\bR)\to \PSL_2(\bC)$ to consider the 
$\prP^1_\tau$ bundle over $X$
with a natural flat connection (algebraic from Serre's GAGA). 
\end{edesc} 

Mochizuki calls this uniformization construction the {\sl canonical indigenous
bundle\/}. Compare to our Prop.~\ref{compTest} with 
several complex conjugation operators
$\hk$, corresponding to a type of degenerate behavior at a cusp (and the
location of a set of branch points over
$\bR$ on $\prP^1_z$). His simplifying assumption ($\phi(x^*)=x$) applied to
Modular Towers when $p=2$ would reduce to considering $X$ that appear as
H-M reps.~(say, in Cor.~\ref{reducedRpts}), just one case we must treat in 
describing the full real locus for
a result like Prop.~\ref{oneOrbitkappa} on a Modular
Tower level.

Mochizuki calls his prime $p$, though in our context it would be $\ell$ prime to 
$|G=G_0|$, since we use $p$ as a
prime for the construction of a Modular Tower. His construction is for any 
family of curves over any base. It extends to
work for
$\ell$-adic uniformization giving a notion of $\ell$-ordinary points on a family 
of curves over $\bQ_\ell$. The Modular
Tower goal  would be to detect  {\sl ordinary\/}  $\bQ_\ell$ (with 
$(\ell,|G|)=1$) points on a Modular Tower, 
generalizing (from
$G_0$ a  dihedral group) ordinary  elliptic curves over $\Spec(\bZ_\ell)$ 
\cite[p.~1089]{MochpUnif}. Only the genus 12 
component at level 1 
of the $(A_5,\bfC_{3^4} )$ Modular Tower supports totally degenerate cusps, 
associated 
with the name {\sl 
Harbater-Mumford\/}  representatives (\S\ref{startHM}). An analysis that 
combines {\sl  
Harbater 
patching\/}  and Mumford's theory of total $\ell$-adic degeneration shows most 
points on 
such components 
will exhibit {\sl ordinary\/}  $\ell$-adic behavior. Nilpotent fundamental 
groups enter 
directly in all Modular 
Tower definitions.  

Wewers \cite{WeFM}
uses the  degeneration behavior type (a tangential basepoint) to consider 
$\ell$-adic ($(\ell,|G|)=1$)
information for level 0 of our main Modular Tower. His goal is to detect  
fields of moduli versus fields of definition (following \cite{DDEm}).  
\S\ref{nofineMod} explains why the geometry of H-M reps.~and near H-M 
reps.~should have an $\ell$-adic
analog  at all levels of our $A_5$ Modular Tower.  His analysis is also a tool 
for our approach to generalizing Serre's
Open Image Theorem (\S\ref{findFrob}). 

\section{Presenting $H_4$ as an extension}

\subsection{Notation and the groups $B_r$, $H_r$ and $M_r$} \label{braidGp} 
Denote 
the Artin Braid group on $r$ strings by $B_r$. It has  generators $\row Q{r-1}$ 
satisfying these relations: \begin{eqnarray} \label{bdgp} 
Q_iQ_{i+1}Q_i=Q_{i+1}Q_iQ_{i+1},& &i=1,\ldots,r- 2; \text{and}\\ 
Q_iQ_j=Q_jQ_i,& &1\le i< j-1\le r-1.  \end{eqnarray}
\newcommand{\bars}{{\bar \sigma}} \newcommand{\bbs}{{\pmb {\bars}}}

\subsubsection{$B_r$ and $H_r$ as automorphisms}
We start from Bohnenblust \cite{Bohnen} (or \cite{Markov}). Let $F_r$ be 
the free group of rank $r$ on generators $\{\row  {\bar\sigma} r\}$, and  denote  
$\Aut(F_r)$ by $\bar A_r$. Then $B_r$ embeds in $\bar A_r$ via this right action 
of  $Q_i$:  
\begin{equation} \label{brAct} 
(\row{\bars}r)Q_i=(\bars_1^{\vphantom{1}},\ldots,\bars_{i- 
1}^{\vphantom{1}},\bars_{i}^{\vphantom{1}}\bars_{i+1}^{\vphantom{1}}\bars_{i} 
^{- 1},\bars_{i}^{\vphantom{1}},\bars_{i+2}^{\vphantom{1}},\ldots, 
\bars_r^{\vphantom{1}}).\end{equation} Also, $B_r$ consists of automorphisms of 
$F_r$ fixing $\bar\sigma_1\cdots \bar\sigma_r$ and mapping $\row  {\bar\sigma} 
r$ to 
permutations of conjugates of these generators. Let $G_r$ be the  quotient of  
$F_r$ by 
the relation  $\Pi(\pmb \bars)=\bars_1\cdots \bars_r=1$. 

The {\sl  Hurwitz monodromy group\/} $H_r$ is the  quotient of $B_r$ by adding 
the 
relation \begin{equation} \label{hwgp} Q_1Q_2\cdots Q_{r-1}Q_{r-1}\cdots 
Q_1=1.\end{equation} It is also the fundamental group of an open subset of 
$\prP^r$ 
(\S\ref{basFGs}).  

Denote inner automorphisms of a group $G$ by $\Inn(G)$.  Then, 
$H_r$ induces automorphisms of $G_r$ mapping the images of $\row  {\bar\sigma} 
r$ to 
permutations of conjugates of these generators, modulo $\Inn(G_r)$. 

If $G\le H$, denote the normalizer of $G$ in $H$ by $N_H(G)$. 
The next lemma is obvious. Suppose $\psi:
G_r\twoheadrightarrow G$. Define an induced action of
$H_r$ on the collection of such $\psi\,$s by
$(\psi(\row{\bars}r))Q\eqdef\psi((\row{\bars}r)Q)$.    

\begin{lem} \label{conjComm} The action of \eqref{brAct} commutes with 
conjugation on an
$r$-tuple by an element of $F_r$. For $G\le S_n$ suppose $\psi:
G_r\twoheadrightarrow G$. The induced action of
$H_r$ on $\psi$ commutes with conjugation by $N_{S_n}(G)$. 
\end{lem}

\subsubsection{A presentation of $M_r$} The {\sl mapping 
class  group\/}  $M_r$ is the quotient by $\Inn(G_r)$ of automorphisms of $G_r$
mapping $\row  {\bar\sigma} r$ to permutations of conjugates. This  
maps $B_r$  to $M_r$ factoring through $H_r$. Further,  $M_r$ is the quotient
of 
$B_r$ by the following relations (\cite[\S3.7]{KMS} or \cite{Mag}): 
\begin{equation} 
\label{mcgp} \begin{array}{rl}\tau_1 = (Q_{r-1}Q_{r- 2}\cdots Q_{2})^{r-
1},\tau_2= 
Q_1^{-2}(Q_{r-1}\cdots Q_{3})^{r-2},\ldots ,\tau_{\ell+1} = &\\ (Q_\ell\cdots 
Q_1)^{-
\ell -1}(Q_{r-1}\cdots Q_{\ell +2})^{r-\ell -1}, \ldots ,\tau_{r-1} = (Q_{r-
2}\cdots 
Q_1)^{1-r}, & \\ \text{and } \tau = (Q_{r-1} \cdots Q_1)^r.& 
\end{array}\end{equation}

This complicated presentation is oblivious to the map $H_r\to M_r$ dominating 
this
paper. Using it conceptualizes the $B_r\to M_r$ kernel, 
$N_r=\ker(B_r\to M_r)$. The switch of emphasis shows 
in Prop.~\ref{Riopers}.  With $Q_1Q_2Q_{3}^2Q_2Q_1=D=R_1$ write generators of 
$N_r$ as follows: $$R_1,\ R_2=Q_1^{-1}R_1Q_1,\ R_3=Q_2^{- 
1}R_2Q_2,\ R_4=Q_3^{-1}R_3Q_3  \text{ and } (Q_1Q_2Q_3)^4.$$   
Thm.~\ref{presH4} presents both $M_4$ and $N_4$ memorably. An appearance of 
$N_4$
in \S\ref{finConc} offers a place for {\sl spin separation\/} to explain
similarities of the two 
$M_4$ orbits in level 1 of our main example.   

\subsubsection{Notation for fundamental groups} \label{notation} The Riemann 
sphere
uniformized by the  variable $z$ is $\prP^1_z=\bC\cup \{\infty\}$. Relation 
\eqref{hwgp}
comes  from the geometry of $r$ distinct points on $\prP^1_z$. We sometimes drop 
the 
notation for uniformizing by $z$.  Similarly, 
copies of $\afA^r$ (affine $r$-space) and $\prP^r$ (projective $r$-space) come 
equipped 
with coordinates suitable for a diagram: 
\begin{equation}\label{fundgpdiag} \begin{matrix}\afA^r\setminus 
\Delta_r&\mapright{}&{(\prP^1)^r\setminus \Delta_r}\\ 
\lmapdown{\Psi_r}&{}&\mapdown{\Psi_r}\\ {\afA^r\setminus 
D_r}&\mapright{}&\prP^r\setminus D_r.\end{matrix}\end{equation} The upper left 
copy 
of $\afA^r=\afA^r_\bz$ has coordinates $(\row z r)$.  Then, $\Delta_r$ is the 
fat diagonal 
of $r$-tuples with two or more of the $z_i\,$s equal.  Embed $\afA^r_\bz$ in 
$(\prP_z^1)^r$ by  $(\row z r)\mapsto (\row z r)$. In $(\prP_z^1)^r$ the 
$z_i\,$s may 
take on the value $\infty$. Regard $\Psi_r$ as the quotient action of $S_r$ 
permuting the 
coordinates of  $(\prP^1)^r$. As with $B_r$, put this action on the right.

Thus, the lower copy of $\afA^r=\afA^r_\bx$ has coordinates $(\row x r)$, with 
$x_i$ 
equal to $(-1)^i$ times the $i$th symmetric function in $(\row z r)$.  Regard 
this as giving 
a monic polynomial of degree $r$ in $z$ (with zeros $(\row z r)$). Then, 
complete the 
commutative diagram by taking $\prP^r$ as all monic (nonzero) polynomials of 
degree at 
most $r$. The map $\Psi_r$ on the right sends each 
$r$-tuple $(\row z r)$ to $\prod_{i=1}^r(z-z_i)$. When $z_i=\infty$ replace the 
factor $(z-z_i)$ by 1. The image of $\Delta_r$ is the discriminant locus of 
polynomials 
with two or more equal zeros. For coordinates around $z_i=\infty$, regard a 
monic 
polynomial of degree $t<r$ as having $r-t$ zeros at $\infty$.

This interprets $U_r=\prP^r\setminus D_r$ as the space of  $r$ distinct 
unordered points 
in $\prP^1$, the image of $(\prP^1)^r\setminus \Delta_r=U^r$.  Thus, $\Psi_r: 
U^r\to 
U_r$ is an unramified Galois cover with group $S_r$. 
The alternating group of degree $r$ is $A_r$. The dihedral group of degree $r$ 
(and order 
$2r$) is $D_r$. Denote the $p$-adic integers by $\bZ_p$, with $\bQ_p$ for its 
quotient 
field. Suppose $p$ is a prime and $\C$ is a conjugacy class in a group $G$. Call 
$G$ a 
$p'$ group if $(|G|,p)=1$. Call $\C$ a $p'$-conjugacy class (or a $p'$ class)
if elements in $\C$ have orders prime to $p$. This applies to
conjugacy classes  in profinite groups.

Each {\sl Hurwitz space\/} (our
main moduli spaces) comes from  a finite group $G$, a set of 
conjugacy classes $\bfC$ in $G$ and an equivalence relation on the corresponding 
{\sl Nielsen
class}. \S\ref{defNielClass} reviews this definition 
 while enhancing traditional equivalences in  \cite{FrMT}, \cite{FrVMS},
\cite{MM} and
\cite{VB}.

\subsection{Geometrically interpreting $M_4$ and $\bar M_4$}
\label{basFGs}  The fundamental group of $\afA^r\setminus D_r$
(resp.~$\prP^r\setminus D_r$) is the Artin braid group
$B_r$ (resp.~the Hurwitz monodromy group   
$H_r$) (\cite[\S4]{FrHFGG}, \cite[Chap. III]{MM}, \cite[Chap. 10]{VB}). 
Embedding
$\afA^r$ into 
$\prP^r$ in \eqref{fundgpdiag} gives the lower row surjective homomorphism from 
$B_r$ 
to $H_r$. The result is a commutative  diagram of fundamental groups induced 
from a 
geometric commutative  diagram.

Fundamental groups in the \eqref{fundgpdiag}  upper row are the {\sl straight\/} 
(or {\sl  
pure\/}) Artin braid and Hurwitz monodromy groups. That is, 
$SH_r=\pi_1((\prP^1)^r\setminus\Delta_r,\bx_0)$ is  the kernel  of 
$\Psi^*_r:H_r\to S_r$ 
mapping $Q_i$ to $(i\,i+1)$, $i=1,\ldots,r$. It increases precision to use 
capitals for 
generators $\row Q  {r-1}$ of $B_r$ and small letters for their images $\row q 
{r-1}$ in 
$H_r$.

\subsubsection{A cover corresponding to $M_r$} Consider $\bz_0$, an $r$-tuple of 
distinct points in $\prP^1_z$. 
Let $\row P r$ be classical generators of $\pi_1(U_{\bz_0},z_0)$ as in 
\S\ref{classGens}.
Suppose $\tilde U_r$ is the universal
covering space of
$U_r$. The fiber $\tilde U_{r,\bz_0}$ over $\bz_0$ consists of isotopy classes 
of
classical generators on $U_{\bz_0}$; it is a homogeneous space for $H_r$. Use
the isotopy class of
$\row P r$  as a designated point $u_0$ in the fiber, so each $u\in \tilde 
U_{r,\bz_0}$ is $(u_0)q$
for some $q\in H_r$. The corresponding isotopy takes a base point for classical 
generators
for $u_0$ to a basepoint for classical generators for $u$. Further, a path in
$\tilde U_{r}$ from
$u$ to
$u'$ induces a canonical isomorphism  
$\mu_{u',u}: \pi_1(U_{\bz_0},z_u)\to \pi_1(U_{\bz_0},z_{u'})$: $z_u$
(resp.~$z_{u'}$) is the basepoint for classical
generators for $u$ (resp.~$u'$). 

\begin{prop} \label{HrtoMr}
Equivalence two points $u,u'\in \tilde U_{r,\bz_0}$ if
some  representing isomorphism $\mu_{u',u}$ (with $z_u=z_{u'}$) induces an
inner automorphism on
$\pi_1(U_{\bz_0},z_u)$. This equivalence gives the fiber over $\bz_0$ for an
unramified cover $V_r \to U_r$ corresponding to $M_r$. \end{prop} 

Prop.~\ref{HrtoMr} interprets $V_r$ as isotopy classes
of (orientation preserving) diffeomorphisms of $U_{\bz_0}$. {\sl Classes\/} here 
means
to mod out by diffeomorphisms deforming (through diffeomorphisms on $U_{\bz_0}$) 
to the
identity map. Then, 
$H_r$ acts through equivalence classes of orientation preserving
diffeomorphisms on
$U_{\bz_0}$. This gives a natural map from it to $M_r$: $V_r\to U_r$ has group 
$M_r$. 

\subsubsection{Adding $\PSL_2(\bC)$ action}
Equivalence $U_{\bz_0}$
and its image under $\PSL_2(\bC)$. Elements of $\PSL_2(\bC)$ that fix $\bz_0$ 
may
equivalence some isotopy classes of diffeomorphisms not previously
equivalenced in $V_r$. Map these new  equivalence classes to
$U_r/\SL_2(\bC)=J_r$ (\S\ref{SL2-act}) giving a ramified cover with group  
$\bar M_r$. When $r=4$: $J_4=\prP^1_j\setminus
\{\infty\}$ has trivial fundamental group while $M_4$ is far from trivial. There 
is a
precise Teichm\"uller space $\text{Teich}(U_{\bz_0})$ for this situation. It is 
isotopy classes of orientation preserving diffeomorphisms $\phi:U_{\bz_0}\to 
U_\bz$
with
$\bz\in U_r$. Here $\phi$ is equivalent to $\phi': U_{\bz_0}\to U_{\bz'}$ if
$\phi'\circ\phi^{-1}$ is isotopic to an analytic map $U_\bz\to U_{\bz'}$ (from
$\PSL_2(\bC)$). 

Then $J_r$, the quotient of this Teichm\"uller space (a ball by a famous
theorem) by $M_r$, is a moduli space. For  all $r\ge 4$, we interpret it as the 
{\sl moduli of $r$ branch
point covers\/} of
$\prP^1_z$ modulo
$\PSL_2(\bC)$ (see the proof of Prop.~\ref{j-Line}). For $r=2g+2$ it would
typically be considered the moduli of hyperelliptic curves of genus $g$.

The following proposition geometrically interprets the group
$\bar M_4$ appearing in Thm.~\ref{presH4}. As elsewhere (see \S\ref{covjvalues} 
and
Lem.~\ref{locj}), change the
$j$ variable linearly so 0 and 1 are the locus of ramification of the degree six 
cover
$\prP^1_\lambda\to \prP^1_j$. \cite[Prop.~6.5]{DFrIntSpec} gave an approximation
to Prop.~\ref{isotopyM4} (as in Prop.~\ref{j-covers}) sufficient for
that paper. Its proof comprises \S\ref{provePSL2Action}. 

\begin{prop} \label{isotopyM4} 
When
$r=4$, Thm.~\ref{presH4} shows $(q_1q_3^{-1})^2$ is the center of $H_4$. Then,
$H_4/\lrang{(q_1q_3^{-1})^2}=M_4$. 
Two points $v_1,v_2\in V_r$, respectively over $\bz_{1}$ and $\bz_{2}$, induce
an isomorphism $\mu_{v_1,v_2}: \pi_1(U_{\bz_2})\to \pi_1(U_{\bz_1})$.  
Equivalence $v_1$ and
$v_2$ if some
$\alpha\in \PSL_2(\bC)$ taking $\bz_1$ to
$\bz_2$ induces $\alpha^*:  \pi_1(U_{\bz_1})\to \pi_1(U_{\bz_2})$ inverse to
$\mu_{v_1,v_2}$. These equivalence classes form a space $V^\rd_4$, 
inducing
$\psi^\rd_4: V^\rd_4\to \prP^1_j\setminus\{\infty\}$ an
analytic (ramified) cover with automorphism group $\bar M_4=H_4/\sQ$
(Thm.~\ref{presH4}). The group $\sQ/\lrang{(q_1q_3^{-1})^2}=\sQ''$ is a Klein
4-group. 

All points in the cover $\psi^\rd_4$ over $j=0$ ramify of order 3. All points 
over
$j=1$ ramify of order 2. Any Hurwitz space cover 
$\sH\to U_4$ induces a (ramified) cover $\sH^\rd \to 
\prP^1_j\setminus\{\infty\}$ through a
permutation representation of $\bar M_4$.  

Assume $\sH\to U_4$ and $\bz\in U_4$ have definition field $K$ and   
$\bz\mapsto j_\bz\in \prP^1_j\setminus \{0,1,\infty\}$. Suppose $\sQ''$ acts
trivially on a Nielsen class. Then, each geometric point of the fiber
$\sH_\bz$ on $\sH$ over $\bz$ goes one-one to the fiber
$\sH_{j_\bz}^\rd$ with both points having exactly the same fields of
definition over $K$. 
\end{prop} 

\subsection{Proof of Prop.~\ref{isotopyM4}} \label{provePSL2Action} Relate the
homotopy classes of classical generators for $u'$ to those for $u$ through
applying an element $Q_{u',u}\in H_r$ (\cite[\S4]{FrHFGG} or
\cite[\S10.1.7]{VB}). The action is that given by
\eqref{brAct}. The map $\mu_{u',u}$
explicitly makes this identification with an element of $H_r$. 

From
Lem.~\ref{conjComm}, if $Q_{u',u}$ acts as an inner
automorphism, then it must be in the center of $H_r$. Apply Thm.~\ref{presH4}: 
$(q_1q_3^{-1})^2$ is in the center of $H_4$ and acts as an inner automorphism of
$F_4$. Further, 
since $(q_1q_3^{-1})^2$ generates the whole center of $H_4$,  
$H_4/\lrang{(q_1q_3^{-1})^2}=M_4$. 

\subsubsection{How $\PSL_2(\bC)$ acts on $\pi_1(U_\bz)$} \label{PSL2comp} Now 
consider
the equivalence from composing equivalence classes defining
$V_r$ with  equivalence of maps on fundamental groups of $\pi_1(U_\bz)$ from 
elements of
$\PSL_2(\bC)$. This lies over the equivalence classes for the action of
$\PSL_2(\bC)$ on $U_r$. The $\PSL_2(\bC)$ equivalence
class of $\bz$ has a representative of form 
$(0,1,\infty,\lambda)$. For $\lambda$ not lying over $j=0$ or 1, we show the 
elements $H_\lambda$
fixing the set $\bz$ form a Klein 4-group. These computations undoubtedly occur 
in the
literature, so we give only an outline. 

First:  $H_\lambda$ contains a Klein 4-group.
Example: The $\alpha\in \PSL_2(\bC)$ switching 0 and 1, and switching $\infty$
and
$\lambda$ is $\frac{-z+1}{-z/\lambda +1}$. 

Second: To see that no $\alpha\in \PSL_2(\bC)$ acts like a 4-cycle on the 
support of
$\bz$ unless
$\lambda$ is special, assume with no loss, $$\alpha(0)=1,\ \alpha(1)=\infty,\
\alpha(\infty)=\lambda,\
\alpha(\lambda)=0.$$ Then, $\alpha(z)=\frac{\lambda z -1}{z-1}$. From the last
condition
$\lambda^2=1$ or $\lambda=-1$. 

Third: To see that no $\alpha\in \PSL_2(\bC)$ acts like a 3-cycle on the support 
of
$\bz$ unless
$\lambda$ is special, assume $\alpha$ is one of the permutations of
$\{0,1,\infty\}$ and $\lambda$ is an outside fixed point of $\alpha$. Example:
$\alpha(z)= \frac{z-1}{z}$ also fixes the two roots of $z^2-z+1$, which are 
sixth roots of
1. 

These computations show for each $j\in \prP^1_j\setminus \{0,1,\infty\}$, the 
fiber of
equivalence classes is the same as the fiber of $V^\rd_4$ over $\bz$ modulo the 
faithful
action of a Klein 4-group. Over $j=0$ the fiber is the same as the fiber of 
$V^\rd_4$ over
$\bz$ modulo the faithful action of a group isomorphic to $A_4$, and over $j=1$ 
it is the
fiber over $\bz$ modulo the action of a group isomorphic to $D_4$, the dihedral 
group of
order 8.  

\subsubsection{Identifying the Klein 4-group in $M_4$} \label{identK4} Complete 
the identification
of the group of $V^\rd_4\to \prP^1_j\setminus
\{\infty\}$ with 
$\bar M_4$ from Thm.~\ref{presH4}. This only requires 
$\sQ/\lrang{(q_1q_3^{-1})^2}=\sQ''$ to induce on $\pi_1(U_\bz)$ the Klein 4-
group in
$\PSL_2(\bC)$ that sits in $H_\bz$ above: Those switching the support of $\bz$ 
in pairs.

This could be a case-by-case identification. Since, however, Thm.~\ref{presH4} 
shows
$\sQ$ is the minimal normal subgroup of
$H_4$ containing $q_1q_3^{-1}$, it suffices to achieve $q_1q_3^{-1}$ through an 
element of
$\PSL_2(\bC)$. It is convenient to take 
$\bz=\{0,\infty,+1,-1\}$, though this $H_\bz$ is larger than a Klein 4-group.  

Crucial to the proof of \cite[Prop.~6.5]{DFrIntSpec}  is
$\mu: z\mapsto -1/z\in 
\PSL_2(\bC)$ and its effect on covers $\phi:X\to \prP^1_z$ in a Nielsen class 
branched 
over $\bz=\{0,\infty,+1,-1\}$. Under any reduced equivalence of covers in this
paper,
$\phi$ and 
$\mu\circ \phi$ are equivalent. The argument uses a set of paths 
$\lambda_0,\lambda_\infty,\lambda_{+1},\lambda_{-1}$ based at $i$. These give 
classical generators (\S \ref{classGens}) of $\pi_1(U_\bz,i)$ with  product 
$\lambda_0 
\lambda_\infty \lambda_{+1}\lambda_{-1}$ homotopic to 1. 

The effect of $\mu^{-1}$  (up to homotopy) is to switch $\lambda_0$ and 
$\lambda_\infty$, and to map $\lambda_{+1}$ to a path $\lambda'$ for which 
$\lambda_0\lambda'$ is homotopic to $\lambda_{-1}\lambda_0$. Conjugating the 
resulting collection of paths by $\lambda_0$ gives the effect of $\mu^{-1}$ as 
follows. It 
takes a representative $\bg$ of the Nielsen class of $\phi$ to $(\bg)q_1q_3^{-
1}$ 
conjugated by $g_1$. This fills in \cite[Prop.~6.5]{DFrIntSpec}  and
corrects a typo  giving $\mu(z)$ as $1/z$ (which doesn't fix $i$).

\subsubsection{The mapping from $V^\rd_r$ to reduced Hurwitz spaces} The points 
of a Hurwitz
space
$\sH$ over
$\bz\in U_r$ correspond to equivalence classes of homomorphisms from 
$\pi_1(U_\bz) \to G$
where $G$ is a finite group attached to the Hurwitz space. While $\sH$ may have 
several
connected components, easily reduce to constructing a map from $V^\rd_r$ to any 
one of them. Start
with a base point
$v_0\in V^\rd_r$ over
$\bz$ and a base point $\bp\in
\sH$ over
$\bz$. Then, relative to the classical generators given by $v_0$, $\bp$ 
corresponds to a
specific homomorphism $\psi: \pi_1(U_\bz) \to G$ up to inner automorphism or 
some stronger
equivalence (\S\ref{defNielClass}).  Any other point $v'\in V_r$ lying over 
$\bz'$ comes
with an isomorphism from
$\pi_1(U_\bz) \to \pi_1(U_{\bz'})$. 

This isomorphism takes classical
generators of $\pi_1(U_\bz)$ to classical generators of $ \pi_1(U_{\bz'})$.
Relative to this new set of generators, form the same homomorphism into $G$. 
Interpret as a canonically given homomorphism $ \pi_1(U_{\bz'})\to G$. This 
homomorphism
determines the image of $v$ in $\sH$ and an unramified cover $V_r\to \sH$. This 
map respects
equivalence classes modulo $\PSL_2(\bC)$ action. So it produces $V^\rd_r\to 
\sH^\rd$. 

\subsubsection{$\sQ''$ action when $\sH^\rd$ has genus 0}
\label{applyA5} For $r\ge 5$, the analog of $\sQ''$ is trivial (\S\ref{rge5}(. 
Then, if a
component $\sH_*^\rd $ of $\sH^\rd$ has a dense set of rational points, it is 
automatically true
for $\sH_*$. 

\begin{rem}[If $\sQ''$ acts faithfully: I] \label{nontrivq''1} Assume $r=4$, 
$\sQ''$ acts faithfully and $K=\bR$. If $\sH^\inn$ or $\sH^\abs$ is a fine 
moduli space
(over $\bR$), then its $\bR$ points correspond to equivalence classes of covers
represented by an actual cover over $\bR$. Prop.~\ref{compTest}, efficiently 
tests for $\bR$ points
on $\sH^\rd$ being {\sl cover\/} or {\sl Brauer\/} points 
(Cor.~\ref{reducedRpts}). 
\end{rem}

\begin{rem}[If $\sQ''$ acts faithfully: II] \label{nontrivq''2} Suppose $r=4$ 
and a component
$\sH_*^\rd$ of
$\sH^\rd$ (image of a component $\sH_*$ of the Hurwitz space $\sH$) has genus 0 
and dense $K$
points. If any conjugacy class in
$\bfC$ is
$K$ rational and distinct from the others, then Lem.~\ref{redCocycle} shows 
$\sH_*(K)$ is also
dense in $\sH_*$. This, however, leaves many cases: all four conjugacy classes 
are distinct,
none
$K$ rational; or each conjugacy class appears at least twice in $\bfC$. Further, 
$\sQ''$ can act neither faithfully, nor trivially, on Nielsen classes through 
$\bZ/2$ (see
Ex.~\ref{A4C34-wstory}). As in Lem.~\ref{redCocycle}, it is useful in all cases 
to distinguish
between
$K$-cover and
$K$-Brauer points on $\sH^\rd$.  
\end{rem}

For the \IGP, $r=4$ is special for the possibility of direct computations, 
especially
when the simplest case of Braid rigidity (for $r=3$) does not apply.  So, 
Prop.~\ref{isotopyM4} is valuable for applications to it. 
Here is an example long in the literature. 

\begin{exmp}[$\Spin_5$ realizations] \label{course-finemoduli} The Hurwitz space 
$\bar \sH(A_5,\bfC_{3^4})^{\inn,\rd}$ (over $\bQ$)  has genus 0 and a rational 
point
\cite[\S3]{Fr3-4BP}. From Ex.~\ref{shA5k=0},  
$\sQ''$ acts trivially. Prop.~\ref{isotopyM4}  
shows $\sH(A_5, \bfC_{3^4})^\inn$ has a dense set of $\bQ$ points. 
Lem.~\ref{A5absR} and Lem.~\ref{A5innR} illustrate this phenomena. 
So, even though the reduced Hurwitz space is not a fine
moduli space, its $K$ points are still realized by $K$-covers. This also
happens when the reduced spaces are modular curves
(as in the {\sl moduli dilemma\/} of \S\ref{compMC}). 

It is more subtle to test when $\bp\in
\sH(A_5,\bfC_{3^4})^{\inn,\rd}(\bQ)$ gives a $\Spin_5$ realization. Any $\bp$  
produces a cover
sequence
$Y_\bp\to X_\bp\to \prP^1_z$ with $Y=Y_\bp\to \prP^1_z$ geometrically Galois 
with group $\Spin_5$
having the following two properties. 
\begin{edesc} \item $X=X_\bp\to
\prP^1_z$ is an $A_5$ cover over $\bQ$. 
\item $Y\to X$ is unramified: from the unique $H_4$ orbit (\S\ref{HMrepRubric}) 
on
$\ni(A_5,\bfC_{3^4})$ having lifting invariant 1 (it contains an H-M rep.~or 
from
Prop.~\ref{serLift}). 
\end{edesc} If
$Y\to X$ is over
$\bQ$, this gives a regular
$\Spin_5$ realization. \cite[100]{Se-GT} summarizes a sufficient condition
for this when the branch cycles have odd order (as here). It is that $X$ has 
a $\bQ$ point. 

See this directly by applying Lem.~\ref{charblem} to $Y\to X$. The resulting
cover $Y'\to X$ has degree 2, is therefore Galois, and $Y'\to \prP^1_z$ produces 
the
$\Spin_5$ realization. Yet, $X=X_\bp$ won't have a $\bQ$ point  for all
$\bQ$ points $\bp$. This fails  even if we replace $\bQ$ by $\bR$. Those $A_5$
covers with a nontrivial complex conjugation operator don't achieve a $\Spin_5$ 
realization;
those with a trivial conjugation operator do. See Lem.~\ref{A5absR} and 
\ref{A5innR} to see 
both happen. 

This is a special case of Prop.
\ref{HMnearHMLevel}: The same phenomenon occurs at all levels of the $A_5$ 
Modular Tower
for $p=2$.  That is, there are real points achieving regular realization of a 
nontrivial central
extension of
$G_k$ (H-M reps.), and real points that don't (near H-M reps.; see 
Prop.~\ref{nearHMbraid}).  
\end{exmp}

\subsubsection{Isotopy classes of diffeomorphisms} An elliptic curve over $K$ 
(genus 1 
with a $K$ point) is a degree two cover of the sphere over $K$, ramified at four 
distinct
unordered points
$\bz$. The equivalence class of $\bz$ modulo the action of 
$\PSL_2(\bC)$ determines the isomorphism class of the elliptic curve. Isotopy 
classes
of (orientation preserving) diffeomorphisms of a complex torus $X$ identifies 
with 
$\SL_2(\bZ)$, appearing from
its action on $H^1(X,\bZ)$. 

Isotopy classes of diffeomorphisms of $U_{\bz_0}$
(as elements in $\bar M_4=\PSL_2(\bZ)$: Prop.~\ref{isotopyM4}) relate
to diffeomorphisms of the torus attached to $\bz_0$. The difference? The complex 
torus
has a chosen origin (a canonical rational point) preserved by diffeomorphisms; 
there is no chosen
point in
$\bz$. This appears in mapping a space of {\sl involution realizations of 
dihedral
groups\/} to a modular curve (\cite{FrGGCM}; the most elementary example from
\S\ref{fineModDil}).

\subsection{Statement of the main presentation theorem} 

\subsubsection{Presenting $H_3$} Basics about $H_3$ expressed in  $q_1$ and 
$q_2$:  $q_1q_2=\tau_1$ and $q_1q_2q_1=\tau_2$ generate $H_3$. From   
\eqref{hwgp}, $q_1^2=q_2^{-2}$. Further:
$$ 
\tau_1^3=q_1q_2(q_1q_2q_1)q_2=q_1q_2(q_2q_1q_2)q_2=q_2^2=q_1^2=\tau_2^2.
$$ This identifies
$H_3$.

\begin{defn} \label{dicyclic} Characterize the {\sl dicyclic (or quaternion)\/}  
group
$Q_{4n}$ of order 
$4n$ 
as having generators 
$\tau_1,\tau_2$ with $\ord(\tau_1)=2n$, $\ord(\tau_2)=4$, $\tau_2^{- 
1}\tau^\sph_1\tau^\sph_2=\tau_1^{-1}$ and  $\tau_2^2=\tau_1^n$.  
\cite[p.~43]{StAlg} has a
memorable matrix representation of $Q_{4n}$ in $\GL_4(\bR)$: $\tau_1=\smatrix 
{R_{2\pi/n}} {0_2}
{0_2} {R_{2\pi/n}}$ and $\tau_2=\smatrix {0_2} \mu {-\mu} {0_2}$ with $R_\theta$ 
the rotation in
2-space through the angle $\theta$, $0_2$  the $2\times 2$ zero matrix and
$\mu=\smatrix{1}{0}{0}{-1}$.  For $n=3$: 
\begin{edesc} \label{ordH3} \item $\ord(\tau_1)=6$ and $\ord(\tau_2)=4$; and  
\item $H_3$
is the  dicyclic group of order 12.\end{edesc} \end{defn}

\begin{lem}\label{quat4n}
For each integer $k$,  $\tau_2\tau_1^k$ in $Q_{4n}$ has order 4: 
$\tau_2\tau_1^k\tau_2\tau_1^k=\tau_2^2$.
For $n=2^t$, every proper subgroup of $Q_{4n}$ contains $\tau_2^2$, and 
$Q_{4n}/\lrang{\tau_2^2}$
is  the 
dihedral
group $D_n$ of order
$2n$. If $u|n$, then $Q_{4u}=\lrang{\tau_1^{\frac{n}{u}},\tau_2}\le Q_{4n}$.  
\end{lem}

\subsubsection{Notation for computations in $H_4$} From this point $r=4$. We 
give a
convenient  presentation of $H_4$ and various subgroups and quotients. 
Computations in both
$H_4$ and in
$M_4$ use the letters
$q_1,q_2,q_3$  for the images of $Q_1,Q_2,Q_3$ in either. Distinguishing between 
their
images in these  two groups requires care about the context. Here is how we 
intend
to do that. 

When $r=4$, the extra relations for $M_4$ (beyond those for $B_4$) are 
$$\tau_1=(Q_3Q_2)^3=1,\ \tau_2=Q_1^{-2}Q_3^2=1,\ \tau_3=(Q_2Q_1)^{-3}=1$$ and 
$\tau=(Q_3Q_2Q_1)^4=1$. Lemma \ref{B4F4} shows adding these relations to $B_4$ 
is 
equivalent to adding $q_1^2q_3^{-2}=1$ to $H_4$. This produces new equations: 
\begin{equation} \label{basM4} 
q_1q_2q_1^2q_2q_1=(q_1q_2q_1)^2=(q_1q_2)^3=1.\end{equation} Then, 
$\gamma_0=q_1q_2$, $\gamma_1=q_1q_2q_1$ and $\gamma_\infty=q_2$  satisfy 
\begin{equation} \label{prodOneM4} 
\gamma_0\gamma_1\gamma_\infty=1.\end{equation} 

The relation $q_1q_3^{-1}=1$ is not automatic from \eqref{basM4}. Crucial, 
however, is 
how $M_4$ acts on {\sl reduced Nielsen classes\/} (\S\ref{j-covers}). This 
action {\sl 
does\/} factor through the relation $q_1q_3^{-1}=1$ (\S\ref{SL2-act} and 
Prop.~\ref{PSL2-quotient}). Therefore it factors through the induced quotient 
$H_4/\sQ=\bar M_4$ of Thm.~\ref{presH4}, from $H_4$ acting on Nielsen classes 
(Prop.~\ref{j-Line}).  The upper half plane appears as a classical ramified 
Galois cover of 
the $j$-line minus $\infty$.  The elements $\gamma_0$ and $\gamma_1$ in $\bar 
M_4$ generate the
local monodromy of this  cover around 0 and 1 (\S \ref{j-covers}).  

Sometimes $q_2\in H_4$ acts significantly different when viewed in $\bar M_4$. 
On these occasions
denote it $\gamma_\infty$: representing the local  monodromy action 
corresponding to the cusp
$\infty$ of  the $j$-line (\S\ref{SL2-act}).  Similarly, denote $q_1q_2q_3$ as 
$\sh$: {\sl the
shift\/}. From the above,  $\sh$ and 
$\gamma_1$ are the same in $\bar M_4$ (\S \ref{shiftOp}).   

\begin{thm}\label{presH4} Let $\sQ =\lrang{(q_1q_2q_3)^2, q_1^\sph q_3^{-1}}$.  
The 
following hold. \begin{edesc} \label{h4det} \item\label{h4deta} $H_4$ has one 
nontrivial involution $z=(q_1q_3^{-1})^2$ generating its center.  
\item\label{h4detb} 
$\sQ\triangleleft H_4$ and $\sQ$ is the quaternion  group $Q_8$ 
(Lem.~\ref{quat4n}). 
\item\label{h4detc} 
$H_4/\sQ=\bar M_4\cong \PSL_2(\bZ)$.  \item\label{h4detd} Exactly two conjugacy 
classes of $H_4$  subgroups $U_1=\lrang{q_1,q_2}$ and 
$U_2=\lrang{q_2,q_3}$ (both containing $\lrang{z}$) are isomorphic 
to $\SL_2(\bZ)$. 
\end{edesc} From \eql{h4det}{h4detb}, 
$\sQ$ is the smallest normal subgroup of $H_4$ containing either $(q_1q_2q_3)^2$ 
or 
$q_1^\sph q_3^{-1}$. So, from \eqref{basM4}, $M_4=H_4/\lrang{z}$ and 
$\bar M_4=\lrang{\gamma_0,\gamma_1}=\lrang{\gamma_1,\gamma_\infty}$ 
is $M_4$ modulo the relation $q_1=q_3$. \end{thm}

\subsection{Start of Theorem \ref{presH4} proof} Use \S\ref{braidGp} 
generators 
$\row {\bar \sigma} 4$ of $F_4$, etc. Since $r=4$, use $A$ for 
$\bar A_4$. As above, reserve the  symbols $Q_1,Q_2,Q_3$ for the generators of 
the $B_4$, 
and $q_1,q_2,q_3$  for their images in $H_4$. For convenience,  arrange   
the 
$Q\,$s  action in a table. The $i$th column has $Q_i$ acting on the 4-tuple of $ 
\bars\,$s: 
$$\begin{matrix}{}&{Q_1}&{Q_2}&{Q_3}\\ {\bars_1}&{\bars^\sph_1 \bars^\sph_2 
\bars_1^{-1}}&{\bars_1}&{\bars_1}\cr {\bars_2}&{\bars_1}&{\bars^\sph_2 
\bars^\sph_3 \bars_2^{-1}}&{\bars_2}\cr 
{\bars_3}&{\bars_3}&{\bars_2}&{\bars^\sph_3 \bars^\sph_4 \bars_3^{-1}}\cr 
{\bars_4}&{\bars_4}&{\bars_4}&{\bars_3}.\end{matrix}$$

Let $i:F_4\to A$ map $ \bars\in F$ to the inner automorphism it induces. 
Initially we work 
with the image of $B_4$ in $A$ where $i(F_4)$ and $B_4$ commute. Denote 
$Q_1Q_2Q_{3}^2Q_2Q_1$ by $D$: it maps $\bbs$ to  \begin{equation} 
\begin{array}{rl} \label{Daction} &(\bars^\sph_1 \bars^\sph_2 \bars_1^{-1}, 
\bars^\sph_1 \bars^\sph_3 \bars_1^{-1}, \bars^\sph_1 \bars^\sph_4 \bars_1^{-1}, 
\bars_1)Q_{3}Q_2Q_1= \\ &(\bars^\sph_ 1\bars^\sph_2 \bars^\sph_3 
\bars^\sph_4\bars_1^\sph\bars^{-1}_4 \bars^{-1}_3 \bars^{- 1}_2 \bars_1^{-1}, 
\bars^\sph_1 \bars^\sph_2 \bars_1^{-1}, \bars^\sph_1 \bars^\sph_3 \bars_1^{-1}, 
\bars^\sph_1 \bars^\sph_4 \bars_1^{-1}).\end{array}\end{equation}

Further Notation: Let $G_4$ be the  quotient of  $F_4$ by the relation  
$$\Pi(\bbs)=\bars_1\bars_2\bars_3\bars_4=1.$$ While $G_4$ is a free group on 3  
generators, $H_4$ computations require $\bars_1,\bars_2,\bars_3,\bars_4$ to 
appear 
symmetrically. So, we use this quotient presentation. Let $\Aut(G_4)$ be the 
automorphisms of $G_4$. Then $\Inn(G_4)$ is the normal subgroup of inner
automorphisms  of $G_4$. Recall: The image of $D$, and all its  conjugates, in
$H_4$ is trivial. 

\begin{lem} \label{B4F4} Both $(Q_1Q_2Q_3)^4$ and $(Q_3Q_2Q_1)^4$ map $\bbs$ 
to $$\bars_1^\sph \bars_2^\sph \bars_3^\sph   \bars_4^\sph \bbs \bars_4^{-1} 
\bars_3^{- 
1} \bars_2^{-1}  \bars_1^{- 1}=(\bbs)i(\bars_4^{-1} \bars_3^{-1} \bars_2^{- 1} 
\bars_1^{-1}).$$ Action of $B_4$ on $F_4$ induces $\alpha_4: B_4\to \Aut(G_4)$ 
with 
the center $\lrang{(Q_1Q_2Q_3)^4}$ of $B_4$ generating the kernel. This  induces 
$\mu_4: H_4\to  \Aut(G_4)/\Inn(G_4)$. The 
kernel of $B_4\to M_4$ is the direct product of the free group 
$K_4^*=\lrang{(Q_3Q_2)^3,Q_1^{-2}Q_3^2, (Q_2Q_1)^{-3}}$ and 
$\lrang{(Q_1Q_2Q_3)^4}$.  Above, $Z=i(\bars_1^\sph \bars_2^\sph \bars_3^\sph 
\bars_4^\sph)$ is identical to an element of $B_4$. Consider the image
$z$  of
$Z$ in
$H_4$: 
$$(q_1q_2q_3)^4z=1=(q_3q_2q_1)^4z=q_1^2q_3^{-2}z\text{ and } z^2=1.$$ In 
particular, combining this with \eqref{basM4} identifies $M_4$ with the image of 
$H_4$ 
in $\Aut(G_4)/\Inn(G_4)$. So, the image of $z$ in $M_4$ is 1. \end{lem}

\begin{proof} The map $i$ induces $i: G_4\to \Inn(G_4)$.  Action of $B_4$ 
preserves 
$\Pi(\bbs)$. Thus, it induces a homomorphism of $B_4$ into $\Aut(G_4)$
where $D$ goes to the automorphism $i(\bars_1^{-1})$, $\bbs\mapsto 
\bars_1^\sph\bbs\bars_1^{-1} $ (as in \eqref{Daction}). Conclude:
Modulo inner automorphisms of
$G_4$, $D$  acts trivially, producing the desired homomorphism 
$$H_4\to  \Aut(G_4)/\Inn(G_4).$$

Now consider the explicit formulas. Most of the calculation is in 
\eqref{Daction}: 
$Q_1Q_2Q_3$  cycles entries of $\bbs$ back 1, and conjugates all entries by the 
first 
entry's inverse. So  $(Q_1Q_2Q_3)^4$ leaves entries of $\bbs$ untouched except 
for 
conjugating 
them by $$(\bars_1^\sph \bars_2\sph \bars_3\sph \bars_4\sph \bars_3^{-1} 
\bars_2^{-1}  
\bars_1^{-1}  \bars_1\sph \bars_2\sph \bars_3\sph \bars_2^{-1} \bars_1^{-1}  
\bars_1\sph \bars_2\sph  \bars_1^{-1} \bars_1)^{-1}=(\bars_1^\sph \bars_2\sph 
\bars_3\sph 
\bars_4\sph)^{-1}.$$ Also, $Q_3Q_2Q_1$  cycles the entries of $\bbs$ forward 1. 
The 
new first 
entry is the old 4th entry conjugated by the inverse of the  product of the old 
first three 
entries. Thus, $(Q_1Q_2Q_3)^4$ and $(Q_3Q_2Q_1)^4$ act the same.  Add that $D$ 
maps to 1 in $H_4$ to see $$1=(q_1q_2q_3)^4(q_3q_2q_1)^4=(q_1q_2q_3)^8=z^{- 
2}.$$

Let $Q'=Q_1Q_2Q_3$. Extending the calculation above gives the next list: 
\begin{eqnarray*} &1= (Q')^{- 1}D Q'(Q_1Q_2)^3i( \bars_1 \bars_2 \bars_3 
\bars_4);&\\ &1= (Q'Q_1Q_3^{-1})^2(Q')^{- 1}D Q'i( \bars_1 \bars_2 \bars_3 
\bars_4);&\\ &1=(Q')^4i( \bars_1 \bars_2 \bars_3 \bars_4); \text{and}& \\ 
&1=(Q_1Q_3^{- 1})^2(Q')^{-1}D(Q_1Q_2Q_3)(Q_1Q_2Q_3)^{- 
2}D(Q_1Q_2Q_3)^{2}i( \bars_1 \bars_2 \bars_3 \bars_4).&\end{eqnarray*} Add the 
relation $Q_1Q_3=Q_3Q_1$ to deduce, in order,  these  relations in $H_4$: 
\begin{equation}\begin{matrix}\label{formz}\text{a)}&(q_1q_2)^3z=1; &\text{b)}& 
(q_1q_2q_1)^2z=1;&\text{c)}&(q_1q_2q_3)^4z=1;\\ \text{d)}& 
(q_1q_3^{-1})^2z=1; &\text{e)}& q_1q_3=q_3q_1.&{}&\end{matrix}\end{equation} 
So, the image of $z$ in $M_4$ is 1. 

The kernel of $\alpha_4$ contains elements of $B_4$ 
inducing inner automorphisms commuting with conjugation by $\bars_1^\sph 
\bars_2\sph 
\bars_3\sph \bars_4=c$. Since $c$ generates conjugations commuting with $c$, 
$\lrang{(Q_1Q_2Q_3)^4}$ generates the kernel of $\alpha_4$.

Generators of $K_4^*$ act on $G_4$: Respectively, $(Q_3Q_2)^3,Q_1^{- 
2}Q_3^2,(Q_2Q_1)^{-3}$ induce conjugation by $g_4, g_1g_2,g_1$. These 
conjugations 
on $G_4$ form a free group. So $K_4^*$ is a free group on these generators. 
\end{proof}

The remainder of the proof of Thm.~\ref{presH4} is in \S\ref{quatProp}. 

\subsection{Properties of $\sQ$; combinatorics of \eqref{h4det}} 
\label{quatProp}
First we show
$\sQ$ is a normal subgroup isomorphic to 
$Q_8$. Then we list  
combinatorics contributing  to \eqref{h4det}.  

\subsubsection{$z$ is the involution in $Q_{16}$} From
(\ref{formz}d),
$z$ is a word in 
$q_1$ and $q_3$. So, (\ref{formz}e) shows  $z$ commutes with $q_1$ and $q_3$. 
From $$q_2(q_1q_2)^3=(q_1q_2q_1)(q_1q_2)^2=(q_1q_2)^3q_2,$$ $q_2$
commutes  with $(q_1q_2)^3$.  Apply (\ref{formz}a): $z$  commutes with
$q_i$, $i=1,2,3$. Thus, Lemma \ref{B4F4} shows   $z$ is a central involution. 
\S\ref{znezero} shows
$z\ne1$ and illustrates its significance. \S\ref{invH4} reformulates property
\eql{h4det}{h4deta}: 
\begin{triv}\label{zinv} $z$ is the only involution of $H_4$.\end{triv}

Consider $\sQ'=\lrang{q_1^\sph 
q_2^\sph q_3^\sph , q_1^\sph q_3^{-1}}$. Set $\alpha=q_1q_2q_3$, $\beta=q_1q_2$ 
and $\gamma=q_1q_3^{-1}$. These generators  simplify presenting $H_4$.

\begin{lem} \label{sQ'lem} Rewrite\/ {\rm (\ref{formz}a--d)} using $\alpha$, 
$\beta$ and
$\gamma$  to  give \begin{equation} \label{formzp} \text{\rm\ a) 
}(\beta)^3=z;\text{\rm\
b)  }(\alpha\gamma)^2=z; \text{\rm\ c) }(\alpha)^4=z; \text{\rm\ and \ d)\ }
(\gamma)^2=z. 
\end{equation} Further, $H_4=\lrang{\alpha,\beta,\gamma}$.  Since $z$ is a 
central 
involution, {\rm (\ref{formzp}b--d)} show $\sQ'=\lrang{\alpha,\gamma}$ is 
$Q_{16}$ 
and $\sQ' \bmod\ z$ is the dihedral group of order 8 (Lem.~\ref{quat4n}).  
\end{lem}

\subsubsection{$\sQ=\lrang{\alpha^2,\gamma}\norm
H_4$ and $\Cusp_4$} 
\label{sQpart2}  From relations for  $B_4$,
$$\alpha^2\gamma=q_1q_2q_3q_1q_2q_1= 
q_1q_2q_3q_2q_1q_2=q_1q_3q_2q_3q_1q_2=q_3\alpha^2q_3^{-1}$$ and 
$q_3^\sph(\gamma)q_3^{-1}=\gamma$. Thus $q_3$ normalizes $\sQ$. Since 
$q_1^\sph q_3^{-1}\in \sQ$,  $q_1$ also  normalizes $\sQ$. We now show $q_2$ 
does 
also. 
Apply (\ref{formz}a) in the form $q_2^{-1}q_1^{-1}q_2^{-1}=zq_1q_2q_1$ to get
$$q_2^{-1}\gamma^{-1} q_2^\sph
=q_2^{-1}q_1^{-1}q_2^{-1}q_2q_3q_2=z(q_1q_2q_1)q_3q_2q_3=z\alpha^2.$$
Also, since $z$ is a central involution, conjugate by $q_2^{-1}$
to get $q_2^\sph\alpha^2q_2^{-1}=z\gamma^{-1}$.  

\begin{defn}[Cusp group] With $\sQ''=\sQ/\lrang{z}$, denote the subgroup $\lrang{\sQ'',q_2}$, the
{\sl cusp  group\/} in $M_4$, 
by $\Cusp_4$. It is
$\lrang{(q_1q_2q_3)^2}/\lrang{z}\times \lrang{q_1q_3^{-1}}/\lrang{z}\xs
\lrang{q_2}$ with $q_2$ switching the two factors on the copy of the Klein 
$K_4$.\end{defn} Nontrivial calculations with 
$\Cusp_4$  figure in computing the genus of reduced Hurwitz space components 
(see
\S\ref{sQorbits2}). 

\subsubsection{Proof of\/ {\rm \eql{h4det}{h4detc}}}
\label{order2and3} From
\eqref{basM4}, the  image of $q_1q_2$ and $q_1q_2q_1$ give elements of
respective orders 3 and 2  generating $H_4/\sQ$.  A well-known abstract
characterization of   $\PSL_2(\bZ)$ is as  the free product of elements
of order 3 and 2. Further, $\SL_2(\bZ)$ is a free amalgamation  of 
cyclic groups of order 6 and  4 along their common subgroups of order 2 
\cite{SeArbres}. Thus, $H_4/\sQ$ is isomorphic to a quotient of  $\PSL_2(\bZ)$. 
The 
proceeding arguments, however, have derived a presentation for $H_4$, and so for 
$H_4/\sQ$. Thus \eql{h4det}{h4detc} holds.

\subsubsection{Conclude properties\/ {\rm \eqref{h4det}}} 
\label{invH4}  \S \ref{sQpart2} and \S\ref{order2and3}  give
\eql{h4det}{h4detb} and
\eql{h4det}{h4detc}.  In any quotient of the group $\lrang{a,b}$ where 
$a2=b^3=1$, all
involutions are conjugate to $a$. Thus,  the only  possible  involutions in 
$H_4$ are in $\sQ$, and 
this contains just one.   From \eqref{formz} (expression c)),  $q_1q_2q_1$ has 
order 4. This gives  
\eql{h4det}{h4deta}. It remains to discuss $z$ being nontrivial (as
above) and  
\eql{h4det}{h4detd}.

From  \eql{h4det}{h4deta}, each subgroup of $H_4$  containing (a
copy of) 
$\SL_2(\bZ)$ contains  $z$. Thus, there is a one-one correspondence between 
subgroups 
$G$ of $H_4$  isomorphic to $\SL_2(\bZ)$ and subgroups of $H_4/\lrang z$ 
isomorphic to $\PSL_2(\bZ)$ via $G \to G/\lrang{z}$. So it suffices to show 
$H_4/\lrang 
z$ contains precisely two conjugacy classes of subgroups   isomorphic to 
$\PSL_2(\bZ)$. 
(These join together by the relation $q_1q_3^{-1}=1$.)

\newcommand{\Cen}{{\text{Cen}}} Expression \eqref{basM4} identifies two such 
groups, the images of $\lrang{q_1,q_2}$ and $\lrang{q_2,q_3}$.  Consider $\sQ/ 
\lrang{z}\eqdef \sQ''$, $\beta\lrang{z}\eqdef \beta''$ and  
$\alpha\gamma\lrang{z}\eqdef 
\alpha''$.  Note that $\sQ''$ is a Klein 4-group.
Denote  the group $\lrang{\beta'',\alpha''}$ by $\Gamma'_1$. It is isomorphic to
$\PSL_2(\bZ)$ and a complement to $\sQ''$ in $H_4/\lrang{z}$. Let
$C$ be the  centralizer, $\Cen_{\Gamma'_1}(\sQ'')$, of $\sQ''$ (regard all as
subgroups of  
$H_4/\lrang{z}$). From (\ref{formzp}), the quotient  $\Gamma'_1/'C$ is 
isomorphic to $S_3\cong \Aut(\sQ'')$. This  identifies $C$ with a well-known 
rank 2 free 
subgroup  of $\PSL_2(\bZ)$. Also,  $\sQ''$ is absolutely irreducible as an 
$\bF_2\bigl[\Gamma'_1\bigr]$ module.

Complements of $\sQ''$ in $H_4/\lrang{z}$ correspond to elements of  
$H^1(\Gamma'_1,\sQ'')$ \cite[p.~244]{NCHA}. Conclude \eql{h4det}{h4detd} if this 
cohomology  group has order 2. Let $B$ be the  largest elementary abelian 2-
group 
quotient of the rank 2 free  group $\Cen_{\Gamma'_1}(\sQ'')$.  Then $B$  and 
$\sQ''$ 
are 
isomorphic as $\bF_2\bigl[\Gamma'_1\bigr]$ modules. For any group $G$ with 
normal 
subgroup $H$ acting on a module $M$, there is  an exact sequence 
\cite[p.~100]{AW}: 
\begin{equation} \label{exCoh}  0\to H^1(G/H,M^H) \mapright{\text{inf}} 
H^1(G,M)\mapright{\text{res}}  H^1(H,M)^G\to H^2(G/H,M^H). \end{equation}  Apply 
this with 
$M=\sQ''$,  $G=\Gamma'_1$ and $H=\Cen_{\Gamma'_1}(\sQ'')$, so $G/H=S_3$. 
Restrict to a 2-Sylow
$P_2$ of $S_3$. As in \S\ref{splitvs2-lift}, if $H^*(P_2,M^H)$ is trivial, so is 
$H^*(G/H,M^H)$.

The action of $P_2$ on  $\sQ''$ is the regular representation on two
copies of
$P_2$, so $\sQ''$ is a free $P_2$ module and the cohomology groups $H^1(G,M)$ 
and $H^2(G,M)$ are
trivial. Conclude: 
$H^1(\Gamma'_1,\sQ'')=\Hom(B,\sQ'')^{S_3}=\bZ/2$.

\subsubsection{Comments on $q_1q_3^{-1}\ne 1$ in 
$M_4$}  The nontrivial action of $q_1q_3^{-1}$ appears in many Nielsen classes. 
Especially in the  action of $\sQ''=\sQ/\lrang{z}$ on
$\ni(G_k,\bfC_{3^4})$, $k\ge 1$:   all orbits have length four on 
$\ni({}_2^1\tilde  A_5,\bfC_{3^4})^{\inn}$ (Lem.~\ref{4HM4NHM}).
We make much, however, of the trivial action of $q_1q_3^{-1}$ when $k=0$. 

Branch cycles for the
$j$-line cover from the Nielsen class 
$\ni(A_5,\bfC_{3^4})^\abs$ and $\ni(A_5,\bfC_{3^4})^\inn$ are in
\cite[(7.7)]{Fr-Schconf}.   Conjugate by $S_5$  to assume $\bg$ representatives  
for $\ni(A_5,\bfC_{3^4})^{\abs}$ have $g_1=(1\,2\,3)$. Further conjugation by 
elements 
of  
$\lrang{g_1,(4\,5)}$ gives the $(g_2,g_3,g_4)$ entries of {\sl absolute\/}  
Nielsen class 
representatives (Table \ref{lista5}).

\begin{table}[h] \caption{List $\ni(A_5,\bfC_{3^4})^{\abs}$} \label{lista5}
\vbox{\halign{&\hfil ${}_{#}\bg$:& $\,\,((#)$,&$\,\,(#)$,&$\,\,(#\ $\cr 
1&1\,3\,2&1\,4\,5&1\,5\,4));&2&1\,4\,5&1\,5\,4&1\,3\,2)); 
&3&1\,4\,5&2\,1\,5&2\,4\,3));\cr 
4&1\,4\,5&3\,2\,1&3\,5\,4)); 
&5&1\,4\,5&4\,3\,2&4\,1\,5));&6&1\,4\,5&5\,4\,3&5\,2\,1));\cr 
7&2\,1\,4&2\,4\,5&5\,3\,2));&8&2\,1\,4&3\,2\,5&5\,4\,3)); 
&9&2\,1\,4&4\,3\,5&2\,4\,5)).\cr}} \end{table} 

Then, $\sQ$ acts trivially on everything in Table \ref{lista5}. Compute: 
\begin{equation}\label{actQs} Q_1=(2\,5\,3\,6\,4)(7\,9\,8),\quad 
Q_2=(1\,4\,9\,8\,5)(3\,6\,7) 
\text{ and } Q_3=(2\,5\,3\,6\,4)(7\,9\,8).\end{equation} Even on 
$\ni(A_5,\bfC_{3^4})^{\inn}$, $\sQ$ is trivial: $q_1q_3^{-1}$ maps 
$\bg=((1\,2\,3),(3\,2\,1),(1\,4\,5),(5\,4\,1))$ to 
$\bg^*=((3\,2\,1),(1\,2\,3),(5\,4\,1),(1\,4\,5))$, conjugate by $(2\,3)(4\,5)$ 
to $\bg$. 

\subsubsection{Comments on $(q_1q_3^{-1})^2\ne 1$} \label{znezero}
From Lem.~\ref{B4F4},
$z=1$ in the
$M_4$ quotient.  Acting on an inner Nielsen class (\S\ref{defNielClass};
conjugating by $G$), $z$ will  always be trivial.  
The following self-contained computation shows it is nontrivial in
$H_4$ (compare with \cite{GilBus} quoted in \cite{Birman}). 


Suppose $z=1$ in $H_4$. Apply this to
$$1=q_1q_2q_3^2q_2q_1=q_1q_2q_1^2q_2q_1=(q_1q_2)^3.$$ Project
$(\prP^1_z)^4\setminus \Delta_4\to (\prP^1_z)^3\setminus \Delta_3$ onto the 
first three
factors (as in \eqref{fundgpdiag}). Then, $(q_1q_2)^3$ goes to
$(q_1'q_2')^3$ with $q_1',q_2'$ the generators of $H_3$. This implies 
$(q_1'q_2')^3=1$
contrary to $\ord(q_1'q_2')=6$ in \eqref{ordH3}.  

Consider $\bg\in
\ni(G,\bfC)/N'$ with
$N'$ a group  from \S\ref{refEquiv} and its  corresponding curve cover 
$\phi_\bg:X_\bg\to
\prP^1_z$.  Denote the subgroup of $H_4$  fixing $\bg$ by $H_\bg$.  
\cite[\S3.1]{FrCCMD}
produces an explicit monodromy action  of $H_\bg$ on $H_1(X_\bg,\bZ)$. We ask if 
$z$
appears nontrivially in a structure from a Nielsen class, as in the
situation of \S\ref{finquotz}. 

\begin{prob} Since $z$ fixes each element of the inner Nielsen classes, when is 
it 
nontrivial 
on each representative cover of the corresponding Hurwitz space? What are its 
$\pm1$ 
eigenspaces on the first homology of curves in the family?\end{prob}

\subsubsection{Finite quotients of $H_4$ where $z$ does not vanish} 
\label{finquotz} Add
the relations
$q_i^4=1$, $i=1,2,3$, to $H_4$. 

\begin{prop} Let $G_4^*$ be the quotient of $H_4$ by the relations $q_i^4=1$, 
$i=1,2,3$.
Denote the quotient of $\bar M_4$ by the (image of) the relations
$q_i^4=1$, $i=1,2,3$, by $\bar M_4^*$. Then, $G_4^*$ is a (finite) 
nonsplit extension of
$S_4=\bar M_4^*$ by $\sQ$. \end{prop}

\begin{proof} The image of $a_{12}=q_1^2$,
$a_{13}=q_1^{-1}q_2^2q_1$ and $a_{14}=q_1^{-1}q_2^{-1}q_3^2q_2q_1$ in $\bar M_4$ 
identify
with classical generators of  $\prP^1_\lambda \setminus\{0,1,\infty\}$ 
(\S\ref{SL2-act}).
Each 
$a_{1j}$ has image an involution $\bar a_{1j}$ in $\bar M_4^*$. The product-one
condition,
$\bar a_{12}
\bar a_{13} \bar a_{14}=1$ implies these elements generate a Klein 4-group $K_4$ 
(Lem.~\ref{quat4n}). The group $\Aut(\prP^1_\lambda/\prP^1_j)$ is $S_3$ (in its
regular representation). To check  how $S_3$ acts on the $K_4$, note 
the element of order $3$ (image of $\gamma_0=q_1q_2$) acts nontrivially on $K_4$ 
(by
conjugation). Here is the effect (in $\bar M_4^*$) of conjugating $q_1^2$ by
$q_1q_2$:  
$$q_1q_2(q_1^2)(q_1q_2)^{-1}=(q_1q_2q_1)(q_1q_2^3q_1^{-1})=q_1^{-
1}(q_1q_2)^3(q_2^2q_1^{-1})
=q_1^{-1}q_2^2q_1^{-1}.$$ Since
$q_1^{-1}q_2^2q_1^{-1}=q_1^{-1}q_2^2q_1^{3}=a_{13}a_{12}$, the action is 
nontrivial and 
$\bar M_4^*$ is $S_4$. 

To see these relations don't kill $z$, verify the normalizing action of
$H_4$ on $\sQ$ gives  $q_i^4$ acting trivially, $i=1,2,3$. 
Use conjugation by $q_i$ as in \S\ref{sQpart2}. 

Finally, to see it is nonsplit, consider the element $q_1q_2q_3\in H_4$ whose 
image in
$\bar M_4^*$ is $q_1q_2q_1$ (image of $\gamma_1$) of order 2. A splitting of 
$S_4$ would
give a lift of this to an element of order 2, contrary to Lem.~\ref{sQ'lem} 
giving
no such element in $\sQ'$. 
\end{proof}

\subsection{The shift $\gamma_1=q_1q_2q_1$  and the middle twist 
$\gamma_\infty=q_2$} \label{shiftOp} The action of $\gamma_1=q_1q_2q_1$ (from 
Prop.~\ref{j-Line}) on $\ni(G,\bfC)^{\inn,\rd}$ is the same as $q_1q_2q_3$. It 
sends 
$(g_1,g_2,g_3,g_4)$ to $(g_2,g_3,g_4,g_1)$ modulo conjugation by $g_1$. 
Sometimes 
we want $\sh=q_1q_2q_3$, the {\sl shift}, instead of $\gamma_1$, on 
$\ni(G,\bfC)^\inn$. 
Here also, use $\sh$. 

\begin{lem} The (middle) twist $\gamma_\infty$ and $\sh$ (on 
$\ni(G,\bfC)^{\inn,\rd}$) 
entwine by \begin{equation} \label{entwine} \sh \gamma_\infty 
\sh=\gamma_\infty^{-1} 
\sh \gamma_\infty^{-1}. \end{equation}  Knowing $\sh$ is an involution and 
\eqref{entwine} holds gives a presentation of $\bar 
M_4=\lrang{\sh,\gamma_\infty}$. 
\end{lem}

\subsubsection{Twist action on products of two conjugacy classes} 
\label{twistAct} Let 
$G$ be a finite group with conjugacy classes $(\C_1',\C_2', \C_3')=\bfC'$.  
Consider: 
$$S(\bfC')=\{(g_1,g_2)\mid g_i\in \C_i,\ i=1,2,\ \text{with } g_1g_2\in 
\C_3\}.$$
Consider the {\sl twist\/}
$\gamma$ mapping $(g_1,g_2)\in S(\bfC')$ to $(g_1g_2g_1^{- 
1},g_1)\in S(\bfC'')$ with $\bfC''=(\C_2',\C_1',\C_3')$. 

\begin{lem} For $g_3\in \C_3'$ let $S(\bfC')_{g_3}$ be those
$(g_1,g_2)\in S(\bfC')$ with $g_1g_2=g_3$. Then, $\gamma$ maps one-one
from   $S(\bfC')_{g_3}$ to $S(\bfC'')_{g_3}$. Also,
$|\C_3'||S(\bfC')_{g_3}|=|S(\bfC')|$. 
\end{lem}

The structure constant formula (\S\ref{congClProd}) calculates $|S(\bfC')|$ 
using  
complex representations of $G$. Calculations in 
\S\ref{2A5Cusps} compute the length of  $\gamma$ orbits on $S(\bfC')_{g_3}$. 
For $g_1,g_2$ in a group, denote the centralizer of $\lrang{g_1,g_2}$ by 
$Z(g_1,g_2)$. 

\begin{prop} \label{gamOrbit}  Assume $g_1\ne g_2$. Let $(g_1,g_2)\in 
S(\bfC')_{g_3}$. The orbit
of 
$\gamma^2$ containing $(g_1,g_2)$ is $(g_3^j g_1 g_3^{-j},g_3^jg_{2} g_3^{-j})$, 
$j=0,\dots, \ord(g_3)\nm1$. So, the  orbit of $\gamma^2$ has length  
$$\ord(g_3)/|\lrang{g_3}\cap Z(g_1,g_2)|\eqdef  o(g_1,g_2).$$ Let 
$x=(g_1g_2)^{(o-1)/2}$ and
$y=(g_2g_1)^{(o-1)/2}$, so $g_1y=xg_1$ and $yg_2=g_2x$. Assume: 
\begin{triv} 
\label{xyswitch} $o(g_1,g_2)=o$  is odd; and $yg_2$ has order 2. \end{triv}  
\noindent Then,  the orbit of $\gamma$ has
length
$o(g_1,g_2)$. Otherwise,   the orbit of $\gamma$ containing $(g_1,g_2)$ has 
length
$2o(g_1,g_2)$.  \end{prop}

\begin{proof} For $t$ an integer, $$(g_1,g_2)\gamma^{2t}=(g_3^tg_1g_3^{-t}, 
g_3^tg_2g_3^{-t}) \text{ and } (g_1,g_2)\gamma^{2t+1}=(g_3^tg_1g_2g_1^{-1}g_3^{-
t}, 
g_3^tg_1g_1g_1^{-1}g_3^{-t}).$$ The minimal $t$ with 
$(g_1,g_2)\gamma^{2t}=(g_1,g_2)$ is $o(g_1,g_2)$. Further, the minimal $j$ with 
$(g_1,g_2)\gamma^j=(g_1,g_2)$ divides any other integer with this property. So
$j|2o(g_1,g_2)$ and if $j$ is  odd, $j|o(g_1,g_2)$.  

From the above, if the orbit of $\gamma$ does not have length $2o(g_1,g_2)$,  it 
has length $o(g_1,g_2)$. Use the notation around \eqref{xyswitch}. The 
expressions  
$g_1y=xg_1$ and $yg_2=g_2x$
are tautologies. If $o$ is odd,
then $(g_1,g_2)q_2^o=x(g_1,g_2)q_2x^{-1}$. Assume this equals $(g_1,g_2)$, which 
is true if
and only if
$xg_1=g_2x=yg_2$. The expression $(g_1g_2)^o=1$ and $xg_1yg_2=1$ are equivalent. 
Conclude
$(yg_2)^2=1$. So long as the order of $yg_2$ is not 1, this shows    
\eqref{xyswitch} holds. If, however, $yg_2=xg_1=g_2x=g_1y=1$, then $g_1=g_2$, 
contrary to
hypothesis. 

This reversible argument shows the converse: $(g_1,g_2)q_2^o=(g_1,g_2)$ follows 
from
\eqref{xyswitch}. This concludes the proof.  
\end{proof} 

\begin{exmp}[$Z(g_1,g_2)$ in Prop.~\ref{gamOrbit}] Prop.~\ref{gamOrbit} applies 
to
compute the $q_2$ orbit length of  
$\bg\in \ni(G, \bfC)$ where $r=4$ and $G$ is centerless (see Lem.~\ref{qi2}).  
With
$\bg=(g_1,g_2,g_3,g_4)$, replace $(g_1,g_2)$ in the proposition by the entries
$(g_2,g_3)$. To give examples where
$\lrang{g_3}\cap Z(g_1,g_2)\ne \{1\}$, suppose $V=\lrang{g_2,g_3}$ is a vector 
space
and
$A$ is an endomorphism fixing no nontrivial element of $V$. Entries of 
$(Ag_2^{-1},g_2,g_3,g_3^{-1}A^{-1})$ generate a centerless group. 
\end{exmp} 
 
\subsubsection{Products of conjugacy classes} \label{congClProd} The {\sl 
structure 
constant formula\/} for $\bZ[G]$ (for example, \cite[\S7.2]{Se-GT}) shows the 
result of 
multiplying elements from conjugacy classes. Let $m$ be a function constant on 
conjugacy 
classes $\row {\C} r$ and $g\in G$. Let the common value of $m$ on  $\C_i$ be 
$m(\C_i)$. Denote $\sum_{i=1}^r\sum_{u_i\in C_i} m(u_1\cdots u_r g)$ by 
$I(m;\bfC,g)$.  If $m=\chi$ is an irreducible character of $G$,  
$$I(\chi; \bfC,g)=\sum_{i=1}^r\sum_{u_i\in C_i} \chi(u_1\cdots
u_rg)=\chi(g)\prod_{i=1}^r 
\chi(u_i)/\chi(1)^r.$$ 

List the irreducible complex characters, $\row \chi s$, of $G$. Write 
$m=m_i\chi_i$. Then, $I(m;
\bfC,g)=\sum_i  m_i I(\chi_i; \bfC,g)$.  Consider 
$\psi_G=\frac{1}{|G|} 
\sum_{i=1}^s \chi_i(1)\chi_i$: 1 at $1_G$ and 0 otherwise. So, $I(\psi_G; 
\bfC,g)$ 
counts solutions of $\row u r g=1$ with $u_i\in \C_i$: $$ N(\row \C r,g)=|G|^{r- 
1}\sum_{i=1}^s \prod_{j=1}^r\chi_i(\C_j)\chi_i(g).$$

\subsection{General context with modular curve illustrations} 
\label{modCurveIllus} We do many computations
of $q_2$ orbits on absolute, inner and reduced Nielsen classes. For
general $r$  an obvious, but valuable, listing of $Q_i$ orbits appears in the 
next
lemma. 

\begin{lem} 
\label{qi2} Let $\bg\in \ni(G,\bfC)$ be a Nielsen class representative. With
$\mu=g_ig_{i+1}$, the orbit  of 
$Q_i$ on $\bg$ is the collection   $$\scases{(\bg)q_i^k}{(g_1,\dots,g_{i- 
1},\mu^lg_i\mu^{-l},\mu^lg_{i+1}\mu^{-l}, g_{i+2},\dots)}{$k=2l$} 
{(g_1,\dots,g_{i- 
1},\mu^lg_ig_{i+1}g_i^{-1}\mu^{-l},\mu^lg_i\mu^{- l},g_{i+2},\dots)} 
{$k=1+2l$}$$ 
\end{lem}

Producing modular curves from Hurwitz spaces illustrates modding out by various 
groups
of automorphisms of $G$ on $\ni(G,\bfC)$ (\cite[\S2]{Fr3-4BP} or 
\cite{FrGGCM}). 

\subsubsection{Modular curve cusps from a Hurwitz viewpoint} \label{HurView} The 
Nielsen class in this case is for s-equivalence (\S\ref{refEquiv}) of {\sl 
dihedral involution 
realizations\/}.

Regard $D_{p^{k+1}}=\bZ/p^{k+1}\xs \{\pm1\}$,  the 
dihedral group of order $2\cdot p^{k+1}$ with $p$ an odd prime, as $2\times 2$
matrices. Start with absolute  Nielsen classes (\S\ref{intEquiv}) from the 
degree
$p^{k+1}$ representation of 
$D_{p^{k+1}}$. Take $\C_2$ the 
conjugacy class of $\smatrix {-1} {\ 0} 0 1$; and $T:D_{p^{k+1}}\to S_{p^{k+1}}$ 
the 
standard representation. The upper-right element $b\in \bZ/p^{k+1}$ determines  
elements of $\C_2$.  Let $N_k$ be the normalizer in $S_{p^{k+1}}$ of 
$D_{p^{k+1}}$.  

List $\ni(D_{p^{k+1}},\bfC_{2^4}, T)^\abs=\ni_k^\abs$ as  $(b_1,b_2,b_3,b_4)\in 
(\bZ/p^{k+1})^4$ modulo $N_k$. Count these using 
representatives with these properties: \begin{edesc} \item $b_1=0$, $b_2- 
b_3+b_4=0$; and \item $b_2=1$ or $p\,|\, b_2$ and $b_3=1$. \end{edesc} Then: 
$|\ni_k^\abs|= p^{k+1}+p^k$. Modding out only by $D_{p^{k+1}}$ \wsp instead of 
$N_k$\wsp gives $|\ni_k^\inn|= (p^{k+1}+p^k)\phi(p^{k+1})/2$ with $\phi$ the Euler 
$\phi$-function.  

Renormalize: Use $b_2-b_3=ap^u$  with $(a,p)=1$ in place of $b_2=1$ ; conjugate 
by 
$\smatrix \alpha 0 0 1$ to take $a=1$. This allows further conjugation with 
$\alpha\equiv 1 
\bmod p^{k+1-u}$. 

We pose some computations as an exercise. The superscript ${}^\rd$ means to add 
the 
quotient action of $\sQ''$. \begin{edesc} \label{modcurve} \item  
\label{modcurveb} 
Show 
$|\ni^{\abs,\rd}_k|=|\ni^\abs_k|$ and $|\ni^{\inn,\rd}_k|=|\ni^\inn_k|$: $\sQ''$ 
action is trivial.
Hint:  Use H-M rep. $b_2=0, b_3=1=b_4$. \item \label{modcurved} Find a length 1 
(resp.~$p^{k+1}$) $q_2$ orbit on $\ni_{k}^{\abs,\rd}$: $(b_2,b_3)=(1,1)$ (resp. 
$(b_2,b_3)=(0,1)$).   \item \label{modcurvee} For $1\le u\le k$ count $q_2$ 
orbits on 
$\ni^{\abs,\rd}_k$ with $b_2-b_3=p^u$. \end{edesc}

Hints on \eql{modcurve}{modcurvee}: Let ${}_u\bg$ have $b_2=1$ and $b_3=1-p^u$. 
Use:
\begin{triv} 
\label{b2list} $({}_u\bg)q_2^{l}$ corresponds to the pair $(b_2,b_3)=(1+lp^u, 
1+(l-
1)p^u)$.
\end{triv}  \noindent Compute the minimal  $l$ so the 2nd entry of 
$({}_u\bg)q_2^{l}$ equals
$\smatrix {-1}  {1} 0 1$, $1\le l\le p^{k+1-u}$. Complete  the 
set of $b_2\,$s: prime to $p$, modulo $p^u$ up to multiplication by 
$H_u=\{\alpha\mid 
\alpha\equiv 1\bmod p^{k+1-u}\}$. List $q_2$ orbits on $|\ni_k^\abs|$.

\begin{lem} Let $S_u\le H_u$ be the subgroup stabilizing \eqref{b2list}. Then, 
$|S_u|=p^{\min(u,k+1-u)}$ and this gives a $q_2$ orbit of length $p^{k+1- 
u}/p^{\min(u,k+1- u)}$. Each $u$, $1\le u\le k$ contributes $\phi(p^u)/|H_u|$ 
orbits of 
$q_2$ of this length. \end{lem}

A traditional way to compute $X_0(p^{k+1})$ cusp data is in 
\cite[p.~25]{ShAutFuncts}.
Note: In  \eql{modcurve}{modcurved}, an {\sl H-M rep.\/}~(from 
$(b_2,b_3)=(0,1)$, see
\S\ref{startHM}) gives the longest $q_2$ orbit. Its shift gives the shortest
$q_2$ orbit. This is  the exact analog of Prop.~\ref{HMreps-count}. 

\subsubsection{Modular curves from \S\ref{HurView}} \label{compMC} 
Any $\bp\in \sH_k^\abs$ represents a cover in the Nielsen class $\ni^{\abs}_k$: 
$\phi_\bp:X_\bp=X\to \prP^1_z$.  Its four branch cycles all have order two, each 
fixing just one integer in the representation. This gives $g(X_\bp)=0$ from 
Riemann-Hurwitz (as in \cite[\S2]{Fr3-4BP} or \cite{FrGGCM}). The {\sl 
geometric\/} (over $\bC$) Galois closure of $\phi_\bp$ is a nondiagonal 
component $\hat 
X_\bp$ of the $p^{k+1}$ times fiber product $X_{\prP^1_z}^{(p^{k+1})}$ 
(\S\ref{sequiv} has 
the fiber product construction): 
$g(\hat X_\bp)=1$. The copy of 
$\bZ/p^{k+1}$ is the (unique) $p$-Sylow of the automorphism group of $\hat 
X_\bp\to 
\prP^1_z$. Let $\Pic^0(\hat X_\bp)$ be the elliptic curve of degree 0 divisor 
classes on 
$\hat X_\bp$.  

Let $Y_0(p^{k+1})$ be $X_0(p^{k+1})$ without its cusps (points over 
$j=\infty$). Points of  $Y_0(p^{k+1})$ are elliptic curve isogeny classes $E\to 
E'$ with kernel
isomorphic to 
$\bZ/p^{k+1}$. So, $\bp$ gives a point in $Y_0(p^{k+1})$: $\Pic^0(\hat 
X_\bp)=E\mapsto E/(\bZ/p^{k+1})=E'$. 

Produce $\sH^\inn_k\to Y_1(p^{k+1})$: If $\bp^*\in \sH^\inn_k$ lies over $\bp\in 
\sH^\abs_k$, 
then $\phi_{\bp^*}$ identifies $\Aut(\hat X_\bp/\prP^1_z)$  with $D_{p^{k+1}}$. 
Pick 
a generator $e\in \bZ/p^{k+1}$. Let ${\pmb o}_E$ be the origin on $E$. This 
gives a $p$-division
point on
$E$:
$e\mapsto  e({\pmb o}_E)$. For the Hurwitz interpretation of $Y(p^{k+1})$ add 
{\sl markings\/} to
the  Nielsen classes \cite{BFr}.

\subsubsection{Moduli dilemma relating $\sH(D_{p^{k+1}},\bfC_{2^4})^{\inn,\rd}$ 
to 
$Y_1(p^{k+1})$} \label{fineModDil} \S\ref{compMC} produces a $\bQ$ rational 
isomorphism of 
$\sH(D_{p^{k+1}},\bfC_{2^4})^{\inn,\rd}$  and 
$Y_1(p^{k+1})$. Their moduli problems, however, are different. 

Points of $Y_1(p^{k+1})$ correspond to pairs $(E,e)$ with $E$ an elliptic curve 
and $e$ 
a $p^{k+1}$ division point. Since $Y_1(p^{k+1})$ is a {\sl fine\/} moduli 
space, a $K$ point $\bp\in Y_1(p^{k+1})$ automatically produces the pair 
$(E_\bp,e_\bp)$ over $K$. This gives $E_\bp\to E_\bp/\lrang{- 1,e_\bp}$: $-1$ is 
the unique involution of $E_\bp$ fixing the origin; and $e_\bp$ is shorthand for 
translating 
points on $E_\bp$ by $e_\bp$. 

By contrast, $\sH(D_{p^{k+1}},\bfC_{2^4})^{\inn,\rd}=\sH_k^\rd$ is the moduli 
space of 
genus 1 Galois covers of $\prP^1$ with these further properties. \begin{edesc} 
\item The covering
group has an isomorphism with $D_{p^{k+1}}$, up to conjugation. \item
Branch cycles for  the  cover are in $\bfC_{2^4}$.   \item Covers $\hat\phi: 
\hat X\to \prP^1_z$ 
and 
$\alpha\circ \hat\phi: \hat X\to \prP^1_z$ for $\alpha\in \PSL_2(\bC)$ are 
equivalent if pullback
by $\alpha$ respects the isomorphism $G(\hat X/\prP^1_z)$ to $D_{p^{k+1}}$ (in 
\S \ref{innEquiv}). \end{edesc} In  \eql{modcurve}{modcurveb}, 
$\sQ''$ acts trivially on inner Nielsen classes. So $\sH_k^\rd$ is not a fine 
(or even b-fine)
moduli space  (Prop.~\ref{redHurFM}). We explain further. Let 
$U_j=\prP^1_j\setminus\{\infty\}$. 

Let $\Phi:\sT\to S$ be a smooth family of curves with an attached algebraic map 
$\Psi: S\to 
U_j$. Denote the fiber of $\Phi$ over $s\in S$ by $\sT_s$. Suppose  
$D_{p^{k+1}}$ acts as a group scheme on $\sT$  preserving each fiber 
$\sT_s$. Write this as $\Gamma: D_{p^{k+1}}\times \sT\to \sT$. Then,  
form the quotient map $\phi_s: \sT_s\mapsto \sT_s/D_{p^{k+1}}$. Assume for 
$s\in S(\bC)$:  \begin{triv} \label{sQuotient}  An isomorphism of 
$\sT_s/D_{p^{k+1}}$ 
with $\prP^1_z$ presents $\phi_s$ in the inner Nielsen class of 
$(D_{p^{k+1}},\bfC_{2^4})$ with branch points in the equivalence class 
$\Psi(s)$. 
\end{triv}   By assumption, $\Phi$ induces a natural map $\psi: S\to\sH_k^\rd$. 
If 
$\sH_k^\rd$ were a fine moduli space, the family $\Phi$ would pull back from a
{\sl unique\/}  total family on $\sH_k^\rd$. Though it isn't, there is at least 
one total family on
$\sH_k^\rd$ (identified with $Y_1(p^{k+1})$).  

\begin{prop} \label{X_1TotFam} For some total 
family $\Phi^*:\sT^*\to Y_1(p^{k+1})$: 
\begin{edesc} \item $D_{p^{k+1}}$ acts on the total space commuting with 
$\Phi^*$; and 
\item  the fiber over 
$\bp\in Y_1(p^{k+1})$ is $E_\bp$ (notation above) with restriction of 
$D_{p^{k+1}}$ producing
the  cover $E_\bp\to E_\bp/\lrang{-1,e_\bp}$. \end{edesc} 

Any $K$ point on $\sH_k^\rd$ produces a 
$D_{p^{k+1}}$ regular realization in the inner Nielsen class of 
$(D_{p^k+1},\bfC_{2^4})$. \end{prop} 

\begin{proof} It only remains to show the last statement. Prop.~\ref{isotopyM4} 
gives
this because \eql{modcurve}{modcurveb} shows $\sQ''$ acts trivially on the 
Nielsen
class. A classical argument, however, suffices in this case by contrast to
\S\ref{applyA5}.  

Given $\bp\in
\sH_k^k(K)=Y_1(p^{k+1})(K)$, consider the
$K$ cover 
$E_\bp\to E_\bp/\lrang{-1,e_\bp}$. If $E_\bp/\lrang{-1,e_\bp}$ is isomorphic to 
$\prP^1_z$ over $K$, this gives a regular $D_{p^{k+1}}$ realization in the inner 
Nielsen class of $(D_{p^k+1},\bfC_{2^4})$. Such an isomorphism occurs if and 
only if 
$E_\bp/\lrang{-1,e_\bp}$ has a $K$ point. The image from  origin in 
$E_\bp$ is such a $K$ point. 
\end{proof}

Suppose $\sH_k^k$ were a fine (resp.~b-fine) moduli space. Then the family 
$\Phi$ (resp.~off of
$j=0$ or 1) would be the pullback  of this family. There would then be a section 
$\gamma: S\to \sT$
for
$\Phi$ from pulling  back the section for the origin on the fibers of $\sT^*\to 
Y_1(p^{k+1})$.
That section  would have the same field of definition as does 
$(\Phi,\sT,\Psi,\Gamma)$. The
structural maps, including $\Gamma$, have field of definition $\bQ$.  We don't 
assume  the genus
zero curve
$\sT_s/D_{p^{k+1}}$ is isomorphic to $\prP^1_z$ over 
$\bQ(s)$; this would give a $\bQ(s)$ regular realization $\sT_s\to \prP^1_z$ of 
$D_{p^{k+1}}$. It is exactly in this way the next proposition says the family 
$\sT^*\to 
Y_1(p^{k+1})$ is special. 

\subsubsection{$D_{p^{k+1}}$ involution covers of forms of $\prP^1$} 
\label{DpkRegReal} 
Prop.~\ref{X_1TotFam} would be most satisfying if there was a diophantine 
consequence 
for $\sH_k^\rd$ not being a fine (or b-fine) moduli space. Best of all would be 
values of $s$ with
$\sT_s/D_{p^{k+1}}$ not isomorphic over
$\bQ(s)$ to $\prP^1_z$. For example, suppose $K$ is a $p$-adic field and
$E$ is an elliptic curve over $K$ with a $p^{k+1}$ division point $e$. From this 
form $\psi_e:E\to
E/\lrang{-1,e}$  as in the Prop.~\ref{X_1TotFam} proof. Let $\mu$ be the $G_K$ 
module
of all roots of 1 and let $E^*$ be the dual abelian variety to $E$. Tate has a 
theory
(\cite{TateWC}, \cite[p.~93]{SeGalCoh}) identifying homogeneous spaces for $E$ 
with
$\Hom_{G_K}(H^0(G_K,E^*),\mu)$. 

We suggest a project. Let $C$ be a nontrivial homogeneous space for $E$, and 
$C^{(2)}$ its 
2-fold symmetric product. Map positive divisors on $C$ of degree 2 to
divisor classes as $s^{(2)}: C^{(2)}\to \Pic(C)^{(2)}$. Fibers of $s^{(2)}$ are 
forms of
$\prP^1$ over $K$ with a degree 2 map from $C$. Assume further that
$\Pic(C)^{(2)}$ has $K$ points; it is isomorphic to $E$. \cite{FrBGJac} shows 
there
is a $K$-fiber $Y$ of $s^{(2)}$ not isomorphic to $\prP^1$ over $K$.

\begin{prob} Assume $j(E)\ne  0$ or 1 (the special values for our normalized
$j$). Does
$s^{(2)}$ have a
$D_{p^{k+1}}$ invariant
$K$ fiber
$Y$ not isomorphic to
$\prP^1_K$? 
\end{prob}   

A positive solution gives a sequence $C\to Y\mapright{\alpha} Y$ (with $\alpha$ 
a degree
$p^{k+1}$ map) that is a
$D_{p^{k+1}}$ involution realization over $\bar K$ ($\alpha$ agrees with 
\S\ref{genus0}).
Despite Ch\^ atelet's
interest in Brauer-Severi varieties and homogeneous spaces for elliptic curves
(Rem.~\ref{chatalet}) he seems not to have come across this problem.  

\subsection{$p$-lifting property and using $\sh$-incidence}
\S\ref{1stOrbComp} looks at 
$\gamma_\infty$ orbits. \S\ref{HMrepRubric} applies to 
$\ni(A_5,\bfC_{3^4})^{\inn,\rd}$ a preferred rubric for computing $\bar M_4$ 
orbits. 

\subsubsection{Orbits for $\gamma$ related to $\ni(A_5,\bfC_{3^4})^{\inn}$} Let 
$\gamma$ be
the twist operator of \S\ref{twistAct}. 
\label{1stOrbComp}  Take $(g_1,g_2)=((1\,2\,3),(1\,4\,5))$. Then,  
 the $\gamma$ orbit of $(g_1,g_2)$ has length 5: 
$$(1\,2\,3\,4\,5)^2(1\,2\,3)(1\,2\,3)(3\,2\,1)(1\,2\,3\,4\,5)^{-2}=(1\,4\,5).$$ 
Similarly, 
take 
$(g_1,g_2)=((1\,2\,3),(1\,3\,4))$. On this element $\gamma$ produces an orbit of 
length 
3: 
$(1\,2\,4)(1\,2\,3)(1\,2\,3)(3\,2\,1)(4\,2\,1)=(1\,3\,4)$. 

\subsubsection{The $p$-lifting property} \label{plifting}  The following lemma 
from 
\cite[Lifting 
Lemma  4.1]{FrKMTIG} appears often in our computations.

\begin{lem} \label{FrKMTIG} For any finite group $G=G_0$, let $G_1$ be the first 
characteristic quotient of the universal $p$-Frattini cover ${}_p\tilde G$ of 
$G_0$ 
(\S\ref{pfrattini}). Suppose $g\in G_0$ has order divisible by $p$. Then, any 
lift of $g$ to 
$G_1$ has order $p$ times the order of $g$. \end{lem} 

Given $G=G_0$, an extension $\psi: 
H\twoheadrightarrow G$ has the {\sl $p$-lifting property\/} if the conclusion of 
Lem.~\ref{FrKMTIG} holds (for $p$). Often a much smaller cover than $G_1\to G_0$ 
will suffice for the $p$-lifting property. For example, $A_5=\PSL_2(\bZ/5)$ has 
one 
conjugacy class of elements of order $2$. That lifts to an element of order 4 in 
$\hat 
A_5=\SL_2(\bZ/5)$. 

\begin{lem} \label{Q8Sp5} Let $g_i$, $i=1,2$, be distinct order 2 elements in 
$A_5$.  
Denote lifts of them to $\hat A_5$ by $\hat g_i$, $i=1,2$. Then, $\lrang{\hat 
g_1,\hat g_2}$ 
is a quaternion group. 

Consider $\hat A_n$, the Spin cover of $A_n$. Write any $g\in A_n$ of order 2 
($n\ge 
6$) as a product $g_1g_2$ of two commuting elements $g_1,g_2$ of order 2  with 
the 
$g_i\,$s having lifts $\hat g_1,\hat g_2$  of order 4. Then, $\lrang{\hat 
g_1,\hat g_2}$ is 
isomorphic to $\bZ/2\times \bZ/4$ if and only if $\hat g_1\hat g_2$ has order 2 
(otherwise 
it is isomorphic to $Q_8$). Further, $\hat g_1\hat g_2$ 
has order 4 if and only if $g_1g_2$ is a product of $2s$ transpositions, with 
$s$ odd. 
\end{lem} 

\begin{proof} The opening statement follows because the quaternion group of 
order 8 is  
only group of order 8 with a Klein 4-group as quotient having 
all its elements lifting to have order 4. The last statement is a variant on 
that using
Prop.~\ref{liftEven}. \end{proof}

Let $G_1$ be ${}_2\tilde A_5/\ker_1$ for $p=2$, the first
characteristic quotient of  the universal 2-Frattini cover of $A_5$ 
(\S\ref{pfrattini}). Suppose
$(g_1',g_2')$ is a lift  of 
$(g_1,g_2)$ to elements of order 3 in $G_1$ in either of these cases. Then, 
$\gamma$ on 
$(g_1',g_2')$ gives orbits of length $2\ord(g_1'g_2')$. For example, let 
$\alpha=(1\,2\,3\,4\,5)^2(1\,2\,3)=(2\,4)(3\,5)$. Suppose the orbit has length 
$\ord(g_1'g_2')$. Then, some lift of $\alpha$ to $\alpha'\in G_1$ conjugates 
$g_1'$ to $g_2'$ and $g_2'$ to $g_1'$. Such an element, however, would have 
order 2. 
Since 2 divides the order $\alpha$, this contradicts Lemma \ref{FrKMTIG}: 
$\alpha'$ has 
order 2 times the order of $\alpha$. Similarly for 
$\alpha=(1\,2\,4)(1\,2\,3)=(1\,3)(2\,4)$. 

We don't know how to compute $\ord(g_1'g_2')$ in general as in Lemma 
\ref{FrKMTIG}. When $G_0$ is a subgroup of an alternating group and $p=2$, 
Prop.~\ref{serLift} 
contributes considerably. Here is an immediate example. \S\ref{onlytwo} uses it 
on 
indirect examples. 

\begin{lem} \label{A5Prod} Take $G_0=A_5$ and $p=2$. Suppose $(g_2',g_3')$ is a 
pair of elements of order 3 lifting to $G_1$ either 
$(g_2,g_3)=((1\,2\,3),(1\,4\,5))$ or 
$(g_2,g_3)=((1\,2\,3),(1\,3\,4))$. Then, $\ord(g_2'g_3')=2\ord(g_2g_3)$. 
Further, 
assume 
$\bg=(g_1,g_2,g_3,g_4)$ is in\/ {\rm Table \ref{lista5}}. Then, if 
$\ord(g_2g_3)=3$ or
5, the  conclusion holds for $\ord(g_2'g_3')$.   \end{lem} 

\subsubsection{Generating orbits from H-M reps.~for 
$\ni(A_5,\bfC_{3^4})^{\inn,\rd}$} 
\label{HMrepRubric} Improve Table \ref{lista5} to get a complete orbit for $\bar 
M_4=\lrang{\sh,\gamma_\infty}$ on $\ni(A_5,\bfC_{3^4})^{\inn,\rd}$ in 2-steps.

First: List $\gamma_\infty$ orbits for the H-M reps. 
$$\bg_1=((1\,2\,3),(1\,3\,2),(1\,4\,5),(1\,5\,4)) \text{ and  } \bg_2= 
((1\,2\,3),(1\,3\,2),(1\,5\,4),(1\,4\,5)).$$ Here are $(\bg_1)\gamma_\infty^j$, 
$j=1,2,3,4$: 
\begin{equation} \label{q2hmrepa5} 
(\!\bullet\!,(2\,4\,5),(1\,3\,2),\!\bullet\!),(\!\bullet\!, 
(5\,1\,3),(2\,4\,5),\!\bullet\!), (\!\bullet\!,(3\,2\,4),(5\,1\,3),\!\bullet\!),  
(\!\bullet\!,(4\,5\,1),(3\,2\,4),\!\bullet\!).\end{equation} A $\bullet$ means 
repeat that  
position 
from the previous element.  Conjugate these by $(4\,5)$ for list 
$(\bg_2)\gamma_\infty^j$, 
$j=2,3,4,5$. The {\sl middle product (order)\/} of a 4-tuple $\bg$ is 
$u=(\bg)\eqdef\ord(g_2g_3)$. 
 
Second: List elements $\bg'=(\bg)\sh$ where $\bg$ is in \eqref{q2hmrepa5} and 
$\mpr(\bg')=1$ ($\bg=\bg_1$) or 3 ($(\bg_1)\gamma_\infty^j$, $j=2,3$) and close 
under 
$\lrang{\gamma_\infty}$ action.  Those with $\mpr(\bg')=1$ are $(\bg_1)\sh$ and 
$((\bg_1)\sh)\gamma_\infty$. Those with $\mpr(\bg')=3$ are in 
$((\bg_1)\gamma_\infty^j)\sh\lrang{\gamma_\infty}$, $j=2,3$: \begin{equation} 
\label{q2shq2a5}\begin{array}{rl} 
&((5\,1\,3),(2\,4\,5),(1\,5\,4),(1\,2\,3)),(\bullet,(2\,1\,4),(2\,4\,5),\bullet) 
, 
(\bullet,(4\,1\,5),(2\,1\,4),\bullet)\\ 
&((3\,2\,4),(5\,1\,3),(1\,5\,4),(1\,2\,3)),(\bullet,(3\,4\,5),(5\,1\,3),\bullet) 
, 
(\bullet,(5\,4\,1),(3\,4\,5),\bullet).\end{array} \end{equation} Note: Each row 
of 
\eqref{q2shq2a5} contains two elements which $\sh$ maps to each of the type 
$(5,5)$ 
orbits (5 is the order of both the middle twist and $\gamma_\infty$ orbits; as 
in
\S\ref{cuspNotat}). 

Example: Let
$\bg^*$ be 
$(\bullet,(2\,1\,4),(2\,4\,5),\bullet)$ as above. 
Then, $$(\bg^*)\sh^{-1}=((1\,2\,3),(5\,1\,3),(2\,1\,4),(2\,4\,5)),$$ so 
$(1\,2\,3)(\bg^*)\sh^{-1}(3\,2\,1)$ 
is in the $\gamma_\infty$ orbit of $\bg_2$. The 3rd elements of each 
\eqref{q2shq2a5} 
row map to each other by $\sh$. Conclude: Two steps produce a list of 4-tuples 
of 3-cycles 
(modulo conjugation by $A_5$) closed under $\sh$ and $\gamma_\infty$. 

\subsection{The $\sh$-incidence matrix} \label{shincidence} For a
general reduced Nielsen class, list $\gamma_\infty$ orbits as $\row
O n$. The
$\sh$-incidence matrix $A(G,\bfC)$ has $(i,j)$ term $|(O_i)\sh\cap O_j|$.
\subsubsection{The $\sh$-incidence matrix when $r=4$}
Since $\sh$ has order two on reduced Nielsen classes, this is a symmetric
matrix. Equivalence $n\times n$ matrices $A$ and $TA{}^\tr T$ running
over permutation matrices $T$ (${}^\tr T$ is its transpose) associated to
elements of $S_n$. List $\gamma_\infty$ orbits as $$\row {O_1,} {t_1},
\row {O_2,} {t_2},\dots, \row {O_u,} {t_u}$$ corresponding to $\bar M_4$
orbits. Choose $T$ to assume $A(G,\bfC)$ is
arranged in blocks along the diagonal.  

\begin{lem} \label{shincBlocks} If $A_j$ is the $j$th block of $A(G,\bfC)$, then 
$A_j$ doesn't 
break into smaller blocks. So, $\bar M_4$ orbits form irreducible 
blocks in the $\sh$-incidence matrix.\end{lem} 

\begin{proof} With no loss assume one $\bar M_4$ orbit and two 
blocks, with orbit listings as $\row O k, O_{k+1}, \dots, O_t$.
As, however, there is one orbit, for some $j\le k$, $|(O_i)\sh\cap
O_j|\ne 0$ for some $i>k$. This contradicts there being two
blocks. \end{proof}

In practice it is difficult to list the $\gamma_\infty$ orbits.
So, we start with the H-M reps., apply $\sh$, then complete the
$\gamma_\infty$ orbits and check $|(O_i)\sh\cap
O_j|$. Sometimes we'll then be done. The case $(A_5,\bfC_{3^4})$
illustrates this. Denote (as above) the 
$\gamma_\infty$ orbits of $\bg_1$ and $\bg_2$ by $O(5,5;1)$ and $O(5,5;2)$;   
$\gamma_\infty$ orbits of $$((5\,1\,3),(2\,4\,5),(1\,5\,4),(1\,2\,3)) \text{ and 
} 
((3\,2\,4),(5\,1\,3),(1\,5\,4),(1\,2\,3))$$ by $O(3,3;1)$ and $O(3,3;2)$; and of 
$(\bg_1)\sh$ 
by $O(1,2)$.

\begin{table}[h] \caption{$\sh$-Incidence Matrix for $\ni_0$} \label{shincni0}
\begin{tabular}{|c|ccccc|} \hline Orbit & $O(5,5;1)$\ \vrule  &
$O(5,5;2)$\ \vrule & 
$O(3,3;1)$\ \vrule &$O(3,3;2)$ \ \vrule& $O(1,2)$\\ \hline $O(5,5;1)$ 
&0&2&1&1&1\\ 
$O(5,5;2)$ &2 &0&1 &1 & 1 \\ $O(3,3;1)$ &1&1&0&1&0\\ $O(3,3;2)$ & 1 &1 &1  
& 0 & 0 \\  $O(1,2)$ & 1 &1 & 0 & 0 & 0\\  \hline \end{tabular} \end{table}

\subsubsection{The $\sh$-incidence matrix for general $r$} \label{genrShInc} 
For general $r$, denote $q_1\cdots q_{r-1}$ at the Hurwitz
monodromy level to be the shift $\sh_r$, so $\sh_4$ is what we call the {\sl
shift\/} above. Ideas for $r=4$ generalize to indicate cusp
geometry for general $r$. 

The element $\sh_r$ plays the role of a shift in two ways. Consider an
$r$-tuple
$\bar\psigma= (\row {\bar \sigma} r)$ of free generators of $F_r$. The effect of $\sh_r$ on
$\bar \psigma$  is to give 
$$(\bar \psigma)q_1\cdots q_{r-1}=({\bar \sigma}_1{\bar
\sigma}_2{\bar \sigma}_1^{-1},\dots,{\bar \sigma}_1{\bar \sigma}_r{\bar \sigma}_1^{-1},
{\bar \sigma}_1{\bar \sigma}_1{\bar \sigma}_1^{-1}).$$ In specializing to a Nielsen class
the effect of $\sh_r$ is to the shift the Nielsen class representative entries
by 1.  Iterate this
$r$ times to see the effect of $\sh_r^r$ is conjugation on  
$\psigma$ by the product $\bar \sigma_1\cdots \bar
\sigma_r$ of these generators. 

Such a conjugation commutes with the action of the
braid group. So we have an interesting interpretation for the action of {\sl
conjugating\/} by  
$\sh_r$ on the generators $q_1,\dots,q_r$.  Define
$q_0$ to be
$\sh_r^{-1}q_1\sh_r$. Then, conjugation by $\sh_r$ on the left of the array
$(q_0,q_1,\dots,q_{r-2},q_{r-1})$ maps its entries to 
$$\begin{array}{rl} \sh_r(q_0,\dots,q_{r-1})\sh_r^{-1}=&
(q_1,q_1q_2q_1q_2^{-1}q_1^{-1},q_1q_2q_3q_2q_3^{-1}q_2^{-1}q_1^{-1},\dots) \\
=&(q_1,q_2,\dots, q_{r-1},q_0).\end{array}$$ To see the effect of conjugation of $\sh_r$ on
$q_{r-1}$ use that $\sh_r^r$ is in the center of $H_r$ (or of $B_r$). Then,
$q_0=\sh_{r}^{r}(q_0)\sh_{r}^{-r}=\sh_rq_{r-1}\sh_r^{-1}$.

Denote $q_{r-1}q_{r-2}\cdots q_1$ by $\sh_r'$. Notice $(\sh_r')^r$ has exactly the same effect on
$\psigma$ as does $\sh_r^r$. In $H_r$ use that $q_1\cdots q_rq_r\cdots q_1=1$ to see
$\sh_r^r(\sh_r')^r=1$, so
$\sh_r^r=z$ has its square equal to 1.  When $r=4$ the group $M_4$ is exactly
$H_r/\lrang{\sh_4^4}=H_r/\lrang{z}$. An especially handy description of $z$ in this case is
$q_1q_3^{-1}$. In general there is a $\sh_r$-incidence matrix. As in
the case
$r=4$, it suffices to choose the image of $q_v$ in $\bar M_r$ for some value of $v$. It doesn't
make any difference which $v$, though for $r=4$ it was convenient to take $v=2$. Call the
resulting element $\gamma_\infty$. List the $\gamma_\infty$ reduced orbits as $\row O t$ and
define
$A(G,\bfC)$ to be the matrix with 
$(i,j)$ term
$|(O_i)\sh_r\cap O_j|$. For general $r$ it won't be symmetric. 

\section{Nielsen classes and the Inverse Galois Problem} \label{arithmetic}

Let $G$ be a centerless transitive subgroup of $S_n$. Assume $\bfC$ is $r\ge 3$ 
conjugacy classes of $G$. This data produces several canonical moduli spaces. 
Four occur 
often: inner and absolute (adding a permutation representation on $G$) Hurwitz 
spaces and  
their reduced versions (quotients by a natural $\PSL_2(\bC)$ action).  We review 
their 
definitions  and  relation to $H_r$. When $r=4$, the reduced versions are 
curves: upper 
half plane quotients by finite index subgroups of $\PSL_2(\bZ)$ covering the 
$j$-line and 
ramified at traditional points.

\subsection{Equivalences on covers of $\prP^1_z$} Consider covers $\phi_i: 
X_i\to \prP^1_z$,
$i=1,2$. There are two natural 
equivalences between isomorphisms $\alpha:X_1\to X_2$.  
\begin{edesc} 
\label{covEq} \item \label{covEqs} Call $\alpha$ an  {\sl s-equivalence\/} if 
$\phi_2\circ
\alpha=\phi_1$. 
\item \label{covEqw} Call $\alpha$ a {\sl w-equivalence} if there exists 
$\beta\in \PSL_2(\bC)$
with $\phi_2\circ\alpha=\beta\circ\phi_1$. 
\end{edesc}

\begin{defn} In the notation above, call $\phi_1$ and $\phi_2$ s-equivalent 
(resp.~w-equivalent)
if \eql{covEq}{covEqs} (resp.~\eql{covEq}{covEqw}) holds for some $\alpha$. 
\end{defn}

Moduli spaces formed from equivalence \eql{covEq}{covEqs} support a natural 
$\PSL_2(\bC)$ action. The quotient is a moduli space for equivalence 
\eql{covEq}{covEqw}. Let $\sH^\rd$ be one of these reduced moduli spaces. If 
$r=4$, its 
projective completion gives a ramified cover $\bar \psi:\bar\sH^\rd\to 
\prP^1_j$. We refine 
\cite[Prop.~6.5]{DFrIntSpec} for computing a branch cycle description of $\bar 
\psi$.  
This clarifies that a geometric equivalence from a Klein 4-group comes exactly 
from the group 
labeled $\sQ''=\sQ/\lrang{z}$ in \S\ref{invH4}.

\subsubsection{Setup for Nielsen classes} \label{setupNC} Consider any smooth 
connected family of $r$  (fixed) branch point covers of $\prP^1_z$. Denote the 
degree
of the covers in the  family
by $n$. It has this attached data:
\begin{edesc}
\item an  associated group $G$; \item a permutation  representation $T:G\to 
S_n$; and 
\item a set of  conjugacy classes $\bfC=(\row \pC r)$ of  $G$.\end{edesc} 

Denote the subgroup of 
$S_{n}$ normalizing $G$ and  permuting the conjugacy classes $\bfC$ by
$N_{S_{n}}(G,\bfC)=N_T(G,\bfC)$. If a cover $\phi:X\to \prP^1_z$ in this family 
has 
definition
field
$K$, then the Galois closure of the cover has Galois ({\sl arithmetic 
monodromy\/}) group  
a
subgroup of
$N_{S_{n}}(G,\bfC)$.  Below use the notation $G(T,1)$ for
the  elements $g\in G$ with $(1)T(g)=1$: $G(T,1)$ is the stabilizer of 1 in the
representation.

\subsubsection{Nielsen class and $H_r$ action} \label{defNielClass} Relate 
covers with  
data 
$(G,\bfC)$ to homomorphisms of fundamental groups using the associated  {\sl 
Nielsen  
class\/}: $$\ni=\ni(G,\bfC)=\{\bg= (\row g r) \mid g_1\cdots g_r=1, 
\lrang{\bg}=G \text{\ 
and\ } \bg\in \bfC\}.$$ Writing $\bg\in\bfC$ means the  conjugacy classes of the 
$g_i\,$s 
in $G$ are,  in some possibly rearranged order, the  same (with multiplicity) as 
those in 
$\bfC$.

Suppose given $\bz$, $r$ branch points, and a choice $\bar \bg$ of classical 
generators for 
$\pi_1(U_\bz,z_0)$ (\S \ref{classGens}). Then, $\ni(G,\bfC)$ lists all 
homomorphisms 
from $\pi_1(U_\bz,z_0)$ to $G$ giving a cover with branch points $\bz$ 
associated to $(G,\bfC)$.  A Galois cover (over $\bQ$) with attached data 
$(G,\bfC)$ 
is a $(G,\bfC)$ {\sl realization\/}, suppressing $\bQ$. Analogous results for 
any 
characteristic 0 field (including non-Galois covers) have easy formulations as 
in  
\cite[Thm.~5.1]{FrHFGG}.

Denote the subgroup of automorphisms $\Aut(G)$ of $G$ 
permuting elements in $\bfC$ by $\Aut(G,\bfC)$. An $\alpha\in \Aut(G,\bfC)$ 
respects the
multiplicity of a conjugacy class in $\bfC$. It acts on $\ni(G,\bfC)$ by
$$(\bg)\alpha=((g_1)\alpha,\dots,(g_r)\alpha)).$$  Given $N^*\le \Aut(G)$ form
$\ni(G,\bfC)/N^*$ in the obvious way. 

Suppose $T:G\to S_n$ is a permutation representation. Let $N'$ be a subgroup of 
$N_T(G,\bfC)$
(\S\ref{setupNC}). Conjugation on $G$ gives a homomorphism $N'\to \Aut(G,\bfC)$.
This is an isomorphism into if and only if $N'$ contains no nontrivial 
centralizer of
$G$. Abusing notation, denote the quotient of $N'$ action by
$\ni(G,\bfC)/N'$. 

The braid group $B_r$ acts from $\row Q {r\nm1}$ acting on $(\row {\bar g} r)\in 
\ni(G,\bfC)$
as in \eqref{brAct}. If $N'\le N_T(G,\bfC)$ contains $G$, this $B_r$ action 
induces an 
$H_r=\pi_1(U_r)$ action (\S\ref{notation}). Assume $G\le N'$. 
Permutation representations of fundamental groups produce covers. For this $H_r$ 
action, denote 
the ({\sl Hurwitz space\/}) cover of $U_r$ by $\sH_{N'}=\sH(G,\bfC)_{N'}$.

\subsubsection{Galois closure and permutation representations} \label{sequiv} 
Use 
notation from \S\ref{setupNC} for a cover 
$\phi:X\to \prP^1_z$ of degree $n=\deg(\phi)$ with  $(G,\bfC)$ and a 
permutation representation $T_\phi=T:G\to S_n$ attached.

The Galois closure of $\phi$ over any defining field for  $(X,\phi)$ has a 
geometric 
formulation. Take the fiber product $X^{(n)}\eqdef X^{(n)}_\phi\eqdef 
X^{(n)}_{\prP^1_z}$ of $\phi$,
$n$ times.  Points on 
$X^{(n)}$ consist of $n$-tuples $(\row x n)$ of points on $X$ satisfying 
$\phi(x_1)=\phi(x_2)=\cdots =\phi(x_n)$. This variety will be singular around 
$n$-tuples 
where $x_i$ and $x_j$ are both ramified  through $\phi$. Replace $X^{(n)}$ by 
its 
normalization to make it now a non-singular  cover. Retain the notation 
$\phi^{(n)}: 
X^{(n)}\to \prP^1_z$. Then,  $X^{(n)}$  has components where at least two of the 
coordinates are identical. These  form the {\sl fat diagonal\/}.

Remove components of this fat diagonal to give $X^*$.  Over the algebraic 
closure $X^*$ 
has as many components as $(S_n:G)$. List one of these components over $\bar K$
as $X^\dagger$. The stabilizer in $S_n$ of $X^\dagger$ is a conjugate of $G$.
Normalize by choosing $X^\dagger$ so the stabilizer is actually $G$. Now,
choose any
$K$  component $\hat X$ of $X^*$ containing $X^\dagger$.  Then, $\hat\phi: \hat 
X\to
\prP^1_z$ is Galois (over $K$)  with group $\hat G\le N_{S_n}(G,\bfC)$ having 
$G$ as
a normal subgroup. Also, $\hat 
\phi$ has the same  conjugacy classes $\bfC$ attached to the branch points $\bz$ 
and  it 
factors  through $\phi$ (project on any coordinate of  $X^{(n)}$). The Galois 
cover  $\hat 
X\to X$ has group $\hat G(1)=\hat G(T,1)$ where  $T$ is the coset representation 
of $G$  
on  $G(1)$. The following is elementary \cite[Lem.~2.1]{FrHFGG}.

\begin{lem}\label{autos} The centralizer of $G$ in $N_{S_{n}}(G,\bfC)$  induces 
the 
automorphisms of $X$ that commute with $\phi_{T}$. \end{lem}

Consider any permutation representation $T':G\to S_{n'}$. 
This provides $\phi_{T'}: X_{T'}\to \prP^1_z$;  $X_{T'}$ is the 
quotient $\hat X/G(T',1)$ (with $G(T',1)$ as in \S\ref{setupNC}). 

\subsection{Motivation from the Inverse Galois Problem}\label{IGMot}  The {\sl 
regular\/} version of the \IGP\ asks if every finite  $G$ is the group of a 
cover of $U_\bz$, 
for some integer $r$, with equations for its automorphisms and the cover over 
$\bQ$. Call 
such an $r$-branch point realization over $\bQ$. To find $r$, $\bz$ and such a 
cover 
needs structure. Form the profinite completion $\pi_1(U_{\bz},z_0)^\alg$ of 
$\pi_1(U_{\bz},z_0)$ with respect to all its finite index subgroups.  Suppose 
$z_0\in 
\bQ\cup\{\infty\}$ and $\bz$ is a $G_\bQ$ stable set. There is a natural short 
exact 
sequence  \begin{equation} \label{exactArithBP} 1\to \pi_1(U_{\bz},z_0)^\alg\to 
\pi_1(U_{\bz},z_0)^\ari \to G_\bQ \to 1\text{ (\cite[p.~57]{Se-GT} or 
Prop.~\ref{splitseq})}.\end{equation} For an analog substituting $K$ for $\bQ$, 
replace 
$G_\bQ$ by $G_K$ in \eqref{exactArithBP}.

\subsubsection{Grothendieck's $G_\bQ$ action} \label{groth} Let 
$\sE(U_\bz,z_0)^\alg$ 
be all meromorphic algebraic functions in a neighborhood of $z_0$ that extend 
analytically 
along each path in $U_\bz$. Regard them as in the Laurent series $\sL_{z_0}$ 
about 
$z_0$. Each is meromorphic in a disk about $z_0$ having the nearest point of 
$\bz$ on its 
boundary. These are the {\sl extensible\/} functions in $U_\bz$. To each such 
$f$ 
associate a {\sl definition field\/} $K=K_f$. This has generators (over $\bQ$) 
the ratios of
coefficients of any (nonzero) polynomial $F(z,w)$ in $z$ and $w$ with 
$F(z,f(z))\equiv 
0$ near $z_0$. The field $K(z,f)$ that $f$ generates over $K(z)$ is the field of 
meromorphic 
functions of a cover $\bar X\to \prP^1_z$ over $K$ (restricting to an unramified 
cover 
$X=X_f\to U_\bz$). 

Suppose $\psi: X\to U_\bz$ is a finite (unramified) cover. Then, $\psi$ is s-
equivalent to 
$\psi_f:X_f\to U_\bz$ for some $f$ with the following properties.  \begin{edesc} 
\item 
$K_f$ is in the algebraic closure of  $\bQ(z_0,\bz)$.  \item  $f$ and each 
analytic 
continuation of it (around a closed path based at $z_0$) have nonzero 
differential at $z_0$.  
\end{edesc} So, these analytic continuations have 
$n=\deg(\psi_f)$ distinct values at 
$z_0$. Their power series
around $z_0$  have coefficients in the algebraic closure of
$\bQ(z_0,\bz)$.

With no loss, when $z_0$ and $\bz$ are over $\bQ$, assume $\sE(U_\bz,z_0)^\alg$ 
consists of power series around $z_0$ with coefficients in $\bbQ$. There is an 
ordering on 
elements of $\sE(U_\bz,z_0)^\alg$: $f\le g$ if $\bbQ(z,g)\supset \bbQ(z,f)$. 
Analytic 
continuation over a path, as an operation on $\bbQ(z,g)$,  determines the path's
analytic  continuation on $\bbQ(z,f)$. So, paths acting on these equivalence 
classes respect this 
ordering. Each equivalence class represents a specific function field, with all 
functions in it 
expanded around $z_0$. It is the exact data one expects from a cover and a point 
on the 
cover over $z_0$. This action by paths in  $U_\bz$ based at $z_0$ on the points 
above 
$z_0$ respects the action on points of a higher covering. For $\gamma\in 
\pi_1(U_\bz,z_0)$, let $f_\gamma$ be the analytic continuation of $f$ to the 
endpoint of 
$\gamma$.

\begin{defn} Define how $\sigma\in G_\bQ$ acts on $\gamma\in \pi_1(U_\bz,z_0)$ 
by how it acts 
on extensible algebraic functions: $f\mapsto f_{\sigma^{-1}\circ \gamma\circ 
\sigma} =  
f_{\gamma^\sigma}$. In words: Apply $\sigma^{-1}$ to the coefficients of $f$, 
continue 
$f$ around $\gamma$ and then apply $\sigma$ to the coefficients of the result. 
The effect 
of $\gamma$ on algebraic functions determines it. This defines $\gamma^\sigma$. 
\end{defn} Unless $\sigma$ is complex conjugation, no actual path gives 
$\gamma^\sigma$. With our right action of $G_\bQ$ replacing a left action, the 
following 
is compatible with \cite[p.~103]{IharIntCong}.

\newcommand{\PP}{\sP\kern-6pt\sP} \begin{prop} \label{splitseq} Fiber compatible 
action on a projective system of points  $\PP_{z_0}$ determines paths in 
$\pi_1(U_\bz,z_0)$. Then, $\{\gamma^\sigma\}_{\gamma\in \pi_1(U_\bz,z_0), 
\sigma\in 
G_\bQ}$ also acts on $\PP_{z_0}$.

Define $\pi_1^\alg$ as the projective completion of this action from finite 
Galois covers. 
Then, $\pi_1^\alg$ is the projective limit of $\pi_1$ quotients from all normal 
subgroups of 
finite index.  The previous action of $\sigma\in G_\bQ$ on $\pi_1(U_\bz,z_0)$ 
extends to 
$\pi_1^\alg$,  giving a splitting of the sequence \eqref{exactArithBP}. 
\end{prop}

\subsubsection{The \BCL\ over $\bQ$} The \BCL\ (\cite[before Thm.~5.1]{FrHFGG} 
or
the  case in Lem.~\ref{bcl}) gives a necessary condition for a $(G,\bfC)$
realization (over 
$\bQ$) through a quotient of $\pi_1(U_{\bz},z_0)^\alg$. Let $\bar \bg$ be 
classical 
generators of $\pi_1(U_{\bz},z_0)^\alg$ (\S\ref{classGens}). Any homomorphism 
$\psi:\pi_1(U_{\bz},z_0)^\alg\to G$  produces a geometrically Galois cover 
$\phi_\psi:X_\psi\to \prP^1_z$. Though construction of $\phi_\psi$ is canonical, 
it 
depends on the base point. Remove this by equivalencing homomorphisms $\psi$ 
differing 
by conjugation of $\pi_1(U_{\bz},z_0)^\alg$. \S\ref{refEquiv} explains this and 
related 
equivalences. 

For $g\in G$ and $k\in\hat \bZ$ (a profinite integer), $k$ acts on $g$ by 
putting
it to the $k$th power. Then, invertible elements $n\in \hat \bZ^*$ map $g$ to 
another generator of
$\lrang{g}$. 

\begin{defn} \label{ratunion} Call a set of conjugacy classes $\bfC$ in $G$ a  
{\sl
rational union\/} if  $\bfC^n=\bfC$ (both sides counted with multiplicity) for 
all
$n\in\hat\bZ^*$. Given conjugacy classes $\bfC$, there is a natural {\sl
rationalization\/} $\bfC'$: The minimal rational collection of conjugacy classes
containing $\bfC$. \end{defn} 

\begin{exmp} Suppose $C_1$, $C_2$ and $C_3$ are respectively the conjugacy 
classes of the 5-cycles
in
$A_5$ given by $g_1=(1\,2\,3\,4\,5)$, $g_2=(1\,3\,5\,2\,4)$ and $g_1$ again. 
Then, $C_1,C_2,
C_3$ is not a rational union because the conjugacy class of $g_1$ appears with 
multiplicity 2,
while its square appears only with multiplicity 1. The collection 
$\bfC'=(C_1,C_2,C_1,C_2)$ is its
rationalization. \end{exmp}

The image $\bg$ of $\bar \bg$ gives conjugacy classes $\bfC$ in $G$, with $\C_i$ 
attached to 
$z_i$, $i=1,\dots,r$.  Suppose the following: \begin{triv} \label{defCond} 
$\phi_\psi:X_\psi\to \prP^1_z$ (with its automorphisms) has equations over 
$\bQ$.\end{triv}

\begin{lem}[Branch Cycle Lemma] \label{bcl} If {\rm \eqref{defCond}} holds then 
$\bz$ 
is a $\bQ$ set. A necessary and sufficient condition for {\rm \eqref{defCond}} 
is that 
$\psi$ extends to a homomorphism $\psi:\pi_1(U_{\bz,z_0})^\ari\to G$ (as in 
\eqref{exactArithBP}).
For an analog  with any $K$ of characteristic 0 (including $\bR$ or $\bQ_p$) 
replace $G_\bQ$ with 
$G_K$ in the definition of $\pi_1(U_{\bz',z_0})^\ari$.

Assume $\sigma\in G_K$  maps to $n_\sigma\in \hat \bZ^*=G(\bQ^\cyc/\bQ)$ and 
$\pi\in S_r$ satisfies $z_i^{\sigma}=z_{(i)\pi}$. Then,  {\rm \eqref{defCond}} 
implies 
\begin{equation} \label{conjPerm}  C_{(i)\pi}^{n_\sigma}=C_i,\ 
i=1,\dots,r.\end{equation}

So, when $K=\bQ$ and \eqref{defCond} holds, $\bfC$ is a  rational union 
(Def.~\ref{ratunion}).  
\end{lem}

Prop.~\ref{3355obst} illustrates using Lem.~\ref{bcl}, to answer question 
\eql{compA5p2r4}{compA5p2r4a}.

\subsubsection{Using Lem.~\ref{bcl} over a general field} \label{nongal} 
Suppose $\phi:X\to \prP^1_z$ over $\bQ$ has 
degree $n$ with attached 
$(G,\bfC)$. Then, $\phi$ produces a permutation representation $G\le S_n$. 
Sometimes 
the goal is to find $(X,\phi)$ over $\bQ$, with less concern for the definition 
field of
the Galois closure of the  cover (as in \cite{DFrIntSpec} and 
\cite{Fr-Schconf}).  The
necessary condition  for this requires a choice of Galois closure group $\hat G$ 
between $G$ and 
$N_{S_n}(G,\bfC)$ as in \S\ref{refEquiv}.  For each such $\hat G$, replace 
\eqref{conjPerm} by \begin{equation} \label{conjPermNG} g_\sigma 
C_{(i)\pi}^{n_\sigma}g_\sigma^{-1}=C_i,\ i=1,\dots,r, \text{for some 
}g_\sigma\in \hat 
G.\end{equation}

A general \BCL, including for non-Galois 
covers is in \cite[Thm.~5.1]{FrHFGG}. 
$K=\bQ_p$, it pays to stipulate orbits of $G_K$ on the branch points precisely. 
For a
general $K$, assume $G_K$ acts on the support of $\bz$ as $\mu:G_K\to S_r$.
Regard $G(K^\cyc/K)$ as a subgroup of $\hat
\bZ^*$.  Consider:
$$ \begin{array} {rl} H_{G,\bfC,\hat G,\mu}=\{\sigma\in G_K \mid g_\sigma 
C_{(i)\mu(\sigma)}^{n_\sigma}g_\sigma^{-1}=& C_i,\ i=1,\dots,r, \text{for some 
}g_\sigma\in
\hat  G \}; \\
\text{ and }  H_{G,\bfC,\hat G}=\{n'\in G(K^\cyc/K) &\mid g_{n'} 
C_{(i)\pi}^{n'}g_{n'}^{-1}=C_i,\ i=1,\dots,r, \\
&\text{for some }g_{n'}\in
\hat  G \text{ and }\pi\in S_r\}.\end{array}$$ Let $\phi:X\to\prP^1_z$ be in
the Nielsen class $\ni(G,\bfC,T)$ have definition field $L\ge K$. Call $\phi$ 
{\sl
attached\/} to
$(K,\bz,G,\hat G,\bfC,\mu)$ if its branch points $\bz$ are a $K$ set with the 
given
action $\mu$ and its arithmetic monodromy group over $L$ is a subgroup of $\hat 
G$. 
Then, the
fixed field $K_{G,\bfC,G,\hat G,\mu}$ of
$H_{G,\bfC,\hat G,\mu}$ is a subfield of $L$ (similar to
\cite[\S1.2]{DFrRRCF}). When $G=\hat G$ drop the appearance of $\hat G$. 
Example:
$K_{G,\bfC,G,\hat G,\mu}$ becomes $K_{G,\bfC,G,\mu}$.

\subsection{Frattini group covers} \label{frattGps} The history of the \IGP\ 
shows it is
hard to find $r$ and $\bz$ producing
$\bQ$  realizations.  To systematically investigate  the diophantine nature of 
these
difficulties we take large  (maximal Frattini) quotients of
$\pi_1(U_{\bz,z_0})^\alg$, rather than finite  groups. 

Start with a finite group $G$. Call  
$\mu: H\twoheadrightarrow G$ a {\sl Frattini\/} cover if for any subgroup 
$H^*\le H$, $\mu(H^*)=G$
implies $H^*=H$. This exactly translates the cover property from 
\S\ref{frattiniPre}.
There is a (uni){\sl versal\/} {\sl  Frattini group\/}  $\tilde G$ covering $G$ 
\cite[Chap.~20]{FrJ}. We review properties of $\tilde G$ by analogy with the 
universal Frattini
cover of the  dihedral group $D_p$ ($p$ odd). 

\subsubsection{Using $G=D_p = \bZ/p\xs \{\pm 1\}$ as a model} \label{fratAnalog} 
Consider 
$$\phi_p: {}_p\tilde D_p=\bZ_p\xs \{\pm 1\}\to D_p \text{ and } \phi_2: 
{}_2\tilde 
D_p=\bZ/p\xs \bZ_2\to D_p$$ as pieces of the universal Frattini profinite cover 
of $D_p$. Patch 
these as a fiber product over $D_p$: $$\tilde D_p=\{(g_1,g_2)\in {}_2\tilde 
D_p\times {}_p\tilde D_p \mid \phi_2(g_1)=\phi_p(g_2)\}.$$ This generalizes: For 
each 
prime $p$, $p\,|\, |G|$, there is a universal $p$-Frattini cover $\psi_p:\tG 
p\to 
G$ with these properties \cite[Part II]{FrMT}.   \begin{edesc} \label{tGp-props} 
\item 
$\tilde G$ is the fiber product of $\tG p$ (over $G$) over $p$ primes dividing 
$|G|$. 
\item \label{tGp-propsb} Both $\ker(\psi_p)$  and a 
$p$-Sylow of $\tG p$ are pro-free pro-$p$ groups, and $\tG p$ is the minimal 
profinite cover of $G$ with this
property. 
\item
\label{tGp-propsc}
$\tG p$  has a characteristic sequence of finite quotients 
$\{G_k\}_{k=0}^\infty$. \item \label{tGp-propsd} Each
$p'$-conjugacy class of $G$ lifts uniquely to a $p'$ class of 
$\tG p$. \item If $H\le G_k$ has image in $G$ all of $G$, then $H=G_k$. 
\end{edesc} 

To simplify notation we suppress the appearance of $p$ in forming the 
characteristic sequence $\{G_k\}_{k=0}^\infty$. Denote $\ker(\psi_p)$ by 
$\ker_0$. For 
any pro-$p$ group, the Frattini subgroup is the closed subgroup that  
commutators and 
$p$th powers generate. Let $\ker_1$ be the Frattini subgroup of $\ker_0$. 
Continue 
inductively to form $\ker_k$ as the Frattini subgroup of $\ker_{k-1}$. Then, 
$G_k=\tG 
p/\ker_k$.  Use of modular representation theory throughout this paper is from 
the
action of $G_k$ on $\ker_k/\ker_{k+1}\eqdef M_k$, a natural $\bZ/p[G_k]$ module. 
Rem.~\ref{extH} notes some points about this.  

We occasionally use $G_K$ with $K$ a field for the absolute Galois group of
$K$. The context and use of capital letters for fields in that rare event should 
cause no
confusion with this $G_k$ notation. 

\subsubsection{$p$-projectivity} \label{pproj} 
Let $\psi: H\to G$ be a cover of profinite groups with the kernel of $\psi$ a 
pro-$p$ group. Such a cover
is {\sl $p$-projective\/} if for any cover  $\psi': H'\to G$ with kernel a pro-
$p$ group, there is a
homomorphism $\alpha:H\to H'$ with $\psi'\circ \alpha=\psi$. 

\begin{defn} With the same setup as above and $k$ a positive integer, we say 
$\psi$ is a $(p,k)$-projective 
cover if for any $\psi'$ with kernel a $p$-group of exponent $p^k$, there exists 
$\alpha:H\to H'$ with 
$\psi'\circ
\alpha=\psi$. \end{defn} 

The next proposition captures the most significant property in the list
\eqref{tGp-props}.  Equivalent to \eql{tGp-props}{tGp-propsc} is that $\tG p$ is 
$p$-projective. 

\begin{prop} \label{pprojchar} The group $\tG p$ is the minimal $p$-projective 
cover of $G$. Also, $G_k$
(with its morphism to $G=G_0$ is the minimal $(p,k)$-projective cover of $G$. 
\end{prop} 

An analogous definition without reference to covers is of $p$-projective. That 
would be a profinite group $H$
such that any cover $\psi:H\to G$ is $p$-projective. Many fields relevant to 
arithmetic
and algebraic geometry have absolute Galois groups that are $p$-projective (call 
the field itself
$p$-projective). In particular that includes function fields over algebraically 
closed fields, and the PAC
fields from \cite{FrJ} that play a role in the exact sequence for $G_\bQ$ in 
\eqref{startGK}. 

Further, many
function fields have geometric fundamental groups that are $p$-projective. If 
$C$ is an algebraic
curve (affine or projective) over an algebraically closed field of 
characteristic $p$, then its
(profinite) fundamental group is $p$-projective. This
$p$-projective result is due to Grothendieck, though there are many recent 
revisitings of this topic
 (for example \cite{JarProj}). We know fundamental groups of projective curves 
in 0
characteristic, though not in positive characteristic; $p$-projectivity is a 
very good hint why not. The
$p'$ quotient of such fundamental groups (affine or projective) are like their 
characteristic 0 counterparts from
Grothendieck's enhancement of Riemann's existence theorem (see 
\S\ref{grassman}). Non $p'$ quotients are not, though we
know the group is
$p$ projective.  Clearly the group is not pro-free for you can't put those two 
properties together in a pro-free group.
We don't know how to characterize profinite groups that do put those properties 
together. 
\subsubsection{Is $\ker_0\,$ characteristic?}  General automorphisms of $\tG 
p$ may not
preserve $\ker_0$. That is, the phrase
{\sl characteristic  quotients\/} of $\tG p$ doesn't imply $\ker_0$ is a
characteristic  subgroup of $\tG p$. Still, the following easy lemma often 
applies. 

\begin{lem} \label{charker0} If $\ker_0$ is a characteristic subgroup of the 
$p$-Sylow of $\tG p$, then by definition all  
$\ker_k\,$s are characteristic subgroups. This holds for the 
universal $p$-Frattini cover of $D_p$. If all $p$-Sylows of $G$ intersect in  
$\{1\}$ (as when $G$ is a simple group), then $\ker_0$ is the intersection of 
all
$p$-Sylows of $\tG p$. So, $\ker_0$ is also characteristic: All automorphisms
of $\tG  p$ preserve the intersection of the $p$-Sylows. 
\end{lem} 

\begin{rem}[Checking when $\ker_0$ is characteristic] Suppose $P'_p$ is a $p$-Sylow of $\tG p$.
Since it is a pro-free pro-$p$ group, its  characteristic subgroups appear from 
the filtrations by the
$p^u$th power  subgroups and the lower central series. In principle this allows 
determining if
$\ker_0$ is a  characteristic subgroup of $P'_p$. \end{rem}

\subsection{One prime at a time} The fiber product characterization of 
\eqref{tGp-props} 
allows dealing with one $p$-Frattini cover $\tG p$ of $G$ at a time.

\subsubsection{Realization of $\tG p$ quotients} Fix $G$ and $p$.  Here is a 
naive 
diophantine goal referencing the groups $G_k$ of \eql{tGp-props}{tGp-propsc}.

\begin{prob} \label{homGk} For each $k$ find the following: \begin{edesc} \item 
$r_k$ 
distinct points $\bz_k$ (possibly varying with $k$);  and \item  
$\pi_1(U_{\bz_k},z_0)^\alg\to G_k\to 1$  factoring through 
$\pi_1(U_{\bz_k},z_0)^\ari\to G_k$. \end{edesc} The \BCL\ limits conjugacy 
classes 
$\bfC_k$ to satisfy \eqref{conjPerm}.  \end{prob}

Assume  $\bfC$ consists of $p'$-conjugacy classes. Suppose $G^*$ is 
simultaneously a quotient of
${}_p\tilde G$ factoring through $G$ and of 
$\pi_1(U_{\bz_k},z_0)^\alg$. Call the corresponding cover of $\prP^1_z$ a {\sl 
$p$-Frattini
$(G,\bfC)$ cover\/} if classical generators (\S\ref{classGens}) of 
$\pi_1(U_{\bz_k},z_0)^\alg$ also
map to
$\bfC$.  According to  
\eql{tGp-props}{tGp-propsc}, if $(G,\bfC)$ passes the \BCL\ test, then so does 
$(G_k,\bfC)$ for  {\sl all\/} 
$k$. This illustrates the groups $G_k$ are  similar. The guiding question asks 
this. 
\begin{quest} \label{GkSimilar} Will all the $G_k$ fall to the \IGP\ 
with a bound, independent of  $k$, on $r_k$ in Prob.~\ref{homGk}? \end{quest} 

\subsubsection{Possible $(\bz,r)$ and the Main Conjecture} \label{usebcl} A {\sl 
Yes!\/}  to 
Question \ref{GkSimilar} is equivalent to the following (Thm.~\ref{thm-rbound}; 
also 
with any number field $K$ replacing $\bQ$).  \begin{triv} \label{prodbfC} There 
is a rational set
of 
$p'$ classes $\bfC$ (with cardinality $r$ for some $r$) so  \eqref{homGk} holds 
with
$\bfC=\bfC_k$  (and $r_k=r$) for all $k\ge 0$. \end{triv} \noindent So,  
with no loss, $\bz_k$ has exactly $r$ points in its support for all $k$. 

\newcommand{\sB}{{\tsp {B}}} 
The result is exactly the same if we replace $\bQ$ by any number field $K$, 
which we now do. 
As in \cite[Main 
Conj.~0.1]{FrKMTIG} (or Prob.~\ref{MPMT}), if $G$ is
$p$-perfect and centerless we conjecture the answer is {\sl No!\/}.  
The Hurwitz spaces $\sH(G_k,\bfC)^\inn$ for  the inner $(G,\bfC)$ Modular Tower 
of covers up
to inner equivalence (\S\ref{innEquiv}) generalize the modular curves
$X_1(p^{k+1})$. Assume there exists $r'$ with $r_k\le r'$ for all
positive integers $k$, and let $\bfC$ be the $p'$ conjugacy class in 
\eqref{prodbfC}. Then,
Prop.~\ref{fineMod} shows the $p$-perfect and centerless hypothesis on $G_0=G$ 
guarantees 
\eqref{homGk} is equivalent to producing $K$ points on each level 
$\sH(G_k,\bfC)^\inn$ of the
inner $(G,\bfC)$ Modular Tower. Given $(G,p)$ with $p\,|\, |G|$, denote the set 
of integers $r$ 
where there is an $r$-tuple of $p'$ conjugacy classes in $G$ with a $\bQ$ point 
on every
level of the $(G,\bfC,p)$ Modular Tower by $\sB_{G,p,K,\inn}$. The goal of the 
Main Conjecture is 
is to show $\sB_{G,p,K,\inn}$ is empty. 

The methods of the paper pretty much restrict to
showing $\sB_{G,p,K,\inn}$ does not contain $r=3$ or 4. When $r=4$, $p=2$ and 
$G=A_5$, we 
succeed in showing this, and in finding considerable about excluding $r=4$ for 
general $(G,p)$. 
\S\ref{indSteps} states the unsolved problems left from \S\ref{ellipFix} in 
excluding $r=4$
from $\sB_{G,p,K,\inn}$. One sees there the effect of increasing
components in going to higher levels of a Modular Tower. This justifies 
concentration 
on analyzing the component structures at level 1 in our main examples.  A rubric 
for including all Modular Tower levels
would regard $\bfC$ as $\tG p$ conjugacy classes. This produces the Nielsen 
class $\ni(\tG  p,\bfC)$. 

\subsubsection{Relating the Main Conjecture to present diophantine literature} 
For Modular Towers defined by
absolute equivalence, rational point conclusions at high levels are likely 
weaker.  For absolute equivalence there is
an analog set
$\sB_{G,p,K,\abs}$. As with the modular curves
$X_0(p^{k+1})$, the case $G=D_p$, many times the answer may also be {\sl No 
rational points at
high levels\/}! To get closer to a problem in 
\cite{FKVo}, let $U\le \ker_0$ be a normal subgroup of $\tG p$. An example would 
be $U=(\ker_0,\ker_0)$ as in
\S\ref{grassman}.  Let $G_{k,U}$ be $\tG p/\ker_k\cdot U$. Then, instead of the 
full Modular Tower (for absolute or
inner equivalence), consider
$\{\sH(G_{k,U},\bfC)^\abs\}_{k=0}^\infty$ (or the inner case). We regard this as 
a quotient Modular Tower. Similarly,
consider values $r$ in $\sB_{G,p,U,K,\abs}$ (or $\sB_{G,p,U,K,\inn})$. The 
phrase Thompson tuples is from
\cite[Thm.~4.16]{FKVo}. 

\begin{exmp}[Thompson tuples and challenges to the Main Conjecture] 
\label{ThomTuples} Consider the
groups $\PGL_n(p)$ with
$p$ a prime. To find regular realizations of this group, take $r=n+1$, and what 
V\"olklein calls {\sl Thompson
triples\/}. These produce regular realizations with $r$ branch points, where the 
attendant Hurwitz space is a nearly 
trivial cover of $\prP^r\setminus D_r$. This is surprising, generalizing {\sl 
Belyi tuples}. 
It works only for Chevalley groups using extremely special conjugacy classes.  
The concept that gives the
realizations they call {\sl weak rigidity}.

\cite[Thm.~4.5]{FKVo} describes the data that displays Thompson tuples in 
$\GL_n(q)$. Let $\bfC$ be the
corresponding conjugacy class for these, and let $N_\bfC$ be the least common 
multiple of orders of element in
$\bfC$. They consider such Thompson tuples where $n$ and $N_\bfC=N$ are fixed, 
yet there is an infinite set $S_{n,N}$ of
corresponding primes $p$ (\cite[Lem.~4.7]{FKVo}; the $\bfC$ changes with $p$, 
though we suppress the extra notation).
The following points are no contradiction to our Main Conjecture. Still, they 
reflect on it with $r=n+1$. 
\begin{edesc} \label{chalMC} \item \label{chalMCa} There is a value of $\bz_0\in 
\prP^r\setminus D_r$ and a finite
extension $K_{n,N}$ of $\bQ$ with
$\PGL_n(p)$ regular realizations over $K_{n,N}$ for all $p\in S_{n,N}$. 
\item There is also a cyclic cover $\phi: Y\to \prP^1_z$, over $K_{n,N}$, 
ramified at $\bz_0$, so pullback of the
regular realizations of \eqref{chalMCa} over $Y$ are unramified Galois (over 
$K_{n,N}$) covers of $Y$.  
\end{edesc} 
Consequence: Using that special point $\bz_0\in \prP^r\setminus D_r$, 
\cite[Thm.~4.8]{FKVo} creates a curve $Y$ over a
number field whose geometric fundamental group has an infinite quotient that is 
a projective limit of Galois unramified
covers over the number field. One of the authors suggested to replace  
$\PGL_n(p)$, $p$ varying,  with $\PGL_n(p^t)$
(for suitable
$p$) with
$t$ varying. Analogous to the above, there would be $S_{n,N,p}$ consisting of 
values of $t$ with corresponding
realizations of $\PGL_n(p^t)$. It appears an argument similar to \cite{FKVo}, 
would then produce $\bfC$, a collection
of $p'$ conjugacy classes, a number field
$K_{n,N,p}$ and a point
$\bz_0\in \prP^r\setminus D_r$ with the following properties. 
\begin{edesc} \label{projSysPGL} \item There is a projective system of covers 
$$\cdots \to \hat X_{k+1}\to \hat
X_{k}\to
\cdots
\to
\prP^1_z:$$ $\hat X_0\to \prP^1_z$ has branch
points
$\bz_0$ and $\hat X_{k+1}\to \hat X_{k}$ is unramified, $k\ge 0$.
\item A  point of $\sH(\PGL_n(p^{k+1}),\bfC)^\inn(K_{n,N,p})$, $k\ge 0$, gives 
$\hat X_{k+1}\to \prP^1_z$.
\end{edesc} 
\cite{VoBC} offers a systematic approach through the {\sl BC Functor\/} to 
discussing projective systems of Chevalley
groups, and so one may consider the whole apparatus as a challenge to the Main 
Conjecture. We conclude, however, by
noting why
\eqref{projSysPGL} does not contradict it. It is true that $\PGL_n(p^{k+1})\to 
\PGL_n(p)$ is a $p$-Frattini cover, if
$n\ge2$, $p$ is odd,
$(p,n)=1$ and if $p=3$, $n>2$ \cite[\S4]{VCycCov}. The group $\PGL_n(\bZ_p)$, 
however, is not the universal $p$-Frattini
cover of $\PGL_n(p)$: The kernel to $\PGL_n(\bZ_p)\to \PGL_n(p)$ is not a pro-
free pro-$p$ group
as is the kernel of ${}_p\widetilde {\PGL_n}(p)\to \PGL_n(p)$. Even the rank of 
the latter is larger than that of
$\ker(\PGL_n(\bZ_p)\to \PGL_n(p))$. Example: $\ker(\PGL_2(\bZ_5)\to \PGL_2(5))$ 
has rank 3, while the
kernel of ${}_5\widetilde {\PGL_2}(5)\to \PGL_n(5)$ has rank 6 
\cite[Rem.~2.10]{FrMT}. 
 \end{exmp} 

\subsection{Refined s-equivalence from $N_{S_n}(G,\bfC)$} \label{refEquiv} Make 
a 
choice
$\bar\bg$  of classical  generators for $\pi_1(U_{\bz},z_0)$ 
(\S\ref{classGens}). Let
$T:G\to S_n$ be a permutation representation. Any surjective homomorphism 
$\psi:\pi_1(U_{\bz},z_0)\to G$ produces a degree $\deg(T)$ cover  $X_\psi\to 
\prP^1_z$  
canonically.  Points of $X_\psi$ are homotopy classes of paths $\gamma:[0,1]\to
U_{\bz}$  with $\gamma(0)=z_0$  modulo this relation: $\gamma':[0,1]\to 
U_{\bz}$ 
and $\gamma$ are equivalent if \begin{triv} 
 $ \gamma(1)=\gamma'(1)$ and $\psi(\gamma'\gamma^{- 
1})\in G(T,1)$.  \end{triv} \noindent To be canonical requires the covering 
space have a 
base 
point. 
For this take the constant path from $z_0$.

\subsubsection{Interpreting equivalence classes of Galois covers} 
\label{intEquiv} To 
avoid running into a base point $z_0$, deforming $\bz$ requires moving $z_0$. So 
the 
cover is independent of $z_0$ only up to inner automorphism.  As in 
\S\ref{setupNC}, 
consider
any group $N'$ between $G$ and $N_T(G,\bfC)$. We now explain different 
Hurwitz spaces as  equivalence classes of covers corresponding to points
$\bp\in \sH_{N'}$. The results we quote in \S\ref{innEquiv} and \S\ref{absEquiv} 
are 
refinements of \cite{FrVMS} suitable for this paper. Throughout we equivalence
only $(G,\bfC)$ covers, with $G$ and $\bfC$ fixed. 

\subsubsection{Case: $N'=G$ and inner equivalence} \label{innEquiv} 
The structure on $\sH_G$ is the moduli space of inner equivalence classes, 
denoted 
$\sH^\inn$.  Interpret
$\bp\in \sH_G=\sH^\inn$ as an  equivalence class of pairs: \begin{equation}
\label{Gpairs} (\phi:X\to\prP^1_z,\  G(X/\prP^1_z)\mapright{\alpha}G),
\end{equation} with
$G(X/\prP^1_z )$ the  automorphisms of a Galois cover and $\alpha$ a group 
isomorphism.
Then, 
$(\phi':X'\to\prP^1_z,\ G(X'/\prP^1_z) \mapright{\alpha'}G)$ is equivalent to 
\eqref{Gpairs} if a continuous one-one map $\psi:X\to X'$ induces the latter 
from the former. 
Example of an equivalent pair: Let $X=X'$ and compose $\alpha$ with conjugation 
by $g\in G$ to get $\alpha'$.  
Composing 
$\alpha$, however, with an outer automorphism of $G$ gives a new equivalence 
class of 
pairs.

\subsubsection{Further explanation of automorphisms} \label{exAutos} For $\bg\in 
\ni(G,\bfC)/G$, automorphisms of $X/\prP^1_z$ identify with the centralizer of 
$G$ in 
$N_R(G,\bfC)$ (Lem.~\ref{autos}) with $R$ the regular representation. A map 
between
two pairs as in
\eqref{Gpairs} is  unique if and only if $G$ has no center. This is equivalent 
to
there being a unique total  family $\sT_{G,\bfC}^\inn=\sT^\inn\to \sH^\inn\times
\prP^1_z$. For
$\bp \in \sH^\inn $,  points $\sT^\inn_\bp$ over $\bp\times \prP^1_z$ form
a Galois cover of 
$\prP^1_z$ with group $G$ representing $\bp$.  This makes $\sH^\inn$ a {\sl 
fine\/} 
moduli space. Then, the (minimal) definition field of $\sT^\inn_\bp\to \prP^1_z$
and the automorphisms given by $G$ is 
$\bQ(\bp)$, generated over $\bQ$ by coordinates of $\bp$ \cite{FrVMS}. The total
family  may exist even if $G$ has a center, though it won't be unique. 

\begin{lem} \label{charblem} Even if $G$ has a center, for any point $\bp\in 
\sH^\inn $, there is a
cover $\phi_\bp: Y_\bp\to \prP^1_z$ with definition field $\bQ(\bp)$. It may 
not,
however, be Galois over $\bQ(\bp)$. 

More generally, suppose $\phi: Y\to X$ is a Galois cover of
nonsingular projective curves with $X$ over $K$ and $x_0\in X(K)$ is unramified 
in $Y$. Assume
$\phi^\sigma: Y^\sigma\to X$ is s-equivalent to $\phi$ for each $\sigma\in G_K$. 
Then there
exists $\phi':Y'\to X$ over $K$, s-equivalent to $\phi$ (over $\bar K$) with a 
rational point
over
$x_0$. 
\end{lem}

\begin{proof} A version of this is in \cite{CoombesHarb}. Consider any Galois 
cover $\phi:Y\to
\prP^1_z$ (over some algebraic closure $\bar K$ of
$K=\bQ(\bp)$) in the equivalence class of
$\bp$. Then, choose any $y\in Y$ over $z_0\in \bQ$ unramified in $Y$. From the 
moduli
property, for
$\sigma\in G_K$, there is an isomorphism
$\psi_\sigma: Y\to Y^\sigma$ commuting with the maps to $\prP^1_z$.
Compose such a map with the unique automorphism of the Galois cover 
assuring $\psi_\sigma$ takes $y$ to $y_\sigma$. 

Apply Weil's cocycle condition to
$(Y,y)$ for an equivalent pair $(\phi':Y'\to \prP^1_z,y')$ over
$\bQ(\bp)$ \cite{WeilField}. We may, however, lose the automorphisms: 
$\phi':Y'\to \prP^1_z$ defines a cover in $\sH_{N_R(G,\bfC)}$ where $R$ is the
regular representation of $G$. The cover, however, is special, for it has a 
rational
point $y'$ over $z_0$. 

The proof works with a general curve replacing $\prP^1_z$. 
With $\afA^r\setminus D_r$ replacing $U_r$  and 
$z_0=\infty$, this construction works uniformly to give a fine moduli space of 
geometrically
Galois covers with a point over $\infty$. 
\end{proof}

\subsubsection{Absolute equivalence:  $N'=N_{S_n}(G,\bfC)$} 
\label{absEquiv} For $G\le N'\le
N_{S_n}(G,\bfC)$, $H_r$ action
on  $\ni(G,\bfC)/N'$ gives $\sH_{N'}$ (as in \S\ref{setupNC}). Refer back to
\S\ref{sequiv} for the notation for Galois closure of a cover. Given $\phi:X\to
\prP^1_z$, we chose $X^\dagger$ to be a geometric Galois closure of $\phi$. This 
depended on 
choosing a coset of 
$G$ in $N_{S_n}(G,\bfC)$. A particular choice determined an isomorphism of $G$
with $\Aut(X^\dagger/\prP^1_z)$. 

A point $\bp\in \sH_{N'}$ corresponds to a cover $\phi:X\to \prP^1$ up to a 
choice
of $X^\dagger$ determined by a coset  of $N'$ in $N_{S_n}(G,\bfC)$.  Formally:  
$(\phi':Y\to\prP^1_z,\ G(Y^\dagger/\prP^1_z) \mapright{\beta}G)$ is $N'$-
equivalent
to the corresponding expression for $(X,X^\dagger)$ if some continuous one-one 
map $\psi:X\to
Y$ induces the latter up to conjugation of $\beta$ by $N'$ from the former. This
generalizes inner equivalence, the case $N'=G$. 

When $N'=N_{S_n}(G,\bfC)$ denote $\sH(G,\bfC)_{N'}$ by $\sH(G,\bfC)^\abs$, 
absolute
equivalence classes of covers $X\to\prP^1_z$ with associated  permutation 
representation 
$T$. From \S\ref{sequiv} this requires no choice in $X^\dagger$.  So, two covers 
are
$N_{S_n}(G,\bfC)$-equivalent if there is a map between them commuting with their
maps to $\prP^1_z$. 

Assume $G$ has no center, so 
$\sH^\inn(G,\bfC)$ has a
unique total representing family. Then, $H_r$ acting on $\ni/G$ and 
$\ni/N'$
produces  $\Psi_{N'}^\inn: \sH^\inn(G,\bfC)=\sH^\inn\to \sH_{N'}$.
Rem.~\ref{nongal} gives the cyclotomic field $K_{G,\bfC,N'}$ 
(resp.~$K_{G,\bfC}$) as
the definition field of $\sH_{N'}$ (resp.~$\sH^\inn$ and $\Psi_{N'}^\inn$). 

If $\bfC$ is a rational union of conjugacy classes, this field is $\bQ$. The
following  interprets the main technical result of \cite{FrVMS}. Recall previous 
notation for
the fibers  of a family $\Phi:\sT\to \sH\times \prP^1_z$ (or $\pr_2\circ\Phi: 
\sT\to \sH$).
If $\bp\in 
\sH$, then $\sT_\bp$ is the set of points of $\sT$ 
over $\bp\times\prP^1_z$ and $\phi_\bp:\sT_\bp\to \prP^1_z$ is restriction of 
$\Phi$.

\begin{thm}\label{FrVMS} Suppose $G$ has no nontrivial centralizer in $N'$. Then 
there 
are total representing families $\Phi_{N'}:\sT_{N'}\to \sH_{N'}\times \prP^1$ 
and 
$\Phi^\inn:\sT^\inn\to \sH^\inn\times \prP^1$. For $\bp\in \sH_{N'}$ and $\hat 
\bp\in 
\sH^\inn$ over $\bp$, the covers $\phi^\inn_{\hat \bp}:\sT^\inn_{\hat \bp}\to 
\prP^1_z$ and $\phi_{N',\bp}:\sT_{N', \bp}\to \prP^1_z$ corresponding to these 
points 
have the following properties. 
\begin{edesc} \item $\phi^\inn_{\hat \bp}$ (resp.~$\phi_{N', \bp}$) has field of
definition 
$K_{G,\bfC}(\hat \bp)$ (resp.~$K_{G,\bfC,N'}(\bp)$). 
\item $\phi^\inn_{\hat \bp}$ is an absolutely irreducible component of a 
$K_{G,\bfC,N'}(\bp)$ (arithmetic) component of the Galois closure of $\phi_{N', 
\bp}$ 
(via the \S\ref{sequiv} construction). \item $K_{G,\bfC}(\hat 
\bp)/K_{G,\bfC,N'}(\bp)$ 
is a Galois extension with group naturally isomorphic to a subgroup of $N'/G$.
\end{edesc}
\end{thm}

\newcommand{\icona}{\otimes}
\newcommand{\iconb}{\oplus}
\newcommand{\iconc}{}
\newcommand{\icon}{{\otimes \oplus}} 

\begin{exmp}[$A_n$ and 3-cycles] \label{full3cycleList} Let $T_n$ be the 
(standard) representation of $A_n$,
$n\ge 5$. 
\cite[Thm.~1]{FrLInv} lists the complete set of 
inner and absolute Hurwitz space components at level 0 of the $(A_n,\bfC_{3^r})$ 
Modular Tower. Table \ref{An3rConst} displays these for inner spaces (the result 
is 
nontrivially the same for absolute spaces). All these components have definition 
field 
$\bQ$. Locations in this diagram have an attached integer pair
$(n,r)$. In each case the $\iconb$ inner component $\sH^{\inn,\iconb}_{n,r}$ 
maps to the 
absolute $\iconb$ component $\sH^{\abs,\iconb}_{n,r}$ by a degree 2 (Galois) map 
with 
group identified with $S_n/A_n$. Similarly for the corresponding $\icona$ 
components. 
In Table
\ref{An3rConst}, the notation for 
components corresponds to lifting invariant values as in
Prop.~\ref{serLift} (or \S\ref{NielSep}, specifically in
\eql{exLiftInv}{exLiftInva}). The genus at $(n,r)$ of a degree
$n$ cover is $g=r-n+1$.  

\begin{table}[h]
\caption{Constellation of spaces
$\sH(A_n,\bfC_{3^r})$} \label{An3rConst}
\begin{tabular}{|c|c|c|c|c|c|c|} \hline 
&&&&&&\\ 
$\mapright{g\ge 1}$& $\icon$&$\icon$&
\dots & $\icon$& $\icon$ &$\mapleft{1\le g}$ 
\\ \hline &&&&&&\\ 
$\mapright{g=0}$ & 
$\icona$ & $\iconb$ & \dots & $\icona$ & $\iconb$ & $\mapleft{0=g}$ 
\\ \hline &&&&&&\\
$n\ge 4$ & $n=4$ & $n=5$ & \dots & $n$ even & $n$ odd & $4\le n$\\
\hline
\end{tabular} \end{table}

\end{exmp}

\subsection{$p$-perfect groups and fine moduli} Start with any equivalence 
between covers of 
$\prP^1_z$. The Hurwitz space representing these equivalences is a fine moduli 
space if it 
has a unique total family representing the equivalence classes of its points. 
Prop.~\ref{fineMod}  
considers only Hurwitz spaces of type $\sH(G,\bfC)^\abs$ or $\sH(G,\bfC)^\inn$.

\subsubsection{$p$-perfect groups} \label{pperfSec} For any $G$ module $A$,
consider a group extension $\phi: \hat G\to G$ with $\ker(\phi)=A$ and the 
lifted 
conjugation action of $G$ on $A$  is that given. Then, $H^1(G,A)$ corresponds to
automorphisms of $\hat G$ trivial on $G$ and $A$ modulo automorphisms induced by 
conjugation by $A$
\cite[p.~239]{NCHA}. Suppose
$A=\bZ/p$, with $G$ acting trivially. Then, $H^1(G,A)$ is just the homomorphisms 
of $G$
into $A$. 

\begin{defn}[$p$-perfect groups] \label{pperfect} For $p$ a prime, a group $G$ 
is 
$p$-perfect if it has no $\bZ/p$ quotient. That is, $H^1(G,\bZ/p)$ is trivial. 
\end{defn} 
Let $U_p$ be the $p$ part of the Schur multiplier of $G$ (this may be trivial). 
That $G$ is $p$-perfect
interprets as $G$ having a central extension ${}_p \hat G$ with this property. 
\begin{triv} $U_p=\ker({}_p \hat G\to G)$ and ${}_p \hat G$ is universal for 
central extensions of $G$ with
$p$-group kernel. \end{triv} \noindent Prop.~\ref{hatP-P} illustrates the 
necessity of this condition. 

Any finite group has a centerless cover \cite{FrVMS}. No cover, however, of
$G$ can be $p$-perfect, unless $G$ is. Here is another characterization of $p$-perfect.

\begin{lem} The $p'$
elements in $G$ generate if and only $G$ is $p$-perfect. \end{lem}

\begin{proof} Let $H$ be the (normal) subgroup of $G$ generated by its $p'$ 
elements. If $H$ is a proper subgroup
of $G$, then $G/H$ is a nontrivial $p$-group, and any $p$-group has a $\bZ/p$ 
quotient. Conversely, given $\phi:
G\to \bZ/p$, the kernel of $\phi$ contains all the $p'$ elements of $G$. 
\end{proof} 

Note: Perfect groups are exactly those 
$G$ that are 
$p$-perfect for every prime dividing $|G|$. Let $\{G_k\}_{k=0}^\infty$ be the 
characteristic
quotients of
${}_p\tilde G$, the universal $p$-Frattini  cover of $G_0=G$ and 
$M_k=\ker(G_{k+1}\to G_k)$ as in
\eqref{tGp-props}. Prop.~\ref{fineMod} uses notation from the
Loewy display of $M_k$ as a $G_k$ module (\S\ref{onesAppearNGP}).  

\begin{exmp}[Centers in quotients of $p$-perfect centerless group] 
\label{cenPerf} Suppose $G_0$ is
$p$-perfect, centerless, and has a nontrivial
$p$ part in its Schur multiplier. Let $G_1$ be the first characteristic Frattini 
cover of $G_0$
as in \eql{tGp-props}{tGp-propsc}. Then the canonical map
$G_1\to G_0$ factors through a nontrivial central extension of
$G_0$. 

Further, this automatically replicates at all levels. For all $k$, $G_k$ is $p$-perfect
and centerless (Prop.~\ref{fineMod}) and $G_{k+1}\to G_k$ factors through a {\sl 
nontrivial\/}
central extension of
$G_k$ (Prop.~\ref{RkGk}). Subexample: The universal exponent 2-Frattini 
extension $G_1$ of
$A_n=G_0$ factors through the spin cover of 
$\hat A_n$. \end{exmp} 

\subsubsection{Fine moduli for a Modular Tower} The next proposition is a 
characterization for all levels
of a Modular Tower having fine moduli. 

\begin{prop} \label{fineMod} Suppose a  Hurwitz space is of type 
$\sH(G,\bfC)^\abs$ 
with associated permutation representation $T$.  It is a fine moduli space if 
$T:G\to S_n$ 
has image with no centralizer in $S_n$ (see \S\ref{centSn}). A Hurwitz space of 
type
$\sH(G,\bfC)^\inn$ is a  fine moduli space if $G$ has no center.

Assume $G_0=G$ is a centerless $p$-perfect group. Then, for each $k$, so is 
$G_k$: $\one_{G_k}$ does not appear at the far left of the Loewy display of 
$M_k$. Let
$\bfC$ be a set of
$p'$ classes of $G$. Then, the Hurwitz spaces 
$\{\sH(G_k,\bfC)^\inn\}_{k=0}^\infty$ are all fine  moduli spaces.  \end{prop}

\begin{proof} The first part is a subset of Thm.~\ref{FrVMS}. Assume $G_0=G$ 
is centerless and $p$-perfect. We inductively show $G_k$ also has 
these properties for all $k$. 

\cite[Lem.~3.6]{FrMT}  shows 
$G_{k+1}$ is centerless if the following hold.  \begin{edesc} \label{loewy} 
\item $G_k$ 
has no center.  \item \label{loewyb} $\ker_k/\ker_{k+1}$ has  no $G_k$ 
subquotient of 
Loewy type $\one\to \one$. \end{edesc} The module in 
\eql{loewy}{loewyb} is distinct from  $\one\oplus\one$. It comes from a 
nontrivial 
representation of $G_k$ of form: $g\in G\mapsto \smatrix 1 {a_g} 0 1$. The map 
$g\in 
G_k\mapsto a_g$ is a homomorphism of $G_k$ into the $\bZ/p$. By hypothesis this 
doesn't exist.

That leaves showing $G_{k+1}$ has no quotient isomorphic to $\bZ/p$, assuming 
$G_k$ 
is centerless and has no such quotient. Suppose $\phi:G_{k+1}\to \bZ/p$ is 
surjective with 
kernel $K$. Consider the map from $K$ to $G_k$ induced by the canonical map 
$G_{k+1}\to G_k$. This is a Frattini cover. So, $K$ is not onto  
$G_k$. 
Then, $K$ has image an index $p$ normal subgroup of $G_k$. This is contrary to 
our 
assumptions.  

Finally, since $G_k$ is centerless, $G_{k+1}$ has a center if and only if $G_k$ 
stabilizes some nontrivial
element of $M_k$; if and only if $\one_{G_k}$ appears to the far left in the 
Loewy display of $M_k$). 
\end{proof}

\subsubsection{Ordering branch points} \label{ordBranchPts} We comment on a 
Hurwitz space topic that
arises in \cite{DFrVarFam}. Let $H$ be a subgroup of $S_r$. Use the diagram of
\eqref{fundgpdiag} for
$\Psi_r:(\prP^1)^r\setminus
\Delta_r 
\to \prP^r\setminus D_r$. So, $H$ defines a quotient $U^H=(\prP^1_z)^r\setminus 
\Delta_r)/H$: $\Psi_H: U^r\to U^H$ is
then the canonical map.  A component $\sH_1$ of a Hurwitz space
$\sH$ has an {\sl $H$-ordering on its branch points\/} if $\Psi: \sH'\to 
\prP^r\setminus D_r$ factors through $U^H$. 
Up to whatever equivalence defines the moduli problem for $\sH$, this means for 
any cover $\phi:\bp: X_\bp\to
\prP^1_z$ representing the equivalence class of $\bp\in \sH'$, the effect of 
$G_{\bQ{\bp}}$ on orderings
$(\row z r)$ of the branch points
$\bz_\bp$ of $\phi_\bp$ is conjugate to a subgroup of $H$. The notion  depends 
only on the
conjugacy class of $H$ in $S_r$. This topic arises naturally in considering the 
\BCL\ (see \S\ref{nongal}). 

\begin{defn}  For $H\le S_r$, and $\sH_1$ a component of a Hurwitz space $\sH$, 
let $\sH'$ be a component of the fiber
product $\sH_1\times_{\prP^r} U^H$. Such an $\sH'$ is an $H$-ordering (of the 
branch points) of $\sH_1$. 
When $H=\{1\}$ this is the traditional meaning of an ordering the branch points.
\end{defn} 

There is a simple Nielsen class interpretation for $\sH_1$ having an $H$-
ordering. Let $O$ be the $H_r$ orbit on 
$\ni(G,\bfC)/N$ (as in \S\ref{setupNC}) corresponding to $\sH_1$. We also say,  
an $H$-ordering of $O$. Recall
$\Psi_r^*: H_r\to S_r$ from
\S\ref{basFGs}.
\begin{lem} For $\bg\in O$, let
$H_\bg$ be the subgroup of $H_r$ stabilizing the equivalence class of $\bg$.  
Then, there is an
$H$-ordering of
$O$ if and only if some conjugate of $\Psi^*_r(H_\bg)\le S_r$ is in  
$H$.\end{lem}

\begin{proof} Suppose $\sH_1\to U_r$ factors through $U^r/H\to U_r$. Then a 
point $\bp\in \sH_1$ has image $\bu\in
U^r/H$.  The (geometric) decomposition group for $\bu$ in the cover $U^r\to 
U^r/H$ is a subgroup of $H$. It must
contain the image $\Psi^*_r(H_\bg)$ since $H_\bg$ is a subgroup of the 
fundamental group of $U^r/H$. The argument is
reversible. 
\end{proof} 

\begin{exmp}[Pairs of conjugacy classes] \label{conjClasspairs} Suppose $r=4$, 
$\C_1=\C_2$, $\C_3=\C_4$, $N=G$ and
$\C_1\ne
\C_3$. Then, the minimal 
 (up to conjugacy in $S_4$) group for which there is an $H$-ordering of an $H_4$ 
orbit on $\ni(G,\bfC)^\inn$ is
$H=\lrang{(1\,2),(3\,4)}$. If, however, $\C_1$ and $\C_2$ are conjugate in a 
group $N$ between $G$ and $N_{S_n}(G)$, then
the minimal $H'$-ordering of an $H_4$ orbit on $\ni(G,\bfC)/N$ is
$H'=\lrang{(1\,2),(3\,4), (1\,3)(2\,4)}$. 
\end{exmp}
  
\subsection{$\SL_2(\bC)$ action on Hurwitz spaces and covers of $\Lambda_r$ and 
$J_r$} \label{SL2-act} Use the notation from \S\ref{notation}: Denote 
$(\prP^1)^r\setminus \Delta_r$ by $U^r$ and $\prP^r\setminus D_r$ by $U_r$. Let 
$N'$ be one of the groups from \S\ref{intEquiv}. Action of $H_r$ on 
$\ni(G,\bfC)/N'$ 
produces an unramified cover $\Psi_{G,\bfC}:\sH(G,\bfC)\to U_r$ \eqref{brAct}.  
Pull 
this cover back  to $U^r=(\prP^1)^r\setminus \Delta_r$ as the fiber product 
$\sH(G,\bfC)'=\sH(G,\bfC)\times_{U_r}U^r$.

\subsubsection{Quotient by $\PSL_2(\bC)$}
Consider $\PSL_2(\bC)$ acting diagonally on $r$ copies 
of $\prP^1_z$. For $\alpha\in \PSL_2(\bC)$ and $\bz\in U^r$, $\alpha(\bz)\mapsto 
(\alpha(z_1),\dots,\alpha(z_r))$. This action commutes with 
$S_r$ permuting coordinates; put $\PSL_2(\bC)$ on the left. So, 
$\PSL_2(\bC)\bsl  U^r=\Lambda_r$ generalizes the $\lambda$-line minus 
$\{0,1,\infty\}$ from 
modular curves. Also, $\PSL_2(\bC)\bsl U_r=J_r$ generalizes the $j$-line minus 
$\infty$.  Both spaces have complex 
dimension $r-3$. 

\begin{lem} \label{presNC} For any strong equivalence of covers (from 
\S\ref{intEquiv},
including absolute or inner equivalence), composing $\phi:X\to \prP^1_z$ with 
$\alpha\in
\PSL_2(\bC)$ preserves the Nielsen class. Further, the equivalence classes of 
$\phi$ and
$\alpha\circ\phi$ lie in the same component of the corresponding Hurwitz space. 
\end{lem} 

\begin{proof} A Nielsen class is given by $(G,\bfC, T)$ with $T:G\to S_n$ a 
faithful
permutation representation. The equivalence depends on some subgroup of $S_n$, 
containing
$G$ and normalizing it. A cover in the Nielsen class is of degree $n$ and has 
$T$ the
natural permutation representation associated to it.  When, it is an inner 
class, we also
attach an isomorphism between the group of the cover and the group $G$. The 
monodromy
groups of the covers $\phi$ and $\alpha\circ\phi$ are exactly the same. If $\bz$ 
are the
branch points for $\phi$, then $\alpha(\bz)$ are the branch points of
$\alpha\circ\phi$, with $\alpha(z_i)$ having the conjugacy class $\C_i$ attached 
to it.
This shows $\alpha$ preserves Nielsen classes. 

The action of $\alpha\in \PSL_2(\bC)$ mapping on s-equivalence classes of covers 
in a Nielsen class
is continuous. Given $\phi:X\to \prP^1_z$, map $\alpha\in \PSL_2(\bC)$ to the s-
equivalence class
of $\alpha\circ\phi$. Since $\PSL_2(\bC)$ is connected, the $\PSL_2(\bC)$ orbit 
of $\phi$  lies in
one  connected component of the Hurwitz space. The orbit contains $\phi$; the 
component is that of
$\phi$. 
\end{proof} 

\begin{defn}[Reduced Hurwitz spaces] \label{redHurSpace}  
\cite[\S6.2]{DFrIntSpec} shows 
$\PSL_2(\bC)$ action extends to $\sH_{N'}$ and to $\sH^\inn$. This
produces affine schemes
$\sH_{N'}/\PSL_2(\bC)$ and $\sH^\inn/\PSL_2(\bC)$ covering (usually ramified) 
$J_r$. These are
reduced Hurwitz spaces. \end{defn}  

\subsubsection{The $j$ and $\lambda$-lines}
Take $r=4$.
Four unordered distinct points,
$\bz\in U_4$, are the branch points of a  unique degree two cover $E_\bz\to 
\prP^1_z$. With the
right choice of inhomogeneous  coordinate on $J_4$, the image of $\bz$ in $J_4$ 
is the classical
elliptic curve 
$j$-invariant. Take the elliptic curve to be degree 0 divisor classes on 
$E_\bz$: 
$\Pic^0(E_\bz)$. Identify $J_4$ with $\prP^1_j\setminus \{\infty\}$. 

Suppose $r=4$ and $\psi^\rd: \sH^\rd\to \prP^1_j\setminus \{\infty\}$ is a 
reduced Hurwitz space cover. Assume also
that  a general point $\bp\in
\sH^\rd$ corresponds to the equivalence class of a cover  $X_\bp\to \prP^1_z$ 
whose Galois closure
maps surjectively to the elliptic curve with invariant $j=\psi^\rd(\bp)$. We say 
$\sH^\rd$ is {\sl $j$-aware}. Many
Hurwitz spaces are $j$-aware. 

\begin{exmp}[Branch cycles in $S_n\setminus A_n$ and $j$-awareness] 
\label{SnAnbcycles} Assume each
class in
$\bfC$  is in $S_n\setminus A_n$.  Let $\bp\in \sH(G,\bfC)^{\inn,\rd}$ lie over 
$\bz$.  The regular representation of $G$ gives a map $G\to S_{|G|}$. The cover 
$\phi_\bp:X_\bp\to \prP^1_z$ naturally factors through $E_\bz\to\prP^1_z$:
Quotient 
$X_\bp$ by $G\cap A_{|G|}$.  (This works for any even $r$; $E_\bz$ is then 
hyperelliptic.)

Suppose $r=4$ and $\bg\in \ni(G,\bfC)$ with $G\le A_n$. Choose $h_1,h_2\in 
S_n\setminus A_n$. Then, $(h_1g_1,g_2h_2,h_2^{-1}g_3,g_4h_1^{-1})$ satisfies the 
product-one condition. It produces a Nielsen class (for some new group) with 
moduli problem 
directly recognizing the $j$-line as parameterizing elliptic curves. \end{exmp} 
In summary
we have the following 
\cite[Prop.~6.3]{DFrIntSpec}. Further remarks on $j$-awareness appear in 
\S\ref{jawareMT}.

\begin{prop} \label{PSL2-quotient} The (unramified) cover $\Psi_4: U^4\to U_4$ 
modulo 
$\PSL_2(\bC)$ produces the classical (ramified) map $\psi: 
\prP^1_\lambda\setminus 
\{0,1,\infty\}\to \prP^1_j\setminus\{\infty\}$. This extends to a (ramified) 
cover 
$\sH_{N'}/\PSL_2(\bC)\to \prP^1_j\setminus\{\infty\}$.  \end{prop}

Prop.~\ref{j-Line} gives the precise action of $\bar M_4$ on {\sl reduced 
Nielsen 
classes\/}. This produces a branch cycle description of the cover 
$\sH_{N'}/\PSL_2(\bC)\to \prP^1_j\setminus\{\infty\}$. Many  
computations 
of this paper depend on this. 

Similar to \S\ref{ordBranchPts}, Let $H$ be a subgroup of $S_3$ (the group of 
the Galois cover
$\prP^1_\lambda\to \prP^1_j$).  A component $\sH_1^\rd$ of a reduced Hurwitz 
space
$\sH^\rd$ has an {\sl $H$-ordering on its branch points\/} if $\Psi: 
\sH_1^\rd\to  \prP^1_j\setminus\{\infty\}$ factors
through $\prP^1_\lambda/H$. Up to the reduced equivalence defining the moduli 
problem for $\sH^\rd$, for
any cover
$\phi:\bp: X_\bp\to
\prP^1_z$ representing the reduced equivalence class of $\bp\in \sH'$, the 
effect of $G_{\bQ(\bp)}$ on the orderings 
$(z_1,z_2,z_3,z_4)$, of 
$\bz_\bp$ of $\phi_\bp$ is conjugate to a subgroup of $H$. \S\ref{extA4C34} is 
an application of the next lemma. 

\begin{lem} \label{hurtoredhurorder} Suppose $H\le S_4$, and for $r=4$ a 
component $\sH_1$ of a Hurwitz space $\sH$ has an
$H$-ordering (as in
\S\ref{ordBranchPts}) of its branch points. Then, the corresponding reduced 
Hurwitz space component $\sH_1^\rd$ has an
$H/(H\cap K_4)$ordering of its branch points. Particularly, if $H\le K_4$, then 
$\sH_1^\rd\to
\prP^1_j\setminus\{\infty\}$ factors through the natural map $\prP^1_\lambda\to 
\prP^1_j$.  \end{lem} 

\begin{proof} The argument of \S\ref{identK4} shows how to identify the group of 
the cover
$\prP^1_\lambda\to \prP^1_j$ with the action of $S_4$ on the points of $\bz$ 
modulo a Klein 4-group. 
\end{proof}

\begin{exmp}[Pairs of conjugacy classes continued] \label{conjClasspairscont} As 
in Ex.~\ref{conjClasspairs}, assume
$r=4$, $\C_1=\C_2$, $\C_3=\C_4$,
$N=G$, $\C_1\ne \C_3$, and  $\C_1$ and $\C_2$ are conjugate in a group $N$ 
between $G$ and $N_{S_n}(G)$. Then, the
minimal groups with the reduced Nielsen classes  $\ni(G,\bfC)^{\inn,\rd}$ and  
$\ni(G,\bfC)^\rd/N$ having an ordering
of the branch points are the same: 
$$\lrang{(1\,2),(3\,4)}/\lrang{(1\,2),(3\,4)}\cap K_4=\lrang{(1\,2),(3\,4),
(1\,3)(2\,4)}/\lrang{(1\,2),(3\,4),
(1\,3)(2\,4)}\cap K_4.$$
\end{exmp}

\subsection{Monodromy groups of $\psi:\sH\to U_r$ and $\psi^\rd: \sH^\rd\to 
J_r$}
Let $\sH^\inn$ (resp.~$\sH^{\inn,\rd}$) be a component of a Hurwitz 
(resp.~reduced Hurwitz)
space covering
$U_r$ (resp.~$J_r$). The geometric monodromy (Galois closure) groups 
$G_{\psi}\le
S_{\deg{\psi}}$ and $G_{\psi^\rd}\le
S_{\deg{\psi^\rd}}$ of these covers are invariants of the component (and of the
Nielsen class) describing this cover. 

We call attention to delicate points useful
outside the area of this paper for investigating rational points on these 
spaces. Example: 
Formula \eqref{caction} for how complex conjugation acts on the branch cycle
description of a reduced Hurwitz space may not determine $c$. Recognition, 
however, of how
to locate one H-M rep.~gives the correct determination.

\subsubsection{Paths for computing complex conjugation} \label{Qpaths} 
Apply complex conjugation directly to paths representing $\row q {r-1}$,
the braid generators of $H_r$. For example, suppose $\bz$ consists of real 
points $z_1<
z_2<\dots < z_r \in \prP^1_z(\bR)$ arranged around the real circle. As in 
\cite[App.~C and 
App.~D]{FrMT} let $B_i$ be a clockwise circle on $\prP^1_z$ with a marked 
diameter on 
the real axis having $z_i$ and $z_{i\np1}$ as endpoints. (One of these has $z_r$ 
and $z_1$ at the endpoints of the directed diameter.) Parametrize the top of 
$B_i$ 
with $t\mapsto B_i^+(t)$ on $[0,1]$, so $B_i^+(0)=z_i$ and 
$B_i^+(1)=z_{i\np1}$. Similarly, parametrize the bottom of $B_i$ with 
$t\mapsto B_i^-(t)$ on $[0,1]$ so $B_i^-(0)=z_{i\np1}$ and $B_i^-(1)=z_{i}$. 
Consider the path \begin{equation} \label{parQa} t\mapsto 
(z_1,\dots,z_{i\nm1},B_i^+(t),B_i^-
(t),z_{i\np2},\dots, z_r),\ t\in [0,1].\end{equation} The range of \eqref{parQa} 
in 
$\prP^r\setminus D_r$ represents the braid group generator $q_i^{-1}$. 
The inverse of  path \eqref{parQa} is
\begin{equation} \label{parQb} t\mapsto (z_1,\dots,z_{i\nm1},\bar B_i^-(t),\bar
B_i^+(t),z_{i\np2},\dots,  z_r),\ t\in [0,1]\end{equation} where the notation 
$\bar B$
indicates complex 
conjugation applied to the coordinate.

Apply this when $r=4$ to compute complex conjugation on the branch
cycles for reduced Hurwitz covers of the $j$-line. Use 
$$\gamma_0=q_1q_2,\ \gamma_1=q_1q_2q_1,\ \gamma_\infty=q_2$$ for paths in 
$\prP^1_j\setminus
\{0,1,\infty\}$,  images by $\PSL_2(\bC)$ reduction from the paths above.  

\begin{lem} Let $r=4$ and $j_0\in \prP^1_j(\bR)$ on the interval $(1,\infty)$.
Assume $\bar \psi: \bar \sH^\rd\to \prP^1_j$ is the cover from an absolutely 
irreducible
component of a reduced Hurwitz with $\bar\psi$ over $\bR$. Let $G_{\bar\psi}$ be
its geometric monodromy group and take $N=\deg(\bar\psi)$. Then, an involution 
$c\in N_{S_N}(G_{\bar\psi})$ gives the effect of complex conjugation on the 
points of
$\bar\sH^\rd$ lying over $j_0$. 

Suppose $(\gamma_0',\gamma_1',\gamma_\infty')$  are the
branch cycles from Prop.~\ref{j-Line} from the action on reduced Nielsen 
classes. Then, 
\begin{equation} \label{caction} c\gamma_1c=\gamma_1 \text{ and } c\gamma_\infty
c=\gamma_\infty^{-1}.\end{equation} \end{lem}

\begin{proof} The hypotheses are for the situation of a 3-branch point real 
cover of the
sphere. The element $c$ will be independent of the value of $j_0\in (1,\infty)$. 
Since $j_0$ is the image of
$\bz$ with
$z_1< z_2<\dots < z_r \in \prP^1_z(\bR)$, Lem.~\ref{locj} says $j_0\in
(1,\infty)$. The formula for computing complex conjugation is the special case 
of
Prop.~\ref{compTest} where $r=3$ and all branch points are real.  

The effect of complex
conjugation on the
$q_i\,$s take them to their inverse. This induces the effect of taking 
$q_1q_2q_1$ to its
inverse
$q_1^{-1}q_2^{-1}q_1^{-1}$. The effect, however, of this image element on 
reduced Nielsen
classes is an element of order 2. Since the permutation effect of complex 
conjugation is to
conjugate by an involution $c$, this gives the formula $c\gamma_1c=\gamma_1$. 
Similarly, $c\gamma_\infty
c=\gamma_\infty^{-1}$. \end{proof}

\subsubsection{Elements of $N_{S_{\deg(\psi^\rd)}}(G_{\psi^\rd})$ that
centralize $G_{\psi^\rd}$} \label{centSn} Suppose $\row g r$ are generators of a 
group $G$, and
$c\in G\le S_n$ is an involution. Let $\Cen_{S_n}(G)$ be the centralizer of $G$ 
in $S_n$. 

\begin{princ}[Centralizer Principle] \label{cenCond} If $\Cen_{S_n}(G)$ contains 
no
involutions, then  conjugation by 
$c$ on $\row g r$ determines it. Assume $G$ is transitive. Then, with no 
assumptions on
$\Cen_{S_n}(G)$, if $c$ fixes 1, then its conjugation effect on $\row g r$ 
determines c. 

Suppose a(n irreducible) cover $\phi:X\to Y$ over a field $K$ with 
$n=\deg(\phi)$ has
monodromy group 
$G_*\le S_{n}$ over $K$. Then, the group $\Aut(\phi)_K$ of $K$ automorphisms of
$X$ commuting with $\phi$ identifies with $\Cen_{S_n}(G_*)$. \end{princ}

\begin{proof} Suppose $\Cen_{S_n}(G)$ contains no involutions. Let $c,c'$ be two 
involutions
with the same conjugation effect on $\row g r$. Then $cc'$ is an involution that 
centralizes
$G$, and so it is trivial. 

Now assume nothing about $\Cen_{S_n}(G)$, that $c$ and $c'$ have the same effect 
on $\row g
r$, both fix 1 and G is transitive. Then, $cc'\eqdef c^*\in \Cen_{S_n}(G)$ fixes 
1. From
transitivity, for any $i\in \{1,\dots,n\}$, there is $g\in G$ with $(1)g=i$. 
Conclude: $(i)c^*=(1)gc^*=(1)c^*g=(1)g=i$, and $c^*$ is the identify.  

The statement on $\Aut(\phi)_K$ comes from identifying it with
$N_{G_*}(G_*(1))/G_*(1)$. List the right cosets $\row {G_*(1)g} n$ of 
$G_*(1)$ in $G_*$. Elements of $G_*$ that permute these by {\sl left\/} 
multiplication on
these cosets are in $N_{G_*}(G_*(1))$. Those acting trivially 
are in $G_*(1)$. Left action commutes with the
right action of $G_*$, thus producing elements in
$\Cen_{S_n}(G_*)$.
\cite[Lem.~2.1]{FrHFGG} has complete details. 
\end{proof} 

\subsubsection{Outer automorphisms of $G$ and computing $c$ using
\eqref{caction}}
\S\ref{realPtEx} uses
\eqref{caction} to detect all real points on $\sH(A_5,\bfC_{3^4})^{\abs,\rd}$ 
and
$\sH(A_5,\bfC_{3^4})^{\inn,\rd}=\sH_0^{\inn,\rd}$. Then, \S\ref{level1} does the 
same for the space
$\sH(G_1,\bfC_{3^4})^{\inn,\rd}$ which is level 1 of the reduced Modular Tower 
for
$(A_5,\bfC_{3^4})$. Lem.~\ref{basicA5cover} shows the arithmetic 
(resp.~geometric) monodromy
group
$\hat G_{0,\inn}$ (resp.~$G_{0,\inn}$) of
$\bar \sH_0^{\inn,\rd}\to \prP^1_j$ over $\bR$ (resp.~$\bC$) is $\bZ/2\wr S_9$ 
(resp.
$\bZ/2\wr A_9$) in $S_{18}$. Further,
$\Cen_{S_{18}}(G_{0,\inn})$ is isomorphic to $\bZ/2$. Here it identifies with 
the element
$(1,\dots,1)\in (\bZ/2)^9$ generating the center of $G_{0,\inn}$. Therefore,
\eqref{caction} does not determine the complex conjugation $c$.

We now show this complication is common, occurring at all levels of many Modular 
Towers.
\S\ref{level1}, however, gives a satisfying remedy for it. Recall the normalizer
$N_{S_n}(G,\bfC)$ from \S\ref{sequiv}. 

\begin{prop} \label{aut-cent} Let $G\le S_n$ be a transitive subgroup and $\bfC$
a collection of conjugacy classes from $G$. Use this
permutation representation for absolute Nielsen classes. Let
$H=N_{S_n}(G,\bfC)/G$. Suppose $\sH''$ (resp.~$\sH'$) is an absolutely 
irreducible component
of $\sH(G,\bfC)^{\inn,\rd}$ (resp.~$\sH(G,\bfC)^{\abs,\rd}$) with 
$\psi^*:\sH''\to \sH'$
(over a field $K$) from the natural map $\psi:\sH(G,\bfC)^{\inn,\rd} \to
\sH(G,\bfC)^{\abs,\rd}$. Then, $\psi^*$ is Galois with group $H_{\psi^*}$ a 
subgroup of $H$. 

Denote the monodromy group of $\psi'': \sH''\to J_r$ ( over $K$) by 
$G_{\psi''}$. Then the centralizer of $G_{\psi''}$ in its natural permutation
representation of degree $\deg(\psi'')$ contains a subgroup isomorphic to 
$H_\psi$. 
\end{prop}

\begin{proof} Once we know that $\psi$ is Galois and $H_\psi$ is a subgroup of 
$H$, the
centralizer statement follows from Princ.~\ref{cenCond}. That identification is 
in
\cite[\S2.1]{FrVMS}. 
\end{proof} 

\begin{exmp}[The outer automorphism of $A_n$] \label{AnouterCont} Consider the 
map of
Hurwitz spaces
$\sH^{\inn,\rd}(A_n,\bfC_{3^r}) \to \sH^{\abs,\rd}(A_n,\bfC_{3^r})$ ($r\ge n-1$) 
from
Ex.~\ref{full3cycleList}. As noted there, each
component of  $\sH^{\inn,\rd}(A_n,\bfC_{3^r})$ maps by a degree 2 map to a 
component
of $\sH^{\abs,\rd}(A_n,\bfC_{3^r})$. This gives cases when complex conjugation 
$c$ on the
monodromy group of the cover does not determine its action on the fibers of a 
Hurwitz
(resp.~reduced Hurwitz) space over
$U_r$ (resp.~$J_r$). 
\end{exmp} 

\subsubsection{Using H-M reps.~to determine $c$} \label{detcHM} Often at level 
0, it is easy
to compute elements  that look like complex conjugation
operators (as in Prop.~\ref{aut-cent}) by inspection. A few tricks, however,  
are needed 
if \eqref{caction} does not determine the effect $c$ of complex conjugation on a
$j$-line cover. 

The simplest remedy is to identify an $\bR$-cover point. Suppose all branch 
points $\bz$
of a cover $\phi:X\to \prP^1_z$ are real. We may choose whatever paths we
desire on $U_\bz$ to compute Nielsen representatives of covers. For the next 
lemma,
choose them to detect covers defined over $\bR$ with $\bz$ as real branch points
using the $\hk$ in Prop.~\ref{compTest}. Call these the $\hk_\bz$-paths. 
For any set $\bz\in U_r(\bR)$ there are such $\hk_\bz$-paths. Let $\ni(\hk_\bz)$ 
be
branch cycles in $\ni(G,\bfC)^\inn$ for covers passing the test of 
Prop.~\ref{compTest}.

\begin{lem} \label{fixRed}Suppose $j_0\in (1,\infty)$ is the image of $\bp\in
\sH(G,\bfC)^\inn(\bR)$ corresponding to a cover $\phi:X\to \prP^1_z$ branched 
over $\bz$ consisting
of four real points on $\prP^1_z$. Then, the complex conjugation operator $c$ 
for the cover fixes
the reduced Nielsen classes coming from any elements of $\ni(\hk_\bz)$. 
\end{lem}

\begin{exmp}[Cont. Ex.~\ref{AnouterCont}: $c$ for
$\sH(A_5,\bfC_{3^4})^{\inn,\rd}=\sH_0^{\inn,\rd}$]  \label{4realpoint} 
Consider  
$\bar\psi: \bar\sH_0^{\inn,\rd}\to \prP^1_j$. It has definition field $\bQ$.
\S \ref{innerHSReal} computes $c$ for it. 
This cover has geometric (resp.~arithmetic
monodromy) $G_{\bar\psi}$ (resp.~$\hat G_{\bar\psi}$) in
$S_{18}$ and  
$$c^\dagger=(1\,10)(2\,11)(3\,12)(4\,13)(5\,14)(6\,15)(7\,16)(8\,17)(9\,18)$$ 
generates
$N_{S_{18}}(\hat G_{\bar\psi})$ (\S\ref{innerHSReal}). The criterion
of  Prop.~\ref{compTest} gives two possible complex conjugation operators 
corresponding to $j_0\in
(1,\infty)$, that labeled
$c_{1,\infty}$ in
\eqref{c1infty} and $c^\dagger c_{1,\infty}$. Only, however, the former fixes an 
H-M rep.
(represented by the integers 1 and 10) as is necessary from Lem.~\ref{fixRed}. 
\end{exmp}  

Using that the real components form a 1-dimensional manifold often is effective 
to handle the
intervals outside $j\in (1,\infty)$. \S\ref{innerHSReal} illustrates this. 

At levels beyond the first it is usually prohibitive to produce the complex 
conjugation
operators directly from the monodromy of reduced Hurwitz space covers. For 
example, \GAP\
couldn't do it for level 1 of our main example. Yet, with reasonable
computation ability with the group $G_k$ of the level, it suffices to check what 
$\hk$ (for
four real branch points) does to the elements of
$\ni(G_k,\bfC)^{\inn,\rd}=\ni(G_k,\bfC)^\inn \mod
\sQ''$. \S\ref{level1} illustrates by showing the genus 12 component of
$\sH(G_1,\bfC_{3^4})^\inn$ has one component of real points, while the genus 9 
component has no
real points. 

\section{Moduli and reduced Modular Towers} \label{defMT} Consider Hurwitz 
spaces $\sH(G_k,\bfC)^\inn$ attached to $\Col{(G_k,\bfC)} k$ as in  
Prop.~\ref{fineMod}. This is the Modular Tower for $(G=G_0,\bfC,p)$ (or for 
$(\tG p, \bfC)$
\cite[Part  III]{FrMT}.

\subsection{Reduced Modular Towers} Reduce elements of $\ni(G_{k+1},\bfC)^\inn$ 
modulo the kernel of 
$G_{k+1}\to G_k$. This induces $\sH(G_{k+1},\bfC)^\inn\to \sH(G_k,\bfC)^\inn$. 
The 
$\PSL_2(\bC)$ 
action is compatible with these maps. This produces the sequence for the {\sl 
reduced 
Modular 
Tower\/} for $(\tG p,\bfC)$: \begin{equation} \label{innSeq} \cdots\to 
\sH(G_{k+1},\bfC)^{\inn,\rd}\to \sH(G_k,\bfC)^{\inn,\rd}\to \cdots \to 
\sH(G_{0},\bfC)^{\inn,\rd}\to J_r.\end{equation}

\begin{defn} \label{defComp} Call a  sequence of representations 
$\Col{T_k:G_k\to 
S_{n_k}} k$  {\sl 
compatible\/} if $G_{k+1}(1)$ goes to a conjugate of $G_k(1)$ by the canonical 
map 
$G_{k+1}\to G_k$. 
\end{defn} A sequence of absolute Hurwitz spaces requires a compatible system 
of representations 
\cite[Part III]{FrMT}. Regular representations of each group $G_k$ give one 
example. 
Another example 
appears when $(|G_0(1)|,p)=1$. Apply Schur-Zassenhaus to the inverse image of 
$G_0(1)$ in $G_k$ to 
conclude $G_0(1)$ embeds compatibly in all the $G_k\,$s. Take $T_k$ the action 
of 
$G_k$ on $G_0(1)$ cosets. Example: With $G_k=D_{p^{k+1}}$ in its standard 
representation with $p$ odd, $G_0(1)$ is cyclic of order 2 (\S\ref{HurView}). 

Quotient by $\PSL_2(\bC)$ to produce \begin{equation} \label{absSeq} 
\begin{array}{rl} 
\cdots\to 
\sH(G_{k+1},\bfC,T_{k+1})^{\abs,\rd}&\to \sH(G_k,\bfC,T_k)^{\abs,\rd}\to \\& 
\cdots 
\to 
\sH(G_{0},\bfC,T_0)^{\abs,\rd}\to J_r.\end{array}\end{equation} We  
suppress the 
appearance of $T_k$ when the representation is obvious.

\begin{rem}[Sequences of other equivalences] We don't know example 
sequences of representations $\Col{T_k:G_k\to 
S_{n_k}} k$  {\sl 
compatible\/} for the characteristic quotients of $\tG p$. For example, suppose 
$G_0=A_5$ and $p=2$, and $G_0(1)$ includes the whole 2-Sylow of $A_5$. Then, 
Prop.~\ref{A5frat}  shows $G_k(1)$ must be the pullback of $G_0(1)$ in $G_k$. 
Such an example (giving not faithful representations) is useless for most 
purposes. 
The faithful representations of $G_1$ in Prop.~\ref{spinSepLst} have $G_1(1)$ 
lying  
over $G_0(1)$, a group with 2-Sylow of order 2. These representations give 
{\sl spin separation\/} (Def.~\ref{spinSepRepDef}). We suspect there are 
extending
compatible $G_k(1)$, one for each 
$k$, giving spin separation at all levels, though we haven't found them yet.  

Compatible sequences of permutation representations, suggest considering 
compatible 
sequences of subgroups $N_k'\le N_{S_{n_k}}(G_k,\bfC)$. Compatibility would 
require that 
$N_{k+1}'$ map to $N_{S_{n_k}}(G_k,\bfC)$. The next lemma notes this is not 
automatic. 

\begin{lem} To induce an action on
the cosets of 
$G_k(1)$ in $G_k$ requires knowing $N_{k+1}'$ normalizes the kernel of 
$G_{k+1}\to G_k$. This holds automatically if $\ker_0$ is a characteristic 
subgroup of 
$\tG p$: a common event (see Lem.~\ref{charker0}). \end{lem}
\end{rem} 

\subsection{$j$-line covers when $r=4$} \label{j-covers} Let $\sH_{O'}$ 
correspond to 
an $H_4$ orbit $O'$ in its action on $\sH(G,\bfC)^\abs$ (resp. on 
$\sH(G,\bfC)^\inn$). 
Thus, $\sH_{O'}$ is an absolutely irreducible component of a Hurwitz space 
$\sH(G,\bfC)$ (over some number field $K$), equivalence classes of $r=4$ branch 
point 
covers.  Prop.~\ref{PSL2-quotient} produces a finite cover $$\beta_4^\rd(O'): 
\sH_{O'}/PSL_2(\bC)=\sH_{O'}^\rd\to \prP^1_j\sem\{\infty\}.$$ Complete this to a 
cover $\bar \beta^\rd(O'):\bar\sH^\rd_{O'}\to \prP^1_z$.

Let $q_1,q_2,q_3$ be the images of $Q_1,Q_2,Q_3$ (or $q_1,q_2,q_3\in H_4$) in  
$M_4=H_4/Q$ as in \eqref{h4det}. Form one further equivalence on 
$\ni(G,\bfC)^\abs$ 
(or  $\ni(G,\bfC)^\inn$). Recall: For $\bg\in \ni(G,\bfC)$, $g_1 g_2g_3 g_4=1$.  
For 
$\bg\in \ni(G,\bfC)$, $Q_1Q_3^{-1}$ has this effect: \begin{equation} 
\label{q1q3} 
\bg\mapsto (g_1^\sph g_2^\sph g_1^{-1}, g_1,  g_4, g_4^{-1} g_3^\sph g_4^\sph).  
\end{equation} Similarly, with $\alpha=Q_1Q_2Q_3$, $\alpha^2$ has this effect: 
\begin{equation} \bg\mapsto (g_2,g_3,g_4,g_1)Q_1Q_2Q_3\mapsto (g_3,g_4,g_1,g_2). 
\end{equation}

\cite[Prop.~6.5]{DFrIntSpec} used  the normal 
subgroup of $H_4$ that $Q_1Q_3^{-1}$ generates acting on $\ni(G,\bfC)$. It
simplifies   computations to make these observations from Thm.~\ref{presH4}.
\begin{edesc} \item  The action of $\sQ$ on $\ni(G,\bfC)$ factors through the
Klein 4-group $\sQ''$.  \item 
$\sQ$ is the minimal normal subgroup of $H_4$ containing either $\alpha^2$ or 
$Q_1Q_3^{-1}$.  \end{edesc} Denote the $\sQ''=\lrang{Q_1Q_3^{-1},\alpha^2}$ 
orbits on 
$\ni(G,\bfC)$ by $\ni(G,\bfC)^\rd$: {\sl reduced\/} classes. Apply \eqref{mcgp}: 
Action of $H_4$ on $\ni(G,\bfC)/N'$ (as in \S\ref{refEquiv}) induces $H_4/Q$ 
acting on 
$\ni(G,\bfC)^\rd/N'$. For $\ni(G,\bfC)^\abs$ (resp. $\ni(G,\bfC)^\inn$) there is 
the 
quotient set $\ni(G,\bfC)^{\abs,\rd}$  (resp. $\ni(G,\bfC)^{\inn,\rd}$). 
Continue using 
$q_1$ and $q_2$ for $Q_1$ and $Q_2$ acting on $\ni(G,\bfC)^\rd$. Then, 
\cite[Prop.~6.5]{DFrIntSpec} computes branch cycles  for $\bar \beta^\rd(O')$.

\begin{prop}[$j$-line branch cycles] \label{j-Line} Consider 
$\gamma_0=q_1 q_2 \text{ and } \gamma_1= q_1 q_2 q_1$, generators of $\bar M_4$, 
with relations  
$\gamma_0^3=\gamma_1^2=1:\ \bar M_4\equiv \PSL_2(\bZ)$ \eqref{basM4}. Further, 
with
$\gamma_\infty=q_2$, the product-one condition $\gamma_0\gamma_1\gamma_\infty=1$
holds \eqref{prodOneM4}.  

Then, $\bar M_4$ orbits on 
$\ni(G,\bfC)^{\abs,\rd}$ (resp. $\ni(G,\bfC)^{\inn,\rd}$) correspond one-one to 
$H_4$ orbits on
$\ni(G,\bfC)^\abs$ (resp.~$\ni(G,\bfC)^\inn$).  Let $O'$ be the orbit in the  
discussion above.
Let $\gamma_0'$, $\gamma_1'$ and $\gamma_\infty'$ be respective  actions of 
$\gamma_0$, $\gamma_1$
and $\bar q_2=\gamma_\infty$ on the image of 
$O'$ 
in $\ni(G,\bfC)^{\abs,\rd}$ (resp.~$\ni(G,\bfC)^\inn$). Then $(\gamma_0', 
\gamma_1',
\gamma_\infty')$ is a  branch cycle description of the cover 
$\bar\beta^\rd(O')$.  
\end{prop}

\begin{proof}[Comments on the proof] Suppose $Q\in H_r$ gets killed in the
$\PSL_2(\bC)$ quotient of all moduli spaces of $r$ branch point covers. This 
happens if for every
$\bz\in U_r$, there exists  
$\alpha\in\PSL_2(\bC)$ fixing $\bz$ and inducing on $\pi_1(U_\bz)$ the same 
effect (modulo inner
automorphisms)  as does $Q$. Prop.~\ref{isotopyM4} identifies the group $\sQ$ as 
the
group of such $Q$, and 
$\bar M_4$ as the quotient of
$M_4$ by a Klein 4-group, $\sQ''=\sQ/\lrang{(q_1q_3^{-1})^2}$. 
 \end{proof}
 
\subsection{Reduced moduli spaces: b-fine and fine} \label{HtoredH} Let 
$\sQ''=\sQ/\lrang{z}$ as in 
Thm.~\ref{presH4}.   We start with results assuring some kind of fine moduli 
condition for reduced
Hurwitz spaces.
\S\ref{nofineMod} then illustrates why we cannot escape considering situations 
where it does not
hold. 

\subsubsection{Fine moduli} We interpret the phrase {\sl fine (resp.~b-fine) 
reduced
Hurwitz  space\/} for inner equivalence (absolute equivalence is similar,
Rem.~\ref{fineRedAbs}). 

Let $\sH^\rd$ be the reduced space for inner equivalence on
Nielsen classes of covers  attached to $(G,\bfC)$. Consider any smooth family 
$\Phi:\sT\to S$ of
curves  with an  
analytic map $\Psi: S\to J_r$. Denote the fiber of $\Phi$ over $s\in S$ by 
$\sT_s$. 
Assume $G$ acts as a group scheme on $\sT$ preserving each fiber $\sT_s$: 
$\Gamma: G\times 
\sT\to \sT$. For the quotient $\phi_s: \sT_s\mapsto \sT_s/G$, assume for $s\in 
S(\bC)$:  
\begin{triv} \label{GactionRed} An isomorphism of $\sT_s/G$ with $\prP^1_z$ 
presents 
$\phi_s$ in 
$\ni(G,\bfC)^\inn$ with branch points in the equivalence class 
$\Psi(s)$.\end{triv} 
\noindent By assumption 
this induces $\psi: S\to\sH^\rd$: $\sH^\rd$ is a  target for such maps. The 
quotient $\sT/G$ is a
geometric $\prP^1$ bundle over $S$. \S\ref{redHtoH} briefly discusses the 
obstruction to fibers
being $\prP^1$ over their definition field. 

\begin{defn}[Fine reduced moduli] \label{fineModRed} Call $\sH^\rd$ a fine 
moduli space (has
fine moduli) if for every such family, there is a {\sl unique\/} family 
$\sT^\rd\to \sH^\rd$
satisfying
\eqref{GactionRed} inducing $\sT$ by pullback from $\psi$. \end{defn}

\begin{defn}[b-Fine reduced moduli] \label{bfineModRed} Consider 
$U_j'=\prP^1_j\setminus
\{0,1,\infty\}$ and $\sH^\rd_{U_j'}$, the restriction of $\sH^\rd$ over $U_j'$. 
The weaker notion
b-fine is that
$\sT$ restricted to
$\psi^*(\sH^\rd_{U_j'})$ is the pullback by $\psi$ of \ $\sT^\rd$ restricted to
$\sH^\rd_{U_j'}$. \end{defn} 
Each notion applies separately to any component of $\sH^\rd$. Further, there is 
an obvious
generalization to $r>4$, though we will not be able to be so precise about 
testing for it. The
action of
$\sQ''$ on Nielsen classes gives an if and only if test for a  reduced Hurwitz 
space being a b-fine
moduli space  (Prop.~\ref{redHurFM} for inner equivalence,  
Rem.~\ref{fineRedAbs} for absolute
equivalence). The
$k$th level of the $(A_5,\bfC_{3^4}, p=2)$ Modular  Tower passes this test for 
b-fine moduli, 
$k\ge1$ (Prop.~\ref{HMreps-count}); even for fine moduli (Ex.~\ref{shA5k=0} and
Lem.~\ref{fixedg0g1}).  It is not even b-fine for level $k=0$. 

\begin{prop} \label{redHurFM} Let $\sH$ be a Hurwitz space with inner Nielsen 
class
$\ni(G,\bfC)$.  A component $\sH_*^\rd$ of the reduced space $\sH^\rd$ has fine 
moduli if
and only if there is a unique  total  space 
$\Phi_*:\sT_*^\rd\to \sH_*^\rd$ with a $G$ action on $\sT_*^\rd$ satisfying 
\eqref{GactionRed} when 
$S=\sH_*^\rd$. For $\bp\in \sH_*^\rd(K)$, $\phi_\bp: 
\sT_{*,\bp}\mapsto\sT_{*,\bp}/G$ has definition field $K$. 

Assume  $r=4$ and $\sH$ has fine moduli (as in Prop.~\ref{fineMod}). Let $O$ 
be an $H_4$  orbit on the Nielsen class corresponding to a component $\sH_*$ of 
$\sH$ with
$\sH_*^\rd$ its image in $\sH^\rd$. Then, $\sQ''$ orbits on 
$O$ have length 4 if and only if $\sH_*^\rd$ has b-fine moduli. Further, 
assuming b-fine moduli,
$\sH^\rd_*$ has fine moduli if and only if all its points over $j=0$ and $j=1$ 
ramify
(Prop.~\ref{j-Line}: $\gamma_0'$ and $\gamma_1'$ have no fixed points). 

If $r\ge 5$, and $\sH$ has fine moduli, then all components of $\sH^\rd$ have
b-fine moduli. 

Suppose $\sT_{*,\bp}$ has a $K$ divisor with odd degree image in  
$\sT_{*,\bp}/G$.
(Examples: a branch point conjugacy class is $K$ rational and distinct from 
other branch cycle 
conjugacy classes; or  
$K$ is a finite field.) 
Then, $\sT_{*,\bp}/G$  is  $K$ isomorphic to $\prP^1_z$. 

\end{prop}

\begin{proof} The field of definition statement follows from $\phi_\bp$  defined 
over $K$.  
Then, the divisor hypothesis produces an  
odd degree $K$ divisor on $\sT_{*,\bp}/G$. 
Since  $\sT_{*,\bp}/G$ has genus 0 and an odd degree $K$ divisor, it is
$K$ isomorphic to 
$\prP^1_z$.  If $K$ is a finite field, a homogeneous space for any Brauer-Severi 
variety 
always has a  rational point (the Brauer group of a finite field being
trivial \cite[p.~126]{SeGalCoh}).

Let $r=4$. Let $\bg\in O$ and $\bg'$ be in the $\sQ''$ orbit of $\bg$. Consider 
a 
given set of 
branch points $\bz$ and classical generators $\bar\bg$ of $\pi_1(U_\bz,z_0)$. 
Assume the
image of $\bz$ in $\prP^1_j$ is different from 0 or 1. Denote by 
$\phi_\bg:X_\bg\to\prP^1_z$ and $\phi_{\bg'}:X_{\bg'}\to\prP^1_z$ the covers 
from the 
homomorphisms 
sending $\bar\bg$ respectively to $\bg$ and to $\bg'$. Suppose for some 
$\alpha\in
\PSL_2(\bC)$ there exists an isomorphism
$\psi: X_\bg \to X_{\bg'}$ for which 
$\alpha\circ\phi_\bg= \phi_{\bg'}\circ\psi$. As in \S\ref{PSL2comp} and 
\S\ref{identK4},
$\alpha$ lies in a Klein 4-group identified with $\sQ''$. The b-fine hypotheses 
implies   
$\psi$ is unique.

So, the subgroup of $\PSL_2(\bC)$ fixing $\bz\in U_4$ extends to a faithful 
action of $\sQ''$ on $\sT$ over $\sH\times U_j'$. Quotient action gives the
total  space  
representing w-equivalence classes of covers for points of $\sH_{*,U_j'}^\rd$. 
The fine
moduli hypothesis say $\alpha$ fixed on $\bz$ over $j=0$ or 1
extends to $\sH_*$ without fixed points over that value of $\bz$. 

For $r\ge 5$, the b-fine property follows from \S\ref{rge5}. \end{proof}

\begin{rem}[Absolute equivalence ] \label{fineRedAbs}  Let $\sH^\rd$ be the 
reduced space for
absolute equivalence on Nielsen classes of covers  attached to $(G,\bfC,T)$. 
Consider any smooth
family
$\sT\mapright{\Phi} \sY\mapright{\pi} S$ with an   analytic map  $\Psi: S\to 
J_r$ with $\pi$
a (geometric) $\prP^1$ bundle. Assume also the fiber 
$\sT_s\to \sY_s$ over
$s\in S$ satisfies this: 
\begin{triv} An isomorphism of $\sY_s$ with $\prP^1_z$ 
presents 
$\phi_s$ in 
$\ni(G,\bfC,T)^\inn$ with branch points in the equivalence class 
$\Psi(s)$.\end{triv} 
\noindent The analogs of b-fine and fine moduli are clear, and the proof and 
conclusion of
Prop.~\ref{redHurFM} hold with little adjustment. 
\end{rem}

\subsubsection{The geometry of not having fine moduli} \label{nofineMod} 
Dilemma: $\sH(D_{p^{k+1}},\bfC_{2^4})^{\inn,\rd}$ is $X_1(p^{k+1})$ without its 
cusps. The former 
space, however, isn't a fine (even b-fine) moduli space according to 
Prop.~\ref{redHurFM}. Yet, it
is  classical the latter is 
a fine  moduli space. Resolution: It is for elliptic curves with a $p^{k+1}$ 
division point, 
though not for (inner)  w-equivalence of dihedral group Galois covers  
(\S\ref{fineModDil}). 

\begin{rem}[Serre's criterion] \label{serFixedPt} An observation of Serre
says an automorphism of an elliptic curve fixing a $p^{k+1}\ge 3$ division point 
must be the
identity. This assures $X_1(p^{k+1})$ is a fine moduli space. An appropriate 
generalization
is to decide for $r=4$, when 
$\gamma_0$ and $\gamma_1$ have no fixed points on the level $k$ Nielsen classes
$\ni(G_k,\bfC)^{\inn,\rd}$ for a Modular Tower, when $k$ is large  (see 
\S\ref{indSteps}).    
\end{rem} 

\begin{defn}[Field of moduli] Let $\phi: X\to \prP^1_z$ represent a cover in a 
Nielsen class $\ni$
for inner or absolute equivalence. Assume it has definition field in the 
algebraic closure of a
field
$K$ (assume char. 0 for simplicity). Apply each $\sigma\in G_K$ to $\phi$, 
denoting the result
$\phi^\sigma$. Let
$G_{K,\phi}$ be the collection of $\sigma$ for which $\phi^\sigma$ is equivalent 
to $\phi$. The
{\sl field of moduli}, $K_\phi$ is the fixed field of $G_{K,\phi}$. For reduced 
equivalence, we have
a corresponding field $K_{\phi,\rd}$. \end{defn} 

Some circumstances might use a cover both as an inner cover, and as an absolute 
cover. Reflecting
this the notation would be $K_{\phi,\inn}\supset K_{\phi,\abs}$.  Having a fine
moduli space for the
\IGP\ assures for
$\bp\in
\sH(G,\bfC)^\inn$  (or $\sH(G,\bfC)^{\abs}$ or $\sH(G,\bfC)^{\inn,\rd}$, etc.) 
there is a cover
representing that point over  $\bQ(\bp)$. 

In both the inner case (when $G$ has a center), or the absolute case (when
$\Cen_{S_n}(G)$ is nontrivial), there has been work to make 
Hurwitz spaces useful. Modular Towers include all information about Frattini 
central
extensions, as in Ex.~\ref{cenPerf}. Lem.~\ref{charblem} reminds
of the Harbater-Coombes argument: In the inner case, even with a center, there 
is
an absolute cover over $\bQ(\bp)$. This representing cover, however, may have no
automorphisms over $\bQ(\bp)$. The extension of this result to any base (in 
place of a point) does
not hold
\cite[\S 5]{DDEm}. 

This recognizes existence of a family over $\sH$ and existence of a
representing cover over $\bp\in \sH$ as part of the same problem, though 
changing the base makes
a difference. Both problems suit the language of stacks or gerbes (as in first 
version of this
question \cite[\S4]{FrHFGG}). \cite[\S2]{DDescTh} gives an exposition on the 
gerbe approach
of \cite{DDEm} and attempts to compute the obstructions to these problems.  

\cite{WeFM} uses an approach like \S \ref{startHM} relating cusps and complex 
conjugation operators.
Consider again our main example $\sH(G_k,\bfC_{3^4})^\inn$, $k\ge 0$. There is a 
nontrivial central
Frattini extension $T_k'\to G_k$. Then,
$\sH(G_k,\bfC_{3^4})^\inn$ contains, among its components, the components of
$\sH(T_k',\bfC_{3^4})^{\inn}$ (the moduli  for inner covers with group $T'_k$).  
For each
$k\ge 1$, Prop.~\ref{HMnearHMLevel} notes there are two kinds of
real points on $\sH(G_k,\bfC_{3^4})^{\inn}$: those for H-M reps.~and those 
for near H-M reps. H-M
rep.~points in $\sH(T_k',\bfC_{3^4})^{\inn}(\bR)$ have representing covers over 
$\bR$.
Near H-M rep. points in $\sH(T_k',\bfC_{3^4})^{\inn}$ have no representing cover 
over $\bR$: for these
points the field of moduli is not a field of definition. This interprets as 
different
degeneration at cusps on 
$\sH(G_k,\bfC_{3^4})^{\inn,\rd}$ attached to H-M reps.~from those for near H-M 
reps. 

The phenomena above also happens at $k=0$, though there are no near H-M 
reps.~there. \cite{WeFM}
starts there to explain a 
$p$-adic theory.  
\cite{IharMat}  provides a Hurwitz space context for the  Drinfeld-Ihara-Grothendieck relations 
(that apply  to elements of the absolute Galois group; call these DIG 
relations). The approach was
through {\sl tangential base points\/}.  This is the \cite{WeFM} approach, 
though he does not take
the exact same tangential base points. For example, he often uses complex 
conjugate pairs of branch
points, while they always used sets of real branch points. Ihara's use of the 
DIG relations has
been primarily to describe the Lie algebra of the absolute Galois group acting 
through various
pronilpotent braid groups, especially on the 3 punctured $\lambda$-line.  
\cite[App.~C]{FrMT}  proposed Modular Towers, though a profinite construction, 
as suitably
like finite representations of the fundamental group to see the DIG relations at 
a {\sl finite
level\/}. There is an analogy with Ihara in that modular curves are close to 
considerations about
the
$\lambda$-line, we know no direct phenomenon for modular curves suggesting the 
DIG
relations. Still, modular curves are just one case of Modular Towers.

In Prop.~\ref{nearHMbraid}, the Modular Tower attached to $A_5$ and four 3-
cycles produces  
a system of Serre obstruction situations from covers of $A_5$, typical  for a 
Modular
Tower. For $(\ell,|A_5|)$, $\ell$-adic points on these Modular Towers levels 
should have
a  similar tangential base point (cusp geometry) analysis to the near H-M and H-M reps.~over
$\bR$. We expect the H-M rep.~analysis  appplied to $\ell$-adic fields to give a tower 
of $T_k'$
realizations. The near H-M rep.~analysis, done
$\ell$-adically adds even more interest. We expect these to provide examples at 
each level $k$,
over $\ell$-adic fields of a function field extension with group a central 
Frattini
extension of
$G_k$ whose field of moduli is not a field of definition. Further, this geometry 
should reveal the DIG relations on  actual
covers, instead of as a Lie algebra relation. Though there are some points
that must be handled to do this, the computations of Cor.~\ref{R1G1} for level 1 
of this
$A_5$ Modular Tower should give a precise analog of \cite[Prop.~2.17]{WeFM} (for 
level 0).

\subsection{Points on a Modular Tower} As in \S\ref{HtoredH}, let $\sH$ be 
a Hurwitz space and 
$\sH^\rd$ its reduced version. We assume $\sH^\rd$ has (at least) b-fine moduli 
(for inner or
absolute equivalence). For
$K$ a field, a 
$K$ point on 
$\sH$ produces a $K$ point on $\sH^\rd$. \S\ref{redHtoH} interprets the subtlety 
of the 
converse. Then, 
\S\ref{ptsMT} applies this to points on a Modular Tower. \S\ref{grassman} and
\S\ref{openIm}  formulate Serre's {\sl Open Image Theorem\/} for $\bar \bQ$ 
points on a
Modular Tower. 

\subsubsection{$\PSL_2$ cocycles} \label{redHtoH} Use the setup for inner 
Hurwitz spaces in
(Prop.~\ref{redHurFM}). Suppose
$\sH^\rd$ is a fine (resp.~b-fine) moduli space, and $\bp^\rd\in
\sH^\rd(K)$  (resp.~also, doesn't lie over $j=0$ or 1). To simplify notation, 
refer to the total
family over $\sH^\rd$ (resp.~$\sH^\rd_{U_j'}$) as $\sT$. To simplify further, 
also assume fine
moduli, for the adjustments to b-fine are obvious. The (geometric) 
$\prP^1_z$ bundle
$\sT^\rd/G$ is algebraic from Serre's GAGA. We may even cover $\sH^\rd$ with $K$ 
affine sets
$\{W_i\}_{i\in I}$ so the restriction of $\sT^\rd/G$ over $W_i$ is a conic 
bundle in
$\prP^2\times W_i$. That is, for each $\bp^\rd\in W_i$, the fiber over $\bp^\rd$ 
is a conic in
$\prP^2$.

We drop the subscript referring to $\bp^\rd$. This gives the reduced class of a 
cover
$\phi:X\to Y=\sT^\rd_{\bp^\rd}/G$ over
$K$, with group
$G$, in the Nielsen  class. To 
simplify,  Then, if $Y$ has a $K$ 
point, 
$\phi$ is a {\sl $K$-cover\/} (of $\prP^1_z$) in the Nielsen class. Compatible 
with definitions
of \S\ref{cover-brauer}, call 
each $\bp^\rd$ a {\sl $K$-cover point\/} (the structure of the Nielsen class is 
over $K$).
Otherwise, 
$Y$ is isomorphic to a conic over $K$. When $Y$ has no $K$ point, call 
$\bp^\rd$ a 
{\sl $K$-Brauer point\/} (of $\sH^\rd(K)$).  

The conic attached to a point $\bp\in \sH^\rd(\bar K)$ defines an element of the 
group of 2-torsion
elements in the Brauer group $B_2(K(\bp))$ of the field $K(\bp)$. 
Any $K$ component $\sH^\rd_*$ of the Hurwitz
space defines an element of 
$B_2(K(\sH^\rd_*))$ by the same argument for a generic point of $\sH^\rd_*$. 
When 
the closure of $\sH^\rd_*$ is isomorphic to $\prP^1_w$, regard the conic bundle 
as an
element $b(\sH^\rd_*)\in Br_2(K(w))$. This case arises often in the
\IGP. Given a rational function $g(w')$ in a new variable $w'$, consider $K(w)$ 
as a subfield of
$K(w')$ by setting $w=g(w')$. This induces $B_2(K(w))\to B_2(K(w'))$. 
\cite[p.~114--116]{SeGalCoh} discusses (what is in our notation)   if for some 
choice of $g$ the
image of $b(\sH^\rd_*)$ in $ B_2(K(w'))$ vanishes. There is a natural notion of 
poles of
$b(\sH^\rd_*)$, and \cite{Me2} shows that if there are at most four poles, then 
such a $g$ exists.
As in Prop.~\ref{isotopyM4}, when $\sQ''$ acts trivially on Nielsen classes, we 
can take $g$ of
degree 1, regarding it as a section for the $j$-invariant. So, in the next lemma 
we emphasize the
crucial case when  $\sQ''$ acts faithfully on the Nielsen classes for the 
component $\sH_*$. 

\begin{lem}[Reduced cocycle Lemma] \label{redCocycle} Let $\sH$ be a fine moduli 
space for the
Nielsen class. 
Assume $\sH_*^\rd$, a component of  $\sH^\rd$, is a fine moduli space and  
$y\in 
Y(\bar K)=Y_{\bp^\rd}(\bar K)$ as above for $\bp^\rd\in \sH_*^\rd$. This 
produces a cover $\phi_y: X\to \prP^1_z$ over $K(y)$ and a unique cocycle class 
in
$H^1(G_K,\PSL_2(\bC))$. These
fit in a 
cocycle of Nielsen class covers. The cocycle class is trivial if and 
only if 
$Y$ is isomorphic to $\prP^1_z$ over $K$. In turn this holds if and only if $Y$ 
has an  
odd degree $K$ divisor. 

Conversely, given a  cocycle of $(G,\bfC)$  Nielsen class covers attached to 
$X$, there is $K$
cover 
$\phi:X\to Y$, 
which over $\bar K$ is in the Nielsen class. 

Suppose $Y$ has a $K$ point. Let 
$j_0$ be the image 
in $J_r$ of $\bp^\rd$. Then, there is a one-one association between $K$ points 
in the 
fiber  $U_{r,j_0}$ of $U_r\to J_r$ and $K$ points $\bp\in\sH$ over $\bp^\rd$. 
\end{lem} 

\begin{proof} The linear system $\sL_y$ attached to $y$ gives an isomorphism of 
the 
genus 0 curve $Y$ 
with $\prP^1_z: \mu_y: Y\to\prP^1_z$. Take $\phi_y$ to be the composition of 
$\phi$ and 
this 
isomorphism. Apply each $\sigma\in G_K$ to $\phi_y$ to get ${}^\sigma\phi_y: 
X\to 
\prP^1_z$. The 
isomorphism here is that given by replacing $y$ by $\sigma(y)$. Since $\sH^\rd$ 
is a fine 
moduli space, 
there is a unique $\alpha_\sigma\in \PSL_2(\bar K)$ and $\psi_\sigma: X \to X$ 
satisfying 
$\alpha_\sigma\circ\phi_y={}^\sigma\phi_y\circ\psi_\sigma$. The cocycle 
condition 
follows from the 
uniqueness conditions. Call this data {\sl a cocycle of Nielsen class covers\/} 
attached to 
$X$. It is standard 
the cocycle is trivial if and only if $Y$ has a $K$ point. Since $Y$ has genus 
0, this is 
equivalent to 
$Y$ having a degree one $K$ divisor. Since the canonical class on $Y$ is a 
degree -2 class over $K$,
this  is equivalent to $Y$ having an odd degree $K$ divisor. 

Now suppose we have such a cocycle of Nielsen class covers. This produces $\phi$ 
as 
$\mu_y^{-1}\circ \phi_y$. We have only to check $\phi$ is well-defined, 
independent of 
$\sigma\in G_K$. The cocycle condition guarantees this. 

Finally, if there is a $K$ point on $Y$, this gives $\bp\in \sH$ lying over 
$\bp^\rd$. 
Composing the cover 
$\phi_\bp:X_\bp\to \prP^1_z$ with elements of $\PSL_2(K)$ gives the 
correspondence 
between the $K$ 
points of $U_{r,j_0}$ and the $K$ points of $\sH$ over $\bp^\rd$. \end{proof}

\begin{exmp} \S\ref{compMC} considers the Nielsen class is 
$\ni(D_p,\bfC_{2^4})$. Then, a point 
$\bp^\rd\in \sH^\rd(K)$ has a $K$ point of $\sH$ over it. This is because an 
elliptic curve 
isogeny $E\to 
E'$ over $K$ represents $\bp^\rd$.  The quotient of $E'$ by $\lrang{-1}$ has a 
rational 
point from 
the image of the rational point on $E'$. This gives $E\to \prP^1_z$ representing 
$\bp\in 
\sH(K)$ lying over $\bp^\rd$. \end{exmp} 

\begin{rem} Suppose $\sH^\rd$ is not a fine moduli space. We can still ask which 
$\bp\in
\sH^\rd$ are $K$-cover points or $K$-Brauer points as in Rem.~\ref{nontrivq''1} 
and
Rem.~\ref{nontrivq''2}.  \end{rem}

\subsubsection{Projective systems of points} \label{ptsMT} 
Let $\Col{\bp_k^\rd\in 
\sH(G_k,\bfC)^{\rd}} 
k$ be a projective system of points on a reduced Modular
Tower. Call this a {\sl point\/} on a reduced  Modular
Tower.  Suppose 
$\sH^\rd_k(\bp^\rd)$ is the absolutely irreducible component containing 
$\bp^\rd_k$. Then 
$\sH^\rd(\bp^\rd)=\Col{\sH^\rd_k(\bp^\rd)} k$ is a projective sequence of 
algebraic varieties. This works with $\bp$ a
projective system of points on  Hurwitz 
spaces $\sH(\bp)=\Col{\sH_k(\bp)} k$ rather
than reduced Hurwitz spaces. 

\begin{princ} \label{betweenMTrdMT} Components of $\sH^\rd$
are manifolds and moduli spaces. Conclude: 
$\bp_k^\rd$ being a $K_k$ point implies $K_k$ contains a definition field for 
$\sH^\rd_k(\bp^\rd)$. 

Assume
$\bp_0\in \sH_0$ lies over the level 0 point $\bp^\rd_0$ of
$\bp^\rd$. Then, this produces a point $\bp$ on the Modular
Tower with its level 0 point equal $\bp_0$. \end{princ}

For $W$ any algebraic variety, and $w_0\in W$, 
denote the pro-$p$ completion of $\pi_1(W,w_0)$ by
$\pi_1(W,w_0)^{(p)}$. It is the closure of $\pi_1(W,w_0)$ in the diagonal of the 
product of all
finite $p$-group quotients of $\pi_1(W,w_0)$. When considering homomorphisms 
involving it, 
defined up to conjugation by an element of this group, with no loss drop the 
$w_0$ decoration.

We concentrate now on inner Hurwitz spaces.  For $z_0\in K$,
not in $\bz$, consider  classical 
generators $\bar \bg$ for $\pi_1(U_\bz,z_0)$ (\S\ref{classGens}). Let 
$\Col{\bp_k^\rd\in 
\sH(G_k,\bfC)^{\inn}} 
k$ be a projective system of points over $\bp_0\in \sH(G_k,\bfC)^{\inn}$ with
$X_{\bp_k}\to \prP^1_z$ a representing cover. Any
projective system $\Col{x_k\in X_{\bp_k}} k$ of points over 
$z_0$  gives a 
compatible system of homomorphisms $\psi_k:\pi_1(U_\bz,z_0)^\alg\to G_k$ 
factoring 
through 
$\pi_1(U_\bz,z_0)^\ari$ \eqref{exactArithBP}. This produces $\tilde \psi\in 
\Hom(\pi_1(U_\bz,z_0)^\alg, 
\tG p)/\tG p$ with $\bar\bg$ mapping into $\bfC$ (\S\ref{usebcl}), not depending 
on  
$\{x_k\}_{k=0}^\infty$. 

Restriction of
$\tilde
\psi$ to the kernel of
$\psi_0$ factors  through 
$\pi_1(X_{\bp_0})^{(p)}$ (Thm.~\ref{thm-rbound}, as in \S\ref{usebcl}). Factor 
$\pi_1(U_\bz,z_0)$ by the kernel of this map, and denote the result by 
$M_{\bp_0}$. 

\begin{prop} \label{tGppts} The group $M_{\bp_0}$ fits in a natural exact 
sequence $$ 
1\to 
\pi_1(X_{\bp_0})^{(p)}\to M_{\bp_0}\to G_0\to 1.$$ Keep the notation $\bar\bg$ 
for the 
image of 
$\bar\bg$ in $M_{\bp_0}$.  Points on the Modular Tower for $(G,\bfC)$ over 
$\bp_0\in\sH(G,\bfC)^\inn$  
correspond one-one with elements of  $\Hom(M_{\bp_0}, \tG p)/\tG p$ mapping 
$\bar\bg$ into $\bfC$. In 
turn these correspond with elements of $\ni(\tG p,\bfC)^\inn$. Then, the action 
of $H_r$ 
on $\bar\bg$ 
induces an action of $H_r$ on $\ni(\tG p,\bfC)^\inn$.

Let $\bp\in \sH$, a point on the projective sequence of components 
$\{\sH_k\}_{k=0}^\infty$ containing 
points over $\bp_0$, as above. Suppose this corresponds to $\psi\in 
\Hom(M_{\bp_0}, 
\tG p)/\tG p$. 
Then, the collection of  points $\bp'\in \sH$ above $\bp_0$ correspond to an 
$H_r$ orbit 
of 
$\bp$. \end{prop}

\subsubsection{A Grassman like object} \label{grassman} Let $g_0=g$ be the genus 
of 
$X_{\bp_0}$. For 
$\pi_1(X_{\bp_0})$, too, there is a notion of classical generators. These are 
topological 
generators, 
$\ba=(\row a g)$ and ${\pmb b}=(\row b g)$, satisfying these properties. 
\begin{edesc} 
\label{poinGens} 
\item The only relation in  
$\pi_1$ is the commutator product $\prod_{i=1}^g [a_i,b_i]=1$.  \item
In 
$H_1(X_{\bp_0},\bZ)$, the cup product pairing maps $(a_i,a_j)$ to 0, $(b_i,b_j)$ 
to 0, 
and $(a_i,b_j)$ to 
$\delta_{ij}$ (Kronecker $\delta$ function) for all $i$ and $j$. \end{edesc} 

\newcommand{\sF}{{\Cal F}} Let $\pi_1(X_{\bp_0})^{(p)}=\pi_1^{(p)}$. As in 
Prop.~\ref{tGppts}, 
consider the collection $\sF_{\bp_0}$ of homomorphisms $\tilde \psi: 
M_{\bp_0}\to \tG 
p$ (up to inner 
action) having this property: \begin{triv} $\tilde \psi$ induces the identity 
map 
$M_{\bp_0}/\pi_1^{(p)}\to 
\tG p/\ker_0$. \end{triv} Since $\tG p\to G_0$ is a Frattini cover, such a 
$\tilde \psi$ is 
necessarily 
surjective.  Denote the surjective homomorphisms in $\Hom(\pi_1^{(p)},\ker_0)$ 
by 
$\Hom^*(\pi_1^{(p)},\ker_0)$. Thus, $\sF_{\bp_0}$ maps into this space. Also, 
let 
$\bT_p(X_{\bp_0})=\bT_p=\pi_1^{(p)}/[\pi_1^{(p)},\pi_1^{(p)}]$ denote the (dual 
of 
the) {\sl Tate module\/} of $X_{\bp_0}$.  
  
\begin{prop} \label{GKaction} Let $K$ be a subfield of $\bar\bQ$. Assume 
$\bp_0\in 
\sH(G_0,\bfC)(K)$ and 
$\sH(G_0,\bfC)$ has fine moduli. Then,  
$G_K$ acts naturally on $\Hom^*(\pi_1^{(p)},\ker_0)/\ker_0$. This permutes 
elements of  
$\sF_{\bp_0}$, compatible (according to \S \/{\rm \ref{groth}}) with acting on 
coordinates of points $\bp\in \sH$. It induces a $G_K$ action on 
$\Hom^*(\bT_p,\ker_0/(\ker_0,\ker_0))$. \end{prop}

\begin{proof} The action of $G_K$ on $\pi_1(X_{\bp_0})^\alg$ (\S \ref{groth}) 
induces 
an action on any 
closed subgroup of $\pi_1(X_{\bp_0})^\alg$. For $\sigma\in G_K$ denote its fixed 
field 
by $K_\sigma$. 
Let $\bp\in \sH$, so it defines a homomorphism $\psi_\bp: 
\pi_1(X_{\bp_0})^\alg\to \tG 
p$.   Then, $\bp$ 
is a projective system of points over $K_\sigma$ of $\sigma$ if and only if 
$\sigma$ 
normalizes the kernel 
of $\psi_\bp$, and the induced action of $\sigma$ on the quotient is trivial. 
Since levels of 
the Modular 
Tower are fine moduli spaces, this produces the desired sequence of Galois 
covers over 
$K_\sigma$. 

The argument reverses. Further, $\sigma$ acts on 
$\pi_1^{(p)}/(\pi_1^{(p)},\pi_1^{(p)})$, a characteristic subgroup of 
$\pi_1^{(p)}$ (on 
which 
$G_K$  acts by 
hypothesis). The image of $\psi_\bp$ in $\tG p$ is into $(\ker_0,\ker_0)$. This 
induces a 
$G_K$ action on 
$\Hom^*(\bT_p,\ker_0/(\ker_0,\ker_0))$. \end{proof}

Thm.~\ref{thm-rbound} says $G_K$ has no fixed points on $\sF_{\bp_0}$. More 
generally, there are no $K_\sigma$ points $\bp\in \sH$ if $\sigma$ induces the 
Frobenius on $\bT_\ell$ for some prime $\ell$ not dividing $|G_0|$. This topic 
continues in
\S\ref{openIm}.

\section{Group theory of an $A_5$ Modular Tower} \label{startA5MT}
 
A Modular Tower has levels corresponding to a sequence of groups  
$$\cdots\to G_{k+1}\to G_k\to \cdots \to G_0\ \text{\rm(Prop.~\ref{fineMod})}.$$   
If $G_0$ is
centerless and $p$-perfect,  each
$G_k$ is a  centerless (Prop.~\ref{fineMod}) Frattini extension of $G_0$ with 
$p$-group as
kernel.  Frattini extensions of perfect groups are perfect: Commutators of the 
covering
group  generate the image, so they generate the covering group. 

\S\ref{frattiniPre} recounts the geometry behind a Frattini cover.
\S\ref{pfrattini} extends the discussion of \cite[Chap.~20]{FrJ} on the 
universal $p$-Frattini cover of $G_0$,
starting with the case when its $p$-Sylow is normal.  \S\ref{allA5} describes 
the groups $G_k$
when $G_0=A_5$. \cite[Part II]{FrMT}  describes some aspects of the universal 
Frattini cover $\tilde A_5$ of
$A_5$, especially the ranks of the Universal $p$-Frattini kernels. Then, $\tilde 
A_5$ has three  pieces 
${}_2\tilde A_5$, ${}_3\tilde A_5$ and ${}_5\tilde A_5$, one for each prime $p$ 
dividing $|A_5|$ (\S\ref{frattGps}). We use \cite[Frattini Principle 
2.3]{FrKMTIG} to enhance and simplify properties of the characteristic modules 
$M_k$ ($\ker_k/\ker_{k+1}\eqdef
M_k$ in \S \ref{fratAnalog}) of ${}_2\tilde A_5$: $\tG p$ when
$G=A_5$ and 
$p=2$. 

\S\ref{onlytwo} shows only two Modular  Towers for $A_5$, $p=2$ and $r=4$ have 
$\bQ$ points at level 1. Then, 
\S\ref{startHM} explains H-M and near H-M representatives. When $p=2$, these 
describe connected 
components of real points on levels 1 and above of a Modular Tower. As a prelude 
for 
\S\ref{compA5MT}, \S \ref{realPtEx} uses $j$-line branch cycles 
(\S\ref{j-covers}) for 
diophantine conclusions about regular realizations. This describes 
real point  components on level 0 of the $(A_5,\bfC_{3^4})$ Modular Tower. This 
paper uses properties of $M_1$ as a $G_1$ module. 

\subsection{Frattini curve covers} \label{frattiniPre} Let $\phi:X\to Z$ be a 
Galois
cover with group $G$, and let $\psi: Y\to X$ be a cover for which 
$\phi\circ\psi$ is Galois
with group $G^*$. Call
$\phi\circ\psi$ a {\sl Frattini extension of $\phi$\/} if the following holds. 
For any sequence
$Y\to W\to Z$, of (not necessarily Galois) covers with $W\ne Z$, there is always 
a proper cover of
$Z$ through which both $W\to Z$ and $X\to Z$ factor. A Frattini extension of 
$\phi$ 
has no differentials and functions that are pullbacks from covers disjoint from
$\phi$. 

For a field theoretic restatement of the property let $K\subset \hat L \subset 
\hat M$ be a chain
of fields with
$\hat M/K$ (resp.~$\hat L/K$) Galois with group $G^*$ (resp.~$G$). This is a 
{\sl Frattini
chain\/} if the only  subfield $K\le T\le \hat M$ for which $T\cap \hat L=K$, is 
$T=K$. Denote by
$\rest: G^*\to G$ the natural map. 

Let $T=\hat M^H$ be the fixed field of a subgroup $H$
of $G^*$. Then, $T\cap
\hat L=K$ is equivalent to  $\rest: H\twoheadrightarrow G$. Hint: $T\cap
\hat L=K$ allows extending any automorphism of $\hat L$ to $T\cdot \hat L$ to be 
the identity on $T$. 
The group theoretic restatement is that  
$H\le G^*$ and $\rest(H)=G$ implies $H=G^*$: $\rest:G^*\to G$ is a {\sl Frattini 
cover}.
A Frattini cover $G^*\to G$ always has 
a {\sl nilpotent\/} kernel. 

\subsection{The normalizer of a $p$-Sylow and Loewy display}\label{pfrattini} 

\newcommand{\bker}[1]{{\ker_0(#1)/\ker_1(#1)}}

\subsubsection{Starting with a $p$-Sylow} \label{normp-Sylow} Let $H=N_G(P_p)$ 
be the
normalizer of  a $p$-Sylow $P_p$ of any finite group $G$. Apply Schur-Zassenhaus 
to write
$H$ as 
$P_p\xs H^*$ with $H^*$ having order prime to $p$. Let ${}_p\tilde F_t$ be the 
pro-free 
pro-$p$ group on $t$ generators, with $t$ the rank (minimal number of 
generators) of 
$P_p$. Then, the universal $p$-Frattini cover of $H$ is ${}_p\tilde F_t\xs H^*$; 
extend 
$H^*$ acting on $P_p$ to ${}_p\tilde F_t$ through the map $\psi: {}_p\tilde 
F_t\to 
P_p$ as in Remark \ref{extH}.

For {\sl any\/} group $G$, $\ker_0(G)$ is the pro-free pro-$p$ kernel of 
${}_p\tilde G\to 
G$. Commutators and $p$th 
powers in $\ker_0$ generate $\ker_1(G)$, the Frattini subgroup of $\ker_0(G)$. 
Iterating
this produces $\ker_k(G)\le \ker_{k-1}(G)$. Let 
$\ker_0(G)/\ker_1(G)=M(G)$ be the first characteristic quotient of the 
universal $p$-Frattini cover of $G$. Then $M(G)$ is a $\bZ/p[G]$ module and 
$G_1=\tG p/\ker_1$ is a Frattini extension of $G$ by $M(G)$.

\begin{lem} \label{fratIsom} Let $W$ be any subgroup of $G$. Then, ${}_p\tilde
W$ embeds in $\tG p$. Further, for each $k\ge 0$, $W_k$ naturally embeds in 
$G_k$. If
$\rk(M(W))=\rk(M(G))$, then ${}_p\tilde W$ appears from the universal $p$-Frattini cover
$\tilde\phi: {}_p\tilde G$ as 
$\tilde{\phi}^{-1}(W)$. For $H=N_G(P_p)$ this applies if  
$\rk(P_p)=\rk(M(G))$. \end{lem}

\begin{proof} A $p$-Sylow of ${}_p\tilde  G$ contains a $p$-Sylow of 
$\tilde{\phi}^{-1}(W)$.
So, the latter is profree.  \cite[Prop.~20.33]{FrJ} characterizes ${}_p\tilde W$ 
as the
minimal cover of
$W$  with  pro-free $p$-Sylow. So, there is a natural map
$\gamma_W: 
\tilde{\phi}^{-1}(W)\to  {}_p\tilde W$ commuting with the map to $W$. As 
${}_p\tilde 
W$ is a $p$-Frattini cover of $W$ the map is  surjective. Since the
natural map $\tilde{\phi}^{-1}(W)\to W$ has a pro-$p$ group as kernel, the 
natural map
${}_p\tilde  W\to W$ produces $\psi_W: {}_p\tilde  W\to \tilde{\phi}^{-1}(W)$ 
commuting
with the projections to $W$. 

The composition $\gamma_W\circ\psi_W$ (commuting with the
projections to $W$) is an endomorphism of ${}_p\tilde  W$. The image of
$\gamma_W\circ\psi_W$ is a closed subgroup of ${}_p\tilde  W$ mapping 
surjectively to
$W$. So, from the Frattini property, $\gamma_W\circ\psi_W$ is onto. An onto 
endomorphism of
finitely generated profinite groups is an isomorphism \cite[Prop.~15.3]{FrJ}. In 
particular,
$\psi_W$ is an injection. The characteristic quotients have maps between them 
induced by
$\psi_W$, and so $W_k$ injects into $G_k$, inducing an injection of
$\ker_k(W)/\ker_{k+1}(W)\to \ker_k(W)/\ker_{k+1}(W)$. If (for $k=0$), 
$M(W)$ and
$M(G)$ have the same dimension, they  are isomorphic. As these groups 
characterize
$\ker_0(W)$ and $\ker_0(G)$, that implies they are equal. This gives an 
isomorphism of
$\ker_0(H)$ and
$\ker_0(G)$ in the special case. 
\end{proof}

\begin{rem}[Extending $H^*$] \label{extH} Since ${}_p\tilde F_t$ is a pro-free 
pro-$p$ 
group, it is easy to create a profinite group $\bar H^*$ whose action on 
${}_p\tilde F_t$ extends 
$H^*$ on $P_p$. Any $h\in H^*$ lifts to an automorphism $\bar h$ of ${}_p\tilde 
F_t$
commuting with the map to $P_p$. The (profinite) automorphism group of 
${}_p\tilde F_t$ is
the projective limit of automorphism groups of finite group quotients (by finite 
index
characteristic subgroups) of ${}_p
\tilde F_t$. 

Let $\{h_i\mid i\in I\}$ be any generators of $H^*$. Choose $\bar
H^*$ to be the closure of $\lrang{\bar h_i\mid i\in I}$ in the automorphism 
group of ${}_p\tilde
F_t$. By construction it maps surjectively to $H^*$.  The problem is to split 
off a copy of
$H^*$. Automorphisms of 
${}_p\tilde F_t$ trivial on its Frattini quotient have $p$-power order 
\cite[Thm. 12.2.2: 
Philip Hall]{HaThGps}. The kernel from $\bar H^*\to H^*$ will be trivial on the 
Frattini
quotient of ${}_p\tilde F_t$ (which equals the Frattini quotient of $P_p$). So 
Schur-Zassenhaus
(for profinite groups, 
\cite[Chp.~10]{FrJ}) always allows splitting off a copy of $H^*$ in $\bar H^*$. 
Since, however, it
depends on
$\psi$, it is an art to do this  explicitly. Again by Schur-Zassenhaus, this 
copy of $H^*$
is unique up to conjugacy.  
\end{rem}

\subsubsection{Appearance of $\one_{G_k}$ in $M_k$} \label{onesAppearNGP} That
$\ker_k/\ker_{k+1}\eqdef M_k$ in \S \ref{fratAnalog} is a $\bZ/p[G_k]$ module 
appears in using modular
representation theory at each level of a Modular Tower.  Significantly, $M_k$ is 
not a $G_{k-1}$ module (though
it may have nontrivial quotients that are). For each $k$, however, $M_k$ is an 
$H^*$ module through this lifted
action. Let
$P_{p,k}$ be a $p$-Sylow of $G_k$ and $M_k(P_{p,k})$ (resp.~$M_k(P_{p,k}\xs 
H^*$) the
restriction of $M_k$ to $P_{p,k}$ (resp.~the normalizer of $P_{p,k}$ in $G_k$). 

Recall the 
{\sl Loewy display\/} of  composition factors for a $G$ module $N$. It derives 
from the 
{\sl radical submodules\/} of $N$, $$1\le N_t\le N_{t-1} \le \cdots \le
N_1\le N_0=N:$$  
$N_i$ is the minimal $G$ submodule of $N_{i-1}$ with $N_{i-1}/N_i$ semi-simple. 
The display
consists of writing the simple module summands of $N_{i-1}/N_i$ at  the $i\nm1$ 
position
with arrows indicating relations between modules at different  levels. 

Warning! A Loewy
display is not a sequence of module homomorphisms; it indicates the relation 
between
modules at different levels in the Krull decomposition. We write this display 
right to left, 
compatible with the quotients from an exact sequence, instead of top to bottom 
(as group theorists
often do). A Loewy display of all  information on subquotients of a module often 
requires several
arrows between layers; there may be several  arrows from a simple level
$i$ module toward level
$i-1$  modules (see the examples of Prop.~\ref{A5Morbs}). Note: Such a display 
uniquely  determines
a module.

The Loewy
display of $M_k(P_{p,k})$ consists of copies of $\one_{P_{p,k}}$, the only 
simple $p$-group module.
Jenning's Theorem
\cite[p.~87]{Ben1} is an efficient tool to figure the dimension of the Loewy 
layers of
$M_k(P_{p,k})$.  (The
proof of Prop.~\ref{hatP-P} has an example of its use.)  Further, since the 
action of $H^*$ respects this construction, it efficiently reveals
how $H^*$ acts on 
the Poincar\'e-Witt basis of the  $P_{p,k}$ group ring's universal enveloping 
algebra. The inductive detection
of $\one_{H^*}$ on Loewy layers of $M_k(P_{p,k})$  often comes through 
$\one_{H^*}$ appearing from a previous
Loewy layer, or from tensors products of representations from previous Loewy
layers (like a representation and its complex conjugate appearing juxtaposed). 
The following argument (a collaboration with D.~Semmen) for the $p$-split case 
shows  $\one_{H^*}$ usually appears
quickly unless the $p$-Sylow is cyclic. Princ.~\ref{obstPrinc} 
already shows $\one_G\,$ appears at infinitely many levels if $G$ has a center. 

\begin{prop}[$p$-split case] \label{psplitOne} Let $G=P\xs H$ be centerless, 
with $P$ a $p$-group and
$(|H|,p)=1$. If
$P$ is not cyclic, $\one_G$ appears infinitely often. Otherwise, $M_k(G)$ is the 
same cyclic $G$ module
for each $k$. So, $\one_G$ does not appear. \end{prop}

\begin{proof} By replacing 
$P$ by
$P/\Phi(P)$, with $\Phi(P)$ its Frattini subgroup, form $P/\phi(P)\xs H=G^*$. 
The universal $p$-Frattini 
cover of $P/\phi(P)\xs H$ is the same as that of $G$. So, Princ.~\ref{obstPrinc} 
shows  
$\one_{G^*}$ occurring at some level for $G^*$ implies it occurs at infinitely 
many levels for
$G$. So, for our question, with no loss take $G=G^*$ to assume $P$ is an 
elementary
$p$-group and a $\bZ/p[G]$ module.  

Then, with $P_\one$ the projective indecomposable for $\one_P$, 
consider it as an $H$ (therefore a $G$) module. Since $((G:P),p)=1$, Higman's 
criterion  says
it is the projective indecomposable for $\one_G$ \cite[p.~64]{Ben1}. Its most 
natural display might by
$\bZ/p[G](\sum_{h\in H} h)$. In its Loewy display
$\one_G$ appears to the far right. Write the next Loewy layer as $S_1\oplus 
\cdots \oplus S_t$, with the
$S_i\,$s irreducible $H$ ($G$) modules. 

Forming $P_\one^*$ (* is the $\bZ/p$ dual) returns $P_\one$. So the dual of the 
$v$th socle
layer of $P_\one$ is the $v$th Loewy layer of
$P_\one$. Conclude that to the far left there is
$\one_G$ (1st socle layer) with  
$S_1^*\oplus\cdots \oplus S_t^*$  immediately to the right of that (2nd socle 
layer). Now we 
apply that $P$ is an elementary $p$-group. The $P$ module $P_\one$ is the group 
ring $\bZ/p[P]$.
Let $\row x m$ be a basis for $P$ and identify $P$ in the vector space 
$\bZ/p[P]$ as the space spanned by
$\{x_1-1,\dots x_m-1\}$. The action of $H$ on $\bZ/p[P]$ preserves this space 
modulo the second power of the
augmentation ideal of $\bZ/p[G]$. Conclude: $P\equiv S_1\oplus \cdots \oplus 
S_t$ as $H$ (or $G$) modules.  

As in
Prop.~\ref{hatP-P} use the notation of \cite[Part II]{FrMT}. So, the kernel 
$\Omega(1,G)$ of $P_\one\to \one_G$
starts with 
$S_1\oplus \cdots \oplus S_t$. Now consider a minimal projective 
$P_{\Omega(1,G)}$
mapping surjectively to $\Omega(1,G)$. Higman's criterion again implies the $G$ 
module $P_{\one_G}\otimes_{G}
(S_1\oplus 
\cdots S_t)$, being projective for $P$, is also projective for $G$. The Loewy 
layers of $T=P_{\one_G}\otimes_G
(S_1\oplus 
\cdots S_t)$ are just the Loewy layers of $P_{\one_H}$ tensored over $\bZ/p[H]$ 
with $S_1\oplus 
\cdots S_t$ modules. Reason:  a Loewy layer of $P_{\one_G}$ tensored with 
$S_1\oplus 
\cdots S_t$ is semisimple. Therefore $T=P_{\Omega(1,G)}$. Then, 
$$M(G)=M_0(G)\equiv
\Omega(2,G)=\ker(T\to
\Omega(1,G))$$ modulo projective summands (as in \cite[Projective Indecomposable 
Lem.~2.3]{FrMT}, called
$\Omega^2(\one_G)$ there). The 2nd socle layer of $\Omega(2,G)$ is 
$(S_1^*\oplus\cdots
\oplus S_t^*)\otimes_G (S_1\oplus \cdots \oplus S_t)$ with one copy of $\one_G$ 
removed (from the end of
$\Omega(1,G)$.  The module
$S_i\otimes_G S_i^*=S_i\otimes_H S_i^*$ has exactly one appearance of
$\one_H$ for each absolutely irreducible factor in $S_i\otimes_{\bZ/p}\bar 
\bZ/p$ (Maschke's Theorem; since $P$
acts trivially on  $S_i$,  the same is true of $\one_G$). Thus, $\one_G$ appears 
in the Loewy display of
$\Omega(2,G)$ as a $G$ module unless
$t=1$,  and $S_1$ is absolutely irreducible. 

Now suppose $t=1$. Instead of looking at $M_0(P)=M_0(G)$, look at $M_1(P)$ as an 
$H$ module. By our
hypotheses, $P$ is not a cyclic module. Therefore, a computation of the rank of 
$M_1(P)$ comes directly from
Schreier's formula for ranks of subgroups of pro-free groups (see 
\S\ref{extA4C34}). The rank of $M_k$, as $k$
increases,   exceeds the degree of any irreducible $H$ module. Replace $P$ with 
a suitable $M_k$ to
revert to the case $P$ is not absolutely irreducible. This  completes showing 
the appearance of
$\one_G$. 
\end{proof}

The structure constant formula
(\S\ref{congClProd}) can detect $\one_{H^*}$ appearing here, as $(|H^*|,p)$ 
means it applies in 
characteristic $p$. This topic continues in \S\ref{genFrattini}.

\subsection{Producing ${}_2\tilde A_5$} \label{allA5} 
 Lemma \ref{fratIsom} hypotheses imply $M(H)$ is a $\bF_p[G]$ module.  

\subsubsection{Converse to Lemma \ref{fratIsom}}
In the proof of Lemma \ref{fratIsom}, $\gamma_H: \tilde{\phi}^{-
1}(H)\to {}_p\tilde H$ induces a surjective $H$ module map $\gamma'_H: M(G)\to 
M(H)$.

\begin{prop} \label{fratExt} Suppose $H=N_G(P_p)$ ($P_p$ a $p$-Sylow of $G$) and 
$M(H)$ is a $G$
module  extending the $H$ action so the following holds. 
\begin{triv} \label{fratExtH} There is an extension $G_1^*\to G$ with kernel 
$M(H)$ so the
pullback of $H$ in
$G_1^*$ is the natural quotient ${}_p\tilde H/\ker_1(H)$. \end{triv} Then, the
Lemma \ref{fratIsom} conclusion holds. This applies with $p=2$ and
$A_5=G$. 
\end{prop}

\begin{proof} The hypothesis says the morphism $\gamma'_H$, a priori an $H$ 
module
homomorphism, is actually a $G$ module homomorphism. Here is why. Suppose 
$G^\dagger$ is a
proper subgroup of $G_1^*$ mapping surjectively to $G$. Then, the pullback $H^*$ 
of $H$ in
$G^\dagger$ is a proper subgroup of the pullback of $H$ in $G_1^*$. Further, 
$H^*$ maps
surjectively to
$H$, contrary to ${}_p\tilde H/\ker_1(H)\to H$ being a Frattini cover. Conclude 
that
$G_1^*\to G$ is also a Frattini cover.  So, there is a natural map from $G_1=\tG
p/\ker_1(G)\to G_1^*$ inducing a surjective $G$ module homomorphism $M(G)\to 
M(H)$. 

Since ${}_p\tilde H/\ker_1(H)$ is 
the 1st characteristic quotient of ${}_p\tilde H$, universal for covers with 
elementary 
$p$-group kernel, there is an $H$ module splitting of 
$M(G)\to M(H)$. Higman's 
Theorem \cite[p.~64, Prop.~3.6.4]{Ben1} says, since $((G:H), p)=1$,  an  $H$ 
splitting of this $G$ map gives a $G$ splitting.  This is contrary to $M(G)$
being an indecomposable $G$  module unless this is an isomorphism (\cite[p.~11,
Exec.~1]{Ben1} or 
\cite[Indecom.~Lem.~2.4]{FrKMTIG}). 

Suppose a $p$-Sylow $P_p$ of $G$
has this property: Either 
$$gP_pg^{-1}= P_p \text{ or } gP_pg^{-1}\cap P_p=\{1\}$$ for each $g\in G$. 
\cite[Cor.~3.6.19]{Ben1} shows $$H^2(G,M(H))
\mapright{\scriptscriptstyle\text{\rm rest}} H^2(H,M(H))$$ is an isomorphism (so 
both have
dimension 1)  guaranteeing
\eqref{fratExtH}. This holds for $A_5$ and its 2-Sylows since they are 
distinguished by
which integer from $\{1,\dots,5\}$ each element in the 2-Sylow fixes. 
\end{proof}

\begin{rem}[Cases and extensions of Prop.~\ref{fratExt}] \label{induced} Suppose 
$G_k$, $k\ge 1$ is a higher
characteristic quotient of $\tG p$. Let $P_{p,k}$ be a $p$-Sylow of $G_k$. Then, 
the hypotheses of
Prop.~\ref{fratExt} automatically hold. 

Even when $k=0$, given $P_p$, suppose  $H_1\le G$ contains  $H=N_G(P_p)$ and is 
maximal for this property.  
\begin{triv} \label{fratExtH2} For some extension $H_1^*\to H$ with kernel 
$M(H)$, the 
pullback of $H$ in $H_1^*$ is the natural quotient ${}_p\tilde H/\ker_1(H)$. 
\end{triv}

As in \cite[Rem.~2.10]{FrMT}, let $\Ind(M(H))_{H_1}^{G}$ be the $G$ module
induced from $H_1$ acting on $M(H)$. Apply Shapiro's Lemma
\cite[Cor.~2.8.4]{Ben1}:
$H^2(H_1,M(H))=H^2(G,\Ind(M(H))^G_{H_1})$. So there is an extension of $G$ with 
kernel
$\Ind(M)^G_{H_1}$ whose pullback over $H_1$ has
$H_1^*$ as a quotient.  
\cite{FrRIMS} uses $\Ind(M(H))_{H_1}^{G}$ to produce the characteristic $p$-Frattini
module
$M_0=M(G)$ from the $p$-split case using indecomposability of $M_0(G)$ 
\cite[Indecom.~Lem.~2.4]{FrKMTIG}. \end{rem}

\subsubsection{Producing ${}_2^1\tilde A_5$} As in \cite[Part II]{FrMT}, 
producing ${}_p\tilde A_5$
for $p$ either 3 or 5 requires using the $G$ module 
induced from the $H$ module $M(H)$ (generalizing Prop.~\ref{fratExt}). 
Prop.~\ref{fratExt}
applies immediately to $A_5$  when $p=2$.

\begin{prop} \label{A5frat} Let $H=A_4=K_4\xs \bZ/3$. Then, $M(H)$ identifies 
with 
the $\bZ/2[H]$ module generated by the six cosets of a $\bZ/2$ in $A_4$, modulo 
the 
module generated by the sum of the cosets. Any $D_5$ in $A_5$ has a unique 
$\bZ/2$ 
lying in $A_4$. So, the action of $A_4$ on $\bZ/2$ cosets extends to an $A_5$ 
action on 
cosets of a dihedral group.  Thus, {\rm Prop.~\ref{fratExt}} gives ${}_2^1\tilde 
A_5$ as an 
extension of $A_5$ by $\ker_0(H)$.  \end{prop}

\newcommand{\mo}{{^{-1}}} \newcommand{\g}{{\bar\alpha}} \begin{proof} Suppose 
${}_p\tilde F_t$ is pro-free pro-$p$ group on $t$ generators. Let  $\psi: 
{}_p\tilde F_t\to 
P_p$ be a surjective homomorphism, with $P_p$ any (finite) $p$ group. Schreier's 
construction gives explicit generators of the kernel of $\psi$ 
\cite[\S15.6]{FrJ}. Apply this 
with $t=2$ and $P_p=K_4$, the Klein 4-group and $\ker_0=\ker(\psi)$. Let $\g$ be 
a generator of
$\bZ/3$. For $u$ and $v$ generators of ${}_p\tilde F_2$ let $\g$ act on 
${}_p\tilde F_2$ by mapping $(u,v)$ to $(v\mo,v\mo u)$.

Use $S=\{1, u, v\mo, u\mo v\}$ as coset representatives for $\ker_0$ in 
${}_p\tilde F_2$. 
Form the set $V$ of elements in $\ker_0$ having the form $tu s\mo$ or  $tv s\mo$ 
with  
$s,t\in S$. Toss from $V$ those that equal 1.  Now consider the images of $\g$ 
and $\g^2$ on
$uu=u^2=m_1\in V$.  This produces $v\mo v\mo=m_2$ and $u\mo v u\mo v=m_3$.

Consider $\g$ on $m_6= u(v\mo)^2 u\mo$. Recall: Modulo $\ker_1$ any two elements 
in 
$\ker_0$ commute. Apply this to get $$v\mo(u\mo v)^2 v=  (v\mo u\mo v)(u\mo)v^2 
\mod 
\ker_1= vu\mo v u\mo =m_5.$$ Apply $\g$ again to get $v\mo u^2v=m_4$. The action 
of 
$K_4\xs \bZ/3=A_4$ on the $m_i\,$s is the same as the action on the six cosets 
of an 
element of order 2.

Denote a commutator of two elements $w_1,w_2$ by $(w_1,w_2)$. Modulo $\ker_1$ 
there are relations among the $m_i\,$s: $m_1m_2m_3=(u,v) \mod \ker_1$; and 
$m_4m_5m_6=(u,v)\mo \mod \ker_1$. So, the product of the $m_i\,$s is 1. The 
proof 
follows from associating a $\bZ/2$ in $A_4$ with a dihedral in $A_5$ as in the 
statement 
of the proposition.  \end{proof}

\subsubsection{Summary of ${}_2^1\tilde A_5$ properties}\label{tA5prop} 
We finish a self-contained treatment of much of \cite[Part II]{FrMT}. The first 
characteristic quotient of the universal 2-Frattini cover of $A_5$ is a 
(nonsplit) extension 
of $\SL_2(\bZ/5)$ by an irreducible module V using $A_5\cong 
\PSL_2(\bZ/5)\cong\PSL_2(\bF_4)$. Prop.~\ref{A5frat} shows this as
follows. As  previously,  let $M(A_5)$  be the kernel of ${}_2^1\tilde
A_5\to A_5$.  Sums of $D_5$ cosets, $\sum_{i=1}^6 a_i T_i$, represent its 
elements. The 
augmentation map sends such an element to the sum $\sum_{i=1}^6 a_i$.

Let  $V$ be the 4-dimensional kernel of the augmentation map. Nonzero elements 
of $V$ 
have representatives $\sum_{i=1}^6 a_i T_i$ with two of the $a_i\,$s nonzero. 
Let 
$\bZ/3$ be the subgroup $\lrang{\g}\le A_5$ in the proof of Prop.~\ref{A5frat}.  
The 
action of $\bZ/3$ on $V$ is two copies of the 2-dimensional irreducible of 
$\bZ/3$. So, no element
of $V$ centralizes $\bZ/3$ and the centralizer of $\bZ/3$ in 
$M=M(A_5)$ is a $\bZ/2$. Denote ${}_2^1\tilde A_5$ by $G_1$. Let $P_1$ be a 
$p$-Sylow of $G_1$ with $P_0$, a Klein 4-group, its image in $A_5=G_0$. Denote 
the 
module for restriction of $M$ to a subgroup $H\le A_5$ by $M_H$. 

\begin{cor} \label{A5Morbs} Besides the origin there are three conjugacy classes 
(orbits 
for $A_5$ action) in $M$. These are (15) elements of $V\setminus \{0\}$;
10  representing sums $\sum_{i=1}^6 a_i T_i$ with exactly three $a_i\,$s nonzero 
(in
$M\setminus V$) ; and  six representing sums $\sum_{i=1}^6 a_i T_i$ with exactly 
one 
$a_i$  nonzero (in $M\setminus V$). Call the second set $M_3'$ and the
third $M_5'$. Elements  of $M_3'$ (resp.~$M_5'$) are exactly
those in $M\setminus \{0\}$ some  3-cycle (resp.~5-cycle)  stabilizes.
Correspond to each the respective 3 or 5-Sylow of $A_5$ that stabilizes
it. 

The action of $P_0$ on $M$ has a module presentation \begin{equation} 
\label{MA5K4} 
0\to J(P_0)\eqdef\one\oplus\one \to M(A_5)_{P_0}\to 
\one\oplus\one\oplus\one=M/J(P_0)\to 
0.\end{equation} Extending this action to $A_4$ gives a Loewy display  
\begin{equation} 
\label{MA5A4} 0\to U_3\to  U_3 \oplus \one\end{equation} with $U_3$ the two 
dimensional $A_4$ module on which $\bZ/3$ acts irreducibly. Restrict
$A_4$  to $M$ for an exact sequence $0\to J(A_4)\to M(A_5)_{A_4}\to
\one\to 0$. The Loewy display of $J(A_4)$ (resp.~ $\lrang{1_H+3_H+5_H, 
J(P_0)}$) is $U_3\to  U_3$ (resp.~$U_3\to \one$). 

If $g\in G_1$ has order 3 or 5 and $m\in M$ then, giving $mgm$ and stipulating 
$m \in M\setminus V$ determines $m$. Suppose $\alpha\in G_1$ has order 4. Then, 
the 
centralizer $Z_\alpha$ of $\alpha$ on $M$ has dimension 3 and on V has dimension
2, $\alpha^2\in M\setminus V$ and $Z_\alpha$ contains two elements each from
$M_3'$ and from
$M_5'$.
\end{cor}

\begin{proof} The action of $A_5$ preserves cosets. So, orbits for conjugation 
by $A_5$ 
are clear if  $V\setminus\{0\}$, $M_3'$ and $M_5'$ each consist of one $A_5$ 
orbit. 
This follows from triple transitivity  of $A_5$ in the standard representation. 

The exact sequence \eqref{MA5K4} comes from writing the right cosets of 
$H=\lrang{(1\,2)(3\,4)}$ in $A_4$. Use $1_H=H$, 
$2_H=\{(1\,3)(2\,4),(2\,3)(1\,4)\}$, 
$3_H=\{(1\,2\,3),(1\,3\,4)\}$, $4_H=\{(4\,3\,2),(1\,4\,2)\}$, 
$5_H=\{(1\,3\,2),(2\,3\,4)\}$ and 
$6_H=\{(1\,2\,4),(3\,1\,4)\}$. Then, $J(P_0)$ has generators $1_H+2_H$ and 
$3_H+4_H$. Images of $1_H$, $3_H$ and $5_H$ generate $M/J(P_0)$.  

Consider $\alpha\in G_1$ of order 4. To be explicit assume $\alpha=\alpha_{12}$ 
lifts 
$(1\,2)(3\,4)$. Then, $\alpha$ fixes $1_H$, $2_H$ and $3_H+4_H$, so it acts like 
$\beta_\alpha=(3\,4)(5\,6)$ on the cosets. The centralizer of $\alpha_{12}$ in 
$M(A_5)$ 
is  therefore $$U=\lrang{(x_1,x_2,a_1,a_1,a_2,a_2)\mid a_1,a_2,x_1,x_2\in 
\bZ/2}/\lrang{(1,1,1,1,1,1)}.$$  

An element in $U$ is in $M\setminus V$ if and only if $x_1\ne x_2$. Note: 
$\alpha$ 
fixes $1_H,2_h \in M_5'$. Sum both with $3_H+4_H$ to get two elements in $M_3'$.
The observation that $\alpha^2\in M\setminus V$ is a special case of
Prop.~\ref{liftEven}.  

Now consider the sequence \eqref{MA5A4}. Write $A_4$ as $K_4\xs\lrang{\beta}$ 
with 
$\beta=(1\,2\,3)$. First: $J(P_0)=\one\oplus\one$ is an $\lrang{\beta}$ module. 
As $\beta$ sends 
$1_H+2_H$ to $3_H+4_H$, it preserves the leftmost layer of the Loewy display. In 
the quotient $\one\oplus\one\oplus\one=M/J(P_0)$, $\one_H$ and $3_H$ generate a 
copy 
of $U_3$. Then, $A_4$ fixes $1_H+3_H+5_H \bmod J(P_0)$ and 
$\lrang{1_H+3_H,3_H+5_H}$ generates a copy of $U_3$. 

The module $V_{A_4}$ has Loewy display $U_3\to U_3$. If there is a subquotient 
of 
$M(A_5)_{A_4}$ having a Loewy display $U_3\to \one$,  it must be 
$T=\lrang{1_H+3_H+5_H, J(P_0)}$. Clearly this is invariant under $\bZ/3$. The 
action 
of $(1\,3)(2\,4)$ is to take $1_H+3_H+5_H$ to $2_H+4_H+5_H$. Since 
$1_H+3_H+5_H$ and $2_H+4_H+5_H$ differ by an element of $J(P_0)$, this shows 
$T$ is invariant under $\lrang{(1\,2\,3), (1\,3)(2\,4)}=A_4$. \end{proof}

\begin{rem}[A $Q_8$ version of Prop.~\ref{A5frat}] We take advantage of the 
special coset
aspect of module $M(A_5)$ throughout the study of our main example. The seed, 
however, is
an explicit presentation of $M(P_p)$ as a $P_p$ module, with $P_p$ a $p$-Sylow 
of $G$. The
proof of Prop.~\ref{hatP-P} shows, with illustration from $Q_8$ (versus $K_4$ 
above) an
effective way to compute the analog of 
\eqref{MA5K4}.
\end{rem} 

\subsection{Cusps and a theorem of Serre} As with modular curves, the 
Prop.~\ref{j-Line} branch cycle description tells much about the reduced cover 
of the $j$-line,  
especially about its points over $\infty\in \prP^1_j$, the {\sl cusps}. Cusp
widths  are the lengths of these orbits. An analysis of possible cusp widths 
starts with a general 
group computation. We refer to computations in \S\ref{twistAct}. 

\subsubsection{Lifting odd order elements} Suppose $g\in S_n$ has odd order (so 
necessarily $g\in A_n$). Write $g$ as a product $\prod_{i=1}^t \beta_i$ of 
disjoint cycles 
with $\beta_i$ of length $s_i$, $i=1,\dots,t$. Define $\omega(g)$ to be 
$\sum_{i=1}^t 
(s_i^2-1)/8$. Serre notes: $\omega(g)\equiv 0 \mod 2$ if and only if 
$\prod_{i=1}^t s_i 
\equiv \pm 1 \mod 8$. Let $\Spin_n\to A_n$ be the universal 
central exponent 2 extension of $A_n$. Regard the kernel of this map as 
$\{\pm1\}$. 


\begin{prop}[\cite{SeLiftAn}] \label{serLift} Assume entries of $\bg=(\row g r)$ 
(in $A_n$) have odd
order. Lift to $\Spin_n$ the entries of $\bg$, preserving their  respective 
orders. Denote the
lifted
$r$-tuple by $\bg^*$.  Suppose the following hold. 
\begin{edesc} \item Product-one condition: $\prod_{i=1}^r g_i=1$.  \item 
Transitivity: 
$\lrang{\bg}$ is a transitive subgroup of $A_n$.  \item Genus 0 condition: 
$\sum_{i=1}^r 
\ind(g_i)= 2(n-1)$.\end{edesc} Then, $s(\bg)=\prod_{i=1}^r g^*_i=(-
1)^{\sum_{i=1}^r\omega(g_i)}$.\end{prop}

The transitivity hypothesis in Prop.~\ref{serLift} is not serious. Restrict 
$\bg$ to each orbit 
of $\lrang{\bg}$. If the genus 0 condition applies in each such restriction, 
then $s(\bg)$ is 
the product of the $s$ values appearing in each restriction. The genus 0 
condition, 
however, 
is {\sl very\/} serious.  Examples of \S\ref{onlytwo} or \cite[Ex.~3.12]{FrMT} 
show how 
to use Prop.~\ref{serLift} to go beyond the genus 0 condition. 

A preliminary example! Take 
$\bg=((1\,2\,3),(1\,2\,3),(1\,2\,3))=(g_1,g_1,g_1)$. With $n=3$, the lifting 
invariant 
$s(\bg)$ is clearly $+1$: $\bg^*=(g^*_1,g^*_1,g^*_1)$ with $g_1^*$ of order 3. 
The conclusion,
however, of  Prop.~\ref{serLift} seems to be $s(\bg^*)=(-1)^3$. It, however, 
doesn't 
apply because the genus 0 condition doesn't hold.  

\S\ref{5-3-3prod} and Princ.~\ref{5-5-3princ} have cases of this 
computation without the 
genus 0 hypothesis; still by using Prop.~\ref{serLift}.

\subsubsection{Lifting even order elements} Suppose $g\in G$ has even order, and 
$\hat 
G\to G$ is a {\sl central\/} extension with kernel an elementary abelian 2-
group. Let $\hat 
g\in G$ be a lift of $g$ to $\hat G$. Then, the order of $\hat g$ is independent 
of the 
choice of the lift. The following applies the technique of proof from 
Prop.~\ref{serLift} to 
analyze the orders of lifts of elements of order 2 to $\hat A_n$. We 
compare this to an analog question with $G_1$, the first characteristic quotient 
of ${}_2\tilde
A_5$ replacing
$A_n$ (Lem.~\ref{nonsplitLem}).  

\begin{prop} \label{liftEven} Assume $n\ge 4$, and $g\in A_n$ of order 2 is a 
product of 
$2s$ disjoint 2-cycles. Any lift $\hat g\in \hat A_n$ of $g$ has order 4 if $s$ 
is odd and 2 
if $s$ is even. \end{prop} 

\begin{proof} We review the {\sl Clifford algebra setup\/} used in 
\cite{SeLiftAn}. 
Let 
$C_n$ be the Clifford algebra on $\bR^n$ with generators $\row x n$ subject to 
relations 
$$x_i^2=1,\ 1\le i\le n, \text{ and } x_ix_j=-x_jx_i \text{ if } i\ne j.$$ In 
the Clifford 
algebra, write $[i\,j]=\frac{1}{\sqrt 2} (x_i-x_j)$. Then, $[i\,j]^2=1$ and
$[i\,j]=-[j\,i]$. The   collection of $[i\,j]$ under multiplication generate a 
subgroup
$\hat S_n$. Characterization:  It is the central nonsplit extension of $S_n$ 
whose
restriction to transpositions splits, and  whose restriction to products of two 
disjoint
transpositions is nontrivial 
\cite[p.~97]{Se-GT}. The map $\hat S_n\to S_n$ appears from $[i\,j]\mapsto 
(i\,j)$.  


So, $\hat A_n =A_n\times \{\pm 1\}$ if $n\le 3$. That $\hat A_n\to A_n$ is 
nontrivial if
$n\ge  4$ shows from lifts of 
certain elements of order 2. Example: $(1\,2)(3\,4)$ lifts to have order 4: 
$$\Bigl(\frac{1}{\sqrt 2}
(x_1- x_2)\frac{1}{\sqrt 2}(x_3-x_4)\Bigr)^2=-[1\,2]^2[3\,4]^2=-1.$$ Of course 
the order of a
lift is conjugacy class invariant. Similarly, with
$n\ge 8$,
$$([\,12][3\,4]\dots [s\nm1\, s])^2=(- 1)^{2(s\nm2)}([\,12][3\,4])^2([5\,6]\dots 
[s\nm1\, s])^2.$$ By induction, the result is
$(-1)^s$: 
$[1\,2][3\,4]\dots [s\nm1\, s]$ has order $2^{1+\frac{1-(-1)^s} 2}$.   
\end{proof} 

\subsection{Further Modular Towers for $A_5$, $p=2$ and $r=4$} \label{onlytwo} 
Prop.~\ref{3355obst} gives a partial answer to \eql{compA5p2r4}{compA5p2r4a} 
about existence of
rational points at level 1 of other $A_5$ Modular Towers. There  are three  $2'$ 
classes: the two
classes of 5-cycles,
$\C_5^+$ (the class of 
$(1\,2\,3\,4\,5)$) 
and $\C_5^-$ (the class of $(1\,3\,5\,2\,4)$) and 3-cycles. From the \BCL, 
excluding 
$\bfC_{3^4}$, to give a Galois cover over $\bQ$ with $r=4$ requires $\bfC$ be 
$\bfC_{5_+5_-3^2}=(\C_5^+,\C_5^-,\C_3,\C_3)$ or $\bfC_{5_+^{2}5_-
^{2}}=(\C_5^+,\C_5^-,\C_5^+,\C_5^-)$. Both give Modular Towers with two 
components at level 0, one obstructed (nothing above it at level 1), one not. 

\begin{prop}\label{3355obst}  The only Modular Towers for 
$A_5$, $r=4$ and $p=2$ with possible $\bQ$ components at level 1 (for 
realizations of ${}_2^1\tilde
A_5$) have 
$\bfC=\bfC_{3^4}$, $\bfC_{5_+5_-3^2}$ or $\bfC_{5_+^{2}5_-^{2}}$. Each of the 
latter two has exactly two level 0 components, one obstructed, one not. There 
are no $\bR$ 
(so, no $\bQ$) points at level 1 of the $\bfC_{5_+5_-3^2}$ Modular Tower. 
\end{prop}

The next four subsections prove this by applying a lift to $\SL_2(\bZ/5)$ in
the two  new cases. For computational reasons, start with
$\bfC_{5_+5_-3^2}$. 

\begin{rem} \label{notinftyA5} Ex.~\ref{A5C54-wstory} notes the two genus 1 
components at level 1
on the Modular Tower for $(A_5,\bfC_{5_+^{2}5_-^{2}})$ over the level 0 
unobstructed component. If
have definition field $\bQ$, they might have infinitely many 
${}_2^1\tilde A_5$ realizations over $\bQ$. This is the only hope for infinitely 
many (up to
$\PSL_2(\bC)$ equivalence) $\bQ$ regular realizations of ${}_2^1\tilde A_5$ with 
$r=4$. Applying
\S\ref{indSteps} to this to see there will be at most finitely many realizations 
${}_2^2\tilde
A_5$ for $r=4$. \end{rem}  

\subsubsection{Setup for Prop.~\ref{3355obst}} \label{3355setup} Suppose 
$\bg\in \ni(A_5,\bfC_{5_+5_-3^2})$ is in the $H_4$ orbit $O_\bg$. Then, there is 
$\bg^*\in \ni({}_2^1\tilde A_5,\bfC_{5_+5_-3^2})$ over $\bg$ if and only if 
there is 
something from  $\ni({}_2^1\tilde A_5,\bfC_{5_+5_-3^2})$ over each element of 
$O_\bg$. If this holds, call $O_\bg$ unobstructed (\S\ref{IGPobstruct}). We list 
inner Nielsen
classes 
$\ni^\inn$ for $\bfC_{5_+5_-3^2}$ restricting to pieces where it is easy to 
demonstrate 
they lie in one $H_4$ orbit. With no loss use representatives for $\bg\in 
\ni^\inn$ with 
$g_1=(1\,2\,3\,4\,5)\in \C_5^+$ and $g_2\in \C_5^-$.

Conjugate the collection of 5-cycles by $g_1$. There are four length one and 
four
length five orbits.  Two length five orbits of $g_1$, $U_1$ and $U_2$, are 5-
cycles in
$\C_5^-$. There are  two cases: \begin{edesc}\label{pm33} \item \label{pm33a} 
$g_2=g_1^2$ or
$g_1^3$, 
$g_3$ and $g_4$ have 1 in their common support; or \item \label{pm33b} 
$g_2$ is a representative from $U_1$ or $U_2$.  \end{edesc} Here are $g_2$ 
representatives for \eql{pm33}{pm33b}: \begin{equation*} g_{2,1}=(2\,1\,3\,4\,5) 
\text{ 
and } g_{2,2}=(1\,2\,5\,4\,3).\end{equation*}

\subsubsection{Separating obstructed and unobstructed $H_4$ orbits} 
\label{5-3-3prod} 
We show the orbit of $\bg$ satisfying \eql{pm33}{pm33a} is obstructed.  Example: 
$\bg=((1\,2\,3\,4\,5),(1\,3\,5\,2\,4),(3\,5\,1),(2\,4\,1))$. Suppose 
$\bg'\in\ni({}_2^1\tilde 
A_5,\bfC_{5_+5_-3^2})$ lies over $\bg$.  Consider the image $\bg^*$ in 
$\SL_2(\bZ/5)$ of $\bg'$, as in proof of Lem.~\ref{cuspWidth}. Denote the 
product of 
the entries of $\bg^*$ by $s(\bg)$. Then, $g_2^*=(g_1^*)^2$. The product-one 
condition, $s(\bg)=1$, holds for $\bg^*$ if and only if it holds for 
$((g_1^*)^3,g_3^*,g_4^*)$. This is a lift of 
$((1\,2\,3\,4\,5)^3,(3\,5\,1),(2\,4\,1))$, by 
Riemann-Hurwitz a branch cycle description of a genus 0 cover.  

Apply  Prop.~\ref{serLift}. Entries of 
$((g_1^*)^3,g_3^*,g_4^*)$ have product $$(-1)^{\frac{25-1}8+2(\frac{9-1}8)}=-
1\in \Spin_5.$$ 
The product-one condition doesn't hold for $\bg^*$. So, it can't hold in 
${}_2^1\tilde A_5 
$. This is necessary for a Nielsen class element, concluding the proof of 
obstruction.

Now consider which elements of \eql{pm33}{pm33b} are obstructed. Check: 
\begin{edesc} \label{2122} \item \label{g21} $g_1g_{2,1}=(2\,4)(3\,5)$ limits 
$(g_3,g_4)$ to be conjugates of $((4\,3\,2),(5\,3\,2))$ by the Klein four group 
centralizer 
of 
$(2\,4)(3\,5)$.  \item \label{g22} $g_1g_{2,2}=(1\,5\,2)$ limits $(g_3,g_4)$ to 
be one of 
two types: \begin{itemize} \item conjugates of $((5\,3\,2),(2\,3\,1))$ by the 
centralizer of 
$(1\,5\,2)$; or \item both are in $\lrang{(1\,5\,2)}$. \end{itemize} \end{edesc} 
As above, 
consider the unique lift $\bg^*$ to $\SL_2(\bZ/5)$ of $\bg$ that could be an 
image from 
$\bg'\in \ni({}_2^1\tilde A_5,\bfC_{5_+5_-3^2})$ for the second of 
\eql{2122}{g22}. 
Since the lifts are unique, $g_4^*\in \lrang{g_3^*}$. Typically: $\bg^*$ would 
lift 
\begin{equation}\label{3repeat} 
\bg^\dagger=((1\,2\,3\,4\,5),(1\,2\,5\,4\,3),(1\,5\,2),(1\,5\,2)).\end{equation} 
Apply 
Prop.~\ref{serLift} to $((1\,2\,5\,4\,3),(1\,5\,2),(3\,4\,5))$: 
$s((1\,2\,5\,4\,3),(1\,5\,2),(3\,4\,5))=-1$. Therefore, $s(\bg^\dagger)=-
s((1\,2\,3\,4\,5),(5\,4\,3),(1\,5\,2))$. One last application of the same 
computation shows 
$s(\bg^\dagger)=(-1)^2$. To conclude $\bg^\dagger$ is unobstructed, use the 
following 
from \cite[\S3]{FrKMTIG}. The Loewy display  of $\ker_k/\ker_{k+1}$ refers to 
the $G_{k+1}$ action 
in the $p$-Frattini cover $G_{k+1}\to  G_k$. 

\begin{princ} \label{obstPrinc} Obstruction can occur from level $k$ to level 
$k+1$ in a
Modular Tower only where $\one_{G_k}$ appears in the Loewy display  of 
$\ker_k/\ker_{k+1}=M_k$. Translate this
to say there is a (Frattini cover) sequence $G_{k+1}\to U_k\to W_k\to G_k$ with 
$W_k$ and $\ker(U_k\to W_k)$ is
a trivial $W_k$ module. The appearance of $\one_{G_k}$ for any value of $k$ 
implies that 
$\one_{G_{k'}}$ appears in $M_{k'}$ for finitely many $k'$. If $G_0$ is $p$-perfect and
centerless, then so is
$W_k$ (argument of Prop.~\ref{fineMod}). 
\end{princ}

\begin{proof} Everything has an explanation already except the automatic 
appearance of $\one_{G_{k'}}$ in
$M_{k'}$ for finitely many $k'$ given its appearance at level $k$. Apply 
Prop.~\ref{RkGk} with $W_k$ replacing $G_k$. Then, use that the universal $p$-Frattini cover of
both $W_k$ and $G_k$ is the same as that of $G_0$. So, the characteristic quotients for the 
universal $p$-Frattini of $W_k$
(in place of $G_0$) are cofinal in the projective system $\{G_k\}_{k=0}^\infty$. 
Therefore the simple modules
(including copies of the identity representation) appearing in the analog for 
$W_k$ of the modules $M_k(G_0)$
also appear as simple modules for a cofinal collection from 
$\{G_k\}_{k=0}^\infty$.\end{proof} 

Since $\one_{A_5}$ appears only at the head
(in the
$A_5$ Schur multiplier) of 
$\ker_0/\ker_1$ (\S\ref{pfrattini}), $\bg^\dagger$ passes the lifting test to 
$\SL_2(\bZ/5)$ 
exactly if there is $\bg'\in \ni_1$ over it.

\subsubsection{Conclusion for $\bfC_{5_+5_-3^2}$} The principle from 
\S\ref{5-3-3prod} is the following, where $s(\bg)$ is as in Prop.~\ref{serLift}.

\begin{princ} \label{5-3-3princ} Suppose for $\bg\in A_n^r$, the product-one 
condition 
holds and there exists $i,j$ with $g_i$ a 3-cycle, $g_j$ a 5-cycle and their 
product (in 
either order) is a 3-cycle. To compute $s(\bg)$ we may assume $i=1$, $j=2$, and 
$s(\bg)=-s(g_1g_2,g_3,\dots, g_r)$. If $n=5$, $s(\bg)=1$ if and only if $\bg$ is 
unobstructed. \end{princ}

\begin{proof} With no loss assume $i< j$. The exists $Q\in H_r$ braiding $\bg$ 
to 
$\dot{\bg}$ in the same Nielsen class with $\dot{g}_1=g_i$ and $\dot{g}_j=g_2$. 
From 
\cite[Part III]{FrMT}, $s(\bg)=s((\bg)Q)$. (More generally, 
$\nu((\bg)Q)=\nu(\bg)$ 
where $\nu$ is the {\sl big\/} braid invariant of a Nielsen class.)  Apply the 
argument of 
\S\ref{5-3-3prod} to $(\dot{g}_1,\dot{g}_2, (\dot{g}_1,\dot{g}_2)^{-1})$ to 
reduce  
computing $s(\bg)$ to computing $-s(\dot{g}_1\dot{g}_2,\dot{g}_3,\dots, 
\dot{g}_r)$. 
\end{proof}

Consider a case of \eql{2122}{g21}: 
$\bg=((1\,2\,3\,4\,5),(2\,1\,3\,4\,5),(4\,3\,2),(5\,3\,2))$. 
Princ.~\ref{5-3-3princ} applies to $i=1$ and $i=3$. The product 
$g_1g_3=(1\,4\,5)$ is a 
3-cycle. Thus, $s(\bg)=-s(\bg^\dagger)$ where $\bg^\dagger$ is a 3-tuple 
consisting of 
two 3-cycles and a 5-cycle. As previously, $s(\bg^\dagger)=-1$ and $s(\bg)=1$. 
Since 
$n=5$, $\bg$ is unobstructed.

Finally, consider the first case of \eql{2122}{g22}. Note: The genus 0 
hypothesis of 
Prop.~\ref{serLift} doesn't hold.

\begin{princ} \label{5-5-3princ} For each ordering of the conjugacy classes 
$\bfC_{5_+5_-3}$, the nielsen class $\ni(A_5,\bfC_{5_+5_-3})^\inn$ 
has exactly one element, for a total of six elements. All representatives $\bg$ 
have $s(\bg)=1$. \end{princ}

\begin{proof} That there is only one $\bg\in\ni(A_5,\bfC_{5_+5_-3})^\inn$ with 
$g_1\in 
\C_{5}^+$, $g_2\in \C_{5}^-$ and $g_3\in \C_3$ follows from the data in 
\eqref{2122}. The 
computation $s(\bg)=1$ comes from 
$s(\bg^\dagger)=1$ in \eqref{3repeat}.  \end{proof}

Let $\bg=((1\,2\,3\,4\,5),(1\,2\,5\,4\,3),(5\,3\,2),(2\,3\,1))$. As 
above, if $s(\bg)=-1$, then this is obstructed. Princ.~\ref{5-5-3princ} shows  
$s(\bg)=s((1\,5\,2),(5\,3\,2),(2\,3\,1))$. Since the entries satisfy the genus 0 
hypothesis, 
Prop.~\ref{serLift} gives this value as -1.

Thm.~\ref{thm-rbound} shows level 1 of the Modular Tower for $({}_2^1\tilde 
A_5,\bfC_{5_+5_-3^2})$ has no real points; that would require H-M reps.~at 
level 0 (see \eqref{HMdef}). It has none, for the inverse of an element of
$\C_5^+$ is  in $\C_5^+$.

\subsubsection{The case $\bfC_{5_+^{2}5_-^{2}}$} This brings a new issue. There 
are 
H-M representatives: 
$((1\,2\,3\,4\,5),(5\,4\,3\,2\,1),(2\,1\,3\,4\,5),(5\,4\,3\,1\,2))$. So 
some 
elements are unobstructed. For an obstructed element in the Nielsen class use  
$\bg'\in 
\ni(A_5,\bfC_{5_+^3})$, three repetitions of the $C_5^+$ conjugacy class. (Note: 
\eqref{2122} shows $\ni(A_5, \bfC_{5_+^25_-})$ is 
empty.) Write the first entry $g_1'$ as $h^2$. Take $g_1^*=g_2^*=h, g_3^*=g_2'$, 
and 
$g_4^*=g_3'$ to produce $\bg^*\in \ni(A_5, \bfC_{5_+^{2}5_-^{2}})$. Clearly 
$s(\bg^*)=s(\bg')$. Example: 
$$\bg'=((1\,2\,3\,4\,5),(3\,5\,1\,4\,2),(4\,5\,2\,3\,1)).$$ 
The next principle gives an obstructed element in $\ni(A_5,\bfC_{5_+^2 5_-^2})$.

\begin{princ} \label{5-5-5princ} The Nielsen class $\ni(A_5,\bfC_{5_+^3})^\inn$ 
has one element. For $\bg$ in this class, 
$s(\bg)=-1$. \end{princ}

\begin{proof} Write $\bg=(g_1,g_2,g_3)$ with $g_1=(1\,2\,3\,4\,5)$ and $g_2$ a 
conjugate of $g_1$ not in $\lrang{g_1}$. Choose $g_2$ in one of the orbits $U_3$ 
or 
$U_4$ (for conjugation by $g_1$; \S \ref{3355setup}) in $C_5^+$.  
If we choose $U_3$ to be the orbit of $g_2'$ (in $\bg'$), then conjugate by 
$\lrang{g_1}$ 
to see $g_3'$ is also in $U_3$. The other $\lrang{g_1}$ orbit of 5-cycles, 
represented by 
$(2\,4\,1\,5\,3)$, doesn't give an element in the Nielsen class: The product 
$(1\,2\,3\,4\,5) 
(2\,4\,1\,5\,3)$ is a 3-cycle. 

Calculate $s(\bg')$ by  
considering $\bg''=((1\,2\,3),\bg',(3\,2\,1))$. Then, $s(\bg'')=s(\bg')$. Use 
Princ.~\ref{5-5-3princ}, then Princ.~\ref{5-3-3princ} to see $s(\bg'')$ equals 
$$s((1\,3\,2\,4\,5),(3\,5\,1\,4\,2),(4\,5\,2\,3\,1),(3\,2\,1))= -
s((1\,3\,2\,4\,5),(3\,5\,1\,4\,2),(4\,5\,1)).$$ One last application of 
Princ.~\ref{5-5-3princ} 
now shows $s(\bg)=-1$.  \end{proof}

\begin{rem}  Noting the position of a 3-cycle gives 
$|\ni(A_5,\bfC_{5_+^23})^\inn|=3$.  A representative 
for the $H_3$ orbit is $((1\,2\,3\,4\,5),(1\,2\,4),(3\,2\,5\,4\,1))$. Put 
$(3\,2\,1)$
on the left and 
$(1\,2\,3)$ on 
the right as in Princ.~\ref{5-5-3princ}. Conclude its lifting invariant is -
1.\end{rem}

\subsection{Modular Towers from $p$-split data} Let $G_0$ be 
$(\bZ/p)^u\xs H$ where $H$ is a $p'$-group acting irreducibly (and nontrivially) 
through 
$\GL_u(\bZ/p)$. This presents $G_0$ as a primitive affine group where 
$(\bZ/p)^u=M$ 
are the letters of the permutation representation. As in Ex.~\ref{HL}, $\tG p$ 
is ${}_p\tilde 
F_u\xs H$ with ${}_p\tilde F_u$ the pro-free pro-$p$ group on $u$ generators. 
(Action 
of $H$ extends to ${}_p\tilde F_u$ \wsp Remark \ref{extH}\wsp uniquely up to 
conjugacy. Still, it will rarely be easy to find.) Let $G_k$ be the $k$th 
characteristic 
quotient of $\tG p$, $k\ge 0$. 

Assume given $p'$ conjugacy classes $\bfC$ in $H$. The next lemma shows how to 
find 
lifts of an element $\bg\in \ni(H,\bfC)$ to $\bg^*\in \ni(G_0,\bfC)$, ensuring 
$\ni(G_0,\bfC)$ is nonempty.  

\begin{lem}\label{HLcont} As above, let $\bg^*\in \bfC$ have $g_i^*=(m_i,g_i)$ 
with 
$m_i\in M$, $i=1,\dots,r$. Then, $\bg^*\in \ni(G_0,\bfC)$ if and only if 
\begin{edesc} 
\label{condLift} \item  \label{condLifta} there is no $m\in M$ with $m-m^{g_i^{-
1}}=m_i$, $i=1,\dots,r$, and \item \label{condLiftb} $m_1+m_2^{g_1^{-1}}+\cdots+ 
m_r^{g_1^{-1}\cdots g_{r- 1}^{-1}}=0$. \end{edesc} So, there is a unique 
solution for $m_r$ given any elements $\row m {r\nm1}$. \end{lem}

\begin{proof} Condition \eql{condLift}{condLiftb} is the product one condition 
for 
$\bg^*$. Since $M$ is an irreducible $H$ module, given $\row m {r\nm1}$, there 
is a 
unique solution for $m_r$. From the same argument, $\lrang{\bg^*}=G_0$ if and 
only if 
the projection of $\lrang{\bg^*}$ onto $H$ has no kernel; otherwise the kernel 
is all of 
$M$ since it is a nontrivial $H$ invariant submodule of $M$. Conclude, 
$\lrang{\bg^*}$ is a splitting of $H$ into $G_0$. From Schur-Zassenhaus, it is 
conjugate 
to the canonical copy of $H$ in $G_0$. Condition \eql{condLift}{condLifta} is 
exactly 
the computation for that. \end{proof} 

Levels 0 and 1 of the $(A_4,\bfC_{3_+^23^2_-},p=2)$
Modular Tower (Ex.~\ref{A4C34-wstory}) show the $p$-split situation, for $u>1$, 
is a serious
challenge for our Main Conjecture. 

\subsection{Criteria for $G_k$ faithful on $M_k$ and appearance of $\one_{G_k}$} 
\label{genFrattini} These
general comments on the universal $p$-Frattini cover appear here for reference 
in \cite{AGR}. Let $O_{p'}(G)$ be
the maximal $p'$ normal subgroup of any finite group $G$. 

\begin{lem} \label{p'normal} Let $H$ be any $p'$ subgroup of $G$ and let 
$M_0(H)$ be restriction of $G$ on
$M_0$ to $H$. Then, this action extends to an action of $H$ on $\ker_0$, unique 
up to conjugation inside
${}_p\tilde G$. Such an $H$ acts trivially on $\ker_0$ if and only if it acts 
trivially on $M_0(H)$. If
$H=O_{p'}(G)$, then this action is trivial. \end{lem}

\begin{proof} The extending action of $H$ to $\ker_0$ is a special case of 
Rem.~\ref{extH} applied to
$M_0(H)\xs H$. The extending group $\ker_0\xs H$ is the universal $p$-Frattini 
of $M_0(H)\xs H$. If the action
on
$M_0(H)$ is trivial, then the trivial extension of $H$ to $\ker_0$ is an 
extending action. 

Consider the case $H=O_{p'}(G)$. We have the natural short exact sequence $1 \to 
O_{p'}(G) \to \tG p\to {}_p
\widetilde {G/O_{p'}(G)}\to 1$. Form the (group) fiber product $${}_p\tilde 
G^\dagger\eqdef {}_p \widetilde
{G/O_{p'}(G)}\times_{G/O_{p'}(G)} G.$$ By the basic property of the $p$-Frattini 
cover,  ${}_p \widetilde
{G/O_{p'}(G)}$ is the minimal
$p$-projective cover of 
$G/O_{p'}(G)$.  Since $G\to G/O_{p'}(G)$ has
$p'$ kernel, ${}_p\tilde G^\dagger$ is the minimal  
$p$-projective cover of $G$. So, by Prop.~\ref{pprojchar}, it is $\tG p$.  
\end{proof}

Warning! Each  $p'$ subgroup $H$ of $G$ does have its action extend to $\ker_0$, 
and that extending action 
is unique up conjugacy. Assuming a group is $p$-perfect then means we can extend 
the action of generators on
$\ker_0$. Those actions, however,  {\sl won't\/} fit together to have $G$ act on
$\ker_0$ (as in \S~\ref{onesAppearNGP}). 

In Prop.~\ref{psplitOne} we found many explicit appearances of $\one_{G_k}$ in 
$M_k$ in the $p$-split case, excluding the case that $M_0$ is cyclic (with 
nontrivial $G_k$) action. 
The following result, basically from \cite{GriessFrat}, characterizes when $G_0$ 
is faithful on $M_0$. It is an ingredient in Thm.~\ref{densityOnes}, giving an 
asymptotic formula for the multiplicity of appearance of any simple
$G_0$ module in $M_k$ for large $k$. 

\begin{thm} Suppose the dimension of $M_0$ exceeds 1. Then,  
$G_k/O_{p'}(G_0)$ acts faithfully on $M_k$, for all $k\ge 0$. \end{thm}

\begin{proof} Once we know the rank of $M_0$ exceeds 1, then so does the rank of 
$G_k$ on $M_k$ for all $k$.
Since $O_{p'}(G_0)$ is a $p'$ normal subgroup, there is a splitting of the 
pullback of $O_{p'}(G_0)$ in $\tG
p$. That identifies $O_{p'}(G_0)$ as the maximal normal $p'$ subgroup of $G_k$ 
for each $k$. Given the
Griess-Schmid result for $k=0$ therefore gives it for all $k$. \end{proof}

With the previous notation, assume $\one_{G_k}$ appears in the Loewy display of
$M_k$. Princ.~\ref{obstPrinc} shows $\one_{G_{k'}}$ appears (even {\sl where\/}) 
in the Loewy display 
for infinitely many integers $k'\ge k$. Thm.~\ref{densityOnes} is in  \cite{Darren}. As a
special case it shows that $\one_{G_k}$ appears with an explicit positive 
density in $M_k$ for $k$ large. The
result is effective, though for small values of $k$ it is subtle to predict the 
appearance of $\one_{G_k}$.
Further, it is imprecise  on the place in the Loewy display of the appearance of 
the $\one_{G_k}\,$s.
So, Princ.~\ref{obstPrinc} and the very effective  $p$-split case remain 
valuable. 

Recall: Over an
algebraically closed field the set of simple $G_0$ modules has the same 
cardinality as the set of $p'$
conjugacy classes. Let $S$ be any simple $G_0$ module. Let $K$ be an 
algebraically closed field and
retain the notation $M_k$ after tensoring with $K$. We use $\lrang{S,M_k}$, and 
related compatible notation, 
for the total multiplicity of $S$ in all  Loewy layers of the $G_k$ module 
$M_k$. 

\begin{thm}[Semmen] \label{densityOnes} The condition $\dim_K(M_0) \ne 1$ is
equivalent to
$$\lim_{n\mapsto
\infty}
\frac{\lrang{S,
M_k} }{\dim_K( M_k)} = \frac{ \lrang{S, K[G/ O_{p'}(G)]}}{\dim_K(K[G/
O_{p'}(G)])}.$$
\end{thm}
 
\begin{rem}[Griess-Schmid characterization that $\dim_{\bZ/p}(M_0)=1$] Let 
$O_p(G)$ be the maximal normal $p$-subgroup of
$G$. Recall that a group $H$ is $p$-supersolvable if $G/O_{p}( G/O_{p'}(G))$ is 
abelian of exponent dividing $p-1$. 
\cite{GriessFrat}
has the following characterization of $\dim_{\bZ/p}(M_0)=1$: G is a $p$-supersolvable group with
a cyclic
$p$-Sylow. Specifically, if $G_0$ is simple, then $\dim_{\bZ/p}(M_0)\ne 1$. 
\end{rem}

\section{Real points on Hurwitz spaces} \label{startHM} Investigations of
modular curves often analyze behavior of functions near cusps. 
Degeneration behavior of the curve covers associated to points on Hurwitz 
spaces, as the
points approach the cusps,  hints at diophantine properties. This section 
illustrates with
detailed analysis of the {\sl easiest\/} case: Degeneration behavior of real 
points on Hurwitz and
reduced Hurwitz spaces. So doing, it points to special {\sl cusps\/} named for 
this degeneration.
The word {\sl cusp\/} implies some compactification of the Hurwitz space. 
Therefore, after
generalities on Hurwitz spaces for any $r$, we concentrate on cusp behavior of 
reduced Hurwitz
spaces for the case
$r=4$. Then the spaces are curves, and the compactification behavior stays 
within the
confines of this paper. In Modular Tower higher levels for 
$p=2$, real points on Hurwitz spaces at those levels associate to branch cycles 
we call H-M
and near H-M reps.  Cusps at the end of the corresponding components of real 
points then
inherit the H-M and near H-M moniker. 

\subsection{Diophantine properties of a Modular Tower} In Thm.~\ref{thm-rbound},
$\sH_k$ is level $k$ of a Modular Tower of {\sl inner\/} Hurwitz  spaces.  
In assuming $G_0$ is centerless and $p$-perfect  (Def.~\ref{pperfect}), these 
are fine moduli
spaces (Prop.~\ref{fineMod}). As a diophantine corollary, over {\sl 
any\/} field $K$ with $(\text{char}(K),|G_0|)=1$ (see \cite[Reduction Hypothesis 
E.3]{FrMT} 
or \cite{WeTh}), $\bp\in \sH_k(K)$ if and only if some cover (in the inner 
equivalence 
class) associated to $\bp$ has definition field $K$.

\subsubsection{H-M and near H-M representatives}
A {\sl  Harbater-Mumford\/} 
(H-M) representative of Nielsen class
$\ni(G,\bfC)$ is an
$r$-tuple
$\bg\in\bfC$  with the following property: \begin{triv} \label{HMdef}
$g_{2i-1}=g_{2i}^{-1}$, 
$i=1,\dots, s$, with $r=2s$. \end{triv} \noindent Call a component of 
$\sH(G,\bfC)^\inn$
an H-M component if it corresponds to an $H_r$ orbit of an H-M rep(resentative). 

The definition of a {\sl near H-M\/} representative applies to any Modular Tower 
at level $k\ge 1$
if $p=2$,  and $r=2s$ is even.  It uses an operator $\hk$ on Nielsen classes, 
the special case
$r_1=0$  for the operator $\hk$ in \eqref{compCong}. Define the effect of $\hk$ 
on 
$\bg\in\ni(G,\bfC)$ by consecutively listing entries of $\hk(\bg)$. List entries 
of $\bg$ as 
$(g_{1,1},g_{1,2},\dots, g_{s,1},g_{s,2})$. As in \eqref{compCong}, denote 
$g_{j,1}g_{j,2}$ by $g_{j,.}$.  \begin{equation}\label{hkCong}\begin{array}{rl} 
\hk(\bg)_{1,1}&= (g^\sph_{2,.}\cdots g^\sph_{s,.})^{-1} g^{-1}_{1,2} ( 
g^\sph_{2,.}\cdots g^\sph_{s,.}),\\ \hk(\bg)_{1,2}&=( g^\sph_{2,.}\cdots 
g^\sph_{s,.})^{- 1} g^{-1}_{1,1} ( g^\sph_{2,.}\cdots g^\sph_{s,.}),\\ 
&\ldots,\\ 
\hk(\bg)_{s,1}&= g^{-1}_{s,2},\ \hk(\bg)_{s,2} = g^{-
1}_{s,1}.\end{array}\end{equation}

Then,  $\bg^*\in \ni(G_k,\bfC)^\inn$, $k\ge 1$, is a near H-M rep.~if the 
following hold.  
\begin{edesc} \label{nearHM} \item $\bg^*$ isn't an H-M representative.  \item  
There 
exists an involution $c\in G_k$ satisfying $c\bg^*c=\hk(\bg^*)$.  \item $\bg^* 
\bmod 
\ker_{k-1}$ is an H-M representative.  \end{edesc}

\subsubsection{A general diophantine statement about Modular Towers} 
\label{diophMT}
We use \S\ref{rboundProof} for direct calculation of real points on 
reduced Hurwitz spaces. Reference to $\bfC_k$ means conjugacy classes varying 
with 
$k$, not necessarily consisting of $p^\prime$ classes. For any prime $\ell$, let 
$\bQ^{(\ell)}$ denote the field with all $\ell^n$th roots of 1, $n=1,2,\dots$, 
adjoined. If 
$Y\subset \sH_0$, denote the complement of $Y$ pulled back to $\sH_k$ by 
$\sH_k\setminus Y$.  

\begin{thm} \label{thm-rbound} Assume $G=G_0$ is centerless and $p$-perfect.
Fix $r^*$ and a subfield $K\le \bar\bQ$.  Suppose there are  $(G_k,\bfC_k)$  
realizations (over $K$)
with $r_k\le r^*$ conjugacy classes in $G_k$ as the entries of $\bfC_k$, for 
each $k\ge  0$. 
\begin{triv}
\label{assumep} Assume $[\bQ^{(p)}\cap K:\bQ]< \infty$. \end{triv} 
\noindent Then, the only possibility for such $K$ realizations is there exists 
$r$, with $r\le
r^*$, $p'$ classes $\bfC$ with $(G_k,\bfC)$ realizations (over $K$) for all $k$.

Assume there is a prime $\ell$, $(\ell,|G_0|)=1$,  satisfying the following.
\begin{triv} 
\label{assumel} At some place $P$ of $K$ over $\ell$ the residue class field 
$F_P$ is 
finite.\end{triv} \noindent Then, reduction of $\sH_k$ at $P$ has no $F_P$ 
points  for 
$k$ large. 

Suppose there is a proper algebraic subset $Y\subset \sH^\rd_0$ with this 
property. 
\begin{triv} \label{propSubset} $\sH_k^\rd\setminus Y(K)$ is finite for each 
$k$. 
\end{triv} Then,  for some $k_{K,G_0,\bfC}$,  $\sH_k^\rd\setminus
Y(K)$  is empty for $k\ge k_{K,G_0,\bfC}$. Assume $r=4$ and all
components of 
$\sH^\rd_{k'}$ have genus at least 2 for some $k'\ge 0$. Then, there is
$k_{K,G_0,\bfC}> 0$ with  no $(G_k,\bfC)$ realizations over $K$ for $k\ge 
k_{K,G_0,\bfC}$. 

Assume $p=2$ and $K\subset\bR$. Then all $(G_k,\bfC)$ realizations over $K$ for 
$k\ge 
1$ appear in components of $\sH_k$ for $H_r$  orbits containing H-M or near H-M 
reps. 
(This holds even replacing $\bfC$ by $\bfC_k$ varying with $k$.) A $\bp\in 
\sH_k(\bR)$ 
corresponds to a $(G_k,\bfC)$ realization over $\bR$. A connected component of  
$\sH_k(\bR)$ associated to an H-M rep.~has a connected component of 
$\sH_{k+1}(\bR)$
above it. A connected  component of $\sH_k(\bR)$ for a near H-M rep.~has
nothing in 
$\sH_{k+1}(\bR)$ above it. \end{thm} 

Explicit formulas (as in \eqref{RHOrbitEq}, Lem.~\ref{RHOrbitLem} and 
\S\ref{indSteps}) 
give us confidence in deciding \cite[Main 
Conj.~0.1]{FrKMTIG} when $r=4$.  
If the \cite[Main  Conj.~0.1]{FrKMTIG} conjecture is correct, that $\sH_k^\rd$ 
has no rational
points for
$k$ large,  then $\bQ$ realizations of $G_k$ require increasingly large sets of 
conjugacy 
classes as $k$ grows. This is subtler than information from  the \BCL. 

\S\ref{ptsMT} succinctly interprets points on a Modular Tower as a statement on 
quotients of a
fundamental group. Thm.~\ref{thm-rbound} says if $p=2$, real points on a Modular 
Tower (versus points at various levels) appear only on projective sequences of 
H-M 
components. In \S\ref{diophimp} we see this gives a tool for progress on the 
Main 
Conjecture. 

The proof of Thm.~\ref{thm-rbound} appears in two subsections. 
\S\ref{Kptbound} gives the argument bounding points at high levels on a Modular 
Tower 
of inner Hurwitz spaces. \S\ref{rboundProof} gives details on $\bR$ points. 
Note: 
Hurwitz spaces $\sH_k$ are affine varieties. The theorem says nothing about 
rational 
points on (any) closure of $\sH_k$ including the boundary. 

\subsection{Bounding levels with $K$ points} \label{Kptbound}  
\cite[Thm.~4.4]{FrKMTIG} gives the first statement of Thm.~\ref{thm-rbound}. 
Bounding 
$k$ with  $(G,\bfC_k)$ having $K$ realizations and $\bfC_k$ not consisting of 
$p'$ 
classes requires only the \BCL. This is effective, though dependent on data 
about $\tG p$ for a result referencing only $r^*$, $K$ and $G_0$.

Consider any suitable Modular Tower for inner equivalence 
where $G=G_0$ is  
$p$-perfect and centerless. Then, the inner Hurwitz spaces are fine moduli 
spaces
(Prop.~\ref{fineMod}). Also, they and their reduced versions have good reduction 
modulo any prime
$\ell$, if 
$(\ell, |G|)=1$ (\cite[p.~167]{FrMT} or with more details \cite{WeTh}). 
Use  $$\cdots \to \sH^\rd_k\to
\sH^\rd_{k-1}\to \cdots \to \prP^1_j\setminus \infty$$ for
the  sequence  of reduced spaces. Assume each 
$\sH_k^\rd$ has a  $K$ point. Let $P$ be a  prime of $K$ over $\ell$ with 
residue class 
field $F_P$. 

Suppose, contrary to the conclusion of the theorem, $\sH^\rd_k$ has an $F_P$ 
rational 
point for each integer $k$. The set $\sH^\rd_k(F_P)$ is finite for each $k$. So, 
these 
finite 
sets form  a nonempty projective system. Conclude: There is a projective 
sequence of 
points $\bp^\rd$ on them.  Let $F'$ be a finite extension of $F_P$ over which 
the
point $\bp_0^\rd$ produces an inner cover $X_{\bp_0}\to \prP^1_z$ in the Nielsen 
class
corresponding to $\bp_0$ on the level 0 Hurwitz space lying over $\bp^\rd_0$. 
Apply
Princ.~\ref{betweenMTrdMT} when
$K$ is a finite field. That translates  to a projective system of $(G_k,\bfC)$ 
regular
realizations with definition field $F'$. With no loss, 
take this finite extension to be 
$F_P$.  

Let 
$T_p(X_{\bp_0})=T_p$ be the Tate  module for $X_{\bp_0}$ formed by taking an 
abelian  quotient of the Modular Tower as in Prop.~\ref{GKaction}. Suppose the 
kernel of 
the universal $p$-Frattini  cover has rank $u$. For a general Modular Tower this 
implies  
the Frobenius for $F_P$  acts trivially on a rank $u$  quotient of $T_p$. 
Contradiction: 
The Frobenius has eigenvalues of absolute value  $|F_P|^{1/2}$. This 
contradiction shows 
there cannot be $F_P$ points on the reduced space $\sH^\rd_k$ for all $k$. 

Now assume \eqref{propSubset} holds. There are finitely many points on 
$\sH_k^\rd\setminus Y(K)$. So, the argument above shows that either
$\sH_k^\rd\setminus  Y(K)$ is empty for some value of $k$ or there is a 
projective
system of $K$ points.   Taking a suitably large prime $\ell$ allows reducing the
covers for the projective system  modulo $\ell$, and the proceeding argument 
completes
the proof of this part of the  theorem. 

\begin{rem}[Thm.~\ref{thm-rbound} and absolute classes] The argument of 
Thm.~\ref{thm-rbound} does not hold for a Modular Tower of absolute Hurwitz 
spaces. 
Modular curves give an example of this. Consider a fixed prime $p$, and all of 
the curves 
$Y_0(p^{k+1})$ (\S\ref{compMC}) as in a Modular Tower. Let $\ell\ne p$ or 2 be a
prime, and reduce modulo $\ell$. 

The following facts are in \cite[Chap.~3]{RoEC}. There are many supersingular  
elliptic 
curves in characteristic $\ell$, roughly $\ell/12$ of them. Further, all have 
field of 
definition $K=\bF_{\ell^2}$. Characterize a supersingular curve (over 
$\bar\bF_\ell$) by 
its having no $\ell$ division points. So, any curve isogenous to it, say by a 
cyclic degree 
$p^{k+1}$ isogeny, is also supersingular.  Since both curves have definition 
field $K$, the 
isogeny has definition field $K$ (though it is not a Galois cover over $K$). 
This gives a 
large number of points on $X_0(p^{k+1})$,  as many as can be expected for $k$ 
large. 
This is Ihara's first example for producing many points on high genus curves 
over a
finite  field of square order. Ihara uses {\sl Shimura curves\/} to do the same 
over
finite fields of  order $\ell^{2t}$ for $t>1$. (\cite{IharSC} has an exposition 
on
this result's scattered  literature.) The Shimura curves he uses are {\sl 
compact\/}
families of abelian varieties, unlike levels of a Modular
Tower which have cusps.
\end{rem}
\newcommand{\sr}{{r_1+r_2}}

\subsection{Remainder of Thm.~\ref{thm-rbound} proof} \label{rboundProof} We
apply \cite{DFrRRCF} using the pattern of \cite[App.~C and 
App.~D]{FrMT}. This explicitly described real points on any of the fine moduli 
Hurwitz 
spaces $\sH$ appearing in this paper.

\subsubsection{Connected components of
$\sH(\bR)$} Let 
$\phi:X\to\prP^1_z$ be a cover over $\bR$, branched over $\bz$, with
$z_0\in \bR$.  
Specific classical generators of $\pi_1(U_\bz,z_0)$ produce an
explicit uniformization 
of any real points on $X$ over $z_0$. This uses a complex 
conjugation operator $\hk$ from \eqref{compCong} on Nielsen classes.
It includes
giving a  combinatorial restatement for $X$ having definition field
$\bR$.  So, 
fixed  points of this complex conjugation operator produce points in
$\sH(\bR)$. With a
fine  moduli assumption these are exactly the points of $\sH(\bR)$.
The process  works efficiently for unramified Frattini extensions of the covers 
$X\to
\prP^1$ appearing in an inner family
 (using the language of \S\ref{frattiniPre}). So, this gives a useful 
description of real
points on a Modular Tower.

\subsubsection{Points in $\bz$ over $\bR$} The approach falls into cases
from 
how the support of $\bz$ behaves under complex conjugation.  As $\bz$
has definition 
field $\bR$, it has $s_1$ complex conjugate pairs and $r_1=r-2s_1$ real
points. To 
simplify, call the latter $\row z {r_1}$, arranged left to right on the
real line and
the  former,  $z_{r_1+1},\dots, z_r$. With no loss (see
\S\ref{realPtEx}) assume
$z_0, 
\row z {r_1}$ appear in that order on the real line (circle on the
Riemann sphere). Then, 
\cite{DFrRRCF} produces paths based at $z_0$,
\begin{equation}\label{compPaths} 
P_1,P_2,\dots,P_{r_1},P_{r_1\np1,1}, P_{r_1\np1,2},\dots,P_{\sr,1},
P_{\sr,2}, 
\end{equation}  with explicit complex conjugation action $c$ on these paths and 
on
points over
$z_0$. 

Represent $c$ as $c_{\bz}=c_{z_0,\bz}$, by its permutation effect on 
points over $z_0$. With these conditions this description is uniform in
$z_0$. Under $c$ 
the paths $P_i$ go to a new set of paths, $P_i'$, $i=1,\dots,r$. From 
\cite{DFrRRCF}:
If $P$ and $P'$ 
are one of the $s_1$ pairs of paths, then $c$ sends $P$ (resp.~$P'$) to
a conjugate of 
$(P')^{-1}$ (resp.~$P^{-1}$). Branch cycles for the complex conjugate
cover 
$\phi^{c_\bz}:X^{c_\bz}\to \prP^1_z$ relative  to the $P_i'\,$s are also
$\bg$. A formula 
(dependent on the particular paths) computes branch cycles for
$\phi^{c_\bz}$ 
relative to the $P_i\,$s: Call this $\hk(\bg)=\hk_\bz(\bg)$ (as in
\eqref{compCong}). 
If  $\phi^{c_\bz}$ is s-equivalent to $\phi$
(Def.~\ref{covEq}), then $c_\bg$ conjugates 
the last branch cycle description to $\bg$.

\subsubsection{Formula defining $c_{\bz}=c\in N'$} Assume $G\le N'\le 
N_T(G,\bfC)$ 
as in \S\ref{intEquiv}. Prop.~\ref{compTest} produces for each $\bg\in 
\ni(G,\bfC)/N'$ a 
test for a connected component of $\sH_{N'}(\bR)$. Let $\phi_\bg: X_\bg\to 
\prP^1_z$ be 
a cover with branch cycle description $\bg$ (relative to the chosen classical 
generators for $\pi_1(U_\bz,z_0)$). The test succeeds if some component has a 
point 
corresponding to $\phi_\bg$.  The basic formula \eqref{compCong} is from 
\cite[Thm.~4.4]{DFrRRCF}.

\begin{prop} \label{compTest} Suppose $\bg\in \ni(G,\bfC)/N'$. Then, $\bg$ 
represents 
$\bp_{\bg}\in \sH_{N'}$ over $\bz$ (using the paths \eqref{compPaths}). Denote 
$g_{r_1+1,.}\cdots g_{\sr,.}$ by $g_{.}$  and $g_{j,1}g_{j,2}$ by  $g_{j,.}$. A 
connected 
component of $\sH_{N'}(\bR)$ goes through $\bp_{\bg}$ if for some 
$c=c_\bg\in N'$, $c\bg c=\hk(\bg)$ with $\hk_\bz=\hk$ as follows. For such a
$c_\bg$, $c_\bg^2$ centralizes $G$. 

\begin{equation}\label{compCong}\begin{array}{rl}  \hk(\bg)_1&= g^{-1}_1,\\ 
\hk(\bg)_2&=( g^\sph_3\cdots g^\sph_{r_1} g_{.})^{-1} g^{-1}_2( g^\sph_3\cdots 
g^\sph_{r_1} g_{.}),\ldots,\\ \hk(\bg)_{r_1\nm1}&=( g^\sph_{r_1} g^\sph_{.})^{-
1} g^{-1}_{r_1-1} g^\sph_{r_1} g^\sph_{.},\\ \hk(\bg)_{r_1}&= g^{-1}_{.} g^{-
1}_{r_1} g^\sph_{.}; \\ \hk(\bg)_{r_1+1,1}&= (g^\sph_{r_1+2,.}\cdots 
g^\sph_{\sr,.})^{-1} g^{-1}_{r_1+1,2} ( g^\sph_{r_1+2,.}\cdots 
g^\sph_{\sr,.}),\\ 
\hk(\bg)_{r_1+1,2}&=( g^\sph_{r_1+2,.}\cdots g^\sph_{\sr,.})^{- 1} g^{-
1}_{r_1+1,1} ( g^\sph_{r_1+2,.}\cdots g^\sph_{\sr,.}), \dots,\\ 
\hk(\bg)_{\sr,1}&= 
g^{-1}_{\sr,2},\quad \hk(\bg)_{\sr,2} = g^{-1}_{\sr,1}.\end{array}\end{equation} 
If a cover over $\bR$ represents $\bp_\bg\in \sH_{N'}$ (using the paths
\eqref{compPaths}), then some
involution $c_\bg$ satisfies $c_\bg\bg c_\bg=\hk(\bg)$. If 
$\sH_{N'}$ is a fine moduli space, such an involution $c_\bg$ exists if and only 
if a
connected  component of $\sH_{N'}(\bR)$ goes through $\bp_{\bg}$. \end{prop}

\begin{proof}[Comments on the proof] The operator $\hk$ acts as an involution:
$\hk^2(\bg)=\bg$ and $\hk$ commutes with conjugation by $c_\bg$. So
$c_\bg^2$ acts trivially on $\bg$. The condition for $\sH_{N'}$ to have fine 
moduli
is that the centralizer of $G$ in $N'$ is trivial
(Thm.~\ref{FrVMS}).  \end{proof} 

\begin{rem}[When $c_\bg$ is an involution] Prop.~\ref{HMnearHMLevel}
includes an archetype example for an inner Hurwitz family with group $T'$ having 
a nontrivial
center: It is a Frattini central extension of a centerless group $G$. If $\bg'$
 is branch cycles  for a $T'$ cover, we can usually decide if 
$c_{\bg'}$  is an involution. So, we can decide if there is a cover over
$\bR$ realizing
$\bp_{\bg'}$. 
\end{rem} 

\subsubsection{When $G=G_{k}$ is the level $k\ge 1$ group} Suppose $\bp\in 
\sH(G_k,\bfC)^\inn(\bR)$ lies over $\bz$ and $z_0$ is a  real base point 
relative to which 
\eqref{compCong} is a branch cycle description of $\phi_\bp:X_\bp\to \prP^1_z$. 
Assume: \begin{triv} \label{HMconds} Complex conjugation $c_\bp=c$ is trivial, 
and 
$r_1=0$ in \eqref{compCong}.  \end{triv} \noindent Then, the branch cycle 
description 
for $\phi_\bp:X_\bp\to \prP^1_z$ is an H-M representative. This is in the $H_r$ 
orbit on 
$\ni(G_k,\bfC)$ corresponding to the component of $\sH(G_k,\bfC)^\inn$ 
containing 
$\bp$. To complete the proof  of Thm.~\ref{thm-rbound} requires showing two 
things: 
\begin{edesc}\item  \eqref{HMconds} is equivalent to having an H-M 
representative.  \item 
If $k\ge1$ and  $c\bg c=\hk(\bg)$ with $c\ne 1$, then $\bg$ is a near H-M 
representative.  \end{edesc}

First assume $r_1>0$. Suppose $c=c_{z_0}$ is  the identity. Then 
\eqref{compCong} 
says $g_1=g_1^{-1}$, or $g_1$ has order 2. This can't hold if this conjugacy 
class is 
$2'$. It also can't hold if $k\ge 1$. The reason is this: The hypotheses remain 
the same for 
the reduction of $\bg$ modulo $\ker_{k-1}$. Apply Lem.~\ref{FrKMTIG}: Any lift 
to 
$G_k$ of an element of order divisible by 2 in $G_{k-1}$ increases its order. 
So, $g_1$ 
can't have order 2.

Let $c'$ be complex conjugation at level 0 (reduction modulo $\ker_0$). The same 
applies 
to it. So, $c'$ has order two (not 1) if $r_1>0$.  For inner classes, $c'\in 
G_0$. Apply 
Lem.~\ref{FrKMTIG} again so $c$ has order exceeding that of $c'$. This 
contradicts that 
$c$, a complex conjugation operator, has order 2 and is a lift of $c'$.  Now we 
know 
$r_1=0$. Further, this shows $c'$ is the identity.

This argument applies in going from level $k$ to level $k\np1$, $k\ge 0$. Two 
possibilities happen. Assume $\bg'\in \ni(G_k,\bfC)^\inn$ lies below $\bg\in 
\ni(G_{k\np1},\bfC)^\inn$, where $\bg$ defines a connected component of real 
points on 
$\sH(G_{k+1},\bfC)^\inn$. Let $c\in G_{k+1}$ (resp.~$c'$) be the complex 
conjugation 
operator associated with $\bg$ (resp.~$\bg'$). Then: Either \begin{edesc} 
\label{hmStyles} \item \label{hmStylesa}$c$ is the identity and $\bg$ is an H-M 
representative; or \item \label{hmStylesb} \begin{itemize} \item $c'$ is the 
identity, $\bg'$ 
is an H-M rep., \item $c\in \ker_k/\ker_{k\np1}$ is not the identity, 
and \item nothing in 
$\sH(G_{k+2},\bfC)^\inn(\bR)$ is over the $\bR$ component corresponding
to $\bg$ in $\sH(G_{k+1},\bfC)^\inn$.  \end{itemize}
\end{edesc}

Near H-M representatives  (satisfying \eql{hmStyles}{hmStylesb}) occur  at level 
$k\ge 1$ of the 
Modular Tower for $(A_5,\bfC_{3^4})$ (\S\ref{othLongOrbits}, especially Prop.
\ref{nearHMreps-count}).

\subsection{Real points when $r=4$} \label{realPtEx} \S\ref{covjvalues} gives an
elementary lemma about $j$ values of cover points over $\bR$ on a Hurwitz space. 
Then,
\S\ref{realCompFit} shows how real components over critical intervals fit 
together. 

\subsubsection{$j$ values of cover points over $\bR$} \label{covjvalues} We note 
the $j$ invariant separates  four branch point covers over $\bR$  according to a 
configuration of their branch points. 

\begin{lem} \label{locj} Suppose $\phi: X\to \prP^1_z$ is a four branch point 
cover over
$\bR$ with either 0 or 4 real branch points. Then,  the corresponding $j$ value 
under the
representative of
$\phi$ on a reduced Hurwitz space is in the interval
$(1,+\infty)$ along the real line. 

If $\phi$ has, instead, two complex conjugate and two real branch points, then 
the
corresponding $j$ value is in the interval $(-\infty,1)$. 
\end{lem}

\begin{proof}
Recall the cross ratio of distinct points $\row z 4$: 
$\lambda_\bz=\frac{(z_1-z_3)(z_2-z_4)}{(z_2-z_3)(z_1-z_4)}$.  The basics are in
\cite[p.~79]{Ahlfors}. Four points in complex conjugate pairs (or on the real
line) lie on a circle and the cross ratio is real. The cross-ratio is
invariant under a transform of the points $\bz$ by $\alpha\in \PGL_2(\bC)$. 
Since
there is an $\alpha\in \PGL_2(\bC)$ that takes two complex conjugate pairs of 
points
to four points in the reals, with no loss assume $\bz$ has either two or four 
real
points in its support. For these cases apply $\beta\in \PSL_2(\bR)$ to assume 
$0=z_1$
and $\infty=z_2$. Then, $\lambda_{\bz}=\frac{z_4}{z_3}$. 

In the former case $\lambda_{\bz}$ runs over the unit circle (excluding $1$) and 
in the
latter case over all real numbers (excluding $0$, 1 and $\infty$).  
The $j_\bz$ value corresponding to $\lambda_\bz$ is 
$J_4(\lambda)=\frac{4}{27}\frac{(1-
\lambda_\bz+\lambda_\bz^2)^3}{\lambda_\bz^2(1-\lambda_\bz)^2)}$
with 
\cite[p.~282]{Ahlfors}. (Classically this is without the $\frac{4}{27}$. We 
chose 
it so the ramified
$j$-values are 0, 1, $\infty$.) 

For $\lambda_\bz\in \bR\setminus \{0,1\}$
the connected range of $j_\bz$ includes large positive values and is bounded 
away from 0. So
the range of
$j_\bz$ for real $\lambda_\bz$ is $(1,\infty)$. For
$\lambda_\bz=e^{2\pi i \theta}=\lambda(\theta)$ in the unit circle (minus 1), 
the range of
$j_\bz$ includes both sides of 0. Also, for $\theta$ close to 1, the numerator 
of $j_\bz$ is
positive and bounded, while the denominator is approximately $(i\theta)^2$. 
Therefore the
range is the interval $(-\infty,1)$.
\end{proof}

\subsubsection{Pairings of real point components} \label{realCompFit}
Suppose $\psi: \bar \sH^\rd\to \prP^1_j$ is a reduced Hurwitz space cover 
defined over $\bR$.
Lem.~\ref{locj} shows the intervals $(-\infty,1)$ and $(1,+\infty)$ on 
$\prP^1_j$ lie
under real points on $\sH$ (the original Hurwitz space) coming from covers over 
$\bR$
with two different styles of branch points. The interval $(-\infty,1)$ goes 
through
$0\in \prP^1_j$, though this is a branch point for the cover $\psi$. This is
because ramification over 0 has order 3 (or 1), and a unique 3rd root of 1 (or
-1) is real. 

The same simple observation gives the next {\sl dessins d'enfant\/} type
lemma. Denote by $S_{(1,\infty)}$ (resp.~$S_{(-\infty,1)}$) the real points of 
$\sH^\rd$
over the interval $(1,\infty)$ (resp.~$(-\infty,1)$) of $\prP^1_j(\bR)$. The 
closure of each
component
$C$ of
$S_{(1,\infty)}$ has endpoints $\bp_1$ over 1 and $\bp_\infty$ over $\infty$. 
Let $e_1(C)=e_1$
(resp.~$e_\infty(C)=e_\infty$) be the ramification order of $\bp_1$ 
(resp.~$\bp_\infty$) over
1 (resp.~$\infty$). Note: $e_1=1$ or 2. The same attachment of endpoints applies 
to
components of $S_{(-\infty,1)}$. 

\begin{lem} \label{square1} Suppose for a given $C$ in $S_{(1,\infty)}$,
$e_1(C)$ (resp.~$e_\infty(C)$) is odd. Then there is a unique $C'$ in
$S_{(-\infty,1)}$ with $e_1(C')=e_1(C)$ (resp.~ $e_\infty(C')=e_\infty(C)$). 

If $e_1(C)$ (resp.~$e_\infty(C)$) is even, then there is a unique  $C'$ in
$S_{(1,\infty)}$ with $e_1(C')=e_1(C)$ (resp.~ $e_\infty(C')=e_\infty(C)$). So,
no  $C''$ in $S_{(-\infty,1)}$ has $e_1(C'')=e_1(C)$
(resp.~$e_\infty(C'')=e_\infty(C)$). 
\end{lem}

\begin{proof} Suppose $t$ is a local uniformizing parameter for $\bar \sH^\rd$ 
in a
neighborhood of $\bp_\infty$. The argument for $\bp_1$ is the same, so we do 
only the former
case. In local analytic coordinates over $\bR$, choose $t$ so 
$t^{e_\infty}=1/j$ and there is a parametrization of the neighborhood of 
$\bp_\infty$ using
power series in $t$ with real coefficients. If $e_\infty$ is odd, then real 
points
around $t=0$ map one-one to real points around $j=\infty$. If $e_\infty$ is 
even, then 
real points around $t=0$ map two-one to the positive number side of $j=\infty$ 
(no points
falling on the negative side of $j=\infty$). That interprets the
lemma's statement in local coordinates. 
\end{proof}

\subsubsection{$\bR$ points on reduced Hurwitz spaces} Assume $r=4$.
Prop.~\ref{isotopyM4} shows, if $\sQ''$ acts trivially on Nielsen classes,
$\bR$-cover points produce all the points on a reduced Hurwitz space
$\sH^\rd=\sH(G,\bfC)^\rd$ (with any equivalence on Nielsen classes) except 
possibly in
the fibers of
$\sH^\rd$ lying over
$j=0$ or
$1$. Our next result extends Prop.~\ref{compTest}, by combining it
with Lem.~\ref{redCocycle}, to consider any $\sQ''$ action. 

\begin{cor} \label{reducedRpts} Assume, as above, $j\in (1,\infty)$
(resp.~$(-\infty,0)$ or
$(0,1)$) and
$\hk$ is the operator of Prop.~\ref{compTest} for complex conjugate pairs of 
points
$\bz$ (resp.~a complex conjugate point and two real points). Let $\bz\in 
U_4(\bR)$ lie
over $j$. Let $\row {\bg} t$ be a listing of representatives from the 
equivalence
classes of Nielsen classes. Then, $\bR$ points of  $\sH^\rd$ over $j$ correspond
one-one with reduced equivalence classes $\mod \sQ''$ that upon
containing $\bg_i$ also contain $\hk(\bg_i)$. 

Suppose $\sH$ has fine moduli, and  
$\sH^\rd$ has b-fine (resp.~fine) moduli (Prop.~\ref{redHurFM}). Then 
$\bp\in \sH^\rd(\bR)$ over $j\in (-\infty,1)\cup (1,+\infty)$ (resp.~$(-
\infty,+\infty)$) is
either an
$\bR$-cover point, or an
$\bR$-Brauer point. When $\sH$ is an inner Hurwitz space, and $X_\bp\to
\prP^1_z$ corresponds to $\bp$, then $X_\bp(\bR)\ne \emptyset$
if and only if $\bp$ corresponds to an H-M rep. 
\end{cor}

This result applies to the main Modular Tower of this paper. 

\begin{prop} \label{HMnearHMLevel} Each level 
$\sH_k^\rd=\sH(G_k,\bfC_{3^4})^{\inn,\rd}$ of the
$(A_5,\bfC_{3^4},p=2)$ Modular Tower has  an absolutely irreducible component 
$\sH_{k,*}^\rd$ with
both H-M and near H-M reps. Each point of $\sH_{k,*}^\rd$ of either type
produces 
$\bR$-cover points in the Nielsen class.  Suppose $\bg$ is a branch cycle 
description of such a
cover with respect to $\bz$ having complex conjugate pairs of branch points. 
Then, the complement
$(\bg)\comp$ of $\bg$ (Def.~\ref{compDef}) corresponds to an $\bR$-Brauer point 
of
$\sH_{k,*}^\rd$. 

Consider the nontrivial central Frattini extension $T_k'\to G_k$ from
Cor.~\ref{R1G1}. The natural one-one (not necessarily onto) map
$\sH(T_k',\bfC_{3^4})^{\inn,\rd}\to
\sH(G_k,\bfC_{3^4})^{\inn,\rd}$ gives a component of
$\sH(T_k',\bfC_{3^4})^{\inn,\rd}$ isomorphic  (equivalent as covers of
$\prP^1_j$) to  $\sH_{k,*}^\rd$. An
$\bR$-cover  point $\bp\in \sH_{k,*}^\rd$ for an H-M (resp.~near H-M) rep.~in 
the Nielsen class
$\ni(G_k,\bfC_{3^4})$ corresponds to a cover in $\ni(T_k',\bfC_{3^4})$ with 
minimal field of definition  
$\bR$ (resp.~$\bC$). So, an 
$\bR$-cover point for a near H-M rep.~(for $G_k$) corresponds to a 
$T_k'$ cover with field of moduli $\bR$, but minimal definition field $\bC$.
\end{prop} 

\begin{proof} Prop.~\ref{nearHMreps-count} gives the component containing both 
H-M and near H-M
reps. The statement on $T_k'$ is from Prop.~\ref{nearHMbraid}. Nontriviality of 
$T_k'$ 
is from the discussion prior to Cor.~\ref{R1G1}, an inductive consequence of 
Prop.~\ref{RkGk} and that $T_0'=\bZ/2=\ker(\Spin_5\to A_5)$. 

The complex conjugation operator for a near H-M rep.~in this case gives 
conjugation by an
element whose lift to $R_k'$ has order 4. From this, Prop.~\ref{compTest} shows 
a near H-M
$\bp\in \sH_{k,*}(\bR)$ (regarding it in $\sH(T_k',\bfC_{3^4})^\inn$) has no 
cover realizing it 
over
$\bR$. If, however, $\bp$ corresponds to an H-M rep., then regarding it as the 
inner class of
an $T_k'$ cover, it has a trivial complex conjugation operator.  

So, we have $X\to \prP^1_z$
with group $G_k$ over $\bR$ and geometrically $Y\to X\to \prP^1_z$ ($Y\to X$ 
unramified) with field
of moduli
$\bR$, and the group of $Y\to \prP^1_z$ is $T_k'$.   Apply Lem.~\ref{charblem} 
to
$Y\to X$; the resulting cover $Y'\to X$ has degree 2. It is therefore Galois, 
and $Y'\to
\prP^1_z$ produces the $T_k'$ realization over $\bR$.
\end{proof}

\begin{rem}[Change conjugate pairs $\bz$ to all real points]
\label{comp-real} Assume
$\bz$ has two conjugate pairs of points. Choose $\beta\in \PSL_2(\bC)$ to map 
$\bz$ to
$\bz'$, four points in $\bR$. If $\phi:X\to \prP^1_z$ is a cover over $\bR$ with
branch points $\bz$, then $\beta\circ\phi:X\to \prP^1_z$ is in the reduced 
equivalence class of
$\phi$. It is not, however, a cover over $\bR$ as it stands, because $\bR$ is 
not
a field of definition of $\beta\circ\phi$. Sometimes, however, $\beta\circ\phi$ 
is strong equivalent
(so it has the same branch points $\bz'$) to a cover over $\bR$. An example 
arising 
in this paper is the absolute and inner reduced equivalence classes for 
$\ni(A_5,\bfC_{3^4})$.
Consider the cover with H-M description 
$\bg=((1\,2\,3),(1\,3\,2),(1\,4\,5),(1\,5\,4))$ and $\bz$
as branch points. With a base point $z_0\in \prP^1_z(\bR)$, the effect of 
complex conjugation
from Prop.~\ref{compTest} is the identity. The
points on $X$ over $z_0$ are all real.  
 
Transform the paths giving $\bg$ above by $\beta$. This gives the same branch 
cycle description
of $\beta\circ\phi$ with respect to these new paths. Since, however, $\bz'$ 
consists of four real
points, use the corresponding complex conjugation operator. As $\beta(z_0)$ may 
not be real,
this may require adjustment in getting a base point $z_0'\in \prP^1_z\setminus
\beta{\bz}$. Still, with this choice, $c_1=(2\,3)(4\,5)$ relative to $\beta$ 
transformed
paths gives complex conjugation on the points above  $z_0'$. Thus, the cover
(either degree 5 or the Galois closure of the degree 5 cover) is equivalent to 
one over $\bR$.
This might happen whenever $\sQ''$ is not faithful on a Nielsen class,  
especially one with an H-M rep., $\bg=(g_1,g_1^{-1},g_2,g_2^{-1})$.  Then, 
$(\bg)q_1q_3^{-1}=(g_1^{-1},g_1,g_2^{-1},g_2)$ and $h\bg h^{-1}=(bg)q_1q_3^{-1}$ 
is equivalent
to the phenomenon above. For, however, H-M reps. in $\ni(G_1,\bfC_{3^4})$, there 
is no such
$h$. So the corresponding cover with real branch points is not equivalent to a 
cover over
$\bR$.   
\end{rem}

\newcommand{\fix}{\textr{fix}}
\newcommand{\br}{\textr{br}}

\subsubsection{How it works for $r\ge 5$} \label{rge5} Suppose $r \ge 5$ and 
$\bz\in U_r$ is
suitably general (lying off an explicit Zariski closed subset of $U_r$, 
dependent on
$r$). Then, there is no nonidentity
$\alpha\in
\PSL_2(\bC)$ with
$\alpha(\bz)=\bz$. The argument of \S\ref{PSL2comp} 
for (fixed $r$), shows the set of $\bz$, with $\alpha$ acting on it as an
element of order exceeding 3, is a Zariski closed subset. If $\alpha$ has order
$3$, assume with no loss $\{z_1,z_2,z_3\}$ is $\{0,1,\infty\}$. This  reduces  
the problem to
noting there are only two other specific fixed points and no
other cycles of order 3 (\S\ref{PSL2comp}). We show only  a proper Zariski
closed subset of 
$\bz$ is fixed by some involution. With no loss, consider four values $0, 1, 
\infty, z_4$
transposed in pairs. This determines $\alpha$, so other elements of $\bz$,
either fixed or transposed in pairs, are special. Denote this exceptional subset 
by $J_{r,\fix}$. 

\begin{prop} Suppose $r\ge 5$, $\bz\in
U_r(K)$ and $\sH$ is a Hurwitz space over $K$.  
Assume
$\bz \mapsto j'\in J_r\setminus J_{r,\fix}$. Then, 
points of the fiber $\sH_\bz$ go one-one to the points of the fiber
$\sH_{j'}^\rd$, so they have exactly the same fields of definition over $K$.

There is a proper algebraic subset $J_{r,\br}$ of $J_r$ so that if $j'\in
J_r\setminus J_{r,\br}$ is a $K$ point, then the fiber of $U_r$ over $j'$ has 
$K$ points. For
$j'\in J_r(K)\setminus (J_{r,\br}\cup J_{r,\fix})$, a $K$ point on $\sH^\rd$ 
over $j'$ is
the image of a 
$K$ point on $\sH$ for some $\bz$ mapping to $j'$.  
\end{prop} 

So, for $r\ge 5$, some problems for $r=4$ don't appear in comparing $\sH^\rd(K)$ 
and $\sH(K)$.
Particularly, from Prop.~\ref{redHurFM}, if $\sH$ is a fine moduli space, each 
component of
$\sH^\rd$ is a b-fine moduli space; and avoiding 
$\bz$ over 
$J_{r,\br}\cup J_{r,\fix}$, we can lift $K$ points from $\sH^\rd$ to $\sH$. 
Still, the
points of $J_{r,\br}\cup J_{r,\fix}$, including the {\sl orbifold\/} points of 
$J_{r,\fix}$
require a refined analysis to extend the full version of Prop.~\ref{redHurFM} 
when $r=4$. 

\subsection{Level 0 $\bR$ points on the  $(A_5,\bfC_{3^4})$ Modular Tower} 
\S\ref{j-covers} produces a ramified cover 
$\psi_{0,\abs}: \sH^{\abs,\rd}_0 \to \prP^1_j\setminus\{\infty\}$ from $\bar 
M_4=H_4/\sQ$
acting on reduced absolute Nielsen classes. Closure of this cover over the
$j$-line  gives $\bar\phi^{\abs,\rd}: \bar \sH^{\abs,\rd}_0 \to \prP^1_j$ with 
branch cycles 
$(\gamma_0',\gamma_1',\gamma_\infty')$ \cite[\S7.4]{Fr-Schconf}:  
\begin{equation} \label{branchCycabs} \gamma_0'=(2\,1\,4)(3\,7\,8)(5\,6\,9),
\gamma_1'=(4\,5)(3\,9)(1\,2)(8\,6),
\gamma_\infty'=(1\,4\,9\,8\,5)(3\,6\,7).\end{equation} 

This is  a cover with branch points $\{0,1,\infty\}$. We explicitly
uniformize the real points on 
this cover, and on the inner version of this space, including the points
over branch 
points. From Lem.~\ref{basicA5cover}, the geometric (resp.~arithmetic) monodromy 
group of the cover
is
$A_9$ (resp.~$S_9$). So, \eqref{caction} uniquely determines the effect of
complex conjugation. 

\subsubsection{Absolute reduced Hurwitz space real points} For a
reduced Hurwitz space (over $\bR$) covering $\prP^1_j$, there are three complex
conjugation operators $c_{-\infty,0}$, $c_{0,1}$ and $c_{1,+\infty}$ 
corresponding to
the three intervals over which the cover is unramified. As in 
\S\ref{realCompFit},
$c_{-\infty,0}$, $c_{0,1}$ will produce the same number of real points lying 
over a
corresponding $j$ value. Apply
\eqref{compCong}. With $u_1,u_2\in \{0,1,\infty\}$, write the operator as 
$c_{u_1,u_2}$
(suppressing the $\pm$).  Its characterization is 
$c\gamma_{u_i}'c=(\gamma_{u_i}')^{-1}$,
$i=1,2$. As $\sQ''$ acts trivially on $\ni(A_5,\bfC_{3^4})^\abs$, 
Prop.~\ref{isotopyM4} shows
all elements of $\sH^{\abs,\rd}_0(\bR)$ are $\bR$-cover points. We see this 
directly in the next
lemma. 

\begin{lem} \label{A5absR} All three points of $\sH^{\abs,\rd}_0(\bR)$ lying 
over $j\in
(1,\infty)$ are $\bR$-cover points. The same is true for each point over the 
intervals
$(-\infty,0)$ and $(0,1)$.  
\end{lem}

\begin{proof} From the above, 
$c_{\infty,0}$ conjugates both 
$\gamma_0'$ and  $\gamma_\infty'$ to their inverses.  The unique element doing
this is 
$c_{\infty,0}=(1\,4)(3\,7)(5\,9)$.  So there are three real
points  above $j_0\in (-\infty,0)$. 
Check: $c_{1,\infty}=(3\,6)(8\,9)(4\,5)$ and  
$c_{0,1}=(3\,8)(9\,6)(1\,2)$. Now we show cover points account for all
$\bR$ points on
$\sH^{\abs,\rd}_0$.

The fixed points of $c_{1,\infty}=(3\,6)(8\,9)(4\,5)$ point to three covers in 
Table
\ref{lista5} labeled
${}_1\bg,{}_2\bg,{}_7\bg$. Each  attaches (by the $\gamma_\infty$ operator) to a
distinct cusp over
$j=\infty$, of respective widths 5, 1 and 3. That, however, comes from choosing
classical generators of $\pi_1(U_\bz,z_0)$ for some  $\bz$ and
$z_0$. It doesn't show what the description of branch cycles will be for a set 
of paths on
$U_\bz$ supporting a complex conjugation operator of Prop.~\ref{compTest}. 
Actual
covers
$\phi:X\to 
\prP^1_z$ having these real structures come from two types of $\bz\,$s using the 
paths 
\eqref{compPaths}. 
 
For $\bz$ with real entries: 
$\bg_1=((1\,2\,3),(1\,3\,2),(1\,4\,5),(1\,5\,4))$ has complex operator
$c_1=(2\,3)(4\,5)$;  $\bg_2=((1\,2\,3),(1\,4\,5),(1\,5\,4),(1\,3\,2))$ has 
$c_2=(1\,3)(5\,4)$ and $\bg_7=((1\,2\,3),(2\,1\,4),(2\,4\,5),(5\,3\,2))$ has 
$c_7=(2\,3)$. 
\end{proof}

\subsubsection{General expectations for genus 0 covers} \label{genus0} Let $K$
be a perfect field with $Y$ isomorphic to $\prP^1_z$ over $\bar K$. For $X$ 
a projective nonsingular curve, assume $\phi:X\to Y$ is over $K$. If $X$ has 
genus 0,
then $X$ defines an element $[X]\in H^2(G_K,\bar K^*)$ of order 2 
\cite[p.~126]{SeGalCoh}. A
simple formula relates these: $n[X]=[Y]$ with $n=\deg(\phi)$ 
\cite[Thm.~1]{FrBGJac}.
Further, for any odd integer
$n$, $K$ and $[Y] \in H^2(G_K,\bar K^*)$, there is such a 
$\phi:X\to Y$ over $K$ with $X$ of genus 0 
\cite[Thm.~2]{FrBGJac}. 

\begin{prob} Given $[Y]$ and $n$ odd, which
Nielsen classes of genus 0 curve covers have $\phi:X\to Y$ over $K$ of
degree
$n$ in the Nielsen class? \end{prob} 

\begin{rem}[Ch\v atalet's Work] \label{chatalet} See \cite[p.~389]{CT} or 
\cite[Ch\^ atalet's
Thesis]{Chat1} for background on these comments. Ch\^ atalet knew that $X$ a 
Brauer-Severi
variety of dimension $n$ defines a central simple algebra over $K$ of dimension 
$n^2$.  We
recognize this as $[X] \in H^2(G_K,\bar K^*)$. If $X(L)\ne \emptyset$ with
$([L:K], n)=1$, then it possesses a $K$ point (from 
the corestriction-restriction sequence of Galois cohomology \cite[p.~105]{AW}).  
Also, if
$X_{K_\bp}= Y_{K_\bp}$ for each completion
$\bp$ of $K$, then Ch\v atelet knew $X \equiv Y$. This gives the easy direction 
of
\cite[Thm.~1]{FrBGJac} when their images in the
completions of a field determine the Brauer-Severi variety. According to
\cite{CT},  Weil \cite{WeilField} gave the modern Galois cohomology 
interpretation. \end{rem}

\subsubsection{Inner reduced Hurwitz space real points} \label{innerHSReal} Now 
consider 
branch cycles for $\bar\phi^{\inn,\rd}: \bar \sH^{\inn,\rd} \to \prP^1_j$ from 
inner classes 
$\ni(A_5\bfC_{3^4})^\inn$ (\cite[\S7.10]{Fr-Schconf}, \cite{Fr3-4BP}):
$$\begin{array}{rl} \label{branchCycinn}
\gamma_0''=&(1\,13\,2) (3\,7\,17)  (4\,11\,10) (5\,15\,9) (6\,18\,14) 
(8\,12\,16),\text{and 
} \\ \gamma_1''=& (1\,11)(2\,10)(3\,9) (4\,14)
(5\,13) (6\,17) (7\,16) (8\,15) (12\,18)\\ \gamma_\infty''=& 
(2\,11)(1\,4\,18\,8\,5)(10\,13\,9\,17\,14)(3\,15\,16)(12\,6\,7).\end{array}$$ 
For 
$j_0\in (1,\infty)$, inspection gives $c_{1,\infty}$ conjugating $\gamma_1''$ 
and
$\gamma_\infty''$ to their inverses: 
\begin{equation} \label{c1infty} 
c_{1,\infty}=(9\,17)(7\,16)(8\,18)(13\,14)(12\,15)(6\,3)(4\,5).
\end{equation} 
The centralizer of the monodromy group ($c^\dagger$ in Ex.~\ref{4realpoint}) 
moves 1. To show
$c_{1,\infty}$ is the correct complex conjugation, it suffices to check the 
complex conjugation
operator for four real points fixes 1 and 10 (the H-M reps.).   

Similarly, for
$j_0\in (-\infty,0)$, by inspection find 
$$c_{-\infty,0}=(10\,13)(1\,4)(9\,14)(7\,3)(12\,16)(11\,2)(15\,6)(5\,18)$$  
conjugates
$\gamma_0''$ and
$\gamma_\infty''$ to their inverses. Note:  $c_{1,\infty}$ fixes 2 and 11. If 
the complex
conjugation operator for $(-\infty,0)$ fixes either, the set of real points 
would not be a
1-dimensional manifold in the neighborhood of the index 2 branch point over 
$\infty$. 
As $c_{-\infty,0}c^\dagger$ fixes 2 and 11, $c_{-\infty,0}$ is the correct 
complex conjugation
operator.  

Applying
\cite[\S4.1]{DFrRRCF} gives the following information on points  over $\infty$ 
\wsp cusps
\wsp\ in each cover. Again, since $\sQ''$ acts trivially on 
$\ni(A_5,\bfC_{3^4})^\inn$,
Prop.~\ref{isotopyM4} shows all elements of $\sH^{\abs,\rd}_0(\bR)$ are $\bR$-
cover points. Our
computations above confirm this directly. (Ex.~\ref{course-finemoduli} has more 
on 
using this example with $\Spin_5$ covers.)

\begin{lem} \label{A5innR} For each $j_0\in (1,\infty)\subset \prP^1_j(\bR)$ 
(resp.~$j_0\in
(-\infty,1)$) there are four (resp.~two) points of $\sH_0^{\inn,\rd}(\bR)$. 
There are
three
$\bQ$ points over
$\infty$ in the cover
$\bar\phi^{\abs,\rd}$  of ramification indices 1, 3 and 5. Denote these 
$x_1,x_3,x_5$. Real
points over $\infty$  in the cover $\bar\phi^{\inn,\rd}$ correspond to disjoint 
cycles in
$\gamma_\infty''$ that  conjugation by $c_{\inn,\rd}$ maps to their inverses. 
These are the
disjoint cycles 
$(2\,11)$, $(1\,4\,18\,8\,5)$ and $(10\,13\,9\,17\,14)$.  There are three real 
points over 
$\infty$, one over $x_1$ and two over $x_5$.  \end{lem}
  
\section{Cusps on a Modular Tower} \label{2A5Cusps} This section and
\S\ref{compA5MT} describe level $k=1$ components of the $(A_5,\bfC_{3^4})$
Modular Tower,  $\sH(G_1,\bfC_{3^4})^{\inn}=\sH_1$ and its reduced
version $\sH(G_1,\bfC_{3^4})^{\inn,\rd}=\sH_1^\rd$. As usual,
$\bfC_{3^4}$  is four repetitions of the class of 3-cycles in
$A_5$. For this noncongruence tower our detail is comparable to
literature describing some  specific modular curves.  

The idea is to apply Thm.~\ref{thm-rbound}
to inspect those cusps  at the end of real components on $\sH_1$. This 
information about the cusps eventually separates the two components of 
$\sH_1^\rd$
according to the cusps lying in each. When $r=4$, $\sQ''$ action produces the 
phenomenon 
$q_2$ {\sl orbit shortening\/} affecting precise genus calculations. 
\S\ref{sQorbits} and
\S\ref{sQorbits2}  develop lemmas anticipating the general case. 

\subsection{Fundamental domains and cusp data} Let $F$
be a standard  fundamental domain for $\PSL_2(\bZ)$ on the upper
half plane. Branch cycle  data from Prop.~\ref{j-Line} allows
computing a fundamental domain for the subgroup of $\PSL_2(\bZ)$
defining $\psi^\rd: \bar\sH^\rd\to \prP^1_j$ (with compatible notation for inner 
and absolute
reduced spaces). The element
$\gamma_0$ (resp.~$\gamma_1$) identifies with the standard element of order 3
(resp.~order 2) in $\PSL_2(\bZ)$. The action of $\gamma_0$ and $\gamma_1$ on
$\ni(G,\bfC)^{\abs,\rd}$ (or $\ni(G,\bfC)^{\inn,\rd}$ depending on the type of 
Hurwitz space) gives
$\gamma_0'$ and $\gamma_1'$, generating a  permutation group 
$\lrang{\gamma_0',\gamma_1'}$ whose
orbits correspond to
$j$-line covers. 

For further  reference, let $O$ be one of these
orbits. Take $\Gamma=\Gamma(O)$ to be the subgroup of $\PSL_2(\bZ)$
stabilizing an element in $O$. A fundamental domain for a cover
for the orbit $O$ comes from translating $F$ around by
coset representatives for $\Gamma(O)$ in $\PSL_2(\bZ)$. 

\subsubsection{Monodromy for $\bar\sH(A_5,\bfC_{3^4})^{\abs,\rd}$
and $\bar \sH(A_5,\bfC_{3^4})^{\inn,\rd}$} \label{A5Cycles} We give 
the geometric and arithmetic monodromy of both $j$-line covers.
Denote $\sH(A_5,\bfC_{3^4})^{\abs,\rd}$ by $\sH^{\abs,\rd}_0$, and
$\sH(G_k,\bfC_{3^4})^{\inn,\rd}$ by $\sH^{\inn,\rd}_k$.  

Consider any sequence of  (separable) absolutely irreducible covers 
$X\to Y\to Z$ over a 
field $K$. Let $G_Y$ (resp.~$G_X$) be the group of the galois closure of 
$Y\to Z$ (resp.~$X\to Y$) over $K$. Then, the group $G=G_{X/Z}$ of the
Galois closure (over 
$K$) of the cover $X\to Z$ is a subgroup of the wreath product $G_X\wr
G_Y$. Further, it has  the following properties \cite[p.~51]{FrSchur}.
\begin{edesc}
\item The projection of $G$ to 
$G_Y$ is surjective. \item The kernel of $G\to G_Y$ contains a group 
isomorphic to 
$G_X$.\end{edesc} 

\begin{lem} \label{basicA5cover} The geometric monodromy 
$G_{\psi^{\abs,\rd}}=G_{0,\abs}$ of $\bar \sH^{\abs,\rd}_0\to \prP^1_j$ is 
$A_9$. Let 
$\alpha=\sqrt{-15}$. The arithmetic monodromy group $\hat G_{0,\abs}$ is 
$S_9$. The two groups fit in a natural short exact sequence $$
1\to G_{0,\abs} \to \hat  G_{0,\abs}\to G(\bQ(\alpha)/\bQ)\to 1.$$ 

The geometric monodromy group $G_{\psi^{\inn,\rd}}=G_{0,\inn}$ of
$\psi^{\inn,\rd}: \bar\sH_0^{\inn,\rd}\to \prP^1_j$ is
$\bZ/2\wr A_9\eqdef (\bZ/2)^9\xs A_9$ with  $A_9$ acting as permutation of the
coordinates of
$(\bZ/2)^9$. The  arithmetic monodromy group
$\hat G_{0,\inn}$ maps to
$\hat G_{0,\abs}$  producing the extension  $$1\to G_{0,\inn} \to \hat 
G_{0,\inn} 
\to G(\bQ(\alpha)/\bQ))\to 1.$$  \end{lem}

\begin{proof} Use $\gamma'_\infty$ in \eqref{branchCycabs} to see that 
$G_{0,\abs}$ contains
a 5-cycle, $(\gamma'_\infty)^3$, so it must primitive. It also contains a 3-
cycle
$(\gamma'_\infty)^5$. A well-known argument says a primitive subgroup of $A_n$ 
containing a
3-cycle is $A_n$: $G_{0,\abs}=A_9$.  Apply \cite[Prop.~2.1 and Irrational Cycle 
Lemma]{FrExtConst}. The gist: As the branch points of the
cover are  rational, the arithmetic monodromy group contains the character field
of the conjugacy class  of any element of form (1)(3)(5). Elements that are  
products of distinct odd disjoint cycle lengths, form two conjugacy classes in 
$A_9$. The
outer automorphism of $A_9$  (conjugation by $(1\,2)$)
permutes these two conjugacy classes. 

For $\psi^{\inn,\rd}$, the degree 18 cover breaks
into a chain  of degree 9 and degree 2 covers. The largest possible
group for the geometric closure is 
$\bZ/2\wr A_9$ and for the arithmetic closure it is $\bZ/2\wr S_9$.
To complete  the proof only requires showing the kernel of the geometric
closure to $A_9$ contains one  of the factors of $(\bZ/2)^9$. Since
$(\gamma_\infty'')^{15}=(2\,11)$ generates such a factor, we are done.
\end{proof}

\begin{exmp} To illustrate,
consider \eqref{branchCycabs} for $\ni(A_5,\bfC_{3^4})^\abs$ (or 
\eqref{branchCycinn} for
$\ni(A_5,\bfC_{3^4})^\inn$). The representation has degree 9 (or 18). An 
explicit
orbit  comes as follows. Label 9 words $w(\gamma_0',\gamma_1')_i=w_i$,
$i=1,\dots,9$, in the elements $$\gamma_0'=(2\,1\,4)(3\,7\,8)(5\,6\,9) \text{ 
and }
\gamma_1'=(2\,1)(4\,5)(6\,8)(3\,9),$$ so these words applied to
1 give the complete set $\{1,2,\dots,9\}$. That is,
they are coset representatives. Example:
$w_1$ the trivial word, $\gamma_0'=w_2$, $(\gamma_0')^2=w_3$,
$w_4=\gamma_0'\gamma_1'$, $w_5=w_4\gamma_0'$, $w_6=w_5\gamma_0'$,
$w_7=w_6\gamma_1'$, $w_8=w_7\gamma_0'$, $w_9=w_8\gamma_0'$.
Then, $$\begin{array}{rlrlrl} (1)w_1=1,& (1)w_2=4,& (1)w_3=2,&
(1)w_4=5,& (1)w_5=6,&\\ (1)w_6=9,& (1)w_7=3,& (1)w_8=7,&
(1)w_9=8.&\quad&\end{array}$$ Finally, apply these, with $\gamma_0$
and $\gamma_1$ replacing $\gamma_0'$ and $\gamma_1'$ to $F$. Take
the union of the images as a fundamental domain for the $j$-line
cover.\end{exmp}

\subsubsection{(2,3) generation}
Every group 
with two generators, $\gamma_0',\gamma_1'$, with respective orders 3  and 2, is 
a
quotient of $\PSL_2(\bZ)$. \cite{Sh} calls these (2,3) generated  groups. 
Finding which such finite simple groups do so occur has been a problem with much 
literature for over 
100 years. The groups $A_n$ have this property for $n\ge 9$, a result of Miller 
from 1901, and all simple (2,3) generated groups are known \cite{Sh}. The group 
$A_8$ does
not so occur, though
$\PSL_2(\bZ/7)$ (a simple group) does give a primitive subgroup of $A_8$ that is 
$(2,3)$
generated. Though to us it is icing on the cake that the main example of this 
paper
starts at the border of this old area, we consider it more evidence  that
there is a cake. 

\subsubsection{Fundamental domains from cusps} \label{fundDomainCusp} Let 
$\bR_\infty=\bR\cup
\{\infty\}$. Recall: Cusps are equivalence classes for $\bQ\cup\{\infty\}$ on
$\bR_\infty$ for the $\Gamma(O)$ action. These equivalence
classes correspond to orbits of $\gamma_\infty'$. For example, 
in \eqref{branchCycabs}, $\gamma_\infty'=(1\,4\,9\,8\,5)(3\,6\,7)$. So,
there are three cusps, with {\sl widths\/} ({\sl
lengths\/}) 1, 3 and 5. By contrast, $\Gamma_0(p)$, with $p$
a prime, has two cusps  of widths 1 and $p$. With
no loss, $\infty$ corresponds to the cycle of  length 1 in
$\gamma_\infty'$. (Let the image of $F$ corresponding to the integer
2 to have $\infty$ as  its cusp.) Then, there exist elements
$g_3',g_5'\in \PSL_2(\bZ)$ as follows:
\begin{equation}\begin{array}{rl} g_3'\Gamma(g_3')^{-1}\cap
\blrang{\smatrix1101} &= \blrang{\smatrix1301}, \text{
and}\\ g_5'\Gamma(g_5')^{-1}\cap \blrang{\smatrix1101} &=
\blrang{\smatrix1501}.  \end{array}\end{equation}

To get explicit $g_3'$ and $g_5'$ from this calculation
choose an  expression in the $w\,$s taking an integer in
the respective 3 and 5 cycles to the integer  2. The following
summarizes this explicitness result.

\begin{prop} Start with any Nielsen class having $r=4$. Prop.~{\rm
\ref{j-Line}} produces a fundamental domain (in the upper
half plane) and basic data on the cusp points of the
cover $\bar\sH^\rd\to \prP^1_j$ coming from that Nielsen
class. \end{prop}

The main result of \S\ref{compA5MT}
shows $\sH_1^\rd$ has two absolutely irreducible ($\bQ$-curve)
components, of respective genuses 12 and 9. The former component
corresponds to a braid orbit ($O_1^+$ below) on $\ni_1$ containing all H-M
and near H-M representatives (\S\ref{startHM}). So, Thm.~\ref{thm-rbound} 
implies
only the genus 12 component contains $\bQ$ points.   

Branch cycles for these components come from applying the
program \GAP\  to compute $\gamma_0,\gamma_1$ on
$\ni({}_2^1\tilde A_5,\bfC_{3^4})^{\inn,\rd}$. As in Prop.~\ref{j-Line},
this gives $\gamma_0',\gamma_1',\gamma_\infty'$ for a cover 
from each $\bar M_4$ orbit. Especially important are the 
lengths of $\gamma_\infty'$ orbits. 

\S\ref{HkH1}  isolates an expected general 
structure of fundamental domains for higher levels of a Modular Tower. It is 
that the
monodromy groups of the higher reduced levels over level 0 should be {\sl 
near\/}
$p$-groups (2-groups in our example). 

\subsection{Cusp notation} 
Denote the image of $q_2\in H_4$ in $\bar M_4$ by 
$\gamma_\infty$. Refined 
calculations differentiate between $q_2$ orbits on inner Nielsen classes
$\ni_1^\inn$ and $\gamma_\infty$ orbits on
$\ni_1^\rd=\ni_1^{\inn,\rd}$. The latter attach to cusps on 
$j$-line
covers. 

\subsubsection{$\gamma_\infty$ and $q_2$ orbit notation}
\label{cuspNotat}  We often denote $\gamma_\infty$ orbits on
$\ni_1^{\rd}$ by  notation
like $O(u,v;a)$. If
$\bg=(g_1,g_2,g_3,g_4)\in O(u,v;a)$, then 
$u=\mpr(\bg)\eqdef\ord(g_2g_3)$ is the middle product
(\S\ref{HMrepRubric}) of 
$\bg$.
Also,
$v=\wid(\bg)\eqdef|O(u,v;a)|$ is the {\sl width\/} of $\bg$. 

The decoration  
$a$ distinguishes orbits with specific
$(u,v)$. If it is cumbersome, or we don't know
$(u,v)$, a briefer notation suffices. Given 
$\bg\in O(u,v;a)$, refer to its {\sl orbit type\/} as
$(u,v)=(\mpr(\bg),\wid(\bg))$. When the context is clear, use $(u,v)$ as 
the
type of a $q_2$ orbit on
$\ni_1^\inn$. The $q_2$
orbit type of $\bg\in \ni_1^\inn$ may be different
from its $\gamma_\infty$ orbit type if $\sQ''$ {\sl orbit shortening\/} reduces 
the value of
$v$ (\S\ref{sQorbits2}). 

Consider $\bg\in \ni_1^\inn$ lifting $\bg'\in \ni_0^\inn$ with $\mpr(\bg')=3$ or 
5.
According to Lem.~\ref{A5Prod}, $\bg$ has a {\sl complementary element\/}
$(\bg)q_2^{\mpr(\bg)/2}$. This is the unique element in the $\gamma_\infty$ 
orbit of 
$\bg$, distinct
from $\bg$, that maps to $\bg'$. Denote this by $(\bg)\comp$. 

\begin{defn}[More complements] \label{compDef} Prop.~\ref{nearHMreps-count} 
generalizes the case $k=1$ to define $(\bg)\comp$ as  $(\bg)q_2^{\mpr(\bg)/2}$ 
for $\bg\in
\ni_k(G_k,\bfC_{3^4})$, $k\ge 1$. 
\S\ref{indSteps} generalizes further. 
\end{defn} 

\newcommand{\hm}{\text{\rm hm}}

\subsubsection{$\sQ''$ orbits} \label{sQorbits} 
Thm.~\ref{presH4} says $\sQ\norm H_4$ acts on an inner Nielsen class through a 
Klein
4-group, $K_4=\sQ''$. 

As in \cite{Fr3-4BP} or \cite[\S7.10]{Fr-Schconf}, the Nielsen
class $\ni_0=\ni(A_5,\bfC_{3^4})^\inn$ contains exactly 18
elements. Just two are H-M representatives. Branch cycles for
these are $$\begin{array}{rl} \bg_{\hm,1}&=((1\,2\,3), (3\,2\,1),
(1\,4\,5), (5\,4\,1)) \text{ and } \\ \bg_{\hm,2}&=((1\,2\,3),
(3\,2\,1),  (5\,4\,1), (1\,4\,5)).\end{array}$$ These
two are equivalent in the absolute Nielsen classes by  $(4\,5)$ conjugation, an 
outer
automorphism of $A_5$. Further, $H_4$
is transitive on $\ni_0$. So, the number of lifts of any $\bg\in
\ni_0$ to $\ni_1$ is independent of the choice of $\bg$.  

\newcommand{\Stab}{\texto{Stab}}
We first
show faithful action of $\sQ''$ on $\ni_1$. As usual, $\{G_k\}_{k=0}^\infty$ are 
the
characteristic quotients of ${}_p\tilde G$, with $\bfC$ a collection of $p'$ 
conjugacy
classes of
$G_0=G$. Write $\sQ^*=\lrang{Q_1Q_3^{-1},(Q_1Q_2Q_3)^2}$  when we regard it as a 
subgroup of $B_4$. So, $\sQ^*$ acts on $\ni(G,\bfC)$ (no equivalence by 
conjugation by $G$). The
notation $\sQ''$ means this Klein 4-group is acting on inner classes. Use
$\Stab_{\sQ''}(\bg)$ for the stabilizer in $\sQ''$ of $\bg\in \ni^\inn$ (for 
absolute
equivalence, mod out further by $G\le N'\le S_n$). For any group $S$, express 
the orbit of $\bg$ under $S$ as $(\bg)S$. 
\begin{lem} \label{4HM4NHM} Let $\ni(G,\bfC)$ be a Nielsen class for $r=4$ with 
$G$
centerless. Then, $\bg\in \ni(G,\bfC)$ is invariant
under
$(Q_1Q_3^{-1})^2\in \sQ^*$ if and only if $\bg$ is an H-M rep. If
$\bg=(g_1,g_1^{-1},g_2,g_2^{-1})$ is an H-M rep., then $\sQ^*$ acts on 
$(\bg)\sQ^*$ as a quotient of the Klein 4-group $K_4$. The $K_4$ action is 
faithful
unless $g_1$ and $g_2$ are involutions. 

Let $\ni_k$ be the $k$th level inner Nielsen class for any Modular Tower with 
$r=4$. Length
of the
$\sQ''$ orbit on
$\bg\in
\ni_k$ depends only on the $H_4$ (or $\bar M_4$) orbit of $\bg$. Faithful action 
of $\sQ''$
on an
$H_4$ level
$k$ orbit  extends to faithful action on an $H_4$ level $k\np1$  orbit above it.
Similarly, if
$\Stab_{\sQ''}((\bg)\lrang{q_2})=\{1\}$, then 
$\Stab_{\sQ''}((\bg')\lrang{q_2})=\{1\}$ for
$\bg'\in \ni_{k+1}$ above $\bg$. 

Assume $p=2$, and $\ni_k$, $k\ge 1$, is the 
Nielsen class for centerless 2-perfect $G_0=G$. Let $O$ be an $H_4$ orbit
on
$\ni_k$ containing an H-M rep. Then all  
$\sQ''$ orbits on $O$ have length four. So, $\sH_O^\rd$, the reduced
component corresponding to $O$ has b-fine moduli (Prop.~\ref{redHurFM}). 

For $G_0=A_5$, $p=2$  and $\bfC=\bfC_{3^4}$, $\sQ''$ is faithful on all $H_4$
orbits in $\ni_k$, $k\ge 1$.  
\end{lem}

\begin{proof} Since $(\bg)(Q_1Q_3^{-1})^2$ is conjugation of $\bg$ by $g_1g_2$, 
if 
$G$
is centerless, and $Q_1 Q_3^{-1}$ fixes $\bg$, 
then $g_1g_2=1$. So, $\bg$ is an H-M rep. Squares of
elements in $\sQ^*$ act trivially on an H-M rep. This shows $\sQ^*$ acts on the 
set
$(\bg)\sQ^*$ as a quotient of $K_4$. 

Assume $\bg$ is an H-M rep.~and 
$(\bg)Q_1Q_3^{-1}=\bg$. Then $g_1=g_1^{-1}$ and $g_2=g_2^{-1}$: Both are 
involutions.
Similarly, if $(\bg)(Q_1Q_2Q_3)^2=\bg$, then $g_2=g_1$ and $G$ is cyclic (in
particular abelian). Also, $(\bg)Q_1Q_3^{-1}(Q_1Q_2Q_3)^2=\bg$ implies 
$g_2=g_1^{-1}$, contradicting $G$ is centerless.  

Assume $\bg\in
\ni_k$. As $\sQ''\norm M_4$ (Prop.~\ref{isotopyM4}), for $\bg\in \ni$,
$$\Stab_{\sQ''}(\bg)=Q^{-1}\Stab_{\sQ''}((\bg)Q)Q.$$ So, $|(\bg)\sQ''|$ depends 
only on 
the $M_4$ orbit of $\bg$. If $\bg\in \ni_k$, then $\Stab_{\sQ''}(\bg) \bmod
\ker_{k-1}=\Stab_{\sQ''}(\bg \bmod \ker_{k-1})$. This gives the statement on 
faithful action at
level
$k-1$. 

Now assume $p=2$ and $G_0$ is
$p$-perfect and centerless, so these hypotheses apply at all levels 
(Prop.~\ref{fineMod}). With
no loss,  on the general statement on an $H_4$ orbit containing
$\bg=(g_1,g_1^{-1},g_2,g_2^{-1})$, assume
$k=1$.  

Suppose for  $q\in \sQ''\setminus \{1\}$, $(\bg)q=h\bg h^{-1}$ for some $h\in 
G_1$. From
the above, $h\bmod\ \ker_0$ is not the identity. As $q^2$ acts
trivially on $\bg$, $(\bg)q^2=h^2\bg h^{-2}$ with $h^2\in
M_0=\ker_0/\ker_1\eqdef M\setminus \{0\}$ (Lem.~\ref{FrKMTIG}). This contradicts 
$G_1$
being centerless. 

Now assume $G_0=A_5$. We show $\sQ''$ is faithful on {\sl any\/} $H_4$ orbit $O$
in $\ni_1$.  Any orbit at level 1 has elements lying over any element of the 
unique  orbit
at level 0. Anything in $\ni_1$ above an H-M rep.~$\bar\bg\in \ni_{0}$ looks 
like 
$\bg=(g_1,mg_1^{-1}m,ng_2n,g_2^{-1})$ with $m,n\in M_0$, and some $g_1,g_2$ 
generating $G_1$. We handled when
$m$ is in the centralizer of $g_1$, so assume it is not. Apply the previous 
argument when 
$k\nm 1=0$ with
$q=q_1q_3^{-1}$. Conclude: If $(\bg)q=h\bg h^{-1}$, then $h\in G_1$ lifts an 
element of
$A_5$ having order 2. From Prop.~\ref{liftEven}, $h^2\in M\setminus V$. Compute 
the first
two entries of
$(\bg)q^2$ to be $(mm^{g_1^{-1}}g_1mm^{g_1^{-1}}, m^{g_1^{-1}}g_1^{-
1}m^{g_1^{-1}})$.

The remainder of the argument uses Cor.~\ref{A5Morbs}. Let $c\in M\setminus V$ 
be the
generator of the centralizer of $g_1$. Multiply $m$ if necessary by  $c$ to 
assume $m\in
M\setminus V$. Then, 
$h^2=mm^{g_1^{-1}}c=m^{g_1^{-1}}$. Conclude $m=c$ contrary to our assumption. 
For other
elements in $\sQ''\setminus \{1\}$, the argument is similar, though no easier. 
\end{proof} 

\begin{rem}[Lem.~\ref{4HM4NHM}, modular curves and $p\ne 2$] Two involutions 
generate
a dihedral group. So, $\sQ^*$ stabilizing an element of $\ni$ in 
Lem.~\ref{4HM4NHM} comes
precisely from the Hurwitz version of modular curves in 
\S\ref{HurView}.  

Suppose
$r=4$, $\ni(G_k,\bfC)^\inn=\ni_k$ is a level $k\ge 1$ Nielsen class for $G_0=G$, 
centerless and $p$-perfect. Assume $O$, an $H_4$ orbit on $\ni_k$,  
contains an H-M rep.  These are the hypotheses of  Lem.~\ref{4HM4NHM},
except $p\ne 2$. Then, there may be H-M reps.~${}^k\bg\in \ni_k^\inn$ with 
\begin{triv} \label{simulConj}
$h_k\in G_k$ (an involution, or the identity if these are involutions) with
$h_k{}^kg_ih_k^{-1}={}^kg_i^{-1}$,
$i=1,2$ and all 
$k\ge 0$.
\end{triv} 
\noindent We don't know if $h_k\,$s and ${}^k\bg\,$s as in \eqref{simulConj} 
exist in
this case: 
$G_0=A_5$,
$p=5$, $\bfC=\bfC_{3^4}$ and 
${}^kg_1,{}^kg_2\in G_k$ lie over the pair 
${}^0g_1=(1\,2\,3),{}^0g_2=(1\,4\,5)$. 
\end{rem} 

\subsubsection{$q_2$ orbit shortening and $\Cusp_4$} \label{sQorbits2} Again 
$r=4$, with Nielsen classes of any type (though we continue to simplify notation 
by using inner
classes). Slightly abusing the notation of \S\ref{sQpart2} use 
$\alpha^2=(q_1q_2q_3)^2$ and
$\gamma=q_1q_3^{-1}$ as generators of $\sQ''$. Recall the $M_4$ subgroup  
$\Cusp_4=
\lrang{\alpha^2}\times \lrang{\gamma}\xs
\lrang{q_2}$ with $q_2$ switching the two factors on the copy of $K_4$. In this 
subsection
$O(u,v;a)=O(u,v)\subset
\ni^\inn$ is a
$q_2$ orbit with $u=\mpr(\bg)$ for $\bg\in \ni^\inn$. Then, $v$ is the $q_2$ 
orbit length: If
$\mpr(\bg)$ is odd and \eqref{xyswitch} holds, then $v=u$, otherwise $v=2\cdot 
u$. 

For $\bg\in O(u,v)$, $\Stab_{\sQ''}(O(u,v))$ is the subgroup of $\sQ''$ 
stabilizing
$O(u,v)$:
$$\Stab_{\sQ''}(O(u,v))/\Stab_{\sQ''}(\bg) \text{ with cardinality } 
|(\bg)\sQ''\cap O(u,v)|.$$ 
Denote the set of $\sQ''$ orbits on $(\bg)\Cusp_4$ by 
$(\bg)\Cusp^*_4=(O(u,v))\Cusp^*_4$. We
speak of the $\gamma_\infty$ orbit type of $(O(u,v))\Cusp^*_4$. 

\begin{lem} \label{sQorbit}  Suppose $0 <k<v$ and
$(\bg)\alpha^2\gamma=(\bg)q_2^k$; or 
$(\bg)\alpha^2=(\bg)q_2^k$ (or with $\gamma$ replacing $\alpha^2$) and
$k$ is even. Then $k=v/2$. 

Suppose $(\bg)\gamma=(\bg)q_2^k$ 
with $k$ odd. Then, $\Stab_{\sQ''}(O(u,v))=\sQ''$,  $k=v/4$ and 
$(\bg)\alpha^2=(\bg)q_2^{-k}$
and
$(\bg)\alpha^2\gamma=(\bg)q_2^{v/2}$. 

In all cases
$(O(u,v))\Cusp^*_4$ is a
$\gamma_\infty$ orbit of type $(u,v/|(\bg)\sQ''\cap O(u,v)|)$.  
 \end{lem}

\begin{proof} We do all computations on inner Nielsen classes $\bg\in \ni^\inn$. 
Suppose $0<k<v$
and
$(\bg)\alpha^2\gamma=(\bg)q_2^k$ or 
$(\bg)\alpha^2=(\bg)q_2^k$ (or with $\gamma$ replacing $\alpha^2$) and
$k$ is even. As $\sQ''$ commutes with even powers of $q_2$, this gives
$\bg=(\bg)q_2^{2k}$: $k=v/2$. 

Suppose
$(\bg)\gamma=(\bg)q_2^k$  with $k$ odd. Apply $\gamma$ to both sides. From  
\S\ref{sQpart2} conclude  
$$\bg=(\bg)q_2^k\gamma q_2^{-k}q_2^{k}=(\bg)\alpha^2 q_2^{k}.$$ So,
$(\bg)q_2^{-k}=(\bg)\alpha^2$. Apply $q_2$ to both sides of 
$(\bg)\gamma=(\bg)q_2^k$. This 
gives $$(\bg)q_2(q_2^{-1}\gamma q_2)=(\bg)\alpha^2=((\bg)q_2)q_2^k.$$ 
Inductively, this shows
$\Stab_{\sQ''}(O(u,v))=\sQ''$. So, $$(\bg)\gamma 
\alpha^2=(\bg)q_2^k\alpha^2q_2^{-k}q_2^k=
(\bg)\gamma q_2^k=q_2^{2k},$$ and  
$(\bg)\alpha^2=(\bg)q_2^{-k}$, or 
$(\bg)\alpha^2\gamma=(\bg)q_2^{v/2}$. 
\end{proof}

\begin{rem}[Shortening types] \label{shortTypes} When $|\Stab_{\sQ''}(\bg)|=2$, 
the
$\gamma_\infty$ orbit type of
$(O(u,v))\Cusp^*_4$ 
is $(u,v/2)$ if $\Stab_{\sQ''}(O(u,v))=\sQ''$ (two-shortening); of type
$(u,v)$ otherwise. When $\Stab_{\sQ''}(O(u,v))=\{1\}$,  the $\gamma_\infty$ 
orbit type  
of $(O(u,v))\Cusp^*_4$ is $(u,v)$ (no shortening). 

Assume $\Stab_{\sQ''}(\bg)=\{1\}$.  When $\Stab_{\sQ''}(O(u,v))=\sQ''$,
the  $\gamma_\infty$ orbit type  of $(O(u,v))\Cusp^*_4$ is $(u,v/4)$ (total-
shortening).
This is equivalent to $(\bg)\gamma=(\bg)q_2^\ell$ 
with $\ell=v/4$ odd. When $|\Stab_{\sQ''}(O(u,v))|=\lrang{q}= 2$, the 
$\gamma_\infty$ orbit
type of $(O(u,v))\Cusp^*_4$ is 
 $(u,v/2)$ (two-shortening; then $(\bg)q=(\bg)q_2^{v/2}$ with $v$
even). 
\end{rem}

\subsection{H-M and near H-M reps.~in $\ni(G_1,\bfC_{3^4})$}  \label{cuspsH1}  
The subsections of this section together comprise the complete the proof of
Prop.~\ref{HMreps-count}, the main goal of this section. The first two
subsections describe the H-M reps. Then
\S\ref{othLongOrbits} describes how the H-M reps. produce the near
H-M reps.   Denote $\ker(G_1\to A_5)$ by
$M$ and  $\ker(G_1\to \SL_2(\bZ/5))$ by  $V$.

\begin{prop} \label{HMreps-count} There are  2304 elements
in $\ni_1=\ni({}_2^1\tilde A_5,\bfC_{3^4})^\inn$. Exactly 16
are H-M reps.; 16 others are near H-M
reps. (as in \eqref{nearHM}).  The $q_2$ orbit of an H-M
(resp.~near H-M) rep.~contains exactly one H-M (resp.~near H-M)
rep.  All orbits of $\sQ''$ on $\ni_1^\inn$ have length four, and
$\sQ''$ maps H-M  (resp.~near
H-M) reps.~among themselves. In particular, $\gamma_\infty$ orbits in 
$\ni_1^{\inn,\rd}$
containing either H-M or near H-M reps.~have type
$(10,20)$. 

From Lem.~\ref{4HM4NHM}, there are four
orbits of $\gamma_\infty$ on $\ni_1^{\rd}$ containing H-M (resp.~near H-M)
reps.~giving eight  $\gamma_\infty$ orbits on $\ni_1^{\rd}$ 
containing
H-M or near H-M reps.
\end{prop}

\subsubsection{Lifting from $\ni(A_5,\bfC_{3^4})=\ni_0$} 
\label{liftA5} Choose $\bg^*\in \ni_1$ by lifting $(1\,2\,3)=g_1$
and $(1\,4\,5)=g_2$ to $g_1^*$ and $g_2^*$ of order 3. Take
$\bg^*=(g_1^*,(g_1^*)^{-1}, g_2^*,(g_2^*)^{-1})$. As $H_4$
is transitive on $\ni_0$, it suffices to count lifts of
$(g_1,g_1^{-1},g_2,g_2^{-1})=\bg$. Multiply by $18=|\ni(A_5,\bfC_{3^4})^\inn|$ 
to
count elements in $\ni({}_2^1\tilde A_5,\bfC_{3^4})^\inn=\ni_1$. Use 
notation from \S\ref{tA5prop}. Let $\tilde g_1$ and $\tilde
g_2$ be the unique order 3 lifts to $\SL_2(\bZ/5)$ of  $g_1$
and $g_2$.

Lifts of $\bg$ to $\ni_1$ correspond exactly to lifts of
$(\tilde g_1,(\tilde g_1)^{-1},\tilde g_2,(\tilde g_2)^{-1})$
to $\ni_1$. Count  these by counting conjugates of an element
of order 3 by the kernel $V$ from  $G_1\to \SL_2(\bZ/5)$.
If $v\in V$ and $g$ has order 3, then $vgv=gv^{g}v$ (as 
in \S\ref{othLongOrbits}). As $V$ is an irreducible $A_5$ module 
(Cor.~\ref{A5Morbs}), the
set
$\{gv\}_{v\in V}$ gives the complete set of conjugates of $g$ by $M$. 

Three entries of a Nielsen class 4-tuple determine the 4th from the product-one
condition by rewriting entries  $gv^{g}v$. Divide by inner automorphisms from
the kernel from $G_1\to A_5$. As $G_1$ is centerless
(Prop.~\ref{fineMod}), there are  $2^7=2^{12}/2^5$ such lifts
of $\tilde \bg$ to $\ni_1$.  This gives $2^7\cdot 18 = 2304$
total inner Nielsen classes.

\subsubsection{Counting H-M lifts of any H-M rep.~of
$\ni(A_5,\bfC_{3^4})$} Continue using $\bg^*\in \ni_1$. Fix any
H-M representative mapping to $\bg\in \ni_0$. Modulo inner action
of $G_1$, select representatives with $g_1^*,(g_1^*)^{-1}$ in
the first two positions.  Other H-M representatives come from
conjugating $(g_2^*,(g_2^*)^{-1})$ by $\ker_0/\ker_1=M(A_5)$. As
in \S\ref{liftA5} or Cor.~\ref{A5Morbs}, take lifts
by  conjugating
$g_2^*$ by elements of $V$. Or should we choose, by conjugating $g_2^*$ by 
elements 
of
$M\setminus V$. The centralizer of
$\lrang{g_1^*}$ in $M(A_5)$ is a $\bZ/2$ acting nontrivially on $g_2^*$. This
cuts from 16 to 8 the inner classes that are H-M representatives
and lifts of $(g_1,g_1^{-1},g_2,g_2^{-1})$.   This concludes the
part of Prop.~\ref{HMreps-count} counting H-M representatives. 

Since the $g_1$ and $g_4$ positions determine an H-M rep.~$\bg$, 
the $q_2$ orbit of $\bg$ can contain only one H-M rep.  
By inspection, $\sQ''$ maps H-M reps.~to H-M reps. 

\subsubsection{Near H-M reps.~in $\ni({}_2^k\tilde A_5,\bfC_{3^4})$, $k\ge 1$}
\label{othLongOrbits} Unless otherwise said, this subsection is about inner and 
inner
reduced Nielsen classes. The definition of near H-M reps.~is in 
\eqref{nearHM}. A modular representation observation produces near H-M
reps.~$\bg^*\in \ni({}_2^1\tilde A_5,\bfC_{3^4})$
by tweaking an H-M rep.~$\bg=(g_1,g_1^{-1},
g_2, g_2^{-1})\in \ni_1$. 

For $m\in M(A_5)$ and
$g\in {}_2^1\tilde A_5$, the notation $m^g$ is shorthand for
$g^{-1}mg$. This is the same as the right action of $g$ 
in $A_5$ acting on $m$. Let $\hk$ be the complex conjugation
operator in \eqref{hkCong} for $z_1,z_2$ and $z_3,z_4$ as complex
conjugate pairs.

\begin{lem}\label{countWTwenty} Assume $\lrang{g_1,g_2}={}_2^1\tilde A_5$ 
and $g_1,g_2\in \C$. 
Then $g_1g_2$ has order 10 (5 $\bmod \ker_0$). It fixes
a unique nontrivial $c=(g_1g_2)^5\in M(A_5)$ (Cor.~\ref{A5Morbs}). 

Let $d,e\in M(A_5)$. Then $\bg^*=(g_1^{-1},dg_1d,eg_2e,g_2^{-1})\in\ni_1^\inn$ 
if
and only if $d^{g_1}dee^{g_2^{-1}}=1$. In this case, the $q_2$ orbit of
$\bg^*$ in $\ni_1^\inn$ has length 20. 
\end{lem} 

\begin{proof} That $g_1g_2$ has order 10  is a special case of
Lemma \ref{cuspWidth}.  An element of order 5 in $A_5$ acts on
$M(A_5)$ by right multiplying cosets of a $D_5$ in $A_5$
(Prop.~\ref{A5frat}). So $g_1g_2$ fixes one nontrivial
element. It must be $(g_1g_2)^5$.

That $\bg^*$ satisfies the product-one condition is exactly that
$d^{g_1}dee^{g_2^{-1}}=1$. Since $\lrang{g_1^{-1},g_2^{-1}}={}_2^1\tilde A_5$ 
has no
center, what a conjugation on $\bg^*$ does to the 1st and 4th elements 
determines it.
Therefore, the $q_2$ orbit of $\bg^*$ has length  2 times the order of $dg_1d
eg_2e$: The orbit has length 20. 
\end{proof}

Let $G_k$ be the $k$th characteristic quotient of ${}_2\tilde A_5$. We use the 
following
proposition in generality.  Its special case with $g_1,g_2\in {}_2^1\tilde A_5$ 
appears in
precise calculations for level 1 of this Modular Tower. The symbol $\hk$ in
Prop.~\ref{nearHMreps-count} is the complex conjugation operator from
Prop.~\ref{compTest} from complex conjugate pairs of branch points. 

\begin{prop} \label{nearHMreps-count} Suppose $g_1,g_2 \in G_k$ lie over 
$g_1',g_2'\in A_5$ with $g_1'g_2'$ having order $y=3$ or 5. Then, $g_1g_2$ has 
order $y\cdot 2^k$.

If $y=5$, let  
$c=(g_1g_2)^{5\cdot 2^{k\nm1}}$. Then,  
\begin{equation} \label{nearHMform} \bg^*=(g_1,
c^{g_2^{-1}}g_1^{-1}c^{g_2^{-1}}, cg_2c,g_2^{-1}) \in \ni_k^\inn \end{equation}
satisfies the product-one condition and  $c\bg^* c=\hk(\bg^*)$
(so is a near H-M rep.~from \eqref{nearHM} at level $k$). For $k\ge 1$, H-M 
reps.~give Galois covers $\phi: X\to \prP^1_z$
(in $\ni_k^\inn$) over $\bR$ with branch points $\bz$ in complex conjugate pairs 
and $X(\bR)\ne 
\emptyset$. For
near H-M reps.~there are such covers $\phi$ over $\bR$, but $X(\bR)=
\emptyset$. 

When $k=1$, $c$ and $c^{g_2^{-1}}$ are in the $A_5$
orbit labeled $M_5'$ in Cor.~\ref{A5Morbs}. Conversely, given $\bg^*$
satisfying
$c\bg^* c=\hk(\bg^*)$, an H-M rep.~$\bg$  exists giving
$\bg^*$ by \eqref{nearHMform}. 
\end{prop}

\begin{proof} That \eqref{nearHMform} holds is a simple check. The order of the
product $g_1g_2$ is from Lem.~\ref{cuspWidth}. Apply Prop.~\ref{compTest} to an 
H-M
rep.~cover
$\phi:X\to\prP^1_z$ with respect to the
$\hk_0$ operator for two pairs of complex conjugate branch points to compute the 
effect of
complex conjugation over $z_0\in \prP^1_z(\bR)$. It is given by $c_{z_0}$ equal 
the identity. So
all points on $X$ over $z_0$ are real. For a near H-M rep.~ the effect of 
complex conjugation is
given by $c$ as in the statement of the proposition. So, $c_{z_0}$ moves all 
points over
$z_0$. There are no real points on $X$. 

Now assume $k=1$ and $\bg^*$ satisfies $c\bg^*
c=\hk(\bg^*)$ (automatically associated with $\hk$ for $z_1,z_2$
and $z_3,z_4$ as complex conjugate pairs).  Write $\hk(\bg^*)$:
\begin{equation} (cc^{g_2^{-1}}v_1g_1v_1cc^{g_2^{-1}} ,
cc^{g_2^{-1}}g_1^{-1}cc^{g_2^{-1}}, g_2,cg_2^{-1}c).\end{equation}
From Cor.~\ref{A5Morbs}, assume $v_1,c\in M(A_5)\setminus V$. Apply $c\bg^*
c=\hk(\bg^*)$ to see
$$\bg^*=(g_1, v_1g_1^{-1}v_1, cg_2c,g_2^{-1}) = (g_1,
g_1^{-1}v_1^{g_1^{-1}}v_1, cc^{g_2^{-1}}g_2,g_2^{-1}).$$

The product-one condition for $\bg^*$ is equivalent to
$v_1^{g_1^{-1}}v_1= cc^{g_2^{-1}}$.  Then, $c\bg^* c=\hk(\bg^*)$
gives four conditions according to the entries of $c\bg^*c$,
with the last two automatic. The second gives 
\begin{equation} cv_1g_1^{-1}v_1c=cc^{g_2^{-1}}g_1^{-1}cc^{g_2^{-1}},
\end{equation} showing $v_1=c^{g_2^{-1}}$ (Cor.~\ref{A5Morbs}): the first
condition is automatic.  The product-one condition says $c=c^{g_1g_2}$:
$c=(g_1g_2)^5$ generates the centralizer of $g_1g_2$.
\end{proof}

To conclude proving Prop.~\ref{HMreps-count} requires two points.
\begin{itemize} \item A unique 
near 
H-M rep.~is in the $q_2$ orbit of a near H-M rep.
\item   $\sQ''$ is
stable on the set of near H-M reps.  
\end{itemize} As with H-M reps.~the 1st and 4th positions determine
them. So the former is clear. Apply $q_1q_3^{-1}$ to $\bg^*=(g_1,
c^{g_2^{-1}}g_1^{-1}c^{g_2^{-1}}, cg_2c,g_2^{-1})$ to get 
$$(cg_1^{-1}c, g_1, g_2^{-1},c^{g_2^{-1}}g_2c^{g_2^{-1}}).$$ Write 
$g_1'=cg_1^{-1}c$ and 
$c^{g_2^{-1}}g_2c^{g_2^{-1}}=(g_2')^{-1}$, and compute that $c'=c^{g_2^{-1}}$ 
centralizes 
$g_1'g_2'$. Thus, $q_1q_3^{-1}$ maps $\bg^*$ to a near H-M representative. 
Similarly 
for
$(q_1q_2q_3)^2$. 

\section{Cusp widths and the genus of components} \label{cuspWidths} Assume
$r=4$, and $\ni(G,\bfC)$ is a Nielsen class. Distinguishing between absolute and 
inner
Nielsen classes is cumbersome, though computations using them both are 
invaluable (as in
Thm.~\ref{FrVMS}) and similar (add the action of some group $N'$ as in 
\S\ref{refEquiv}). For simplicity, assume inner Nielsen classes so the 
decoration 
$\ni(G,\bfC)^\rd$ means reduced inner classes (unless said otherwise). 

Following two preliminary subsections, for the $(A_5,\bfC_{3^4},p=2)$ Modular 
Tower this section
lists cusps of a given
width from $\gamma_\infty$ action on 
$\ni_1^{\inn,\rd}=\ni_1^\rd$. Possible widths are 2, 4, 6, 8, 10, 12 and 20. The 
aim is to relate all $\gamma_\infty$ orbits to orbits 
of width
20 (especially to H-M reps.).  Lem.~\ref{4HM4NHM} and Prop.~\ref{HMreps-count} 
report
precisely on  H-M and near H-M cusps. Prop.~\ref{expInv} gives the {\sl spin 
separation\/}
ingredient that establishes the  distribution of 
cusps between two $M_4$ orbits $O_1^+$ and $O_1^-$ on $\ni_1^{\inn,\rd}$. 

When no further ramification occurs from
level $k=0$ to level 1 over the elliptic points $j=0$ and $j=1$, 
the genus of level 1 components comes just from the story of $\gamma_\infty$ and
$\lrang{\sh, \gamma_\infty}=\bar M_4$ orbits. \S\ref{ellipFix}  observations on 
this
continue in \S\ref{indSteps}. \S\ref{subComp} has subtle conjugation 
computations for 
$\sQ''$ shortening of 
$q_2$ orbits to $\gamma_\infty$ orbits.  

\subsection{Orbit genus and fixed points of $\gamma_0$ and $\gamma_1$} 
\label{ellipFix} 
For much of this subsection, the prime $p$ is arbitrary. Suppose $O$ is an orbit 
of $\bar M_4$
acting on
$\ni(G,\bfC)^\rd$. Use
\S
\ref{cuspNotat} notation for $\gamma_\infty$ orbits: For $\bg\in O(u,v;a)$, 
$u=\mpr(\bg)$ and $v=\wid(\bg)$. Let $\tr_O(\gamma_i)$ be the number of fixed 
points of
$\gamma_i$ on $O$, $i=1,2$. From Prop.~\ref{j-Line}, the Riemann-Hurwitz formula 
gives the genus
$g_O$ of the reduced Hurwitz space component $\sH_O^\rd$: 
\begin{equation} \label{RHOrbitEq} 2(|O|+g_O-1)=\frac {2 (|O|-\tr_O(\gamma_0))}3 
+
\frac{|O|-\tr_O(\gamma_1)} 2 + \sum_{O(u,v;a)\subset O} v-1. \end{equation} 

Use this to rephrase for $r=4$ our Main Problem~\ref{MPMT} on Modular Towers. 
For $G$ centerless and $p$-perfect, show for each
$\bar M_4$ orbit
$O$ on $\ni(G_k,\bfC)^\rd=\ni_k^\rd$, $g_O\ge 2$ if $k$ is large. If this holds, 
we say $k$ is in
the {\sl hyperbolic\/} range. Thm.~\ref{thm-rbound} says this implies there are 
no $K$ points ($K$ a
number field) on
$\sH(G_k,\bfC)^\rd$ if $k$ is large (possibly larger than the beginning of the 
hyperbolic
range). 

Assume $O_k$ (resp.~$O_{k+1}$) is a $\bar M_4$  orbit 
in $\ni_k^\rd$ (resp.~$\ni_{k+1}^\rd$ over $O_k$). So, cusps of 
$\sH_{O_{k+1}}^\rd$ lie over cusps
of
$\sH_{O_k}^\rd$: $\gamma_\infty$ orbits on $O_{k+1}$ lie over
$\gamma_\infty$ orbits on $O_k$.

\subsubsection{Two helpful assumptions}  
\begin{edesc} \label{bestGenAssump} \item
\label{bestGenAssumpa} $\tr_{O_k}(\gamma_0)+\tr_{O_k}(\gamma_1)=0$; and for 
$\bg'\in O_{k+1}$ over $\bg\in
O_k$,
\item \label{bestGenAssumpb}  if $(\bg')\sh$ not an H-M rep., then 
$\mpr(\bg')=p\cdot\mpr(\bg)$. 
\end{edesc} 
The condition $(\bg')\sh$ not an H-M rep.~in \eql{bestGenAssump}{bestGenAssumpb} 
is equivalent to
$\mpr(\bg')\not=1$, so is absolutely necessary. For a $\gamma_\infty$ orbit
$O(u',v';a')\subset O_{k+1}$ use the notation
$O(u,v;a)\subset O_k$ for cusps of $O_k$ under $O(u',v';a')$. 

\begin{lem} \label{fixedg0g1} If  \eql{bestGenAssump}{bestGenAssumpa} holds (for 
$O_k$), then it
holds with $O_{k+1}$ replacing $O_k$. 
For $\ni_0^\rd=\ni(A_5,\bfC_{3^4})^{\inn,\rd}$ ($p=2$), both 
\eql{bestGenAssump}{bestGenAssumpa}
and 
\eql{bestGenAssump}{bestGenAssumpb} hold for any $\bar M_4$ orbit $O'\subset 
\ni^\rd_k$, $k\ge 1$.
\end{lem}

\begin{proof} Suppose $\bg\in \ni_{k+1}^\inn$ and
$(\bg)Q_1Q_2Q_1=\alpha(\bg)Q'\alpha^{-1}$ for some $\alpha\in
G_k$ and $Q'\in \sQ''$. Reduce all expressions modulo $\ker_k$ to
conclude $\gamma_1$ fixes $\bg \bmod \ker_k$. This is a contradiction. The
same argument works for $\gamma_0$. 

By inspection \eqref{branchCycinn} shows \eql{bestGenAssump}{bestGenAssumpa} 
holds for
$(A_5,\bfC_{3^4},p=2)$, and so at all levels in this $(A_5,\bfC_{3^4})$ Modular 
Tower. 
If $\mpr(\bg)=3$ or 5, \eql{bestGenAssump}{bestGenAssumpb} follows from the 
opening statement of
Prop.~\ref{nearHMreps-count}. The only other possibility is that $(\bg)\sh$ is 
an H-M rep., but
$(\bg')\sh$ above it is not. Lem.~\ref{cuspWidth} says if $p\,|\,\mpr(\bg_k)$, 
then 
$\mpr(\bg_{k+1})=p\cdot
\mpr(\bg_k)$ for $\bg_{k+1}$ over $\bg_k$. That completes the proof.  
\end{proof}

Assume at level $k$ of a Modular Tower there are still components of genus 0 or
1. Lem.~\ref{RHOrbitLem} inspects the contribution of cusp ramification toward 
the genus of
components at level $k+1$. Three phenomena play a role in this contribution for 
each cusp in
going from level $k$ to level $k+1$.  
\begin{edesc} \item Detecting if condition \eqref{xyswitch} in
Prop.~\ref{gamOrbit} changes.  
\item Deciding if there is a multiplication by $p$ factor as in
\eql{bestGenAssump}{bestGenAssumpb}. 
\item Computing $q_2$ orbit shortening changes from one level to another. 
\end{edesc} As above, $O_k$ is a $\bar M_4$ orbit in $\ni_k^\rd$; $O_{k+1}$ is a 
$\bar M_4$ 
orbit in
$\ni_{k+1}^\rd$ above it. 
 
\begin{lem} \label{RHOrbitLem} For $\bg_k\in O(u,v;a)$, assume $Z(g_2,g_3)\cap
\lrang{g_2g_3}=\{1\}$ so the $q_2$ orbit type is $(u,2u)$ unless
$u$ is odd and \eqref{xyswitch} holds  (Prop.~\ref{gamOrbit}). For each
$\gamma_\infty$ orbit
$O(u,v;a)\subset O_{k}$ let $\alpha(u,v;a)=2$ in the former case, 1 in the 
latter. If
$\alpha(u,v;a)=2$, then $\alpha(u',v';a')=2$ for $O(u',v'; a')$ over $O(u,v;a)$. 
Let
$\ind(u',v';a')$ be the index of ramification of $O(u',v'; a')$ over $O(u,v;a)$. 
Also, let
$\beta(u',v';a')=\frac{u'}u$ (1 or $p$; automatically $p$ if $p\,|\,u$ from
Lem.~\ref{FrKMTIG}). Depending on the amount of $q_2$ orbit shortening 
(\S\ref{sQorbits2}), let
$\mu(u,v;a)=1$ (no shortening), 2 (two-shortening) or 4 (total shortening). 
Then,
\begin{equation}\label{indFactor} 
\ind(u',v';a')=\frac{\alpha(u',v';a')\beta(u',v';a')\mu(u,v;a)}
{\alpha(u,v;a)\mu(u',v';a')}.\end{equation} 

Suppose
$g_{O_k}=1$. Then $g_{O_{k+1}}\ge 2$ if for some 
$\gamma_\infty$ orbit $O(u',v';a')\subset O_{k+1}$ over $O(u,v;a)\subset O_k$, $ 
\ind(u',v';a')>1$.
 
Suppose $g_{O_k}=0$. Then, 
\begin{equation} \label{cuspRamEst} 2\Bigl(\frac{|O_{k+1}|}{|O_k|}+g_{O_{k+1}}-
1\Bigr)\ge
\sum_{O(u',v';a')\subset O_{k+1}}
\ind(u',v';a')-1.\end{equation} If
$\tr_{O_k}(\gamma_0)+\tr_{O_k}(\gamma_1)=0$, then equality holds in
\eqref{cuspRamEst}. 
\end{lem} 

\begin{proof} The references explain most of this lemma. Given $O(u',v';a')$ 
over $O(u,v;a)$ the
ramification index of the respective cusps is exactly $\frac{v'}v$ which by 
previous comments is  
$\ind(u',v';a')$ as in \eqref{indFactor}. It is well-known that if a $X\to Y$ is 
a covering of
projective nonsingular curves with the lower curve of genus 1, then the genus of 
$X$ is 1 if and
only if the cover is unramified. Formula \eqref{cuspRamEst} expresses the 
Riemann-Hurwitz formula
applies to the relative curve covering $\sH^\rd_{O_k+1}\to \sH^\rd_{O_k}$, when 
the latter has
genus 0. \end{proof} 

We use Lem.~\ref{RHOrbitLem} to show what goes into computing the genus of the
two $\sH(G_1,\bfC_{3^4})^\rd$ components.  

\begin{cor} \label{orbitGenus} There are two orbits $O_1^+$ and $O_1^-$ of $\bar 
M_4$ on
$\ni(G_1,\bfC_{3^4})^\rd$, each of degree 16 over the unique $\bar M_4$ orbit 
$\ni_0$. The genus
of the orbit $O_1^+$ containing H-M reps.~is 12. The other orbit $O_1^-$ has 
genus 9.\end{cor} 

\begin{proof} Apply \eqref{cuspRamEst} to $O_1^+$. From Prop.~\ref{nearHMbraid}, 
one $\bar M_4$
orbit on
$\ni_1^\rd$ contains all H-M and near H-M reps. Call this orbit $O_1^+$. 
Prop.~\ref{expInv} says there are exactly eight  $\gamma_\infty$ orbits with the 
following
properties: 
\begin{edesc} \item 5 divides 
$\mpr$; 
\item they are in the $\bar M_4$ orbit of a near H-M rep.
\end{edesc} 
Prop.~\ref{HMreps-count} says the $\gamma_\infty$ cusp type of an H-M or near H-M rep.~is
$(10,20)$. So, a cusp of $O_1^+$ lying over an H-M rep.~ of $\ni_0^\rd$ has
ramification index 4. Let $\bar \sH_1^+\to\prP^1_j$ be the component of $\bar 
\sH_1^\rd$
corresponding to $O_1^+$. Conclude: Each cusp (eight total) of $O_1^+$ with 
$\mpr$ divisible by 5
has ramification index 4 over the cusp below them at level 0. Together they 
contribute
$2\cdot 4\cdot 3=24$ to the right side of \eqref{cuspRamEst}. 

Similarly, consider cusps of $O_1^+$ over cusps of $\ni_0^\rd$ with $\mpr$ equal 
3.
Prop.~\ref{bgv} gives  a similar conclusion about cusps at
level 0 with middle product 3. Cusps above them on $O_1^+$ have ramification 
index 4 ($\sQ''$
does not shorten them). They also contribute $2\cdot 3\cdot 4=24$ to the right 
side of
\eqref{cuspRamEst}. 

Now apply Prop.~\ref{sqOnM} for the contribution of cusps on $O_1^+$ over $\sh$ 
applied
to  level 0 H-M reps. This contributes $2$ for the shift of H-M reps. and 
$2\cdot 2$
for the others to the right side of \eqref{cuspRamEst}. So, the right side of
\eqref{cuspRamEst} is 54. The expression $2(16+g_{O_1^+} -1)=4\cdot 3\cdot 
4+2\cdot 3$
gives $g_{O_1^+}=12$.  

Let $O_1^-$ be the collection of cusps at level 1 not in $O_1^+$. First assume 
they
all lie in one $\bar M_4$ component for the computation of the genus of this 
orbit. From
Prop.~\ref{bgv} there are eight cusps in
$O_1^-$ with width 6, and four with width 12. Similarly, there eight cusps in 
$O_1^-$ with
width 10, and four with width 20. Finally, Lem.~\ref{tough2-4} gives 8 type 
(2,4) $\gamma_\infty$
orbits in $O_1^-$.  To complete the calculation above for $g_{O_1^-}$, list 
respective
contributions to the right of \eqref{cuspRamEst}: Type (2,4) contribute 8; type 
(6,6)
contribute 8; type (6,12) contribute $4\cdot 3$; type (10,10)
contribute 8 and type (10,20) contribute $4\cdot 3$. So, 
$2(16+g_{O_1^-} -1)=48$ gives $g_{O_1^-}=9$.  

Since it is true at level 0, it is also
true at level 1 that every
$\bg$ not in
$O_1^+$ is in the $\bar M_4$ orbit of an element $\bg'$ with
$5\,|\,\mpr(\bg')$. From Prop.~\ref{bgv} (and its notation) any cusp of $L_{20}$ 
(resp.~$L_{10}$) connects to each cusp of
$L_6$ (resp.~$L_{12}$). So, to prove all elements of
$O_1^-$ lie in one
$\bar M_4$, it suffices to connect some cusp of $L_6$ (resp.~$L_{20}$) to some 
cusp of $L_{12}$
(resp.~$L_{10}$). This shows a component $\bar\sH'$ containing a cusp of $O_1^-$
has degree sixteen over a cusp at level 0 with
$\mpr=3$. So $\bar\sH'\to \bar\sH_0^\rd$ has degree sixteen everywhere. For the 
degree over
every cusp to be 16 forces including all cusps in
$O_1^-$. To join something in $L_{10}$ to something in $L_{20}$ consider a 
$\gamma_\infty$ type
(2,4) orbit in $O_1^-$. A representative for such an orbit
has $\bg=(g_2^{-1},g_1,ag_1^{-1}a,bg_2b)$. Lem.~\ref{tough2-4} shows each type 
(2,4)
$\gamma_\infty$ orbit in $O_1^-$ is such an element. This concludes the proof of 
the corollary. 
\end{proof}

\subsubsection{The mystery of $\gamma_0$ and $\gamma_1$ fixed points} 
 Consider orbits of length one for
$\gamma_0$ and
$\gamma_1$ (of respective orders 3 and 2) acting on
$\ni(G,\bfC)^\rd$.

Suppse $G\le N'\le S_n$ is a situation for absolute equivalence, as in 
\S\ref{absEquiv}.  Knowing
the length one orbits tells what is the contribution
of $\gamma_0$ and $\gamma_1$ to the genus of components of
$\sH(G,\bfC)^\rd/N'$ in \eqref{RHOrbitEq}. Consider $\gamma_1$; $\gamma_0$ is 
similar. 

Let $\bz=\{1, -1, i, -i\}$ be a set of branch points representing the elliptic 
point for
$\gamma_1$.  Then $\alpha: z\mapsto iz\in \PSL_2(\bC)$ cycles the set $\bz$. A
fixed set of classical generators of $\pi_1(U_\bz,z_0)$ produces a list
of covers (up to $N'$-equivalence) $\phi_i: X_i\to \prP^1_z$,
$i=1,\dots,t$. Composing the $\phi_i\,$s with $\alpha$ permutes
them. For some choice of classical generators, this action is 
$q_1q_2q_1$ modulo the action of $\sQ''$. 

\begin{prob}[Elliptic fixed points] \label{ellFixedPtr} Generalize 
Rem.~\ref{serFixedPt} to the general case of inner reduced Hurwitz spaces when 
$r=4$. Further, if
$\sH(G_k,\bfC)^{\inn,\rd}$ is a Modular Tower, with $p$-perfect and centerless 
$G_0$, when can there
be a projective system of Nielsen classes fixed by $\gamma_i$ for $i=0$ or 1? 
\end{prob} The
Modular Tower version of this question makes sense for any $r$ applied to the 
orbifold points in
$J_r$: When can a Modular Tower have a projective system of Nielsen class 
representatives fixed by a
nontrivial element in $\bar M_r$ associated with an orbifold stabilizer? 

\begin{exmp}[Action on $\ni(A_5,\bfC_{3^4})$] \label{shA5k=0}
According to \S\ref{A5Cycles}, $\sQ''$ acts
trivially on the list of 
Table \ref{lista5} (\S\ref{znezero}).  Note: $\gamma_1$ fixes ${}_7\bg$
from
$\ni(A_5,\bfC_{3^4})^\abs$. Yet, $\gamma_1$ fixes no item of
$\ni(A_5,\bfC_{3^4})^\inn$. Since $\sQ''$ acts trivially,
this means, from the list $\phi_i: X_i\to \prP^1_z$,
$i=1,\dots,9$, of degree 5 covers, exactly one of $\alpha\circ
\phi_i$ is equivalent to $\phi_i$. Here is a fact about the
list: Suppose an item from it has two 3-cycles (not necessarily
consecutive) with exactly two integers of common support. Up to
conjugation by $S_5$ these 3-cycles are $((1\,2\,3),(2\,1\,4))$
where the common support integers (1 and 2 here) appear in {\sl
opposite\/} order in the second 3-cycle. Further,  for ${}_7\bg$,
this is true for all consecutive pairs of 3-cycles, especially
including the 4th and 1st, taken in that order. Now, apply
$\gamma_1$ as the {\sl shift\/} of \S\ref{shiftOp}.  \end{exmp}

\subsection{Subtle conjugations} \label{subComp} 
Conjugation by an element of $G_1$ is determined by its action
on two generators $g_1,g_2$ of order 3. Several computations require the precise 
effect of
that conjugation given the generators and what the conjugation does modulo 
$\ker_0$. Use
the notation of Lem. \ref{sQorbit} for elements of $\sQ''$. For $\bg\in 
\ni_1^\inn$, let
$O_\bg$ be its $q_2$ orbit. Let $\delta\in \sQ''\setminus 0$. If $(\bg)\delta 
\in
O_\bg$, we say $\delta$ {\sl shortens\/} $\bg$ (as in \S\ref{sQorbits2}). 
Variants: $\sQ''$
shortens 
$\bg$ or
$\sQ''$ shortens $O_\bg$.   Use 
$\bh_+=((1\,2\,3),(1\,3\,2),(1\,4\,5),(1\,5\,4))$ and
$\bh_-=((1\,2\,3),(1\,3\,2),(1\,5\,4),(1\,4\,5))$ as representatives of the two 
H-M reps.~in $\ni_0^{\inn}$. \S\ref{cuspNotat} has notation
for the type of  orbits. 

\begin{prop} \label{sqOnM} Let
$g_1,g_2\in G_1$ be order 3 lifts of  $(1\,2\,3),\ (1\,4\,5)\in 
A_5$. Denote
$(2\,3)(4\,5)$, $(2\,4)(3\,5)$ and $(2\,5)(3\,4)$ respectively by
$\alpha_{23}$, $\alpha_{24}$ and $\alpha_{25}$. Each
$\bg\in \ni_1^\inn$ with $\mpr(\bg)=2$ has the form 
$\bg=(g_1,g_2,cg_2^{-1}c,dg_1^{-1}d)$ with $c,d\in M(A_5)$. (So,
$c^{g_2^{-1}}c=dd^{g_1}$.) If $(\bg)q_1q_3^{-1}$ (resp.~$(\bg)(q_1q_2q_3)^{2}$,
$q_3(q_1q_2q_3)^{2}q_3^{-1}$) is conjugate to $(\bg)q_2^2$,  then the
conjugation  is 
by a
lift of $(1\,2\,3)\alpha_{24}$
(resp.~$\alpha_{25}$, $(1\,4\,5)^{-1}\alpha_{23}$). 

Suppose $\bg'$ is a near H-M rep. Then, $\sQ''$ shortens $\bg=((\bg')\comp)\sh$.  
Thus, shifted complements of near H-M reps.~fall into two pairs of 
$\gamma_\infty$ type (2,2) orbits. All other $\bg\,$s with $\mpr(\bg)=2$ fall 
into
type (2,4) $\gamma_\infty$ orbits.  
\end{prop} 

Three subsections cover the proof of Prop.~\ref{sqOnM}.
The first  establishes
the relevance of 
$\alpha_{2j}$, $j=3,4,5$, to the existence
of appropriate conjugations. 
\S\ref{a24} shows neither  $q_1q_3^{-1}$ nor
$q_3(q_1q_2q_3)^{2}q_3^{-1}$ shorten any $\bg$ with
$\mpr(\bg)=2$. \S\ref{a25} shows $(\bg)(q_1q_2q_3)^{2}$ shortens the
complement of a near H-M rep. As this is an even power of $q_2$, 
Rem.~\ref{shortTypes} notes
this suffices to determine exactly the shortening type: It is two-shortening.   

\subsubsection{Preliminaries} With $c,d\in M$,  
$\bg=(g_1,g_2,cg_2^{-1}c,dg_1^{-1}d)\in \ni_1^\inn$ as above:  
\begin{equation} \label{bgq2} (\bg)q_2^2= (g_1,cc^{g_2^{-1}} g_2
c^{g_2^{-1}}c,c^{g_2^{-1}} g_2^{-1}
c^{g_2^{-1}},dg_1^{-1}d).\end{equation} Then, $g_1^{-1}(\bg)q_1q_3^{-1}g_1= 
(g_2,g_1,d^{g_1}g_1^{-1}d^{g_1}, d^{g_1}dcg_2^{-1}d^{g_1}dc)$. If this is 
conjugate to
\eqref{bgq2}, it is by a lift of an $A_5$ conjugation switching $(1\,2\,3)$ and
$(1\,4\,5)$. The element $\alpha_{24}$ gives this conjugation. 

Consider if $(\bg)(q_1q_2q_3)^{2}$ is conjugate to \eqref{bgq2}. Such a
conjugation is by a lift of an $A_5$ conjugation mapping $(1\,2\,3)$ to 
$(1\,5\,4)$ and
$(1\,4\,5)$ to $(1\,3\,2)$. Thus, the conjugation is a lift of $\alpha_{25}$. 

Now consider if  $g_2(\bg)q_3(q_1q_2q_3)^{2}q_3^{-1}g_2^{-1}= 
(d^{g_1}g_1^{-1}d^{g_1},c^{g_2^{-1}}g_2^{-1}c^{g_2^{-1}} ,g_2,
g_1)$ is conjugate to $(\bg)q_2^2$. If so, it is by a lift of an $A_5$ 
conjugation mapping
each of 
$(1\,2\,3)$
and $(1\,4\,5)$ to their inverses. The element $\alpha_{23}$ gives this 
conjugation. 

\subsubsection{$q_1q_3^{-1}$ shortens no $\bg$ with $\mpr(\bg)=2$} 
\label{a24} Suppose $g\in G$ normalizes a subgroup $H$. Denote the
centralizer of 
$g$ in $H$ by  $Z_H(g)$. Denote  $Z_M(g_i)$ by $\lrang{m_i}$, $i=1,2$. 
From the last statement in Cor.~\ref{A5Morbs}, with no loss 
$(g_1,g_2,cg_2^{-1}c,dg_1^{-1}d)$ uniquely determines $c$ and $d$
knowing also $c,d\in M(A_5)\setminus V$.  Suppose
$g_1^{-1}(\bg)q_1q_3^{-1}g_1$ is conjugate to
\eqref{bgq2}.  Then some
lift $\alpha$ of $\alpha_{24}$  conjugates as follows:
\begin{equation} \label{propAlpha} 
\begin{array}{rl} &g_1\mapsto g_2,\ c^{g_2^{-1}} g_2^{-1}
c^{g_2^{-1}}\mapsto d^{g_1}g_1^{-1}d^{g_1}, \\
&\ cc^{g_2^{-1}} g_2
c^{g_2^{-1}}c\mapsto g_1 \text{ and } dg_1^{-1}d\mapsto
d^{g_1}dcg_2^{-1}d^{g_1}dc.\end{array}\end{equation}

The 1st, 2nd and 3rd expressions give the effect of conjugating
$\alpha^2\in M(A_5)\setminus V$ (Lemma \ref{FrKMTIG}) on $g_1$ two ways: 
\begin{equation} \label{alsqm1}
\alpha^2=c^{g_2^{-1}\alpha}d^{g_1}m_1=(cc^{g_2^{-1}})^\alpha
m_1.\end{equation} Conclude 
$d^{g_1}=c^\alpha$.  Similarly, figure the effect of $\alpha^2$ conjugating 
$g_2$:  
\begin{equation} \label{alsqm2}
\alpha^2=(d^{g_1}c^{g_2^{-1}\alpha})^\alpha
d^{\alpha}d^{g_1}dc m_2=(d^{g_1}c^{g_2^{-1}\alpha})^\alpha
m_2.\end{equation} 
Now we show $\alpha$ satisfying \eqref{propAlpha} does not exist.

From \eqref{alsqm2} and the product one condition
($c^{g_2^{-1}}c=dd^{g_1}$)  $c^{g_2^{-1}}=d^\alpha$ and  
\begin{triv} \label{alphaInv} $d^{\alpha}c=d^{g_1}d=c^\alpha d$, an expression
invariant under
$\alpha$. 
\end{triv}

From \eqref{alsqm1} and \eqref{alsqm2}, $m_1=m_2$. 
This is false: The nontrivial element
$m_2$ would centralize $G_1=\lrang{g_1,g_2}$. As $G_1$
has no center (Prop.~\ref{fineMod}), $\alpha$ doesn't exist.  

\subsubsection{$(q_1q_2q_3)^{2}$ shortens the shift of a near H-M 
rep.~complement} 
\label{a25} The following shows there is a natural braid taking the
reduced class of a complement of a near H-M rep. over $\bh_+$ to one over $\bh_-
$. So, for
reduced classes, complements of near H-M reps. are similar to a pair of H-M 
reps. $\bg$ and
$((\bg)q_1$ naturally paired as being over $\bh_+$ and $\bh_-$. 

Suppose
$(\bg)(q_1q_2q_3)^{2}$ is conjugate to
\eqref{bgq2}. A  lift
$\alpha$ of
$\alpha_{25}$  conjugates as  
follows:
\begin{equation} \label{a25eq} \begin{array}{rl} g_1\mapsto cg_2^{-1}c,\ 
&c^{g_2^{-1}} g_2^{-1} c^{g_2^{-1}}\mapsto g_1, \\
cc^{g_2^{-1}} g_2
c^{g_2^{-1}}c\mapsto dg_1^{-1}d &\text{ and } dg_1^{-1}d\mapsto
g_2.\end{array}\end{equation} 

Here are the analogs of \eqref{alsqm1} and \eqref{alsqm2}:
\begin{equation} \label{shsq} 
\alpha^2= c^{g_2^{-1}\alpha}c^\alpha m_1= c^{g_2^{-1}\alpha}dm_1 \text{
and }
\alpha^2=c^{g_2^{-1}} cm_2= c^{g_2^{-1}}d^\alpha m_2. \end{equation}
Conjugate the latter by $\alpha$ to conclude $m_1=m_2^{\alpha}$ and 
$c^\alpha=d$ together are equivalent to  $(\bg)(q_1q_2q_3)^{2}$
shortening $\bg$. The expression $m_1=m_2^{\alpha}$ is automatic from
$g_1^{\alpha}=cg_2^{-1}c$ and conjugating $m_1m_1^{g_1}=1$ by $\alpha$. 

As in the proof of Cor.~\ref{A5Morbs}, compute the effect of $\alpha_{25}$ on
$1_H,\dots, 6_H$. By explicit computation $g_1$ acts as
$\beta_{g_1}=(1_H\,3_H\,5_H)(2_H\,4_H\,6_H)=(1\,3\,5)(2\,4\,6)$. Similarly,
$g_2$ acts as $\beta_{g_2}=(1\,3\,2)(4\,5\,6)$. Finally,  $\alpha_{25}$
acts as $\beta_{\alpha}=(3_H\,1_H)(2_H\,5_H)=(3\,1)(2\,5)$. The remaining list 
of possible $c$ values 
is in the following lemma.
\begin{lem} \label{liftorder2} With $d=c^{\alpha_{25}}$, $(q_1q_2q_3)^{2}$ 
shortens the
$\gamma_\infty$ orbit of $\bg$ if and only if
$u=c^{g_2^{-1}\alpha_{25}}c^{\alpha_{25}} m_1$ is invariant under
$\alpha_{25}$. Given a value of $u$ (or $c$) with this property, all
others arise by running over all possible lifts of $\alpha_{25}$. 
\end{lem}

\begin{proof} There is an orbit shortening for $(q_1q_2q_3)^{2}$ only if
there is  an
$\alpha$ (lifting $\alpha_{25}$) with $\alpha^2=u$. This implies $u$ is
invariant under $\alpha_{25}$. Given one such lift $\alpha$,
Cor.~\ref{A5Morbs} shows you get all
others by multiplying this $u$ by $mm^{\alpha_{25}}$ as $m$ runs over $M$.
Given a lift $\alpha$ giving $u$, multiplying $\alpha$ by $m$ produces $u
mm^\alpha=(m\alpha)^2$. 

Producing one such $\alpha$ is the final step. 
The Cor.~\ref{A5Morbs} proof gives explicit action of $g_1,g_2, \alpha_{25}$ on 
the
cosets of the $D_5=\lrang{(1\,3\,4\,2\,5),(1\,2)(3\,4)}$. 
Respectively:  
$$\beta_{g_1}=(1\,3\,5)(2\,4\,6), \beta_{g_2}=(1\,3\,2)(4\,5\,6), 
\beta_{\alpha_{25}}= (3\,1)(2\,5).$$
Consider a near H-M rep.~as in \eqref{nearHMform}: 
$$\bg'=(g'_1, (c')^{(g'_2)^{-1}}(g'_1)^{-1}(c')^{(g'_2)^{-1}},
c'g'_2c',(g'_2)^{-1}) \in
\ni_1^\inn, $$ with $c'$ centralizing $g_1'g_2'$. For reasons coming up, we
take  
$g_1'=g_2$, a lift of $(1\,4\,5)$ and 
$(g_2')^{-1}=g_1$, a lift of $(1\,2\,3)$. 

Let $t$ be the centralizer of $(g_1')^{-1}g_2'=(g_1g_2)^{-1}$. Then, compute
the shift of the complement of $\bg'$ to get 
$$\bg=(g_1,g_2,t(c')^{g_1}m_2g_2^{-1}m_2(c')^{g_1}t, tc'm_1g_1^{-1}m_1c't)$$
with $c=t(c')^{g_1}m_2$ and $d=tc'm_1$ and $c'$ centralizing $g_2g_1^{-1}$. 

Since $m_2^\alpha=m_1$, the desired shortening amounts to showing
$t^\alpha(c')^{g_1\alpha}=tc'$ and (using that by definition
$m_2^{g_2^{-1}}=m_2$) 
$cc^{g_2^{-1}}m_2=tt^{g_2^{-1}}(c')^{g_1}c'm_2$ centralizes
$\alpha$. Explicit computations are reassuring:
$\beta_{g_1}\beta_{g_2}=(1\,2\,5\,3\,6)$ and so $t=(0,0,0,1,0,0)$.
Similarly, $\beta_{g_2}\beta_{g_1}^{-1}=(2\,5\,4\,3\,6)$, so
$c'=(1,0,0,0,0,0)$. Thus, $t^\alpha=t$ and $(c')^{g_1\alpha}=c'$. The other
check works as easily.  
\end{proof}

\subsection{Length two and four cusp widths} From \S\ref{a24} and 
\S\ref{a25} there is a concise relation between
$\gamma_\infty$ orbits of length two and H-M and near H-M reps. As
previously, use $u$ for $\mpr(\bg)=\ord(g_2g_3)$. 

\subsubsection{Types of width 2 and 4 cusps} Continue the notation from 
\S\ref{subComp} for
$\bh_+$,
$\bh_-$ and the types of cusps. 

\begin{prop} \label{countWTwoFour} There are $16$ total  $\bg\in
\ni_1^\inn$  with $u=1$; all have $(\bg)\sh$ an H-M rep. This 
gives  two $\gamma_\infty$ orbits
$O(1,2;1)$ and 
$O(1,2;2)$ of type $(1,2)$. 

$Sixteen$ total $\bg\in \ni_1^\inn$ have these properties:
\begin{edesc} \item $g_2g_3$ has order $u=2$; and
\item $\bg$ is in a $\gamma_\infty$ orbit of type $(2,2)$. 
\end{edesc}
Such $\bg$ have $(\bg)\sh$ a complement of a near H-M rep. These account
for the $\gamma_\infty$ orbits
$O(2,2;1)$ and
$O(2,2;2)$ of type
$(2,2)$, giving all width 2 cusps. 

There are $16\cdot 16$ total $\bg\in\ni_1^\inn$ with $g_2g_3$ having 
$u=2$.
From these there are 14 $q_2$ orbits in
$\ni_1^{\inn,\rd}$ of type $(2,4)$. All $q_2$ orbits $O$ in this
proposition have  $(O)\sh$ modulo $\ker_0$, mapping surjectively to 
$\{\bh_+,\bh_-\}$.  All H-M and near H-M reps.~lie in one $M_4$ orbit
containing the width two cusps.  
\end{prop} 

\begin{proof} 
Each $\bg$ with $\mpr(\bg)=1$ is the shift of an H-M 
rep.~$(h_1,h_1^{-1},h_2,h_2^{-1})$.  Further, applying $q_2$ to such a $\bg$ 
gives 
another
such element. By inspection  
$\sQ''$ is stable on this set. So, according to Prop.~\ref{HMreps-count}, the 
shift
applied to H-M reps. contributes a total of two length 2 orbits 
for $q_2$ on reduced classes. 

Now, consider the  case $\bg$ has 
$u=2$: $g_2g_3\in
\ker_0/\ker_1=M(A_5)\setminus \{0\}$. With no loss,
$g_1=(1\,2\,3)$ and $g_4=(3\,2\,1) 
\bmod \ \ker_0$. Further, assume $g_2=(1\,4\,5)$ or $(1\,5\,4)
\bmod \ker_0$. Fix
$g_1$. There are $2^3$ choices of $g_2$, all lifts of 
$g_2 \bmod\ \ker_0$ modulo 
conjugation by the centralizer of $g_1$.  Then, there are $2^4$ lifts
of $g_4 \bmod\ 
\ker_0$, now determining $g_3$. 

So, there are $16\cdot 16$ total
$\bg\in\ni_1^\inn$ with $g_2g_3$ having $u=2$. Since $\sQ''$ acts faithfully
(Prop.~\ref{HMreps-count}), there are $16\cdot 4$ elements of $\ni_1^{\inn,\rd}$ 
in 
$q_2$
orbits of length 2 or 4. \S\ref{w2cusps} completes showing all H-M and
near H-M reps.~fall in one $M_4$ orbit. \end{proof}

\newcommand{\cc}{\bar c} 

\subsubsection{Cusp widths containing H-M reps.} Lem.~\ref{qi2} and 
Lem.~\ref{cuspWidth} 
explain the eight length 20 orbits of $\gamma_\infty$ on
$\ni(A_5,\bfC_{3^4})^{\inn,\rd}$  containing
H-M and near H-M reps. Part of the next lemma applies to any Modular
Tower with $r=4$. 

\newcommand{\jmp}{\text{\rm jmp}}

\begin{lem} \label{cuspWidth} Let $G_k$ be the $k$th
characteristic quotient of ${}_2\tilde A_5$. For an H-M
rep.~$\bg=(g_1,g_1^{-1},g_2,g_2^{-1})\in \ni(G_k,\bfC_{3^4})^\inn$, 
$\mpr(\bg)$ is $2^k\cdot 5$. The $q_2$ orbit of
of  $\bg$ has length $2^{k+1}\cdot 5$. Further, only one H-M
rep.~is in a given $q_2$ orbit. 

Generally, let $\{\bg_k\}_{k=0}^\infty$ be a projective system of Nielsen class
representatives in the $(A_5,\bfC_{3^4}, p=2)$ Modular Tower. If $\mpr(\bg_0)\ne 
1$, then 
$\mpr(\bg_k)=2^{k}\mpr(\bg_{0})$,
$k\ge 0$. For $\bg\in
\ni(G_1,\bfC_{3^4})^\inn$ with $\mpr(\bg)\ne 1$, its $q_2$ orbit 
has length 4, 12 or 20.

Let $\ni(G_k,\bfC)$ be the $k$th level Nielsen class  in any
Modular Tower (any $p$).  For $\tilde \bg=\{\bg_k\}_{k=0}^\infty$, a projective 
system
of Nielsen class representatives, either: 
\begin{edesc} \label{projDivisible} \item $\mpr(\bg_k)=\mpr(\bg_0)$ for all $k$; 
or 
\item some smallest $k_0\eqdef\jmp(\tilde \bg)$ satisfies $p\, |\,
\mpr(\bg_{k_0})$ and $$\mpr(\bg_k)=p^{k-k_0}\mpr(\bg_{k_0}),\ k\ge k_0.$$ 
\end{edesc} 
\end{lem} 

\begin{proof} We show the conclusion of \eqref{projDivisible} first. Consider 
the
projective system
$\tilde
\bg=\{\bg_k\}_{k=0}^\infty$ in the Modular Tower for $(G_0,\bfC,p)$. Let 
$g_{2,k}, g_{3,k}$
be the second and 3rd entries of $\bg_k$.  Suppose $k$ exists so $p| g_{2,k}
g_{3,k}$. By assumption, $g_{2,k\np1} g_{3,k\np1}$ is a lift of $g_{2,k} 
g_{3,k}$
to
$G_{k+1}$.  Thus, Lem.~\ref{FrKMTIG} says $\mpr(\bg_{k+1})=p\mpr(\bg_{k})$. 
Inductively
apply this for the conclusion of the lemma.  

For $\bg\in \ni(A_5,\bfC_{3^4})$ the value of $\mpr(\bg)$ is either 1, 3
or 5. If $\mpr(\bg)=1$, then $\bg$ is an H-M rep. If $\mpr(\bg)=3$, then
$\lrang{g_2,g_3}=A_4$, and the 3-tuple $(g_2,g_3,(g_2g_3)^{-1})$ satisfies the
product-one and genus 0 ($2(4+g-1)=6$ implies $g=0$) conditions of 
Prop.~\ref{serLift}.  Conclude that $s(g_2,g_3,(g_2g_3)^{-1})=-1$.
Equivalently: If $(\hat g_2,\hat g_3)$ are the lifts of $(g_2,g_3)$ to 
$\Spin_4$, then
$\ord(\hat g_2\hat g_3)=6$.  

Similarly, if $\mpr(\bg)=5$, then $\lrang{g_2,g_3}=A_5$,
the genus 0 condition holds, and the conclusion from Prop.~\ref{serLift} is that  
$\hat g_2\hat g_3$ has order 10. From Lemma~\ref{qi2}, if $2|\mpr(\bg)$, then 
the length of
the orbit of
$q_2$ on $\bg$ is twice  $\ord(g_1^{-1}g_2)$ with $\bg$ written as a 
perturbation of an H-M
rep.~$(g_1,ag_1^{-1}a,bg_2b,g_2^{-1})$, $a,b\in M$. 
 \end{proof}

\begin{rem} \label{mpr1} When $r=4$, $\mpr(\bg)=1$ is equivalent to $(\bg)\sh$ 
is an H-M
rep. Suppose 
$\tilde
\bg$ is a projective system of Nielsen class representatives as in 
Lem.~\ref{cuspWidth}
with $\mpr(\bg_0)=1$ and $p=2$. Then, $\jmp(\tilde \bg)$ is the smallest integer 
$k_0$
with  $(\bg_{k_0})\sh$ not an H-M rep. If all the
$\bg_k\,$s are H-M reps., then $\jmp(\tilde \bg)=\infty$. 
\end{rem} 

\begin{defn} \label{hmorbit} If $\bg$ is an H-M (or near H-M) rep., refer to
its $\gamma_\infty$ orbit as an H-M (or near H-M) rep.~orbit. Call a
$\gamma_\infty$ orbit of an element in $(O)\sh$ with $O$ an H-M (or
near H-M) orbit the shift of an H-M (or near H-M) orbit. \end{defn}

\subsubsection{$\bar M_4$ braids H-M to near H-M reps.}
\label{w2cusps} 
Prop.~\ref{nearHMreps-count} braids a near H-M rep., $$\bg^*=(g_1,
c^{g_2^{-1}}g_1^{-1}c^{g_2^{-1}}, cg_2c,g_2^{ -1}),$$ to
the H-M rep.~$(g_1,g_1^{-1},g_2,g_2^{- 1})$. That
$\bg^*$ satisfies the product one condition is equivalent to  $c$
centralizes $g_1g_2$.
As in \S \ref{fratAnalog}, denote $\ker_k/\ker_{k+1}$ by $M_k$ and $M_0=M(A_5)$. 
The central
extension
$T_k'\to G_k$ appears above Cor.~\ref{R1G1}. Its characterization is that  
$\bZ/2=\ker(T_k'\to
G_k)$ and any lift to
$T_k'$ of the nontrivial element in $\ker(\Spin_5\to A_5)$ has order $2^{k+1}$. 
 
\newcommand{\da}{{{\scriptscriptstyle\dagger}}}  

\begin{prop}\label{nearHMbraid} Suppose
$\bg'=(g_1, \cc g_1^{-1}\cc, \cc g_2\cc,g_2^{ -1})$ with $\cc\in
M(A_5)$ (and $g_1^{-1}g_2$ of order 5). Then $\bg'$ satisfies
the product one condition if and only if  $\cc$ centralizes
$g_1^{-1}g_2$. If $\cc\ne 1$, it is $(g_1^{-1}g_2)^5$. Then, 
$(\bg')q_2^{10}=(\bg')\comp$
is the H-M rep. $(g_1, g_1^{-1},
g_2,g_2^{-1})$.  Also, $(\bg^*)q_1^{-1}q_2^{10}=(\bg^*)q_1^{-1}\comp$ is 
a different H-M rep.~lying over the same H-M rep. at level 0 as does 
$(\bg')\comp$ 

Then, if $\bg$ is an H-M rep.,   
$(\bg)\comp \gamma_\infty^{-1}\sh \gamma_\infty^{-1}$ represents the reduced 
class  of a 
near H-M rep. Each type (2,2) cusp has
Nielsen class  representatives consisting of $\sh$ applied to
elements lying over either $\bh_+$ or $\bh_-$. Conclude: All
near H-M reps.~and elements in width 2 cusps fall in one $M_4$
orbit. 

More generally, consider $g_1^\da,g_2^\da\in G_k$, and suppose the image of 
$(g_1^\da)^{-1}g_2^\da$ in
$A_5$ has order 5. So, $(g_1^\da)^{-1}g_2^\da$ has order $5\cdot 2^k$
(Lem.~\ref{cuspWidth}) and $((g_1^\da)^{-1}g_2^\da)^{5\cdot 2^{k\nm 1}}=c^\da$ 
has
order 2. Further, any lift of $c^\da$ in the cover $T_k'\to G_k$ has order $4$. 
There is a
braid from $(g_1^\da, (g_1^\da)^{-1}, g_2^\da,(g_2^\da)^{ -1})$ to a near H-M
rep.~$\bg''$. Complex conjugation $\hk$ (for two pairs of complex conjugate 
branch points)
applied to $\bg''$ gives $c''\bg''c''$ with
$c''$ lifting to
$T_k'$ to have order $4$.  
\end{prop}

\begin{proof} The statement on $\bg'$ comes to noting the
product-one condition is equivalent to $\cc^{g_1^{-1}}\cc
\cc\cc^{g_2^{-1}}=\cc^{g_1^{-1}}\cc^{g_2^{-1}}=1$. Apply $g_2$
to both sides to see $g_1^{-1}g_2$ fixes $\cc$. 
Prop.~\ref{nearHMreps-count} shows $\cc$ is the unique nontrivial
element of $M(A_5)$ in $\lrang{g_1^{-1}g_2}$. Lem.~\ref{qi2}
shows $q_2^{10}$ braids $\bg'$ to $(g_1, g_1^{-1}, g_2,g_2^{-1})$.

Apply $q_1^{-1}$ to $\bg^*$ to
get $$(c^{g_2^{-1}}g_1^{-1}c^{g_2^{-1}},
c^{g_2^{-1}g_1^{-1}}c^{g_2^{-1}}g_1c^{g_2^{-1}}c^{g_2^{-1}g_1^{-1}},
cg_2c,g_2^{-1}).$$ Use $c^{g_1g_2}=c$ to rewrite this
as $(g_1',c(g_1')^{-1}c,cg_2c,g_2^{-1})=\bg''$ with
$g_1'=c^{g_2^{-1}}g_1^{-1}c^{g_2^{-1}}$. Apply the first part of
this lemma with $\cc=c$: An element of $H_4$ takes $\bg''$
to the H-M representative $(g_1',(g_1')^{-1},g_2,g_2^{-1})$. 

As in Lem.~\ref{cuspWidth}, $c=(g_1g_2)^5\in
M(A_5)\setminus\{0\}$ centralizes $g_1g_2$. Use that $c$
and its conjugates don't centralize $g_1$ (Cor.~\ref{A5Morbs}) to see
$$(g_1',(g_1')^{-1},g_2,g_2^{-1})q_1=((g_1')^{-1},g_1',g_2,g_2^{-1})$$ and 
$(g_1, g_1^{-1},
g_2,g_2^{-1})$ give two H-M reps.~at level 1 in one $H_4$ orbit  lying over the 
same level 0 H-M
rep.   

We already know how to braid from
$\bh_+$ to $\bh_-$. Since, there are only two H-M reps.~over
$\bh_+$, it suffices that we have this braid from one to the other. So,
\S\ref{a25} gives the braiding between all the width two cusps. 
The last expression relating H-M reps. to complements of near H-M reps. comes 
from  writing $q_1$ as
$$q_2^{-1}(q_1q_2q_1)q_2^{-1}=q_2^{-1}(q_2q_1q_2)q_2^{-1}$$ to express this in 
standard generators of $\bar M_4$. 

The relation between H-M reps.~$\bg$ and near H-M reps.~ given by 
$(\bg)\comp \gamma_\infty^{-1}\sh \gamma_\infty^{-1}$ works at any level (with 
$p=2$).
Prop.~\ref{nearHMreps-count} has already established that near H-M reps.~have 
the complex
conjugation properties stated here. The element $c^\da$  
has the form
$((g_1^\da )^{-1}g_2^\da )^{5\cdot 2^{k\nm 1}}$. Lift it to $T_k'$ by 
lifting $g_1^\da$ and $g_2^\da$ to (respectively) $g_1$ and $g_2$ of order
3 in $T_k'$. Then form $(g_1)^{-1}g_2)^{5\cdot 2^{k\nm 1}}$. Let $h=((g_1)^{-
1}g_2)^{5}$.
The characterizing property of $T_k'$ implies $h$ has order $2^{k\np1}$ or 
$h^{2^{k\nm1}}$ has
order 4. So, any lift of $c^\da$ to $T_k'$
has order 4. Since $c''$ is just a conjugate of $c^\da$, it applies to $c''$ as 
well.    
\end{proof}

By now it is clear  H-M reps.~figure in many
computations. Denote the collection of these by $S_{\text{HM}}$.

\subsection{Cusp widths of length 6 and 10} 
The effect of $q_2^2$ on $\bg\in \ni_1^\inn$ is conjugation of 
$(g_2,g_3)$ by $g_2g_3$
(Lem.~\ref{qi2}) leaving $g_1$ and $g_4$ fixed. 

\subsubsection{$\gamma_\infty$ orbits of $\bg\in \ni_1^{\inn,\rd}$ with
length six}  A $\gamma_\infty$ orbit of length six implies a length
three orbit of
$\gamma^2_\infty$ on
$\bg$. The group $H=\lrang{g_1,g_4}$ is the pullback in
$G_1={}_2^1\tilde A_5$ of a copy of
$A_4$ in $A_5$ (Table \ref{lista5}, items 3,6,7 or
one of the conjugates of these by $(4\,5)$). So, $H$ has
no nontrivial centralizer in $G_1$.  Conclude:   
Three divides the order of $g_2g_3$, and with $v=(g_2g_3)^3\in M(A_5)$, 
$\bg_v=(g_1,v g_2 v, v g_3 v,g_4)=(\bg)q_2^6$ does not equal $\bg$ in 
$\ni_1^{\inn}$ (Prop.~\ref{serLift}).  
Note: $v\in M(A_5)\setminus V$ generates the centralizer
in $M(A_5)$ of $g_2g_3$.  

So, the $\gamma_\infty$ orbit of $\bg\in \ni_1^{\inn,\rd}$ 
has length six if and only if 
\begin{triv} $(\bg_v)Q'$ is conjugate to $\bg$ for some 
$Q'\in \sQ''\setminus\{0\}$. 
\end{triv}
\noindent 
Prop.~\ref{HMreps-count} shows no element of $\sQ''\setminus \{1\}$  fixes any 
inner classes. 

\begin{prop} \label{bgv} There are eight $\gamma_\infty$ orbits on 
$\ni_1^{\inn,\rd}$ of length six. Let $\bh\in \ni_0^{\inn,\rd}$ be in a
length three $\gamma_\infty$ orbit, $\bh\lrang{\gamma_\infty}$. Then there
are exactly six length twelve and four length six $\gamma_\infty$ orbits on
$\ni_1^{\inn,\rd}$ lying over $\bh\lrang{\gamma_\infty}$. They have these
further properties. For each H-M and near H-M orbit $O'$:
\begin{edesc} \label{interHM} \item four length twelve orbits $O$ over
$\bh\lrang{\gamma_\infty}$ satisfy
$|(O)\sh\cap O'|=1$; and  
\item the remaining length twelve or six over $\bh\lrang{\gamma_\infty}$
satisfy 
$|(O)\sh\cap O'|=0$. \end{edesc} Denote the collection of length 6 orbits and 
those of length 
12 not meeting H-M or near H-M reps., respectively by $L_6$ and $L_{12}$. 

Similarly, there are eight $\gamma_\infty$ orbits on $\ni_1^{\inn,\rd}$ of 
length 10 and four of
length 20 that are not H-M or near H-M orbits. Denote the former by $L_{10}$ and 
the latter by
$L_{20}$. For $O\in L_6$ and $O'\in L_{20}$ or for $O\in L_{12}$ and $O'\in 
L_{10}$, $|(O)\sh\cap
O'|=1$. This accounts for all entries in the $\sh$-incidence matrix for orbits 
with 3
dividing $\mpr$ paired with orbits with 5 dividing $\mpr$.  
\end{prop} 

Subsections \S\ref{lsix} and \S\ref{ltwelve} report on  respective
length six and twelve orbits to conclude the proof of Prop~\ref{bgv}. The 
calculation 
shows that various $\sh$-incidence intersections are 0 or 1 (at most 1). 
These incidence intersections divide according to how the orbits relate to H-M 
and near H-M
reps. This is a consequence of the lifting invariant from Prop.~\ref{expInv}. 
Calculations for
\eqref{interHM} are similar to those in the last half of the statement. So, we 
only give details
for  the former. 

\subsubsection{Length six orbits over $\bh\lrang{\gamma_\infty}$}
\label{lsix}
Use the previous notation for $\bg_v$ with $\bg$ lying over
$\bh\lrang{\gamma_\infty}$. Then $\bg\lrang{\gamma_\infty}$ has length six
exactly when $(\bg_v)Q'$ is inner equivalent to $\bg$ for some
$Q'\in\sQ''$. Modulo $\ker_0$,
$\bg=\bg_v$. Since nontrivial elements of
$\sQ''$ have order 2, an element of order 4 must give any conjugation of
$\bg$ to $(\bg_v)Q'$. 

The only possibility is that $\bg$ is conjugate
in $G_1$ to one of the following: 
$$\begin{array}{lr} (vg_3v,g_4,g_1,vg_2v),\ &
(g_1vg_2vg_1^{-1},g_1,g_4,g_4^{-1}vg_3vg_4),\ \text{or }\\
& (g_3vg_4vg_3^{-1},g_3,g_2,g_2^{-1}vg_1vg_2).\end{array}$$

Try these by cases. Suppose $\alpha$ conjugates $\bg$ to 
$ (g_1vg_2vg_1^{-1},g_1,g_4,g_4^{-1}vg_3vg_4)$. Then,
$g_2=\alpha^{-2}(g_2)\alpha^2=g_1g_2g_1^{-1} \bmod \ker_0$. 
There is a contradiction: When $\mpr(\bg)$ divides three, the
images of $g_1$ and $g_2$ in $A_5$ don't commute (check
Table \ref{lista5}). 

If $\alpha$ conjugates $\bg$
to  $(g_3vg_4vg_3^{-1},g_3,g_2,g_2^{-1}vg_1vg_2)$, then $\alpha^2$
commutes with $g_2$ and $g_3$. These generate the
pullback of $A_4$ in $G_1$ which has no center (Cor.~\ref{A5Morbs}). This
implies $\alpha^2=1$, contrary to  $\alpha$ having order 4. 

This leaves the serious case. Consider if $\bg$ is
conjugate to $$(vg_3v,g_4,g_1,vg_2v)=(\bg_v)Q',\ 
Q'=(Q_1Q_2Q_3)^2.$$ That is, for some $\alpha\in G_1$:
$$\alpha g_1\alpha^{-1}=vg_3v,\ \alpha g_2\alpha^{-1}=g_4,\
\alpha g_3\alpha^{-1}=g_1,\ \text{and } \alpha g_4\alpha^{-1}=vg_2v.$$  
Apply conjugation by $\alpha$ to these expressions with $v'=v^{\alpha}$. 
Conclude:  Conjugation by $\alpha^2$ on $\bg$ gives 
$(v'g_1v',vg_2v,vg_3v,v'g_4v')$. 

Opening arguments in \S\ref{a24} imply $\alpha^2=v'=v$. 
Yet, $v$ generates the centralizer of
$g_2g_3$. So, $H=\lrang{g_2g_3,\alpha}$ stabilizes $v$. There are two
cases. If the image of $H$ in $A_5$ is transitive on four
letters, then the image is isomorphic to $A_4$:  its pullback in
$G_1$ (as above) has no center. Otherwise, the image of
$H$ in $A_5$ is isomorphic to $S_3$. 

\subsubsection{Counting length six orbits} \S\ref{lsix} gives the only
allowable shape for elements $\bg\in \ni_1$ in a length six orbit of
$\gamma_\infty$. We prove there are in total
eight such orbits. To be explicit, assume $h_2h_3\equiv (1\,4\,2) \bmod\ \ker_0$
(with $\bh={}_3\bg$ in Table \ref{lista5}). If $\bg\in \bfC_{3^4}$ lies over
$\bh$, conjugating by $\alpha
\equiv (1\,2)(3\,5) \bmod\ \ker_0$ has the effect of switching $h_1$
with $h_3$ (resp.~$h_2$ with $h_4$). 

Apply Lem.~\ref{liftorder2} to find $\alpha$ having an additional property:
$\alpha^2$ generates the centralizer $\Cen_M(h_2h_3)$ of $h_2h_3$ in
$M(A_5)$. Fix some
$g_1$ in its conjugacy class, lifting $h_1$. Let $g_2$ be arbitrary in its
conjugacy class lifting $h_2$. Define $g_3$ to be $g_1^\alpha$ and $g_2$ to
be $g_4^\alpha$. The resulting $\bg$ gives a $\gamma_\infty$ orbit of
length six if 
$\Pi(\bg)=1$. 

If $\Pi(\bg)=v'$, then $s(\bh)=1$ (as in
Prop.~\ref{serLift}) implies
$v'\in V$. For $m\in M(A_5)$, $(\alpha m)^2=\alpha^2$ if and only if
$mm^\alpha=1$: $\alpha$ fixes $m$. Conjugate $g_4$ by
$u\in V(A_5)$ and replace $\alpha$ by $\alpha m=\alpha'$ with $m^\alpha=m$
to construct $\bg'$ from $\bg$. Compute $\Pi(\bg')$ as 
\begin{equation} \label{comp-utom} g_1m(ug_4ug_1)^\alpha mug_4u
= g_1g_4^\alpha
(u^{g_4g_1}u^{g_1})^\alpha
uu^{g_4^{-1}} mm^{(g_4g_1)^\alpha}g_1^\alpha g_4.
\end{equation}  

\begin{lem} \label{sixCent} For each $m\in \Cen_M(\alpha)$, there is a unique 
$u\in V$ giving product 1 in \eqref{comp-utom}. Representatives
$\bg$ for distinct length six orbits of
$\gamma_\infty$ over $\bh$ correspond to the four elements $m\in V$ that
centralize $\alpha$ (Cor.~\ref{A5Morbs}). \end{lem}

\begin{proof} The first sentence follows if 
$u\in V\mapsto (u^{g_4g_1}u^{g_1})^\alpha uu^{g_4^{-1}}$
is one-one. Consider
$A=1+(2\,3\,4)+(1\,2)(3\,5)((1\,2\,3)+(2\,4\,1))(1\,2)(3\,5)$ as an
element in the group ring $\bZ/2[A_5]$. It suffices to show
$A$ is invertible on $V$. 

Use notation from the proof of Lem.~\ref{liftorder2}:
$g_1,g_2$ and $\alpha_{25}$ on $D_5$ cosets act by
$\beta_{g_1}=(1\,3\,5)(2\,4\,6)$, $\beta_{g_2}=(1\,3\,2)(4\,5\,6)$ and
$\beta_{\alpha_{25}}=(3\,1)(2\,5)$. Since
$(1\,2)(3\,5)=(g_1g_2)^{2}(2\,5)(3\,4)(g_1g_2)^{-2}$, compute
$\beta_\alpha$ to be  
$(\beta_{g_1}\beta_{g_2})^2\beta_{\alpha_{25}}(\beta_{g_1}\beta_{g_2})^{-2}=
(1\,6)(2\,3)$. Then $A$, in the group ring, acts on
the six cosets generating $M(A_5)$ as follows:
$$\beta_A=1+(5\,4\,1)(6\,2\,3)+(1\,6)(2\,3)((1\,3\,5)(2\,4\,6)+
(1\,6\,4)(2\,5\,3))(1\,6)(2\,3)$$ or as
$1+(5\,4\,1)(6\,2\,3)+(1\,3\,4)(2\,5\,6)+(1\,4\,6)(2\,3\,5)$. 

Recall $\sum_{i=1}^6 i_H=0$ (in $M(A_5))$. Check $A$ on a basis
for $V$: $$1_H+2_H\mapsto 5_H+6_H,\ 2_H+3_H\mapsto 1_H+5_H,\ 
3_H+4_H\mapsto
3_H+5_H,\ 4_H+5_H\mapsto 2_H+5_H.$$ The range vectors form a
basis. So, $A$ is invertible on $V$. The remainder of the proof follows
from the setup to this lemma.  
\end{proof} 

\subsubsection{Length twelve orbits over $\bh\lrang{\gamma_\infty}$}
\label{ltwelve}  We apply the higher lifting invariant of Prop.~\ref{expInv}  
to compute the $\sh$-incidence entries for H-M and near H-M orbits
(Def.~\ref{hmorbit}) and
$\gamma_\infty$ orbits of elements with $\mpr$ six.  
 
\begin{lem} Suppose $O_\bg$ is the $\gamma_\infty$ orbit of $\bg \in
\ni_1^{\inn,\rd}$ where $\mpr(\bg)=6$ and
$(\bg)\sh$ is in the $\gamma_\infty$ orbit of an H-M or near H-M rep. Then,
for each $\gamma_\infty$ orbit $O$ of an H-M or near H-M rep.,
$|(O_\bg)\sh\cap O|=1$. Suppose $O'$ is the $\gamma_\infty$ orbit of 
$\bg'\in \ni_1^{\inn,\rd}$ with $\mpr(\bg')=10$ with $O'$ not the orbit
of either an H-M or near H-M rep. Then, $|(O_\bg)\sh\cap O'|=0$. \end{lem}

\begin{proof} Suppose $\bg$ is an H-M or near H-M rep., or a complement of
one of these. Then, $((\bg)\gamma_\infty^j)\sh$ has $\mpr$ six if only
if $j\equiv 2, 3 \bmod\ 5$. So, for values of $j \bmod\ 20$, there are
eight distinct elements $((\bg)\gamma_\infty^j)\sh$ with $\mpr$ six. We show
no two are on the same $\gamma_\infty$ orbit. 

First consider level 0 from Table \ref{lista5} and 
$(({}_1\bg)\gamma_\infty^2)\sh$ and
$(({}_1\bg)\gamma_\infty^3)\sh$. Conjugation by $(4\,5)$ takes one of the two
$\gamma_\infty$ orbits to the other. So, it
suffices to show 
$((\bg)\gamma_\infty^{2+m5})\sh$, $m\in \{0,1,2,3\}$ are in distinct
$\gamma_\infty$ orbits. That requires showing
$((\bg)\gamma_\infty^{2})\sh$ and $((\bg)\gamma_\infty^{2+10})\sh$ are in
distinct $\gamma_\infty$ orbits. This is saying the complements shift to
distinct $\gamma_\infty$ orbits. 

The computations are similar for H-M and near H-M reps., so we do the former
only. For simplicity write $\bg$ as $(g_1^{-1},g_1,g_2,g_2^{-1})$.  
With $u=g_1g_2$ and $m=u^5$, we show $(ug_1u^{-1}, ug_2u^{-1},
g_2^{-1}, g_1^{-1})=\bg'$ and  $(mug_1u^{-1}m, mug_2u^{-1}m,
g_2^{-1}, g_1^{-1})=\bg''$ are in distinct $\gamma_\infty$ orbits. By
reducing modulo $\bmod \ker_0$, the only possibility is 
$(\bg')\gamma_\infty^6$ is in the same reduced equivalence class as $\bg''$.
Check easily this isn't so. 

The level one lifting invariant of Prop.~\ref{expInv} has value $+1$ on 
$\sh$ applied to the $\gamma_\infty$ orbits of elements
$((\bg)\gamma_\infty^j)\sh$ with
$\mpr$ six described above. From Prop.~\ref{expInv}, such an element
$\bg'$ is in the $\bar M_4$ orbit of an H-M or near H-M rep.
Prop.~\ref{expInv} also says that orbits $O'$ in the statement of the lemma
have lifting invariant $-1$. Since the lifting invariant is an $\bar M_4$
invariant, they cannot intersect the set from $\sh$ applied to the 
$\gamma_\infty$ orbit of
$((\bg)\gamma_\infty^j)\sh$.  
 \end{proof}
 
\subsection{The $\sh$-incidence matrix at level 1} \label{shinc1} Table 
\ref{shincni0}
has the
$\sh$-incidence matrix for the reduced Nielsen classes
$\ni(A_5,\bfC_{3^4})^{\inn,\rd}=\ni_0$. We have all the data 
in place to compute the $\sh$-incidence matrix for $\ni_1$ (replacing
$G_0=A_5$ by $G_1$). As at level 0, a 2-stage process gives the block for orbit
$O_1^+$: List H-M and near H-M reps.~and compute their $\gamma_\infty$ orbits; 
then, apply
$\sh$ to these, and compute their $\gamma_\infty$ orbits. What results is closed 
under
$\sh$, so the process completes after two stages.  Use the notation for
$\bh_+$ and
$\bh_-$ in the above discussion. 

Label the type (10,20) $\gamma_\infty$ orbits as $O(10,20;j)$, $j=1,\dots,4$
for those containing H-M reps.~and $O(10,20;j)$, $j=5,\dots,8$ for those
containing near H-M reps. Refinement: $O(10,20;j)$ with $j$ odd
(resp.~even) lies over the $\gamma_\infty$ orbit at level 0 containing $\bh_+$
(resp.~$\bh_-$). Further refinement: $|(O(1,2;1)\sh\cap O(10,20;j)|=1$,
$j=1,2$;  $|(O(1,2;2)\sh\cap O(10,20;j)|=1$,
$j=3,4$. According to Prop.~\ref{sqOnM}, complements of near H-M reps. produce   
$|(O(2,2;1)\sh\cap O(10,20;j)|=1$,
$j=5,6$; and  $|(O(2,2;2)\sh\cap O(10,20;j)|=1$,
$j=7,8$. 

\begin{lem} \label{easy2-4} Two type (2,4)
$\gamma_\infty$ orbits, $O(2,4;i)$, $i=1,2$,  contain
$\sh$ applied to two near  H-M reps.~and two
complements of H-M reps. So, $\sh$ applied to $O(2,4;i)$ covers 
$\{\bh_+,\bh_-\}$ evenly, $i=1,2$. 

Let 
$\bg$ be an H-M rep. Let $\bg'$ be the near H-M rep.~attached to $\bg$ by
Prop.~\ref{nearHMreps-count}. Let
$\bg''$ be  the near H-M 
rep.~attached to the H-M rep.~$(((\bg)\sh) q_2)\sh$. Then,
$(\bg')\sh$ is different from
$((\bg'')\sh) q_2$. 
\end{lem}

\begin{proof} Let $\bg$ be an H-M rep. The relationship 
$(\bg)\comp \gamma_\infty^{-1}\sh \gamma_\infty^{-1}$ already says 
$\gamma_\infty^{-1}\sh \gamma_\infty^{-1}$ takes the complement of an H-M rep. 
to a 
near H-M rep. Now use $\sh \gamma_\infty^{-1}\sh \gamma_\infty^{-1}\sh 
=\gamma_\infty$ in $\bar M_4$. Let $\bg'$ be $(\bg)\comp$, and apply $\sh$ to
it, to get  
$$(\bg')\sh(\sh \gamma_\infty^{-1}\sh \gamma_\infty^{-1}\sh)$$ is the shift of a 
near H-M rep. This shows why a type (2,4) orbit containing $\sh$
applied to a complement of an H-M rep.~also contains $\sh$ applied to a near
H-M rep.  The remaining conclusions are similar. \end{proof} 

As in Prop.~\ref{countWTwoFour} use  
$\bh_+$ and $\bh_-$ for the two 
H-M reps.~in $\ni_0^{\inn}$. 

\begin{lem} \label{tough2-4} Beyond Lem.~\ref{easy2-4}, there are four 
type (2,4) $\gamma_\infty$ orbits
$O(2,4;j)$, $j=3,\dots,6$, with lifting invariant +1 (Prop.~\ref{expInv}).
For $j=3,\dots,6$, $$|(O(2,4;j))\sh\cap O(10,20;k)|=0 \text{ or }1.$$ For $\bg$ 
an H-M rep.,  $(\bg)q_2^5$ and $\bg$ lie over the same element of
$\ni_0^\rd$. Each of these $(O(2,4;j))\sh\,$ sets contains exactly one reduced
equivalence class of form $((\bg)q_2^5)\sh$.   

That leaves eight type (2,4)  $\gamma_\infty$ orbits with lifting
invariant -1 (from Prop.~\ref{expInv}). A representative for such an orbit
has $\bg=(g_2^{-1},g_1,ag_1^{-1}a,bg_2b)$. The shift applied to the 
$\gamma_\infty$ orbit of $\bg$ 
intersects exactly two type (10,20)
$\gamma_\infty$ orbits, with these covering $\{\bh_+,\bh_-\}$. This
accounts for all 16 of the elements with shifts of type (2,4) in the four 
(10,20) orbits with
lifting invariant -1. Shifts of the remaining two elements from each (2,4) orbit  
meet type (10,10) $\gamma_\infty$ orbits, covering $\{\bh_+,\bh_-\}$.  \end{lem}

\begin{proof} Let $\bg$ be an H-M rep. of the form
$(g_1^{-1},g_1,g_2,g_2^{-1})$.  Then, 
$(\bg)q_2^5$ and $\bg$ lie over the same element of $\ni_0^\rd$. 

In each width 20 $q_2$ orbit, there are four elements to which $\sh$ gives 
an element with
$\mpr(\bg)=1$ or 2. Apply the shift to the remaining (from 16) elements in the 
orbits 
containing H-M and near H-M reps. that lie over $\bg\in \ni_0$ with 
$\mpr((\bg)\sh)=3$. 
There are $8\times 8=64$ elements in the width 20 cusps with shifts  of
this  type. 

The remaining elements with $\mpr$ 2 must have lifting
invariant -1. A representative has the form
$\bg=(g_2^{-1},g_1,ag_1^{-1}a,bg_2b)$ with
$(a,b)\in M_3'\times M_5'\cup M_5'\times M_3'$ as in Prop.~\ref{expInv}. With 
$\lrang{c_1}$ the
centralizer of $g_1$ in $M$, the 
$\gamma_\infty$ orbit of $\bg$ includes
$$\begin{array}{rl} &(\bg)q_2=(g_2^{-1},a^{g_1^{-1}}g_1^{-1}a^{g_1^{-
1}},g_1,bg_2b),\\
&\bg_a=(g_2^{-1},c_1aa^{g_1^{-1}}g_1c_1aa^{g_1^{-1}},a^{g_1^{-1}}g_1^{-
1}a^{g_1^{-1}},bg_2b),
\text{and}  (\bg_a)q_2.\end{array}$$ An analog the Lem.~\ref{sixCent} proof 
shows that if 
$(\bg)\sh\in L_{10}$, then $(\bg_a)\sh\in
L_{20}$. As there, this argument uses the special presentation of $M(A_5)$ from
Cor.~\ref{A5Morbs}. \end{proof}

Let $O$ be one of the (6,12) type
$\gamma_\infty$ orbits with lifting invariant +1. From \S\ref{ltwelve},
$|(O)\sh\cap O(10,20;j)|=1$, $j=1,\dots,8$. Further, $(O)\sh$ hits two other
$(6,12)$ $\gamma_\infty$ orbits with multiplicity two. The argument is similar
to that of Lem.~\ref{meetO-O}: Use the $\sh\gamma_\infty\sh$ operator applied to 
$\bg$
with $\mpr$ six with $(\bg)\gamma_\infty$ also having $\mpr$ six. This gives the 
$8\times 12=96$ elements in (6,12) type orbits with lifting invariant
$+1$. 

The only unaccounted nonzero entries in the $\sh$-incidence matrix for orbit 
$O_1^+$ have the
form $|(O)\sh\cap O'|$ with $O$ and $O'$ running over H-M and near H-M orbits. 
From
symmetry, it suffices to fix $O$ an H-M rep.~orbit, and account for the 
intersections as
$O'$ varies. 

\begin{lem} \label{meetO-O} As $O'$ runs over H-M and near H-M orbits, 
$\sum_{O'}
|(O)\sh\cap O'|=8$. 
If $\bg$ is an H-M rep., then the H-M orbit $O_{\bg'}$ of $\bg'=(((\bg)\sh)
\gamma_\infty)$ meets it with multiplicity four. Two near H-M orbits $O'_1$ and 
$O'_2$ meet
$O$ with multiplicity two. \end{lem}

\begin{proof} Lemmas \ref{easy2-4} and \ref{tough2-4} account for all of the 
intersections of
$\sh$ applied to orbits with $\mpr=2$ with $O$. There are four such 
intersections.
\S\ref{ltwelve} shows that $\sh$ applied to  orbits with $\mpr=6$ contributes 
eight
intersections with $O$. Since there are 20 elements in $O$, 
$\gamma_\infty$ orbits $O_{\bg'}, O_1', O_2'$ account for all the remaining 
intersections
with $O$. Now we show why $O_{\bg'}, O_1', O_2'$ give the indicated 
intersections. 

Such intersections lie above corresponding level 0 intersections. 
\S\ref{shiftOp}
gives the formula  $\sh \gamma_\infty \sh=\gamma_\infty^{-1} 
\sh \gamma_\infty^{-1}$ acting on reduced inner Nielsen classes. Suppose $\bg\in 
O$ lies
over $\bh_+$ and
$(\bg)\sh$ has middle product 2. For example, this would hold if $\bg$ is an H-M 
rep. Then  
(Lem.~\ref{easy2-4}),
$\bg'=(((\bg)\sh)\gamma_\infty)\sh$ lies over $\bh_-$, on an H-M or near H-M 
rep.
orbit, and its shift also has middle product 2. So 
$((\bg')\gamma_\infty)\sh=(\bg)\gamma_\infty^{-1}$, $\sh$ applied to 
$(\bg')\gamma_\infty$
is  an intersection of $O_\bg$ and $(O_{\bg'})\sh$. 
If we use instead  
$\sh \gamma_\infty^{-1} \sh=\gamma_\infty \sh \gamma_\infty$, this  gives 
another
intersection of $O_\bg$ and $(O_{\bg'})\sh$:
$((\bg')\gamma_\infty^{-1}) \sh=(\bg)\gamma_\infty$. 

If $\bg$ is the  H-M rep.~on $O$, then $\bg'$ is an H-M rep. Suppose 
$$(\bg)\gamma_\infty^{10}=(g_1^{-1},dg_1d,dg_2d,g_2^{-1})$$ ($d=(g_1g_2)^5$; the 
complement of
an H-M rep.), then $\bg'$ is a near H-M rep.~as given by 
Prop.~\ref{nearHMreps-count}.
Similarly, replacing $\bg$ by $(\bg)\gamma_\infty^{15}$ and 
$(\bg)\gamma_\infty^{5}$
replaces  $\bg'$ by the complement of $\bg'$ and the other near H-M rep.~over 
$\bg_-$.  
\end{proof} 

Lem.~\ref{meetO-O} fully accounts for the $\sh$-incidence matrix restricted to
the block corresponding to the pairings of $\gamma_\infty$ orbits in $\bar M_4$ 
orbit $O_1^+$
having lifting invariant +1. The $32\times 32$ block of pairings of 
$\gamma_\infty$ orbits
in $\bar M_4$ orbit $O_1^-$ having lifting invariant -1 is more complicated. 
\S\ref{B4Act}
discusses why there may be a representation theory connection between the two 
orbits, though
it would not be numerically simple. To conclude this section we illustrate the 
lemmas
above with a contribution to the complete
$\sh$-incidence matrix at level 1 for  having $+1$ lifting invariant.

\newcommand{\lngorb}{{\scriptstyle\dagger}}
From here to the end of this subsection, $\gamma_\infty$ orbits refer to 
suborbits of
$O_1^+$. The part of the
$\sh$-incidence matrix pairing $O$ and $O'$ with respective middle products 6 
and 10 is
significant. Still, this $8\times 8$ block consists of just one's 
(Prop.~\ref{bgv}).
More interesting is the $8\times 8$ block of intersections from (10,20) type 
$\gamma_\infty$
orbits. For economy, replace $(10,20)$ as the $\gamma_\infty$ orbit type in 
Table
\ref{10-10shinc} by the symbol
$\dagger$. {\small\begin{table}[h] \caption{\rm $\mpr\,$s $10\times10$ part
of $\sh$-Incidence Matrix for $O_1^+$} \label{10-10shinc}
\begin{tabular}{|c|cccccccc|} \hline Orbit & $\!\!\!O(\lngorb;1)$\vrule  &
$\!\!\!\!\!O(\lngorb;2)$\vrule & 
$\!\!\!\!\!O(\lngorb;3)$\vrule &$\!\!\!\!\!O(\lngorb;4)$ \vrule&\!\! 
$\!\!\!\!\!O(\lngorb;5)$\vrule  &
$\!\!\!\!\!O(\lngorb;6)$\vrule  &
$\!\!\!\!\!O(\lngorb;7)$\vrule  & $\!\!\!\!\!O(\lngorb;8)$ \\ \hline \!\! 
$\!\!\!O(\lngorb;1)\!\!$
&0 &0& 0 &4 &0& 0 &2 &2\\ 
$\!\!\!O(\lngorb;2)\!\!$ &0& 0 &4& 0 &2 &2& 0 &0 \\  $\!\!\!O(\lngorb;3)\!\!$  
&0& 4& 0 &0 &0&
0 &2 &2\\
$\!\!\!O(\lngorb;4)\!\!$  &4& 0 &0 &0 &2 &2& 0 &0\\ $\!\!\!O(\lngorb;5)\!\!$ &0& 
2& 0& 2 &0
&0 &0 &4\\  
$\!\!\!O(\lngorb;6)\!\!$ &0 &2 &0& 2 &0 &0 &4 &0\\ 
$\!\!\!O(\lngorb;7)\!\!$ &2 &0& 2 &0 &0 &4 &0 &0\\ 
$\!\!\!O(\lngorb;8)\!\!$ &2 &0 &2 &0 &4& 0 &0 &0\\
\hline \end{tabular} \end{table}}

\subsection{Real points on level 1 of the $(A_5,\bfC_{3^4})$ Modular Tower} 
\label{level1} We show the genus 12 component of
$\sH(G_1,\bfC_{3^4})^\inn$ has one component of real points, while the genus 9 
component has no
real points. 

We check the effect of $\hk$ at level 1 extending the use of the complex
conjugation operator at level 0 as Lem.~\ref{A5absR} reports. Note: $\hk^2$ acts
trivially on
$\bg$ satisfying the product-one condition. For $Q\in \sQ''$, while $Q^2$ acts 
trivially
on a Nielsen class, it usually is a nontrivial conjugation dependent on the 
$\bg$ representing a
Nielsen class. 

First, check the effect over $j_0\in
(1,\infty)$. Corresponding to $j_0$, we know there are eight cover points 
corresponding to the 
four H-M reps. and the four near H-M reps.~lying on $\sH_1^{\inn,\rd}$ 
(Prop.~\ref{HMreps-count})
and the complex conjugation operator $\hk_0$. The Nielsen classes for the
corresponding components all lie in distinct
$\bar M_4$ orbits. For each, moreover, their complementary Nielsen class gives 
another
point over
$j_0$ from Lem.~\ref{square1}. These 16 real components over the interval 
$(1,\infty)$ all lie on
the component $\sH_1^+$ (as in Prop.~\ref{finTwoOrbits}). Denote the complete 
set of real points
(corresponding to reduced Nielsen classes) lying on $\sH_1^{\inn,\rd}$ over 
$j_0$ by $\ni_{j_0}$.  

\begin{prop} \label{oneOrbitkappa} For $j_0\in (1,\infty)$, elements of 
$\ni_{j_0}$ lie on exactly 
16 components of  $\sH_1^+(\bR)$ lying over $(1,\infty)\subset \prP^1_j(\bR)$. 
The closure in
$\bar \sH_1^+(\bR)$ of these components consists of one connected component of 
real points. There
are no $\bR$ points on $\sH_1^{\inn,\rd}$ over $(-\infty,1)\subset 
\prP^1_j(\bR)$, and no
$\bR$ points at all on $\sH^-_1$. 
\end{prop}

The proof of Prop.~\ref{oneOrbitkappa} takes up the next five subsections. Its 
arrangement starts
with an equivalent equation defining a point on $\sH_1^{\inn,\rd}(\bR)$ for a 
reduced equivalence
class of covers represented by a particular cover with branch points all in 
$\bR$ or in complex
conjugate pairs. Assume a representing cover has branch cycles $\bg\in 
\ni(G_1,\bfC_{3^4})$.
Then,  $\hk (\bg) = h(\bg)Q h^{-1}$ with $\hk$ the complex conjugation operator
corresponding to the branch point configuration, $h\in G_1$ and $Q\in \sQ''$. 
Also,
$\mod \ker_0$, $Q$ acts trivially modulo conjugation. 

\subsubsection{$\hk_4$ for four real branch points} \label{kappa4} The level 0 
real point
data for the  absolute Hurwitz space appears at the end of the proof of 
Lem.~\ref{A5absR}. 

For four points $\bz$ in $\bR$: 
$\bg'_1=((1\,2\,3),(1\,3\,2),(1\,4\,5),(1\,5\,4))$ has complex
operator $c_1=(2\,3)(4\,5)$;  $\bg'_2=((1\,2\,3),(1\,4\,5),(1\,5\,4),(1\,3\,2))$ 
has
$c_2=(1\,3)(5\,4)$. The involution $c_{1,\infty}=(3\,6)(8\,9)(4\,5)$ acts on 
Nielsen classes in
Table \ref{lista5}, revealing $\bg_1'$ and $\bg_2'$ as the first two items in 
that table. Since the
complex conjugation operator for $\bg_7'$ is $(2\,3)$, the Galois closure of 
that cover is not over
$\bR$. So, points corresponding to these branch cycle descriptions on the 
absolute space have no
real points on the inner space above them.   

\S\ref{innerHSReal} has the level 0 inner Hurwitz space data summarized in the
complex conjugation operator from \eqref{c1infty}:
$c_{1,\infty}=(9\,17)(7\,16)(8\,18)(13\,14)(12\,15)(6\,3)(4\,5)$. The fixed 
points correspond to
real points, with branch cycle descriptions 
$\bg'_1$ and $\bg_2'$ and their conjugates by $(4\,5)$. 

Let $\hk_4$ be the complex conjugation
operator for four real branch points. It suffices to handle the defining 
equation, when 
$\hk_4 (\bg) = h(\bg)Q h^{-1}$. With no loss, take $\bg \mod \ker_0$ equal 
either $\bg'_1$
or $\bg'_2$. 

\subsubsection{$\hk_0$ for complex conjugate pairs of branch points} 
First consider the case $\bg_1'$, which is an H-M rep.~at level 0. Here and 
below we use the
normalization that writes $\bg=(g_1,ag_1^{-1}a,bg_2b, g_2^{-1})$, a {\sl 
perturbation\/} of an H-M
rep.  With no loss $(g_1,g_2)$ lies over $((1\,2\,3),(1\,4\,5))$ or 
$((1\,2\,3),(1\,5\,4))$
with $a,b\in M\setminus V$.  It is convenient to use past data by replacing 
$\hk_4$ by
$\hk_0$, the complex conjugation operator for complex conjugate pairs of branch 
points. As 
in Rem.~\ref{comp-real}, this doesn't change the real points representing 
reduced
classes. 

Assume $\bg$ covers $((1\,2\,3),(1\,2\,3)^{-1},(1\,4\,5),(1\,4\,5)^{-1})$. 
Double the number of
components corresponding to using  $(1\,5\,4)$ in place of 
$(1\,4\,5)$. Use the notation around Prop.~\ref{expInv} to distinguish the $\bar 
M_4$ orbit of $\bg$
by whether
$a$ and $b$ are in the same $G_1$ orbit (acting on $M\setminus V$). 

Suppose $\bg$ corresponds to $(a,b)$ in different $G_1$ orbits. We run through 
nontrivial elements $Q\in \sQ''$, showing that $\hk_0(\bg)=h(\bg)Qh^{-1}$ is not 
possible
for any $h\in G_1$. The product of 3rd and 4th entries of $\bg$ is  
$v=bb^{g_2^{-1}}$. So,
$\hk_0$ maps $\bg$ to $(vag_1av, vg_1^{-1}v, g_2^{-1}, bg_2b)$. The components
$\sH_1^+$ and  $\sH_1^-$ have field of definition $\bQ$ and so $\hk_0(\bg)$ 
gives branch
cycles for a cover in $\sH_1^-$. Prop.~\ref{finTwoOrbits} implies it is attached
to $(a',b)$ as a perturbation of an H-M rep. One case is easy: 
$Q=(q_1q_2q_3)^2$. On the
right side these produce $\bg^*$ conjugate to
$((bg_2b, g_2^{-1}, g_1,ag_1^{-1}a)$. Though $\bg^*$ is still a perturbation of 
an H-M rep., it
corresponds to $(a^*,b^*)$ with $a^*$ (resp.~$b^*$) in the opposite $G_1$ orbit
to $a$ (resp.~$b$). Similarly for $Q=(q_1q_2q_3)^2q_1q_3^{-1})$. 

Yet, when $Q=q_1q_3^{-1}$, 
$(\bg)Q=(ag_1^{-1}a,g_1,g_2,bg_2^{-1}b)$. A conjugation of $(\bg)Q$ that gives 
$\hk_0(\bg)$
is by a lift $\alpha$ of $(2\,3)(4\,5)$ with square the centralizer $c_2$ of 
$g_2$ in
$M\setminus V$, and
$(\alpha v)^2$ the centralizer $c_1$ of $g_1$ in $ M\setminus V$. The last is 
equivalent to
$\alpha^2 v^\alpha v=c_1$, or $c_1=v^\alpha vc_2$. 
Cor.~\ref{A5Morbs} notes that $\alpha$ fixes exactly four elements in $V$ and 
two elements of
$M\setminus V$ in the conjugacy class $M_3'$.  They are $c_1$ and $c_2$. So, 
$v^\alpha v=c_1c_2$ is
possible for exactly four elements in $V$. (The proof of Cor.~\ref{A5Morbs} has 
explicit
computations, for a different $\alpha$, corroborating this.) The expression 
$w^{g_2^{-1}}w=v$
then has a unique solution $w\in V$ from which we derive values of $b=c_2w$. We 
already know four
values of $b$ that satisfy this. These correspond to the complements of H-M and 
near H-M reps.
lying above the real component of $\sH_0$ (over $(1,\infty)\subset \prP^1_j$) 
corresponding to the
cover points with branch cycles $((1\,2\,3),(1\,3\,2),(1\,4\,5),(1\,5\,4))$ and 
complex conjugation
operator $\hk_0$. Thus, Prop.~\ref{HMreps-count} shows for 
$(a,b)$ in different $G_1$ orbits on $M\setminus V$,
$v^\alpha v=c_1c_2$ is not possible. 

\subsubsection{Action of $\sQ''$ on $\bg$ over $\bg'_2$} \label{pertHMs} Adopt 
the idea of using a
perturbation of an H-M rep.~here by writing $\bg=(ag_1a,bg_2^{-1}b, g_2,g_1^{-
1})$. The product-one
condition gives
$a^{g_1}a=bb^{g_2}$. We show $\hk_4(\bg)=h(\bg)Qh^{-1}$ is impossible for $Q\in 
\sQ''$. All
computations are similar. We do just $Q=q_1q_3^{-1}$.  This will complete that 
there are exactly 16 components of $\sH^{\inn,\rd}(\bR)$ over the $j$-interval 
$(1,\infty)$. 
Compute $\hk_0(\bg)$: 
\begin{equation}   
\begin{array} {rl} &(ag_1^{-1}a, g_1g_2^{-1}bg_2bg_2g_1^{-1}, g_1g_2^{-1}g_1^{-
1},
g_1)\\
=&h(a(ab)^{g_1^{-1}}g_1g_2^{-1}g_1^{-1}a(ab)^{g_1^{-1}},ag_1a,g_1^{-1},
g_1g_2g_1^{-1})h^{-1}.\end{array}\end{equation} The effect of $h$ on the 3rd and 
4th entries shows
$h^2$ is the identity. Since, however, $h$ is a lift of $(1\,4)(3\,5)$, it has 
order 4
(Lem.~\ref{FrKMTIG}). 

\subsubsection{Connecting 16 components of $\sH_1^+(\bR)$ over 
$(1,\infty)\subset \prP^1_j$}
Now we show these 16 components of real points close up to one connected 
component on $\bar\sH_1^+$.
That is, if we add the endpoints over $\infty$ (cusps \S \ref{cuspNotat}) and 1 
to these components,
they form one connected set. Branch cycles for an H-M rep. and its complement 
give real components
on $\sH_1^+$ meeting (each pair) at a cusp of type $(10,20)$ 
(Prop.~\ref{HMreps-count}). The same is
true for near H-M reps. 

Apply $\sh$ to these to get another real component (over $(1,\infty)$). This is 
the result of
reflection off the
$j=1$ boundary of the entering real locus. Apply Prop.~\ref{sqOnM} to
decipher the labeling of branch cycles to components on applying $\sh$ to the H-M, near H-M
reps.~and their complements. The result is a cusp of type (2,2) or (2,4). In the 
former case it 
reflects off the $j=\infty$ boundary to give another real component. In the 
latter case
an application of $\gamma_\infty^2$ gives a real component reflection. Applying 
Prop.~\ref{sqOnM}
to these reflections shows there is one connected component of real points.  

\subsubsection{$\sH_1^{\inn,\rd}(\bR)$ is empty over $j\in (-\infty,1)$} To 
finish
Prop.~\ref{oneOrbitkappa} we show there are no points on $\sH_1^{\inn,\rd}(\bR)$ 
over
$j\in (-\infty,1)$.  According to Lem.~\ref{locj}, this is a check that $\hk_2 
(\bg) =
h(\bg)Q h^{-1}$ is impossible under the following conditions.
\begin{edesc} \item $\hk_2$ is the complex conjugation operator for
two real and a complex conjugate pair of branch points. 
\item $h\in G_1$ and $Q\in \sQ''$. 
\end{edesc} 

As in
\S\ref{pertHMs}, consider branch cycles for real points on $\sH^{\inn, \rd}_0$. 
Write $\hk_2(\bg)=(g_1^{-1}, (g_3g_4)^{-1}g_2^{-1}g_3g_4, g_4^{-1},g_3^{-1})$.
Running through entries of Table \ref{lista5} gives the existence of $h$ in
$\hk(\bg')=h\bg' h^{-1}$ with $\bg'\in \ni_0$ only for ${}_{1}\bg$ with 
$h=(2\,3)$,
${}_{3}\bg$ with $h=(1\,3)(4\,5)$ and
${}_{9}\bg$ with $h=(2\,3)$. As in \S\ref{kappa4}, only the real component 
corresponding to
${}_{3}\bg$ lies below real components of $\sH^{\inn,\rd}_0$. This is compatible 
with
$c_{-\infty,0}=(10\,13)(1\,4)(9\,14)(7\,3)(12\,16)(11\,2)(15\,6)(5\,18)$ from 
\S\ref{innerHSReal}
fixing two integers. It suffices to consider $\bg\in \ni_1^{\inn,\rd}$ lying 
over
${}_{3}\bg=((1\,2\,3), (1\,4\,5),(2\,1\,5),(2\,4\,3))$. Write $\bg$ as
$(g_1,g_2,\beta g_2^{-1}\beta^{-1},\beta g_1^{-1}\beta^{-1})$ with 
$\beta=g_1g_2$ or a
perturbation of this of the form $(g_1,g_2,a\beta g_2^{-1}\beta^{-1}a,b\beta
g_1^{-1}\beta^{-1}b)$ with $a,b\in M\setminus V$. The calculation to show that for 
no $Q$ and $h$
can this give an $\bR$ point on $\sH^{\inn,\rd}_1$ is similar to that of 
\S\ref{pertHMs}. We
leave it to the reader.

\section{Completing level 1 of an $A_5$ Modular Tower}
\label{compA5MT} Cor.~\ref{orbitGenus} computes genuses of the two components 
$\sH_1^+$ and
$\sH_1^-$ of level of the $(A_5,\bfC_{3^4})$ Modular Tower. 
Prop.~\ref{oneOrbitkappa} shows 
only $\sH_1^+$ has real points (one component of them on $\bar \sH_1^+$). 
\S\ref{gaplist}  describes data we got from  \GAP.  Following that, all but the 
last subsection 
is group theory for the {\sl spin structure\/}
explanation that nails the two components in Cor.~\ref{orbitGenus}. The
\GAP\ data was reassuring at several stages. Still, it was crude compared to our 
final arguments.
Despite the detail developed for our special case, there remains the mysterious 
similarity, amidst the
differences, of the two level 1 orbits. \S\ref{finConc} gives ideas for that 
appropriate
for testing the speculative Lie algebra guess in \S\ref{GKS}. The paper 
(excluding
appendices)   concludes with lessons on  computing (genuses of) components for 
all Modular Towers
(\S\ref{indSteps}). 

Orbits of
$\gamma_\infty$ correspond to cusps of a $j$-line cover. Orbits of 
$\lrang{\gamma_1,\gamma_\infty}$
correspond to components of the $j$-line cover. Finding {\sl discrete\/} 
invariants
of a Modular Tower means analyzing the $\lrang{\gamma_1,\gamma_\infty}$ and 
$\gamma_\infty$
orbits on reduced Nielsen classes. That includes analyzing arithmetic
properties of the cusps. Since components are moduli spaces this means analyzing
degeneration of objects in the moduli space on approach (over $\bR$ or $\bQ_p$) 
to the
cusps. This generalizes analysis of elliptic curve degeneration in approaching a 
cusp of a
modular curve (\cite{FrLInv} introduces using $\theta$ functions to make this 
analysis; \S\ref{hcan}). Even for the one Modular Tower for $(A_5,\bfC_{3^4})$, 
it would be a
major event to prove the analog of Prop.~\ref{oneOrbitkappa} for $p$-adic points
(see \S\ref{openIm}). 

\subsection{Diophantine implications} \label{diophimp}
Apply Falting's Theorem \cite{FaMordCon} to these
genus 12 and 9 components to conclude the following from Thm.~\ref{thm-rbound}. 

\begin{thm} \label{twoOrbits} There are only finitely many
$({}_2^1\tilde A_5,\bfC_{3^4})$ realizations over any number field $K$.  For $k$
large there are no $(G_k,\bfC_{3^4})$ realizations over $K$.  \end{thm}

\begin{exmp}[The $(A_4,C_{3^4},p=2)$ Modular Tower] \label{A4C34-wstory} If 
$G_k$ is the level
k characteristic quotient for ${}_2\tilde A_5$, then the pullback of $A_4\le 
A_5$ in $G_k$ is
$G_k'$, the level $k$ characteristic quotient of ${}_2\tilde A_4$. The notation 
$(A_4,C_{3^4})$ is
ambiguous. There are two conjugacy classes of 3-cycles in $A_4$. Compatible with
Ex.~\ref{A5C54-wstory}, we  mean each pair of conjugacy classes appears twice, 
or what we
might label as $(A_4,C_{3^2_+3^2_-})$. From the \BCL\ this is equivalent to the 
levels of the
Modular Tower have definition field $\bQ$.  

Level 0 is a special case of
\cite{FrLInv} (Ex.~\ref{full3cycleList}). There are two components, $\sH^+_0$
and
$\sH^-_0$, corresponding to two
$H_4$ orbits
$O^+_0$ and $O^-_0$ of the lifting invariant from $A_4$ to $\Spin_4$ (analogous 
to
Prop.~\ref{expInv}). Note: $\sQ''$ acts through $\bZ/2$ on these inner classes 
(\S\ref{extA4C34} has
details). For absolute equivalence these spaces are families of genus 1 curves.  
It is a
special case of
\cite{FrKK} that the map of {\sl each\/} component to the moduli space of curves 
of genus 1 is
generically surjective (see \S\ref{hcan}). As in Prop.~\ref{expInv}, there is 
nothing above the component
$\sH^-_0$ at level 1. 

Applying this paper as with the $A_5$ Modular Tower, the first
author's thesis will show {\sl level 1\/} for $(A_4,C_{3^4})$, $p=2$ has six 
components. From 
\S\ref{liftA5}, some  differences between this case and $(A_5,C_{3^4})$ are 
clear. 
For $A_4$, 
$\{gv\}_{v\in V}$ is not the complete  set of conjugates of $g\in G_1'$, though 
one can still
choose conjugation by
$v\in V$ or $m\in M\setminus V$. Also, there are level 0 cases where 
$\mpr(\bg)=2$. 

The same lifting invariant 
from Prop.~\ref{expInv} applies. At level 1 there are components that together
comprise a moduli space $\sH^+_1$, and another $\sH_1^-$. Unlike, however, for 
$A_5$, each
of these has several (three) connected (reduced, inner) components. 

For $\sH^+_1$ the two H-M components have genus 1, and the 3rd, which contains 
near H-M reps., has
genus 3. So, all these components have real points. Two distinct components of 
the inner classes at
level 1 contain H-M reps. Further, the two orbits differ by an outer 
automorphism of
${}_2^1\tilde A_4$. We don't know if the H-M rep.~components are conjugate over 
$\bQ$. They come
from an outer automorphism (as in Thm.~\ref{FrVMS}). So, if they are
conjugate, it must be from an extension of $\bQ$ in the Galois closures of the 
absolute covers in
either one of these H-M families. 

For $\sH^-_1$, there are no real points. It
has two genus 0 components, with corresponding $\bar M_4$ orbits differing by an 
outer
automorphism of ${}_2^1\tilde A_4$ (different from that binding the two H-M 
components). The
genus 0 components are conjugate over $\bQ$ by the complex conjugation operator. 
The
last component of $\sH^-_1$ has genus 3. 

Since the genus 0 components are on $\sH^-_1$, they are obstructed: there is 
nothing
above them at level 2. The Main Conjecture for this case (Prob.~\ref{MPMT}) 
therefore 
follows by assuring any component at level 2 above the genus 1 H-M components 
has
genus at least 2.  This, however, follows easily from Lem.~\ref{RHOrbitLem}, the 
situation where one of
the orbits (at level $k=1$) has genus 1. Use that at level 1, 2 divides all 
$\mpr$ values,
except for H-M reps. Then, use that an H-M rep. orbit contains elements with 3 
dividing their
$\mpr$, so it contains non-H-M reps.~above which are elements with larger 
$\mpr$.  
\end{exmp} 

\begin{exmp}[Level 1 for $(A_5,\bfC_{5_+^{2}5_-^{2}}, p=2)$]
\label{A5C54-wstory} Prop.~\ref{3355obst} discusses level 0 for this Modular 
Tower, while  
the level 1 components and genuses  are exactly as for 
Ex.~\ref{A4C34-wstory}. Here automorphisms join the two H-M and two genus 0
orbits at level 1, generating the full group of ${}_2^1\tilde A_5$ outer 
automorphisms. 
There are two puzzles. 
\begin{edesc} \item Why are Exs.~\ref{A5C54-wstory} and \ref{A4C34-wstory} so 
alike? 
\item If the two H-M components in either case have definition field $\bQ$? 
\end{edesc} \end{exmp}

\subsection{Working with \GAP} \label{gaplist} We had \GAP\ list each
inner Nielsen class, computing the action of the  $Q_i\,$s as elements of
$S_{2304}$. Let
$O_1^+$ and
$O_1^-$ be the braid orbits in
$\{ 1,
\dots, 2304 \}$.  They each have cardinality 1152.  For $j =
1, 2$, mod out on $O_j$ by the action of $\sQ''$. In each $O_j$
there are 288 orbits of $\sQ''$,  each of size 4. \GAP\ automatically computes
the action of $\bar M_4$ on $O_j$. This  gives the branch
cycles of both degree 288 covers of $\prP^1_j$. \GAP\ doesn't
seem to have facility with {\sl split extensions\/}, which meant
we interpreted $G_1$ as embedded in $S_{n}$ for some integer $n$.
We started with a representation of degree   192, and eventually
used a degree 80 representation $T_{80}$ allowing \GAP\ faster
computation.  In $A_4$ let $\alpha$ have order 3, and let $K$
be a Klein 4-group in $\ker_0/\ker_1=M(A_5)$ on which $\alpha$
acts irreducibly. The group $K$, $\alpha$ and the centralizer of
$\alpha$ generate has order 24. Its cosets give $T_{80}$.

\subsubsection{Branch cycle description on $O_1^+$} By inspection
below: \begin{edesc} \item $\gamma_0$ has ninety-six 3-cycles;
$\gamma_1$, one hundred and forty-four 2-cycles;   and \item
$\gamma_\infty$, four 2-cycles, six 4-cycles, eight 12-cycles
and eight 20-cycles.  \end{edesc} From Riemann-Hurwitz,
the genus $g(O_1^+)$ of a cover with these branch cycles satisfies
$2(288+g(O_1^+)-1)=96\cdot2+144+4+6\cdot 3+8\cdot 11+ 8\cdot 19=598$:  
$g(O_1^+)=12$.

$$\!\!\!\!\scriptsize{\begin{array}{rl}&\gamma_0 = (1\, 53\, 42)(  2\,229\, 72)(
3\, 86\,190)(  4\, 55\,200)(  5\, 94\,205)(6\,125\, 81)
(  7\,127\,109)(  8\,132\, 91)(  9\,191\, 48)\\&( 10\,210\,
46)(11\,208\, 52)( 12\, 41\, 45)( 13\, 54\,199)( 14\,126\,
73)( 15\, 95\,197)(16\, 90\,274)( 17\,213\, 71)\\ &( 18\,214\,108)( 19\, 
58\,203)(
20\,225\,107)(21\, 65\,262)( 22\, 68\,263)(
23\,193\, 49)( 24\,275\, 50)( 25\,202\, 51)\\ &(26\,206\,
47)( 27\, 63\,278)( 28\, 64\,282)( 29\, 69\,239)( 30\,
70\,246) (31\, 57\,238)( 32\, 62\,245)( 33\,258\,144)\\ &(
34\,280\,234)( 35\,257\,221)(36\,255\,223)( 37\,270\,148)(
38\,271\,233)( 39\,284\,146)(40\,260\,141)(43\,149\,
93)\\ &( 44\, 85\,171)( 56\,189\,152)( 59\,204\,159)(60\,194\,153)
(61\,273\,162)( 66\,150\,220)( 67\,131\,288)( 74\,219\,249)\\
& (75\,128\,237)(76\,230\,287)( 77\,129\,251)( 78\,155\,134)( 79\,163\,232)(
80\,227\,265)(82\,217\,243)( 83\,160\,139)\\ &( 84\,130\,266)(
87\,151\,192)( 88\,169\,207)(89\,173\,209)( 92\,166\,215)(
96\,226\,267)( 97\,174\,212)( 98\,170\,211)\\ 
&(99\,158\,276)(100\,172\,196)(101\,195\,165)(102\,198\,154)(103\,201\,161)(104\
,277
\,175)
(105\,281\,177)\\
&  (106\,157\,133)(110\,164\,228)(111\,156\,140)
(112\,136\,252)(113\,218\,242)(114\,231\,286)\\
&(115\,135\,244)(116\,248\,176)(117\,264\,168)(118\,261\,180)
(119\,241\,178)(120\,247\,179)\\ &(121\,240\,167)(122\,138\,268)
(123\,137\,285)(124\,216\,250)(142\,185\,259)(143\,186\,269)\\
&(145\,187\,283)(147\,184\,254)(181\,253\,222)(182\,272\,236)(183\,256\,224)
(188\,279\,235)\\ &\end{array}}$$

$$\!\!\!\scriptsize{\begin{array}{rl}&\gamma_1 = (1\, 72)(  2\,234)(  3\,109)(
4\, 71)(  5\,107)(  6\,148)(  7\,146)(  8\,141)(  9\,190)(
10\,199)( 11\,197)( 12\, 42)\\ &( 13\, 81)( 14\,144)( 15\, 91)(
16\,108)( 17\,221)( 18\,223)( 19\, 73)( 20\,233)( 21\,
45)( 22\, 49)( 23\,203)( 24\,274)\\ &( 25\,205)( 26\,200)(
27\, 52)( 28\, 51)( 29\, 47)( 30\, 46)( 31\, 48)( 32\,
50)( 33\,239)( 34\,278)( 35\,238)( 36\,246)\\ &( 37\,262)(
38\,263)( 39\,282)( 40\,245)( 41\, 53)( 43\,171)( 44\,
93)( 54\,206)( 55\,210)( 56\,191)( 57\,237)( 58\,193)\\ &(
59\,202)( 60\,208)( 61\,275)( 62\,249)( 63\,287)( 64\,288)(
65\,266)( 66\,152)( 67\,139)( 68\,265)( 69\,243)(
70\,251)\\ &( 74\,213)( 75\,134)( 76\,229)( 77\,126)(
78\,153)( 79\,159)( 80\,232)( 82\,220)( 83\,162)( 84\,125)(
85\,178)( 86\,192)\\ &( 87\,177)( 88\,167)( 89\,179)( 90\,276)(
92\,161)( 94\,207)( 95\,209)( 96\,228)( 97\,168)( 98\,180)(
99\,175)(100\,176)\\ 
&(101\,196)(102\,211)(103\,212)(104\,285)(105\,286)(106\,154)
(110\,165)(111\,149)(112\,140)(113\,214)\\
&(114\,225)(115\,132
)(116\,244)(117\,268)(118\,267)(119\,250)(120\,242)(121\,252)(122\,127)(123\,133
)\\& 
(124\,215)(128\,259)(129\,254)(130\,269)(131\,283)(135\,258)(136\,260)(137\,284
)(138\,270)(142\,156)\\&(143\,160)(145\,157)(147\,155)(150\,224)(151\,189)(158\,
273
)(163\,235)(164\,236)(166\,222)(169\,204)\\&(170\,201)(172\,195)(173\,194)(174\,
198
)(181\,241)(182\,264)(183\,247)(184\,248)(185\,240)(186\,261)\\&(187\,281)(188\,
277
)(216\,255)(217\,253)(218\,257)(219\,256)(226\,271)(227\,272)(230\,279)(231\,280
).
\end{array}}$$

$$\!\!\!\!\scriptsize{\begin{array}{rl}&\gamma_\infty =\! ( 12\, 53)( 23\,
58)( 93\,171)(101\,172)(9\,86\,151\, 56)( 10\, 54\, 26\,
55)( 11\,95\,173\, 60)\\ &( 24\, 90\,158\,61) ( 25\,
94\,169\, 59)(102\,170\,103\,174)(  6\,270\,122\,
7\,284\,123\,157\,283\, 67\,160\,269\, 84)
\\&(8\,260\,112\,156\,259\, 75\,155\,254\, 77\, 14\,258\,115)(
17\,257\,113\, 18\,255\, 124\,166\,253\, 82\,150\,256\, 74)
\\&( 33\, 69\,217\,181\,119\,216\,
36\, 70\,129\,184\,116\,135)(34\, 63\,230\,188\,104\,137\, 39\,
64\,131\,187\,105\,231) \\& ( 35\, 57\,128\,185\,121\,136\, 40\, 
62\,219\,183\,120\,218)
(37\, 65\,130\,186\,118\,226\, 38\, 68\,227\,182\,117\,138)
\\& ( 2\,280\,114\, 20\,271\,96\,164\,272\,80\,163\,279\, 76)\\&( 1\,229\,287\, 
27\,208\,153\,134\,237\,
31\,191\,152\,220\,243\, 29\,206\, 13\,
 125\,266\, 21\, 41)\\ &(
 3\,127\,268\,168\,212\,161\,215\,250\,178\,
44\,149\,140\,252\,167\,207\,  5\,
 225\,286\,177\,192)\\
&(  4\,213\,249\, 32\,275\,162\,139\,288\, 28\,202\,159\,232\,265\,
22\,193\, 19\,
 126\,251\, 30\,210)\\
&( 15\,132\,244\,176\,196\,165\,228\,
 267\,180\,211\,154\,133\,285\,175\,276\, 16\,214\,242\,179\,209)\\
&( 42\, 45\,262\,148\,81\,199\,
 46\,246\,223\,108\,274\, 50\,245\,141\, 91\,197\, 52\,278\,234\,
 72)\\
&( 43\, 85\,241\,222\,  92\,201\, 98\,261\,143\, 83\,273\,
99\,277\,235\, 79\,204\, 88\,240\,142\,111) \\ &( 47\,239\,144\,
73\,203\, 49\,263\,233\,107\,205\,51\,282\,146\,109\,190\,
48\,238\,221\,  71\,200)\\ &( 66\,189\, 87\,281\,145\,106\,198\,
97\,264\,236\,110\,195\,100\,248\,147\, 78\,
 194\, 89\,247\,224).\end{array}}$$

\subsubsection{Branch Cycle Description on $O_1^-$} By
inspection as above: \begin{edesc} \item $\gamma_0$
has ninety-six 3-cycles; $\gamma_1$, one hundred and
forty-four 2-cycles;   and \item $\gamma_\infty$, eight
4-cycles, eight 6-cycles, eight 10-cycles, four 12-cycles and
four 20-cycles.  \end{edesc} Riemann-Hurwitz gives $g(O_1^-)$:
$2(288+g(O_1^-)-1)=96\cdot2+144+8\cdot 3+8\cdot 5+8\cdot 9+ 4\cdot
11+4\cdot 19=592$ or $g(O_1^-)=9$.

$$\!\!\!\!\!\!\!\!\!\!\scriptsize{\begin{array}{rl}&\gamma_0\!\! =\! (1\,129\, 
63)(2\,218\,172)(79\,237\,169)( 
5\,151\,142)(6\,131\,150)( 7\,211\,166)(89\,216\,165)( 18\,178\,192)(  9\, 33\, 
143)\\ &( 10\,277\,174)
(122\,155\, 254)( 12\,152\,187)( 13\, 40\,190)( 14\,157\,182)( 15\,175\, 66)( 
16\,132\,146)( 17\,177\,148) (3\,75\, 57)\\&
( 19\,196\,242)( 20\, 67\,286)( 21\,271\,227)( 22\,
64\,257)( 23\,160\,241)( 24\,164\,273)( 25\,252\, 88)
( 26\,158\,239)( 27\,248\, 86)\\ &( 28\,269\, 87)( 29\,287\,224)(
30\,247\,281)( 31\,206\,135)( 32\, 91\,139)( 34\,115\,180)(
35\, 94\,184)( 36\,105\,133) ( 37\,116\,179)\\ &( 38\,130\,101)(
39\,236\,121)( 41\,183\,113)( 42\,138\,108)( 43\,188\, 93)(
44\,255\, 96)( 45\,272\,117)( 46\,235\,207)( 47\,283\,
95)\\ &( 48\,280\,246)( 49\,100\,221)( 50\, 80\,245)(
51\, 92\, 72)( 52\,240\,124)( 53\, 81\,244)( 54\,71\,259)(
55\,214\,285)( 58\,112\,278)\\ &( 59\,109\,217)( 60\,213\,275)(
62\,176\,111)( 68\,107\,222)(4\, 61\,134)( 69\, 74\,284)( 70\,219\,260)( 
76\,171\,106)( 78\,253\,198)
\\ &( 82\,200\,205)(
83\,210\,250)( 84\, 98\,288)( 85\,125\,270)(102\,156\, 136)(
90\,154\,137)( 97\,268\,228) ( 99\,229\,161)(  8\, 77\, 65)\\&(103\,141\,145)
(104\,185\,149)(110\,153\,181)
(114\,140\,191)(118\,163\,262)(119\,276\,223)(120\,195\, 261)(11\, 73\, 
56)\\&(123\,251\,234)
(126\,265\,225)(127\,266\,226)(128\,249\,282)\!(144\,202\,193)
(147\,203\,186)(159\,256\,208)\!(162\,215\,258)
\\ &(167\,279\,243)(168\,231\,263)(170\,220\,264)(173\,212\,274)(189\,204\,194)
(197\,201\,232)\!(199\,230\,238)
\!(209\,267\,233)\end{array}}$$

$$\!\!\!\!\!\!\!\!\!\scriptsize{\begin{array}{rl}&\gamma_1 = (  1\,142)(  
2\,224)(  3\,
86)(  4\, 65)(  5\,172)(  6\,143)(  7\,227)(  8\, 88)(
9\, 57)( 10\,281)( 11\, 87)( 12\,166)( 13\, 56)( 14\,174)\\ &(
15\,190)( 16\,134)( 17\,187)( 18\,182)( 19\,192)( 20\,
66)( 21\,273)( 22\, 63)( 23\,146)( 24\,150)( 25\,239)(
26\,148)( 27\,242)\\ &( 28\,257)( 29\,286)( 30\,241)( 31\,207)(
32\, 95)( 33\,139)( 34\,124)( 35\, 96)( 36\,121)( 37\,117)( 38\,135)( 39\,245)( 
40\,179)\\ 
&
( 41\,184)(
42\,133)( 43\,180)( 44\,259)( 45\,275)( 46\,244)( 47\,285)(
48\,278)( 49\,113)( 50\, 73)( 51\, 93)( 52\,246)( 53\,
72)\\ &( 54\, 75)( 55\,221)( 58\,108)( 59\,101)( 60\,217)(
61\,129)( 62\,175)( 64\,260)( 67\,284)( 68\,111)( 69\, 77)( 70\,222)( 71\,269)\\ 
&(
74\,288)( 76\,169)( 78\,270)( 79\,250)( 80\,252)( 81\,248)( 82\,198)( 83\,205)( 
84\, 92)(
85\,106)( 89\,228)\!( 90\,165)\!( 91\,137)\\ &( 94\,188)( 97\,255)( 98\,283)(
99\,234)(100\,226)(102\,171)(103\,138)\!(104\,181)\!(105\,136)
(107\,225)(109\,223)\\&(110\,161)(112\,282)(114\,130)(115\,183)(116\,176)(118\,1
49
)(119\,272)(120\,191)(122\,145)(123\,236)(125\,262)\\ 
&(126\,254)(127\,261)(128\,240
)(131\,154)(132\,151)(140\,193)(141\,156)(144\,206)(147\,208)(152\,185)(153\,177
)\\ 
&(155\,263)(157\,186)(158\,237)(159\,258)(160\,243)(162\,216)(163\,264)(164\,274
)(167\,277)(168\,229)(170\,218)\\&(173\,211)(178\,194)(189\,203)(195\,253)(196\,
238
)(197\,202)(199\,232)(200\,204)(201\,233)(209\,256)(210\,235)\\ 
&(212\,276)(213\,271
)(214\,287)(215\,266)(219\,268)(220\,265)(230\,251)(231\,267)(247\,280)(249\,279 
).
\end{array}}$$

$$\!\!\!\!\!\!\!\scriptsize{\begin{array}{rl}&\gamma_\infty = (1\,151\, 16\, 
61)(  6\,
33\, 91\,154)( 15\, 40\,116\, 62)( 42\,105\,156\,103)(
17\,152\,104\,153) ( 18\,157\,203\,194)\\ &(
41\, 94\, 43\,115)  ( 38\,206\,193\,114)( 7\,271\,
60\,109\,276\,173)( 10\,247\, 48\,112\,249\,167)(
21\,164\,212\,119\, 45\,213)\\ & ( 99\,251\,199\,201\,267\,168)  
(106\,270\,198\,205\,250\,169)( 28\, 64\,219\, 97\, 44\,
71) ( 29\, 67\, 74\, 98\, 47\,214)\\ & (30\,160\,279\,128\,
52\,280)   (4\, 77\,284\, 20\,175\,111\,222\,260\, 22\,129)(9\,75\,259\,
96\,184\,113\,221\,285\, 95\,139)\\ &
( 36\,236\,234\,161\,181\,149\,262\, 85\,171\,136)
( 56\,190\, 66\,286\,224\,172\,142\, 63\,257\, 87)\\&(
32\,283\, 84\, 51\,188\, 35\,255\,228\,165\,137) (
19\,178\,204\, 82\,253\,120\,140\,202\,232\,238)\\& (
26\,177\,110\,229\,263\,122\,141\,102\, 76\,237) (25\,158\,
79\,210\, 46\, 81\, 27\,196\,230\,123\, 39\, 80) \\ &(31\,235\, 83\,200\,189\, 
147\,256\,233\,197\,144)
(2\,287\,55\,100\,266\,162\, 89\,268\,70\, 107\,265\,170)\\ &(3\,248\,
53\, 92\,288\, 69\,  8\,252\, 50\, 11\,269\, 54)(
78\,125\,163\,220\,126\,155\,231\,209\,159\,215\,127\,195)
\\
&( 5\,218\,264\,118\,185\,12\,211\,274\,  24\,131\,
90\,216\,258\,208\,186\, 14\,277\,243\, 23\,132) \\ &( 13\,
73\,245\,121\,133\,108\,278\,246\,124\,180\, 93\,72\,
 244\,207\,135\,101\,217\,275\,117\,179)\\
&( 34\,240\,282\, 58\,138\,145\,254\,225\, 68\,176\, 37\,272\,223\,
59\,130\,191\,261\,226\,
  49\,183)\\
&( 57\,143\,150\,273\,227\,166\,187\,148\,239\, 88\,
65\,134\,146\,241\,281\,174\,182\,
 192\,242\, 86).\end{array}}$$

From Lem.~\ref{HMreps-count},
there are 16 H-M representatives in $\ni({}_2^1\tilde
A_5,\bfC_{3^4})^\inn=\ni_1$.  \S\ref{sQorbits} shows all
$\sQ''$ orbits on $\ni_1$ have length four.

\begin{res} \GAP\ identifies the H-M representatives in $O_1^+$ as
corresponding to the integers $S_{\text{HM}}=\{41,42,43,44\}$
in the permutations $\gamma_0$, $\gamma_1$, and
$\gamma_\infty$. Further, $41$ and $44$ lie over
$((1\,2\,3),(1\,3\,2),(1\,5\,4),(1\,4\,5))\in
\ni(A_5,\bfC_{3^4})$ while $42$ and $43$ lie over
$((1\,2\,3),(1\,3\,2),(1\,4\,5),(1\,5\,4))$. Also, it also shows
that $S_{\text{nHM}}$, the near H-M representatives satisfying
\eqref{nearHM}, lie in this orbit. \end{res}

\subsection{Schur multipliers detect two $H_4$ orbits on
$\ni_1=\ni(G_1, \bfC_{3^4})$} \label{schMultInv}  Let $O$ be an
$H_4$ orbit on a Nielsen class $\ni$. 
\cite[Part III]{FrMT} uses an invariant $\nu(O)$, the {\sl big\/} invariant of
$O$, a union of conjugacy classes in $\ker_0=\ker(\tG p\to G)$. We adapt
this for an $H_4$ orbit on $\ni_1$. For $k\ge 1$ define $$\nu_k(O) =
\{\Pi(\tilde \bg)\mid \tilde \bg \in G_{k+1}^4\cap \bfC_{3^4} \text{ and } 
\tilde \bg
\bmod 
\ker_{1} \in O\}.$$
We use a quotient invariant of $\nu_1(O)$. Prop.~\ref{RkGk}, with $p=2$,
$k=1$ and
$G_1={}_2^1\tilde A_5$ produces
$R=T_1'$ with $\ker(R\to G_1)=\bZ/2$, and any lift of $m\in M\setminus V$ to $R$ 
has order 2
(Cor.~\ref{R1G1}).  

\begin{prop} \label{compSerre} There are two $H_4$ orbits on
$\ni({}_2^1\tilde A_5, \bfC_{3^4})$. Denote these by $O_i$; use
$\sH_{O_i}$ for their corresponding Hurwitz space components,
$i=1,2$. Then, $1\in \nu_1(O_1^+)$, but $1\ne \nu_2(O_1^-)$. Therefore,
$\sH_{O_1^-}$ is an obstructed component, having nothing above it
in the Modular Tower.  Suppose $\bp\in \sH_{O_1^+}(K)$ is a $K$
point for any number field $K$. Then, this corresponds to a $K$
cover $\hat\phi_\bp: \hat X_\bp\to \prP^1_z$ with the following
properties.  \begin{edesc}\label{liftProp} \item \label{liftPropa}
$\hat\phi_\bp$ is in the Nielsen class $\ni({}_2^1\tilde
A_5,\bfC_{3^4})$.  \item \label{liftPropb} $\hat\phi_\bp$
extends to a $K$ sequence of covers $\hat Y_\bp\mapright{\hat
\psi_\bp} \hat X_\bp\mapright{\hat \phi_\bp:} \prP^1_z$ with
group $R$.  \item \label{liftPropc} $\hat \psi_\bp$
is unramified of degree 2.  \end{edesc} Each $\bp\in \sH_{O_1^-}(K)$
satisfies \eql{liftProp}{liftPropa}, but not the combination
of  \eql{liftProp}{liftPropb} and \eql{liftProp}{liftPropc}.
\end{prop}

Prop.~\ref{expInv} gives the essential ingredient of the proof: The two orbits
separate by the values of invariants called $s(O_1^+)$ and $s(O_1^-)$, quotients 
of
$\nu_1(O_1^+)$ and $\nu_1(O_1^-)$. 

\subsubsection{Lifting invariants}  
Prop.~\ref{RkGk}  generalizes  \cite[\S4]{FrKMTIG}. Assume
$G=G_0$ is $p$-perfect and centerless with $G_k$ the $k$th 
characteristic quotient of $\tG p$. Let $\psi_k: R_k\to G_k$ be the universal 
exponent
$p$ central extension of $G_k$; this exists from Def.~\ref{pperfect}. Write 
$R_k$ as $\tG
p/\ker_k^*$, and the closure of $(\ker_k,\ker_k^*)(\ker_k^*)^p$ in
$\ker_{k+1}$ as $\ker_{k+1}'$.  The canonical
morphism $\phi_{k+1,k}: G_{k+1}\to G_k$ factors through $\psi_k$. 

\begin{prop} \label{RkGk} Denote $G_{k+1}/ker_{k+1}'$ by $R_{k+1}'$.
Then,  
$G_{k+1}$ acts trivially on $\ker_{k+1}/\ker_{k+1}'$. For $k\ge 0$, the $\bZ/p$ 
module map of
$\ker_k/\ker^*_k$ to $\ker_{k+1}/\ker_{k+1}'$ by the $p$th power is injective. 
So, a lift of $m\in
\ker(\phi_{k+1,k})$ to 
$\ker(R_{k+1}'
\to  G_{k})$ has order $p^2$ if and only if the image of $m$ in 
$\ker(R_{k}\to G_k)$ is nontrivial. This characterizes $R_{k+1}'$ as a quotient
of $R_{k+1}$. 
\end{prop}

\begin{proof} Since $\psi_k$ is a Frattini cover with kernel of exponent $p$,
$\phi_{k+1,k}$ factors through $\psi_k$. There are two types of generators of
$\ker_{k+1}/\ker_{k+1}'$:
$(v,v')$ and $v^p$ with $v,v'\in \ker_k$. Let $g\in G_{k+2}$. By assumption, 
there are
$h,h'\in \ker_k^*$ with $g^{-1}vg=vh$ and
$g^{-1}v'g=v'h'$. From $vh=hv \bmod \ker_{k+1}'$, 
$$g^{-1}v^pg=(vh)^p \equiv v^ph^p \bmod \ker_{k+1}' \equiv v^p.$$ Similarly,
$g^{-1}(v,v')g=(vh,v'h') \equiv (v,v')(h,h') \bmod \ker_{k+1}' \equiv (v,v')$. 
Conclude:
$G_{k+1}$ acts trivially on $\ker_{k+1}/\ker_{k+1}'$. 

Finally, suppose $m\in \ker_k$. Since $\ker_k$ is a pro-free, pro-$p$ group, 
$m^p \bmod
\ker_{k+1}'$ is nontrivial if and only if $m$ is not in  $\ker_k^*$. This 
translates the
last sentence. 
\end{proof}

\S\ref{cliffAlgLift} combines the next corollary with the
computation of Prop.~\ref{liftEven}. 
Prop.~\ref{RkGk}, with $p=2$, $k=0$ and $G_1={}_2^1\tilde A_5$ produces
$T_0'\eqdef R_0'$. In this case $T_0'\to G_0=A_5$ is the spin cover of $A_5$. 
Use
Prop.~\ref{liftEven} to inductively define a quotient $T_k'$ of $R_k'$, 
$k=1,2,\dots$ so that
$\ker(T_k'\to G_k)=\bZ/2$; all lifts of the nontrivial element of $\ker(T_0'\to 
G_0)$ 
to $T_k'$ have order $2^k$.  Let $T_1'=R$, so $\ker(R\to G_1)=\bZ/2$. As 
previously, 
$M(A_5)=M=\ker(G_1\to G_0)$ and $\hat M$ is the pullback of $M$ to $R$. 

\newcommand{\nze}{{\dot{-}}}
\begin{cor} \label{R1G1} Let $R_1\to G_1$ be the representation cover of $G_1$ 
and 
$U_1=\ker(R_1\to G_1)$.  Then,
$U_1$ is a natural $G_1$ (acting trivially) quotient module of $\ker_1/\ker_2$. 
With
$s_1(O)$  the image of $\nu_1(O)$ in $U_1$, $s_1(O)$ is a single element.  

Let $\hat a_i$ be a lift of $a_i\in M\setminus V$ to $\hat
M$ where
$a_i^{g_i}\ne a_i$, $i=1,2$ and $a_1\ne a_2$. Then, each of 
$H_{a_i}=\lrang{\hat a_i^{g_i^{-1}},
\hat a_i}$, $i=1,2$, and $H_{a_1,a_2}=\lrang{\hat  a_1,\hat a_2}$ is isomorphic 
to
$\bZ/4\times 
\bZ/2$.  
\end{cor}

\begin{proof} Since $R_1$ is a Frattini cover of $G_1$ it is a quotient of 
${}_2\tilde
A_5$. This defines $U_1$ as the maximal quotient of $\ker_1/\ker_2$ on which 
$G_1$ 
acts
trivially. So, there is Frattini cover $\psi^*: G^*\to G_1$, a central extension 
of
$G_1$ with kernel
$U_1$. Each $2'$ element of $G_1$ has a unique $2'$ lift to $G^*$. 

Suppose $\tilde \bg \in
G_{2}^4\cap \bfC_{3^4}$ and $\tilde \bg =\bg
\bmod \ker_{1} \in \ni_1$. Then, the unique lifts of the entries of $\bg$ to 
$(G^*)^4\cap
\bfC_{3^4}$ determine the image of $\Pi(\tilde \bg)$ in $U_1$.
Further, this element only depends on the $H_4$ orbit of $\bg$. 

The
two given generators in the groups $H_{a_i}$, $i=1,2$, or 
$H_{a_1,a_2}$ have product of  order
2. The resulting groups have order 8, the Klein 4-group as a quotient and 
a pair of generators with one of order 4 and one of order 2. Only the group 
$\bZ/4\times
\bZ/2$ has these properties: elements  of order 4 don't generate the dihedral 
group, 
and the quaternion group's only element of order 2 is in its center.
\end{proof}

\subsubsection{Identifying lifting invariants from $G_1$ to $R$}
The next proposition gives the key invariant separating the two orbits $O_1^+$ 
and $O_1^-$.
Together with the $\sh$-incidence matrix of \S\ref{shinc1} it explains precisely 
the 
two $\bar M_4$ orbits on the reduced Nielsen class
$\ni(G_1,\bfC_{3^4})^{\inn,\rd}$. 

\begin{prop} \label{expInv} Let $( g_1, g_1^{-1}, g_2, g_2^{-1})$ be an H-M 
representative. 
Consider $$( g_1, a g_1^{-1} a, b g_2 b, g_2^{-1})=\bg^*.$$ With no loss,  
assume $a,b\in M\setminus V$ (Cor.~\ref{A5Morbs}).  Let $\hat
\bg^*$ have as entries the unique lifts of $\bg^*$ entries 
to $(R,\bfC_{3^4})$. Then, $\Pi(\hat
\bg^*)=s(\bg^*)$ equals 
$\hat a^{g_1^{-1}}\hat a\hat b\hat b^{g_2^{-1}}$. 
With no loss, take
$\bg^*$ as a representative of a given $\bar M_4$ orbit $O$.
Then, if $O=O_1^+$:
\begin{triv} \label{conjSep} $s(\bg^*)=\hat
a^{g_1^{-1}}\hat a\hat b\hat b^{g_2^{-1}}=+1$; and  $(a,b)\in
M_3'\times M_3'$ or $M_5'\times M_5'$. 
\end{triv} The corresponding statement for $\bg^*\in O_1^-$
is 
\begin{triv} \label{conjSep2} $\hat a^{g_1^{-1}}\hat a\hat
b\hat b^{g_2^{-1}}=-1$, with $(a,b)\in M_3'\times M_5'$ or $M_5'\times M_3'$. 
\end{triv}
\noindent A $\bar M_4$ orbit satisfying either \eqref{conjSep} or 
\eqref{conjSep2} contains 8
reduced classes of such $g^*\,$s. 
\end{prop}

\begin{proof} If $g\in G_1$ has order 3, and $a\in M$,  both $\hat a\hat
g\hat a$ and
$\widehat{aga}$ are order 3 lifts to $R$ of $aga$. So, they are equal.
The first computation is almost trivial:
$$\hat
\bg^*=(\hat  g_1,\hat
g_1^{-1}\hat a^{g_1^{-1}}\hat a,\hat b\hat b^{g_2^{-1}}\hat g_2,\hat g_2^{-
1}).$$ The 
lifting
invariant is a braid invariant. The formulas of \eqref{conjSep} and 
\eqref{conjSep2} come to computing 
$\hat a^{g_1^{-1}}
\hat a
\hat b^{g_2^{-1}}
\hat b$  when without the hats it is 1. 

Each element in the pullback $\hat V$ of $V$ in $R$ has 
order 2. Therefore,
$\hat V$ is an elementary 2-group and a $G_1$ module. Suppose $g\in G_1$ has
order 3. For
$v\in V$,
$\hat  v^{g}\hat v$ is the trivial element if and only if
$v=0$. So, as
$v\in V$ varies, 
$\hat v^{g}\hat v$ (or $\hat v^{g^{-1}}\hat v$) runs over an $\bF_2[G_1]$ module
$V_{g}=V_{g^{-1}}$ of $R$ mapping to $V$ by the canonical map $\phi_1:R\to
G_1$. Further,
$m\mapsto
\hat m^{g}\hat m$ has the same image as $v$, where $m$ and
$v$ differ by the element generating the centralizer of $g$. In
particular, $M$ also maps to 
$V_{g}$
under this map. 

Now suppose $V_{g_1}\ne V_{g_2}$. Then, there are $m_1,m_2\in M$ with
$\hat
m_1^{g_1^{-1}}\hat m_1\ne \hat
m_2^{g_2^{-1}}\hat m_2$ having the same image in $V$.
Conclude:  There 
exists
$\bg\in
\ni_1$ with
$s(\bg)=-1$ if and only if
$V_{g_1}\ne V_{g_2}$. 

If, however,  $V_{g_1}= V_{g_2}=\tilde V$, 
then $\tilde V$ is a proper $\bF_2[G_1]$ module 
of $\hat V$ mapping surjectively to $V$.  So, as a $\bF_2[G_1]$ module, $\hat V$
is the direct sum of $\tilde V$ and $\one$.  Consider $\hat
M=\phi_1^{-1}(M)\le R$.  Cor.~\ref{R1G1} shows this is a $\bZ/4[G_1]$ module.
Then, $\hat M/\tilde V$ is a $\bZ/4[G_1]$ module of order 4. Thus, it is
either $\bZ/4$ or $\bF_2\oplus\bF_2=C_1\oplus C_2$. 

Suppose  $\hat M/\tilde V=\bZ/4$, and $G_1$ acts 
nontrivially here. This would force $G_1$ to have a nontrivial $\bZ/2$ quotient, 
a
contradiction to $G_1$ being 2-perfect. Now assume $\hat M/\tilde V=\bZ/4$ and 
$G_1$ acts 
trivially.  Then, $\hat M/\tilde V\to R/\tilde V\to A_5$ is a central Frattini 
cover with kernel
$\bZ/4$. Exclude this because the Schur multiplier of
$A_5$ is $\bZ/2$. 

So $\hat M/\tilde V=C_1\oplus C_2$ holds. We may choose $C_1$ to map trivially  and $C_2$
to map  nontrivially by
$\phi_1$ to
$M/ V$. Thus, as an abelian group we can can write $\hat M$ as $C_1\oplus \lrang{C_2,\tilde
V}$. Since $ \lrang{C_2,\tilde V}/C_1$ maps surjectively to $M$ and has the same order, It is
isomorphic to
$M$.  Contrary to how we formed $R$, each element in $\hat M$ would have
order 2. This concludes showing $V_{g_1}\ne V_{g_2}$:  Both
$\bar M_4$ orbits are nonempty. 

We 
show there are eight values of $(a,b)$ with $\hat a^{g_1^{-1}}
\hat a \hat b^{g_2^{-1}} \hat b=+1$ (resp.~-1). Consider pairs $(\hat
v_1,\hat v_2)\in V_{g_1}\times V_{g_2}$ with $\hat v_1\hat v_2$ over
the trivial element in $V$. These form an $\bF_2$ module $\hat W$ consisting of
the two cosets of the submodule $\hat U$ where $\hat v_1\hat v_2=+1$. So, 
each of the two invariants contribute eight elements. 

Now we show the distinguishing conditions of
\eqref{conjSep} and \eqref{conjSep2} hold. 
\begin{edesc} \item If $\bg^*=(g_1, a g_1^{-1} a, b g_2 b,
g_2^{-1})$ is in an H-M rep.~orbit, then $$(a,b)\in
M_3'\times M_3'\cup M_5'\times M_5', \text{ so }s(\bg^*)= +1.$$
\item Otherwise $(a,b) \in M_3'\times M_5'\cup
M_5'\times M_3'$ and $s(\bg^*)=-1$.
\end{edesc} Toward the former, consider $(a,b)$ when
$\bg^*$ is an H-M, a near H-M rep.~or a complement of such. 
Let
$c_i\in M\setminus V$ centralize 
$g_i$, $i=1,2$. Then,
$(a,b)=(c_1,c_2)\in M_3'\times M_3'$. Let $d\in M\setminus \{0\}$ 
centralize $g_1^{-1}g_2$. For the complement $\bg'$ of
$\bg^*$, $(a',b')=(d,d)\in M_5'\times M_5'$.

Next, suppose $\bg^*= ( g_1, c^{g_2^{-1}}
g_1^{-1}c^{g_2^{-1}},cg_2c, g_2^{-1})$ where $c\in M\setminus \{0\}$ 
centralizes $g_1g_2$. As in Prop.~\ref{nearHMreps-count}, this is a near
H-M rep. Here $(a,b)=(c^{g_2^{-1}},c)\in M_5'\times M_5'$. Let
$(a',b')=(dc^{g_2^{-1}}c_1, dcc_2)$ be the $(a,b)$ for the
complement $\bg'$ of $\bg^*$. To show $(a',b')\in M_3'\times
M_3'$, use  the Lem.~\ref{liftorder2} proof. 
Assume $g_1$ is a lift of $(1\,2\,3)$ and $g_2$ is a lift of
$(1\,4\,5)$. As usual write $M$ as
$(\bZ/2)^6/\lrang{(1,1,1,1,1,1)}$: $ g_1$ (resp.~$g_2$)
permutes the coordinates as
$\beta_{g_1} =(1\,3\,5)(2\,4\,6)$
(resp.~$\beta_{g_2}=(1\,3\,2)(4\,5\,6)$). 

Since $\beta_{g_1}
\beta_{g_2}=(1\,2\,5\,3\,6)$ (resp.~$\beta_{g_1^{-1}}
\beta_{g_2^\sph}=(1\,6\,5\,2\,4)$), $c=(0,0,0,1,0,0)$
and
$d= (0,0,1,0,0,0)$. Also, $c_1=(1,0,1,0,1,0)$ and
$c_2=(1,1,1,0,0,0)$. Then,
$(dc^{g_2^{-1}}c_1, dcc_2) = ((1,0,0,0,1,1), (1,1,0,1,0,0))\in
M_3'\times M_3'$ and $(a, a^{ g_1^{-1}},
b, b^{g_2^{-1}})$ is $$((1,0,0,0,1,1), (0,0,1,1,1,0),
(1,1,0,1,0,0), (0,1,1,0,0,1)).$$

Let $U^*\le V$ be the subgroup mapping isomorphically to $\hat U$ by the action
of $(u_1,u_2)\in U^*\mapsto (\hat u_1\hat u_1^{g_1}, \hat u_2\hat u_2^{g_2})$.
To simplify exponents replace $u_i$ by $u_i^{g_i}$.
Similarly, form $W^*$. The above shows that for all $(\hat u_1,\hat u_2)\in
\hat U$, $c_1u_1$ and $c_2u_2$ are  in the same conjugacy class of 
$M\setminus V$. To finish it suffices that the remaining eight pairs,
with $(w_1,w_2)\in W^*\setminus U^*$, $(c_1w_1,c_2w_2)$ are in different
conjugacy classes in $M\setminus V$. 

Do an explicit calculation using $\beta_{g_1}$ and $\beta_{g_2}$,  
$c_1$ and $c_2$. Elements in $V$: 
$$\scriptsize{\begin{array}{lllll} 1_v=(1,1,0,0,0,0),& 2_v=(1,0,1,0,0,0),&
3_v=(1,0,0,1,0,0),& 4_v=(1,0,0,0,1,0),&  \\
5_v=(1,0,0,0,0,1),& 6_v=(0,1,1,0,0,0),& 7_v=(0,1,0,1,0,0),& 8_v=(0,1,0,0,1,0),&
\\
9_v=(0,1,0,0,0,1),& 10_v=(0,0,1,1,0,0),& 11_v=(0,0,1,0,1,0),&
12_v=(0,0,1,0,0,1),&\\ 13_v=(0,0,0,1,1,0),& 14_v=(0,0,0,1,0,1),&
15_v=(0,0,0,0,1,1).\end{array}}$$  

Apply $1+\beta_{g_1}$ (resp.~$1+\beta_{g_2}$) to this list, to get two new
lists $L_i$, $i=1,2$. For example, the 6th item in $L_1$ (resp.~the 10th
item in List $L_2$) is
$6_v+(6_v)\beta_{g_1}= (1,0,0,0,0,1)=5_v$ (resp.~$10_v+(10_v)\beta_{g_2}=
(1,0,0,0,1)=5_v$). 

Treat these two lists as giving a permutation
$P_v$ of $1_v,\dots, 16_v$, where $1+\beta_{g_i}$ fixes
$16_v=(0,0,0,0,0,0)$, $i=1,2$. The rules for $P_v$: $(i_v)P_v=j_v$ if the ith
element of $L_1$ equals the $j$th element of $L_2$. The result:
$$P_v=(16_v)(1_v\,14_v\,4_v\,9_v\,13_v)(2_v\,12_v\,11_v\,6_v\,10_v)(3_v\,5_v\,
15_v\,8_v\,7_v).$$  
The final step is to check the collection of pairs $(c_1i_v,c_2j_v)\in
M\setminus V\times M\setminus V$ where $(i_v)P_v=j_v$. Eight of these pairs
are in $M_3'\times M_5'$ or in $M_5'\times M_3'$. Examples:
$(c_16_v,c_210_v)\in M_3'\times M_3'$  while $(c_15_v,c_215_v)\in M_3'\times
M_5'$. 
\end{proof} 

\begin{prop} \label{finTwoOrbits} 
There are two $\bar M_4$ orbits on
$\ni_1=\ni(G_1,\bfC_{3^4})^{\inn,\rd}$. These correspond to two 
components
$\sH_1^+$ and $\sH_1^-$ of $\sH(G_1,\bfC_{3^4})^{\inn,\rd}=\sH_1$, and 
in turn
to the two values achieved by
$s(\bg)$ as $\bg$ runs over $\ni_1$. Points on the genus 12 component $\sH_1^+$ 
have branch cycle descriptions $\bg$ for which $s(\bg)=1$ 
(Prop.~\ref{compSerre}).
\end{prop}

\subsection{Spin separation and two $\bar M_4$ orbits on
$\ni_1^\rd$}
\label{cliffAlgLift} 
Alternating groups come with a standard faithful permutation
representation. Some characteristic Frattini
covers have faithful (less obvious)  permutation
representations. This subsection has our main illustration of {\sl spin
separation\/} (Def.~\ref{spinSepRepDef}). It is that the lifting invariant
$s(O)$ for an $H_4$ orbit
$O$ on
$\ni(G_1,\bfC_{3^4})^\inn$ (from
$R\to G_1$) comes from the spin cover of 
$A_{40}$  (or of either
$A_{60}$ or $A_{120}$). 

\newcommand{\that}[1]{{}^{{\scriptscriptstyle T}} \hat #1}

\subsubsection{Embedding groups in alternating groups}
Suppose $k\ge 4$ and $T:G\to S_k$ is a permutation representation of a finite 
group $G$ of degree $k$. Assume the image of $G$ is in $A_k$, and let $\hat A_k$ 
be the central nonsplit extension of $A_k$ with kernel $\bZ/2$. 
Finally, let $\that G$ be the fiber product of $T$ and $\hat A_k\to A_k$. 

\begin{lem} \label{nonsplitLem} If elements of odd order generate 
$G$, then any permutation representation $T:G\to S_k$ has range
in $A_k$. Also, 
$\that G\to G$ is nonsplit if some $g\in G$ of order $2$ has a lift of 
order 4 in $\that G$ (Prop.~\ref{liftEven} applies; involution lifting 
property).  
\end{lem}

\begin{proof} Elements of odd order in $S_k$ have a presentation as a product of 
an even number of disjoint 2-cycles, so they are always in $A_k$. If $\that G$ 
is 
split, then any element $g\in G$ has a lift $\hat g\in \that G$ having the same 
order. Since $\that G\to G$ is a central extension with order 2 kernel, the 
order of a lift of an element of even order is independent of the lift. This 
concludes the proof. \end{proof}

Let $G$ be a finite group, with $\one_p$ denoting the trivial $G$ module of 
dimension 1
over $\bZ/p$. Consider
$\alpha\in H^2(G,\one_p)$. This always defines a central extension of $G$, and 
it does so
canonically if $G$ is $p$-perfect (Def.~\ref{pperfect} and the following 
comments). 

\begin{defn} \label{spinSepRepDef} Suppose $p=2$. Call a faithful $T:G\to A_k$ 
an {\sl
$\alpha$-spin separating\/} representation if $\that G$ realizes $\alpha$. 
\end{defn}  

When reference to $\alpha$ is clear refer to the
representation as spin separating. In our application to $G_1$ the extension
$\alpha$ is from Prop.~\ref{expInv2}, spin separating representations of
$G_1$ interpret that $\that {G_1}$ is $R$. Even if each involution of $G$ lifts 
to an
involution of $\that G$ (as in Lem.~\ref{nonsplitLem}), $\that  G\to G$ may
still not split (\S\ref{splitvs2-lift}). 

\subsubsection{Embedding $G_1$ in alternating groups} 
Use the notation of Cor.~\ref{R1G1}. 
\begin{prop} \label{expInv2} For $T:G_1\to S_k$, a faithful permutation 
representation, $T(G_1)\le A_k$. If a lift of 
$m\in M\setminus V$ to $\that G_1$ has order 4,  
then the following hold. 
\begin{edesc} \item Each $m\in M\setminus V$ lifts to have
order 4. 
\item Each $v\in V$ lifts to have order 2. 
\item $\that G_1\to G_1$ is the extension $R\to G_1$.  
\end{edesc} Therefore, with $a,b\in M\setminus V$, compute the lifting invariant 
$s(\bg)$ of the
$\bar M_4$ orbit of $$\bg=( g_1,  a g_1^{-1} a, b g_2 b, g_2^{-1})\in 
\ni(G_1,\bfC_{3^4})$$ (from 
Prop.~\ref{expInv}) as  
${\that a}^{g_1^{-1}}\that a\that b{\that b}^{g_2}$. 
\end{prop}

\begin{proof} Apply Lem.~\ref{nonsplitLem}. Elements of order 3 generate $G_1$. 
Now suppose $m'\in M\setminus V$ lifts to an element of order 4. Then, the 
nonsplit central extension $\that G_1$ must be the quotient of the universal 
central extension of the perfect group $G_1$ that characterizes orders of lifts 
from $M$ in Prop.~\ref{RkGk}. 
The criterion for this is that $T(m')$ is a product of $2l'$ disjoint 2-cycles 
with $l'$ odd 
(Prop.~\ref{liftEven}). 
\end{proof}

\subsubsection{Spin separating representations of $G_1$}  
Consider $M=\ker(G_1\to A_5)$ as an $A_5$ module from Prop.~\ref{A5frat}. 
Let $H\le G_1$, and denote the
corresponding (right) coset representation by $T=T_H:G_1\to A_{k}$ with
$k=(G_1:H)$. We list faithful 
spin separating representations
$T_H$. This is equivalent to  
$m'\in M\setminus V$  lifts to order 4 elements of $\that G_1$.  

Prop.~\ref{liftEven} characterizes this    
with two properties.
\begin{edesc} \label{G1prop} \item $H$ contains no nontrivial subgroup normal in
$G_1$.
\item \label{G1propb} $m'$ acting on $H$ cosets has $t$ fixed points with
$k-t=4l'$  and $l'$ odd. 
\end{edesc} 

A faithful representation of $G_1$ is not compatible
(as in Def.~\ref{defComp}) with the standard degree 5
representation on
$A_5$.  Reason: Any lift of $A_4$ (which contains a 2-Sylow
of $A_5$) contains a 2-Sylow of  $G_1$. 
So, it does not give a faithful representation.
Denote a $2$-Sylow of $G_1$ (resp.~$A_5$)
by $P_1$ (resp.~$P_0$).  Notation for the centralizer in $H$ of the set $S$ is 
$Z_H(S)$. 

\begin{lem} \label{2Sylspin} Let $T:G\to S_n$ be any
permutation representation with $H=H_T$ defining the cosets.
For $h\in H$, denote by $U_h\le G$ the subset of $G$ (not necessarily subgroup)
whose inverses conjugate
$h$ into $H$. Since multiplication by $H$ on the left maps $U_h$ into $U_h$, 
$U_h$ consists of right cosets of $H$. Use $(U_h:H)$ as the number of these.
Then
$T(h)$ fixes
$\tr(T(h))=(U_h:H)$ integers. Denote the orbit of $h$ under $U_h$
by $O$. Let
$\row O u$ be the orbits of $H$ on $O$. Choose $h_i\in O_i$,
$i=1,\dots,u$. Then,
$\tr(T(h))= \sum_{i=1}^u(Z_G(h_i):Z_H(h_i))$. 

Suppose $G=G_1$ and $T_H$ is spin separating. Then, some $m'\in
M\setminus V$ is in $H$. 
Further, the 2-Sylow of $H$ has order at least $2^4$ and at most $2^5$. 
\end{lem} 

\begin{proof} If $Hgh=Hg$, then $gh(g)^{-1}\in H$, and $g\in
U_h$. From our previous notation $g^{-1}$ conjugates $h$ into $H$. Right cosets
of $H$ in $U_h$ are the 
$H$ cosets that $h$ fixes. Suppose $g_1^{-1}$ conjugates $h$ to $h_i\in O_i$.
Running over $Hg_1$ gives elements whose inverses conjugate
$h$ into $O_i$. If $g_2^{-1}$ conjugates
$h$ to $h_i$, then $sg_1=g_2$ with $s\in Z_{G_1}(h_i)$. So, the cosets
with representatives conjugating $h$ to $h_i$ are in
one-one correspondence with $Z_{G_1}(h_i)/Z_H(h_i)$.  

The test for spin separation does not depend the choice of 
$m'\in M\setminus V$. If $H$ contains no 
$m'\in M\setminus V$, then $V$ contains the 2-Sylow of $H$ and
$m'$ fixes no cosets:
$t=0$. This implies $k=4l$ with $l$ odd. Also, the 2-Sylow of $H$ can't contain 
$V$, or else the representation won't be faithful. So the 2-Sylow $U$ of $H$ has
order at most 8.  This implies $k=(G_1:H)$ is divisible by $2^4$, a
contradiction to \eql{G1prop}{G1propb}.

Now assume the 2-Sylow of $H$ has order at most 8. As above, let $t$ be the
number of cosets $m'\in H\cap M\setminus V$ fixes. Then,
$2^4|k$ and $k-4l=t$ with $l$ odd, implies $t=4t'$ with $t'$ odd. The
centralizer of $m'$ has a 2-Sylow of order $2^6$, and this centralizer is in
$U_{m'}$. This implies $8|t$, a contradiction. 

Now suppose the 2-Sylow $U$ of $H$ has order at least $2^6$. Assume first that 
$U$
surjects onto a 2-Sylow $P_0$ of $A_5$. Then  $U$ is the whole 2-Sylow $P_1$ of
$G_1$, for $P_1\to P_0$ is a Frattini cover. If, however, $U$ does not surject
onto a 2-Sylow of $A_5$, then the kernel of the map has order $2^5$. So, $U$
contains $V$, and $T_H$ is not faithful.  Conclude: $|U|\le 2^5$. 
\end{proof}

\subsubsection{Degrees of spin separating representations of $G_1$}
Prop.~\ref{spinSepLst} lists properties of the spin separating
representations of $G_1$ of degree 120. Prop.~\ref{spinSepLst2} does the same
for the degree 60 and 40 representations. 

\begin{prop} \label{spinSepLst} Let $\alpha'\in P_1$ have order 4.  Its
image $\alpha\in A_5$  and 
$(\alpha')^2$ determine $\alpha'$. Given $\alpha$, possible 
$(\alpha')^2$ are elements of $M(A_5)\setminus V$
fixed by $\alpha$. 

Let $K\le V$ be a Klein 4-group on which $\alpha'$ acts. With
$H_{\alpha',K}=H=\lrang{\alpha',K}$,
$T_{H}$ spin separates if and only if $|M_3'\cap H|$ is odd:
$\alpha'$ is nontrivial on $K$. 

If $|M_3'\cap H|=1$ (resp.~$|M_5'\cap H|=1$), then  $|M_5'\cap
H|=3$ (resp.~$|M_3'\cap H|=3$),  $|V\setminus \{1\}\cap H|=3$. Further,
corresponding to the two choices:
\begin{edesc}  \item 
$\tr(T_H(m))=12$ (resp.~$\tr(T_H(m))=36$) if $m\in M_3'$; 
\item $\tr(T_H(m))=60$ 
(resp.~$\tr(T_H(m))=20$) if $m\in M_5'$; and 
\item $\tr(T_H(v))=24$ if $v\in
V\setminus \{1\}$.  \end{edesc} 

Four conjugacy classes of subgroups
$H\le G_1$ correspond to these choices. They give all spin separating
representations of degree 120. 

Suppose $|M_3'\cap H|=1$. Then, $T_H$ acts imprimitively
on ten sets of integers of cardinality 12. Each set consists of the
integers fixed by an element of $M_3'$. This induces the degree 10
representation of
$A_5$ on ordered pairs of distinct integers from
$\{1,2,3,4,5\}$. If $g\in G_1$ has order 10, then $T_H(g)$ is a product of
twelve 5-cycles and six 10-cycles. If $g\in G_1$ has order 6, then $T_H(g)$ is
a product of four 3-cycles and eighteen 6-cycles.  
\end{prop}

\begin{proof} Use the notation for $\alpha'\in P_1$ of order 4 in the
statement. To be explicit assume
\wsp  up to conjugacy \wsp $\alpha'$ lifts
$\alpha_{24}=(2\,4)(3\,5)$. Given one lift  $\alpha_0$, all others have the 
shape $m\alpha_0=\alpha'$. Map $M\to V$ by $m\mapsto mm^{\alpha_{24}}=\psi(m)$.
Cor.~\ref{A5Morbs} says the kernel of $\psi$ has rank 3, so the image $\psi(M)$ 
of
$\psi$ is a homogeneous space for the squares of lifts of $\alpha_{24}$. 
Conclude: 
$\alpha_0$ and $(\alpha')^2$ determine $\alpha'$; and  $(\alpha')^2$ runs over 
the subset $W_{24}$ of $M(A_5)\setminus V$ that $\alpha_{24}$ fixes.

Suppose $H$ has order 16 and $T_H$ is spin separating.
Then, some $\alpha'\in H$ has order 4, $m'=(\alpha')^2\in M\setminus V$
and $\alpha'$ stabilizes $K=V\cap H$. Write 
$H_{\alpha',K}=H=\lrang{K,\alpha'}$. Eight divides the degree of the 
representation. As in
Lem.~\ref{2Sylspin}, spin separation is equivalent to 
$t=|U_{m'}|/|H|$ being four times an odd number. 

First case: Suppose 
$\alpha'$ is trivial on the Klein 4-group $K$. 
Cor.~\ref{A5Morbs} says $H$ contains exactly two conjugates of $m'$,
and all of $H$ centralizes
$m'$. Let $c$ be 3 (resp.~5) if $m'\in M_3'$ (resp.~$m'\in M_5'$). So, $m'$
fixes $$t=2\cdot |Z_{G_1}(m')|/|Z_{H}(m')|=2c\cdot 2^6/2^4$$
integers. Therefore, 8 divides $k-t$, and
$T_H$ is not spin separating. 

Second Case: Suppose $\alpha'$ is nontrivial on $K$. If $(\alpha')^2=m'$ is the
only element in $H_{\alpha',K}$ in its conjugacy class, then the computation
just concluded gives $$\tr(T_H(m'))=
|Z_{G_1}(m')|/|Z_{H}(m')|=c2^6/2^4=4c.$$ So, 4 exactly divides $k-t$ and $T_H$
{\sl is\/} spin separating. If three elements in $H$ are conjugate to
$(\alpha')^2=m'$, then $\tr(T_H(m'))=12c$. Now suppose $v\in
V\setminus\{1\}$. As $\alpha'$ is nontrivial on $K$, conjugates (under $G_1$) of 
$v$ in $H$ fall
in two $H$ orbits.  Let $v'$ be a
representative of one of these orbits. Then, $Z_G(v')/Z_H(v')=4$, so 
$(U_{v'}:H)=8$. 

The argument is similar for 3 conjugates of
$m'$ in $H$.  Now we count subgroups $H$, up to conjugacy in $G_1$,
giving spin separating representations of degree 120.

Giving $H$ is equivalent to giving pairs $W^*_{24}=\{(w,m)\in
W_{24}\times M(A_5)\setminus  V\}$
with $m$ subject to these conditions. 
\begin{equation} m^{\alpha_{24}}\ne m, \text{ and } |M_5'\cap
\{w,m,m^{\alpha_{24}},wmm^{\alpha_{24}}\}| \text{ is odd}. \end{equation}
Note: $m$ and $m^{\alpha_{24}}$ are in the same conjugacy class of $M\setminus 
V$;
this asks exactly that $w$ and $wmm^{\alpha_{24}}$ are in different conjugacy 
classes
of $M\setminus V$.  

Count elements of
$W^*_{24}$ from the proof of Cor.~\ref{A5Morbs}.  With the cosets of $M(A_5)$ 
labeled
as there,
$\beta_\alpha=(1\,6)(3\,5)$ represents
$\alpha_{24}$. Thus, its centralizer in $M(A_5)$ is
$$U=\lrang{(a_1,x_1,a_2,x_2,a_2,a_1)\mid a_1,a_2, x_1,x_2\in
\bZ/2}/\lrang{(1,1,1,1,1,1)}.$$ Calculate: $W_{24}=\{(w_1,0,w_2,0,w_2,w_1)\mid
w_1,w_2\in \bZ/2\}$. Uniquely determine a representative for $w$ by taking
$x_1=1$ and $x_2=0$. The elements $w$ run over all choices $a_1,a_2\in
\bZ/2^2$. So, $|W^*_{24}|$ equals the number of pairs
$(a_1,1,a_2,0,a_2,a_1)$ and $(a_1+w_1,1,a_2+w_2,0,a_2+w_2,a_1+w_1)$ with one
having three nonzero entries, and the other not. Given any pair $(a_1,a_2)$, 
there are
two choices for
$(w_1,w_2)$ satisfying this condition. Conjugate the set of  groups $H_{w,m}$ 
attached
to these choices by any lift of $(2\,3)(4\,5)\in A_5$ (centralizing
$\alpha_{24}$). The result is a  distinct subgroup. Conclude there are four such 
subgroups,
falling in pairs according to $w\in M_3'$  or $w\in M_5'$.   

Now, suppose $|M_3'\cap H|=1$. As above, for $m\in M_3'$,  $T_H(m)$ fixes
twelve integers. Running over the ten elements of $M_3'$, transitivity
of the representation guarantees some element of $M_3'$ fixes any given
integer in $\{1,2,\dots,120\}$. As $|M_3'|=10$, these sets of fixed integers
are disjoint, and they form a set of imprimitivity. With no loss assume $m\in
H$, so $H\le Z_{G_1}(m)$. Elements of $Z_{G_1}(m)=H^*$ are exactly those
permuting the fixed integers for $T_H(m)$.  Acting on  $H^*$ cosets 
gives the degree 10 representation of $A_5$ on pairs of distinct
integers from
$\{1,2,3,4,5\}$. 

Continue the hypothesis $|M_3'\cap H|=1$. Suppose $g\in G_1$ has order 10.
Then, $g^2$ has order 5, and $T_H(g)$ fixes no integers. On the other hand,
$g^5\in M_5'$. From above, $T_H(g^5)$ fixes exactly sixty integers. Conclude:
$T_H(g)$ is a product of twelve 5-cycles and six 10-cycles. For $g\in G_1$ of 
order 6, similarly deduce $T_H(g)$ is a product of four 3-cycles and eighteen
6-cycles. 
\end{proof} 

The following is similar to Prop.~\ref{spinSepLst}. We leave details to
the reader. 
\begin{prop} \label{spinSepLst2} Continue Prop.~\ref{spinSepLst} notation (for 
$\alpha'$). 
Assume $\beta\in G_1$ has order three,
$\lrang{\beta,\alpha'}$ has order $3\cdot 2^4$ and the image of
$\lrang{\beta,\alpha'}$ is an
$S_3\le A_5$. Then, with
$H_{\alpha',\beta}=H'$, $|M_3'\cap
H'|=1$ and $T_{H'}$ spin separates. There are
two conjugacy classes of such groups
$H'$ corresponding to these choices giving all spin separating representations
of degree 40. If $m,m'\in M_3'$ are distinct, the four fixed integers of
$T_{H'}(m)$ are distinct from those fixed by $T_{H'}(m')$. 

Let $H$ be one of the groups with $T_H$ of degree 120, and let
$\alpha^*\in V\setminus H$. Then, $H''=\lrang{H,\alpha^*}\cap M_3'$ has 4 or 6
elements. Exactly in the former case $T_{H''}$ spin separates. There are 
two conjugacy classes of such groups
$H''$ corresponding to these choices giving all spin separating representations
of degree 60.
\end{prop}

\begin{proof}[Comments on the proof] If $|M_3'\cap H|=1$, as in
Prop.~\ref{spinSepLst}, the conditions on $\beta$ are that it centralizes this
element of $M_3'$ and it acts on $V\cap H$ in the standard Klein 4-group
representation. For this case, $T_{H'}(m)=4$ (resp.~20 and 8) if $m\in M_3'$
(resp.~$m\in M_5'$ and $m\in V$). 
\end{proof}
 
\begin{cor} Prop.~\ref{spinSepLst} and Prop.~\ref{spinSepLst2} include the
complete list of spin separating permutation representations of the $\bar M_4$
orbits on $\ni(G_1,\bfC_{3^4})^{\inn,\rd}$. 
\end{cor} 

\begin{proof} Assume for $H\le G_1$, $T_H$ gives a faithful spin separation.
According to Lem.~\ref{2Sylspin}, the 2-Sylow $U$ of $H$ has order $2^4$ or
$2^5$ and it contains an element $m'$. Prop.~\ref{spinSepLst} lists all cases 
where
$k=(G_1:H)$ is 120, 40 or 60. The only degrees left as dangling possibilities 
are
where $H$ has order divisible by 5 or where it has order $2^5\cdot 3$. In the
former case, since $H$ contains nontrivial elements of $V$, it must contain all 
of
$V$. So it is not faithful. Similarly for the latter case: The action of an 
element
of order 3 on $V\cap H$ forces this to be all of $V$.  
\end{proof} 
 
\subsubsection{Splitting versus the involution lifting property}  
\label{splitvs2-lift}

Restriction of an element from $H^2(G,A)$ to $H^2(P_p,A)$ is injective if $P_p$ 
is the
$p$-Sylow of $G$ \cite[p.~105]{AW}. Sometimes this allows taking 
$G$ to be a $p$-group. We use this and produce a 2-group extension
$\psi: \hat P\to P$ with the following properties. 
\begin{edesc} \label{2splitprop} \item The kernel of $\psi$ is $\bZ/2$. 
\item $\psi$ is nonsplit and spin separating: $\that P=\hat P$, pullback from an 
embedding
$P\le A_n$ for some $n$. 
\item  \label{2splitpropc} Each involution of $P$ lifts to an involution of 
$\hat P$. 
\end{edesc}  The quaternion group presentation is in 
Def.~\ref{dicyclic}: 
$Q_{4n}$ of order 
$4n$ 
has generators  
$\tau_1,\tau_2$ with $\ord(\tau_1)=2n$, $\ord(\tau_2)=4$, $\tau_2^{- 
1}\tau^\sph_1\tau^\sph_2=\tau_1^{-1}$ and  $\tau_2^2=\tau_1^n$. Let $\hat 
Q_{4n}$ be the
group from dropping the condition $\tau_2^2=\tau_1^n$. Then, the only involution
$\tau_2^2$ in $\hat Q_{4n}$ lifts to an involution in $\hat Q_{4n}$. This 
extension,
however, does not split. 

The calculations below apply for $n$ even to all the groups $\hat
Q_{4n}$. For simplicity we do the case $n=2$ and $Q_8$, using $\tau_1=a$ and
$\tau_2=b$; $\hat Q_8$ has  generators $\hat a$ and $\hat b$, with relations 
$\hat a^4=
\hat b^4=\hat a^{-1}\hat b\hat a\hat b=1$. For $n$ even,  $Q_{4n}$ doesn't have
an extension satisfying
\eqref{2splitprop}, though $\hat Q_{4n}$ does. (This observation should start 
with the
Klein 4-group whose $Q_8$ extension in Lem.~\ref{Q8Sp5} is an example of spin
separation.)

\begin{prop} \label{hatP-P} The extension $\hat P=\hat Q_8\to P=Q_8$ has no
spin separating extension. Still, $\hat P$ has an extension $\hat{\hat
P}$ satisfying all properties of \eqref{2splitprop}. \end{prop} 

\begin{proof} The group $Q_8$ has only one faithful transitive permutation 
representation
as a subgroup of $S_8$. In this representation elements of order 4 have the 
shape
$(4)(4)$. So these generators are in $A_8$: $Q_8\le A_8$. This is the only
embedding we must test to check if $Q_8$ has a spin separating representation of 
any
kind. 

Form a central extension of $Q_8$ by mapping $\hat a$ to $a$ and $\hat b$ to 
$b$.  In the extension, the square of $\hat a$ is cleaved away from the squares $\hat 
b^2= (\hat a\hat b)^2$ (they remain equal).  Let $\one_2$ be the trivial $Q_8$ 
module. Recall: $Q_8$ has rank $e=2$ and order $n=8$. So,  the 1st characteristic 2-Frattini
module $M(Q_8)=\ker_0(Q_8)/\ker_1(Q_8)$ (notation of \S \ref{normp-Sylow}) of
$Q_8$  has dimension
$9=1+(e-1)n$ \cite[\S15.6]{FrJ}: Schreier's formula for the number of generators 
of a subgroup of a free group of rank $e$ and index $n$. As a $Q_8$ module it is
indecomposable \cite[Indecom.~Lem.~2.4]{FrKMTIG}. The maximal quotient of it on
which $Q_8$ acts trivially is $H^2(Q_8,\one_2)$. 

We ouline our computation that  
$H^2(Q_8,\one_2)$ has dimension two.  Use the notation of
\cite[Part II]{FrMT}. The augmentation map gives $\bF_2[Q_8]\to \one_2$. Denote 
its
kernel by $\Omega(1)$. The version of Jenning's Theorem (\cite{Je} or 
\cite[p.~89]{Ben1})
in \cite{Qu} gives an effective tool for computing the Loewy layers of any $p$-group ring
using the Poincar\'e-Witt basis of its universal enveloping algebra. The Loewy
layer at the top of
$\Omega(1)$ (the second radical layer in $\bF_2[Q_8]$) is $\one_2\oplus\one_2$. 
The
projective module $\bF_2[Q_8]\oplus \bF_2[Q_8]$ maps naturally and surjectively 
to
$\Omega(1)$ extending the direct sum of two augmentation maps 
$\psi:\bF_2[Q_8]\oplus
\bF_2[Q_8]
\to
\one_2\oplus\one_2$. The kernel of this map is $\Omega(2)=M(Q_8)$
\cite[Projective Indecomposable Lem.~2.3]{FrMT}. The explicit basis of 
$\bF_2[Q_8]$ from
the group elements produces a natural interpretation of $\psi$ as a matrix. Row 
reduce
and compute explicitly a basis of $\Omega(2)$ from that of $\bF_2[Q_8]$. This 
gives 
$\Omega(2)$ a 
$Q_8$ module structure. Rational canonical form of the matrix action
$M_a$ of $a$ (resp.~$b$ and $M_b$) on this module determines the maximal 
quotient on which
$M_a$ (resp.~$M_b$) acts like the identity matrix, giving the maximal quotient 
on which
$Q_8$ acts like the identity. We let {\bf GAP} do this computation.  

All non-split extensions of $Q_8$  
are isomorphic as groups though not as extensions. As in \S\ref{pperfSec}, the
isomorphisms between them come from $H^1(G,\bZ/2)$ which is isomorphic to a 
Klein
4-group. Extensions correspond to which order 4 subgroup of
$Q_8$ cleaves from the other order 4 subgroups.   

The proof of Prop.~\ref{liftEven} reviews the generators and relations for $\hat 
S_n\to
S_n$.  In that notation, symbols $[i\,j]$, $1\le i, j\le n$ (in the 
multiplicative group of
units in the Clifford algebra),
generate subject to the relations
$[i\,j]^2=1$ and
$[i\,j]=-[j\,i]$. The map $\hat S_n\to S_n$ appears from $[i\,j]\mapsto (i\,j)$. 
Finally, $[i\, j][j\, k][i\, j] = [k\, i]$.  

Like $Q_8$,  $\hat Q_8$ has no non-trivial coreless subgroup. So, only  its
regular representation is faithful and its generators $\hat a$ and $\hat b$, of 
order 4,
have the shape $(4)(4)(4)(4)$; they are in $A_{16}$. Involutions in $\hat Q_8$ 
are the
products of 8 disjoint two cycles. So Prop.~\ref{liftEven} implies they lift to 
the
same order on pullback in $\Spin_{16}$.  

Now we show $\hat{ {\hat Q_8}}$ (the pullback of $\hat Q_8$ to $\Spin_{16}$) 
does not split
off
$\hat Q_8$. If the extension splits, the lifts would retain
the relation $\hat b^2=(\hat a\hat b)^2$. So, it suffices that the square
of any lift of
$\hat b$  differs from the square of any lift of $\hat a \hat b$.  With no loss:
$$\begin{array}{rl}\hat b&=(1\, 5\, 9\, 13)(2\, 6 \,10\, 14)(3\, 7\, 11\, 
15)(4\, 8\, 12\,
16), \text{ and}\\
\hat a\hat b&= (1\, 6\, 9\, 14)(2\, 7\, 10\, 15)(3\, 8\, 11\, 16)(4\, 5\, 12\,
13).\end{array}$$  Write each 4-cycle's lift to $\Spin_{16}$ as the product of 
three
obvious generators. Example:
$(1\, 5\, 9\, 13)$ lifts to $[1\, 5][5\, 9][9\, 13]$.  Square the lifts of $\hat 
b$ and
$\hat a\hat b$ using the relations from the previous paragraph. Then multiply 
the two
squares together to get -1, the generator of the kernel of $\Spin_{16}\to 
A_{16}$.  This
also shows -1 is in the Frattini subgroup, so the extension is Frattini, and
therefore nonsplit. 

The same procedure for $Q_8$ shows the squares of the
lifts of $a$, $b$, and $ab$ are the same.  Every non-split extension of $Q_8$
has one of these cleaved off from the rest. This shows  the spin
extension of $Q_8$ splits. \end{proof}

\subsection{Understanding spin separation} 
\label{finConc}  Representations in Prop.~\ref{spinSepLst} give spin
separation of elements in $\ni_1=\ni({}_2^1\tilde A_5,\bfC_{3^4})^{\inn,\rd}$.
\S\ref{expSpin} gives the lifting invariant that separates the two $\bar M_4$
orbits on $\ni_1$. \S\ref{orbMyst} and
\S\ref{B4Act} together suggest why the two orbits have the same
$\bar M_4$ monodromy groups. \S\ref{indSteps} concludes the paper with
lessons for similar computations at higher levels. 

\subsubsection{Explicit lifting invariants} \label{expSpin} First we
characterize commuting pairs of elements of order 2 in $S_n$. Denote the
elements of $S_n$ that are of products of $u$ disjoint 2-cycles by
$C_u=C_{u,n}$. In the situation where $g\in C_u$ and $g'\in C_{u'}$, and $g$
and $g'$ commute, we will identify the product.  The
most complicated standard situation is where $(g,g')$ is
a {\sl $K4$-pair\/}: 
\begin{triv} $g=(i_1\,i_2)(i_3\,i_4)$, $g'=(i_1\,i_3)(i_2\,i_4)$ and 
$gg'=(i_1\,i_4)(i_2\,i_3)$. 
\end{triv}

Given $(g,g')$ as above, let $a(g,g')$ be the number of disjoint 2-cycles 
appearing in $g$ having no support in $g'$. Note: $a(g,g')$ may not equal
$a(g',g)$. Let $b(g,g')$ by the disjoint 2-cycles appearing in both $g$ and
$g'$. Suppose $(h,h')$ is a $K4$-pair, $g_1$ and $h$ (resp.~$g_1'$ and $h'$)
have no common support, and $g=g_1h$ and $g'=g_1'h'$ commute. We say the
$K4$-pair $(h,h')$ appears in $(g,g')$. Finally, let
$c(g,g')$ be the number of $K4$-pairs appearing in $(g,g')$.   
\begin{lem} \label{2-commute} Suppose $g\in C_u$ and $g'\in C_{u'}$ and
$gg'=g'g$. Then,
$gg'\in C_{u''}$ with $a(g,g')+a(g',g)+2c(g,g')=u''$. \end{lem}

\begin{proof} Write $g$ as $g_{a(g,g')}g_{b(g,g')}g_1$ with $g_{a(g,g')}$
(resp.~$g_{b(g,g')}$) the disjoint cycles in $g$ that don't (resp.~do) appear
in $g'$. Similarly, write $g'$ as $g'_{a(g,g')}g'_{b(g,g')}g'_1$. Since
$g$ and $g'$ commute, so do $g_1$ and $g_1'$. Proving the lemma amounts to
extracting  a $K4$-pair
$(h,h')$ that appears in $(g_1,g_1')$. Then, an induction characterizes
$(g_1,g_1')$ as given by products of disjoint $K4$-pairs. 

Suppose $i_1$ appears in a disjoint cycle $(i_1\,i_2)$ in $g_1$. By assumption, 
$i_1$ appears in a disjoint cycle $(i_1\,i_3)$ in $g_1'$. Since
we have extracted $g_{b(g,g')}$, $i_3\ne i_2$. Suppose $i_3$ is not in a
disjoint cycle of $g_1$. Then, $g_1g_1'$ applied to $i_2$ has the effect
$i_2\mapsto i_3 \mapsto i_1$. Thus, $g_1g_1'$ has order larger than 2. 
So, contrary to a previous deduction, $g_1$ and $g_1'$ don't commute.

Therefore, $i_3$ appears in a disjoint cycle $(i_3\,i_4)$ of $g_1$. Calculate:
$g_1g_1'$ applied to $i_4$ has the effect $i_4\mapsto i_3\mapsto i_1$. Now
compute the effect of $g_1'g_1$ on $i_4$ to conclude $g_1'$ contains the
disjoint cycle $(i_2\,i_4)$. That is, $(h,h')=(i_1\,i_2)(i_3\,i_4),
(i_1\,i_3)(i_2\,i_4)$ appears in $(g,g')$. 
\end{proof} 

\begin{exmp}[Applying Lem.~\ref{2-commute}] Suppose $|M_3'\cap H'|=1$, as
in Prop.~\ref{spinSepLst2}.  Suppose $m,m'\in M_3'$ are distinct. For this
case, $\tr(T_{H'}(m))=4=\tr{T_{H'}(m')}$ and $\tr(T_{H'}(mm'))=8$: $m,m' \in C_{18}$ and
$mm'\in C_{16}$. Integers not in the support of $m$ appear in 2-cycles in $m'$,
and give
$a(m,m')=a(m',m)=2$. Conclude $16=4+2c(m,m')$: $c(m,m')=6$. Use that
$a(m,m')+b(m,m')+2c(m,m')=18$ to conclude $b(m,m')=4$. 
\end{exmp}

\begin{prop} \label{10shape} Let $g\in A_5$ have order 5, and 
denote by $\hat g$ a lift of it to an element of order 10 in $G_1$. With $H'$
giving the degree 40 representation used above, let $H^*$ be its image
(isomorphic to $S_3$) in $A_5$. Then, we can speak of the 20 cosets of $T_{H'}$
lying above the cosets $H^*g^j$, $j=0,1,2,3,4$, of $T_{H^*}$. Analogous to
previous calculations, $T_{H'}(\hat g)$ consists of four 5-cycles and two
10-cycles. The integers moved (resp.~fixed) by $T_{H'}(\hat g^5)$ are in the
support of the 10-cycles (resp.~5-cycles). Suppose $\hat g^5\in H'$. Then,
multiplying $\hat g$ on the right of these cosets gives four 5-cycles of
$T_{H'}(\hat g)$. Two 10-cycles comprise the remaining action of $T_{H'}(\hat
g)$. If $\hat g^5\not\in H'$, then the description of the 5-cycles and
10-cycles of $T_{H'}(\hat g)$ switches. If $m,m'\in M_5'$, $m\not=m'$, then
$a(m,m')=6$, $b(m,m')=0$, $c(m,m')=2$. \end{prop}

\begin{proof} The number of integers in the support of 2-cycles from 
$T_{H'}(\hat g^5)$ is the same as in the 10-cycles of $T_{H'}(\hat g)$. From
Prop.~9.12 this means there are twenty such integers, and therefore two
10-cycles in $T_{H'}(\hat g)$.  Since this element fixes no integers there are
also four 5-cycles in $T_{H'}(\hat g)$. 

Without loss, take the 20 cosets of $T_{H'}$ above the 
cosets $H^*g^j$, $j=0,1,2,3,4$, of $T_{H^*}$ to be $H'm_i\hat g^j$,
with $m_0=1$, $m_i\in M$, $i=1,2,3$, and $j=0,\dots,4$. The $m_i\,$s are four
representatives of $H'\cap M$ cosets in $M$. First assume $\hat g^5\in
H'$. For each $i$, $\hat g$ cycles $H'm_i\hat g^j$, 
$j=0,\dots, 4$ (use that $\hat g^5$ commutes with $m_i$). The argument is the
same if $\hat g^5\not \in H'$: Write the distinct cosets
above as $H'm_i\hat g^j$ with $m_0=1$, $m_i\in M$, $i=1$ and
$j=0,\dots,9$.

Now inspect the relation between two elements of $M_5'$.  
There are six elements in $M_5'$. An element $g\in A_5$ of order 5 that
stabilizes $m$ is transitive on the remaining five elements of $M_5'$. Since
$a(m,m')=a(m^g,(m')^g)$,  $b(m,m')=b(m^g,(m')^g)$ and $c(m,m')=c(m^g,(m')^g)$,
these values are independent of their arguments.  Use that
$a(m,m')+b(m,m')+2c(m,m') =10$ and $a(m,m')+a(m',m)+2c(m,m') =16$. Thus,
$a(m,m')-b(m,m')=6$, and $b(m,m')=2$, 1 or 0. 

If $T_{H'}(m)$ and
$T_{H'}(m')$ have eight common fixed integers, then they also move eight
common integers. So we show the former and conclude $b(m,m')=0$ and
$c(m,m')=2$. With no loss, assume $m\in H'$ and find the cosets from among
$H'm_i\hat g^j$ fixed by  $m'$. Since $m'\not=m$, as $j$ varies, $ (m')^{\hat
g^{-j}}$ runs over all five other elements of $M_5'$. Use $\hat
g^j m'= (m')^{\hat g^{-j}}\hat g^j$ to see $m'$ fixes exactly two of these
cosets as $j$ varies, and eight cosets in all.
\end{proof}

\subsubsection{$\bar M_4$ orbit mysteries} \label{orbMyst}
The  two
$\bar M_4$ orbits on $\ni({}_2^1\tilde A_5,\bfC_{3^4})^{\inn,\rd}$ have the same 
image 
groups, and degrees. This is despite their being considerably different
as representations. They are clearly permutation inequivalent
representations, or else $\gamma_\infty$ would have the same shape on each
orbit. More strikingly, however, the permutation representations are
inequivalent as representations (have different traces). 
The number of $\gamma_\infty$ orbits on each of the
$\bar M_4$ orbits $O_1^+$  and $O_1^-$ (resp.~26 and 32
$\gamma_\infty$ orbits) is different. This shows the two
representations are different (see for example \cite[Lem.~5.3]{Fr-Schconf})  
despite their similarities. 
\S\ref{GKS} speculates on representation theory  (based on \cite{GKS})  that 
might explain this
and also measure the difference between these two representations. 

\subsubsection{$B_4$ as evaluation homomorphisms} \label{B4Act}
From Lem.~\ref{B4F4}, 
the direct product of the free group
$K_4^*=\lrang{(Q_3Q_2)^3,Q_1^{-2}Q_3^2, (Q_2Q_1)^{-3}}$ and
$\lrang{(Q_1Q_2Q_3)^4}$ equals $N_4\eqdef \ker(B_4\to M_4)$. Denote 
$Q_1Q_2Q_{3}^2Q_2Q_1$ ($D$ in \eqref{Daction}) by $R_1$. The following
presentation of
$N_4$ is superior for our purposes. It easily follows from the notation 
and
proof of Lemma \ref{B4F4}. Recall: The effect of
$R_i$ on $G_4$ is conjugation by $\bar\sigma_i$, $i=1,\dots,4$. We state this 
only for
$r=4$, though it works for any value of $r$. 

\begin{prop} \label{Riopers} These 
are generators of $N_4$: 
$$\pmb R=(R_1,\ R_2=Q_1^{-1}R_1Q_1, 
R_3=Q_2^{-1}R_2Q_2, R_4=Q_3^{-1}R_3Q_3)  \text{ and } (Q_1Q_2Q_3)^4.$$ Any
representative element $\bg$ in a Nielsen class $\ni(G,\bfC)=\ni$ produces an 
{\sl
evaluation\/} homomorphism $\psi_{\bg}: N_4\to G/\Cen(G)$ mapping
$\prod_{i=1}^4 R_i$ to 1. Conversely, any homomorphism $\psi:N_4\to H$ mapping 
$\prod_{i=1}^4
R_i$ to 1 produces an associated Nielsen class representative. 

For $Q\in B_4$, act on $\psi_\bg$ by applying evalution of $Q^{-1}{\pmb R}Q$ to 
$\bg$. An
orbit on $\ni$ is equivalent to a $B_4$ orbit on the homomorphisms
$\psi_\bg$ (up to conjugation by $G$).  
\end{prop} 

For the computation of orbits of $H_4$, consider the induced action of
$N_4$ on $G_4$.   Denote images of $\row R 4$ by $\row r 4$, so these generate a
free group with one relation $r_1r_2r_3r_4=1$. An obvious reformulation of
Prop.~\ref{Riopers} has a rephrasing in $H_4$ orbits using $\bar N_4=\lrang{\row 
r 4}$. 
Use
the same notation $\psi_\bg$ for the corresponding evaluation homomorphisms. 
Referring 
to the
generators $\row r 4$ allows us give the Nielsen class of a homomorphism by what 
conjugacy
classes the generators hit. 

\begin{defn} Let $\mu$ be an automorphism of $\bar N_4$. We say $\mu$ preserves 
a 
Nielsen
class $\ni(G,\bfC)=\ni$ if $\psi_\bg\circ \mu$ is in the same Nielsen class as 
$\psi_\bg$.
\end{defn}

\begin{lem} \label{muorbit} If $\mu$ preserves the Nielsen class $\ni$, then it 
permutes 
$H_4$
orbits. Any two $H_4$ orbits under $\mu$ have the same length. Further, the 
monodromy 
groups
(permutation groups for the action of $H_4$) on these orbits will be the same. 
\end{lem}

Note: While the group on two $\mu$ orbits in Lem.~\ref{muorbit} are the same, 
the 
permutation
representations may be inequivalent. Denote the automorphism of $\bar N_4$ by
$(r_1,r_2,r_3,r_4)\mapsto (r_4^{-1}, r_3^{-1},r_2^{-1},r_1^{-1})$ by $\mu_0$. 

\begin{exmp}[Cont.~Ex.~\ref{A4C34-wstory}] Among the six $\bar M_4$ orbits on 
$\ni({}_2^1 \tilde
A_4,
\bfC_{3^4})^{\inn,\rd}$) the automorphism $\mu_0$ fixes four of them, but 
permutes the other two. 
\end{exmp}

\begin{quest} Is there a Nielsen class preserving automorphism interchanging the 
two 
$H_4$
orbits on $\ni({}_2^1 \tilde A_5, \bfC_{3^4})^{\inn,\rd}$?
\end{quest}

\subsection{Inductive steps in going to higher Modular Tower levels} 
\label{indSteps} Explicit formulas 
(\eqref{RHOrbitEq} and Lem.~\ref{RHOrbitLem}) list properties of ${}_p\tilde G$ 
characteristic
quotients  helpful to decide Prob.~\ref{MPMT} and other properties of Modular 
Towers when $r=4$.

\subsubsection{The starting point and Schur multipliers} By example, the paper 
establishes it
is difficult to predict component genuses at level 0 of a Modular Tower. As in 
all our $A_5$ (and
$A_4$) examples, genus 0 components are common. For many applications this is
excellent; rational points from such a moduli space produce Galois
realizations. So, establishing that at some level in the Modular Tower all 
components will have
large genus requires exploiting the module representation
theory of the universal $p$-Frattini cover. Whatever the genus of the components 
at level 0, it
is properties of the universal $p$-Frattini cover that force sufficient 
ramification at higher
levels to have the genus go up (as in \eqref{RHOrbitEq}). 

Lem.~\ref{RHOrbitLem} exploited the hypotheses \eqref{bestGenAssump}, valid  in 
the main
Modular Tower $(A_5,\bfC_{3^4}, p=2)$. A more general hypothesis, with 
$\{O_k\}_{k=0}^\infty$ a
projective system of $\bar M_4$ orbits, is the following: 
\begin{triv}  $\frac{\tr_{O_k}(\gamma_0)+\tr_{O_k}(\gamma_1)}{|O_k|} \mapsto 0$; 
and
 $\frac{I_k}{|O_k|}\mapsto 1$. \end{triv}   These two assumptions give a 
positive
conclusion to Prob.~\ref{MPMT} (when $r=4$) from Lem.~\ref{RHOrbitLem} in that 
there is no bound on
the genus of components corresponding to the projective system 
$\{O_k\}_{k=0}^\infty$. 

Rem.~\ref{serFixedPt} and Prob.~\ref{ellFixedPtr} view 
fixed points of $\gamma_0$ and $\gamma_1$ on projective systems of Nielsen 
classes in a Modular
Tower as generalizing a classical observation about fine moduli for the modular 
curve spaces
$X_1(p^{k+1})$.  

Ex.~\ref{A4C34-wstory} for $(A_4,\bfC_{3^4},p=2)$ is a valuable test case. The 
genus
situation from the view of Prob.~\ref{MPMT} is even worse at level 0 and level 1 
than for our main
example. There are components, some with low genuses. Yet, the universal 2-
Frattini situations
produce the desired ramification.   

\subsubsection{Jumps in $\mpr\,$s} \label{jumps} 
In all examples, perturbation of level 0 H-M reps.~starts the
story. Spin separating representations are the most intriguing element in 
structured 
description of Modular Towers levels, like listing real points in
Prop.~\ref{oneOrbitkappa}. Serre's Prop.~\ref{serLift} is indispensable, 
with \S\ref{onlytwo}  using it indirectly even when the initial hypotheses don't 
apply.
Still, this is only about the prime $p=2$, and only about universal 2-Frattini 
covers that
have an initial relation with alternating groups. 

An elementary approach to producing the desired ramification comes from the 
conjugacy class
counting formula of \S\ref{congClProd}. Suppose $G_k$ is the $k$th 
characteristic quotient of
the universal $p$-Frattini cover of a finite $p$-perfect group $G=G_0$. Assume
$g_1, g_2\in G_{k+1}$ have $p'$ order. Princ.~\ref{obstPrinc} reminds that 
obstruction (from
\cite[\S3]{FrKMTIG}; as in \S\ref{5-3-3prod}) occurs only when $\one_{G_k}$ 
appears in the Loewy
display of $M_k=\ker(G_{k+1}\to G_k)$:  $G_{k+1}\to U_k\to W_k\to G_k$ with 
$W_k$ centerless and $\ker(U_k\to
W_k)$ the trivial $W_k$ module. With no loss for the question of concern, assume 
$g_1,g_2\in W_k$. Suppose also that $p$ does not
divide the order of $g_1g_2$. Then, consider the triple of conjugacy classes 
$\C_1,\C_2,\C_3$ defined by
$(g_1,g_2,(g_1g_2)^{-1})$. 
Consider (any) lifts $\hat g_1,\hat g_2\in U_k$. Let $U_k'$ be the pullback of 
$\lrang{g_1,g_2}$ in $U_k$. Asking if $p$ divides
the order of
$\hat g_1\hat g_2$ is equivalent to asking if element $(g_1,g_2, (g_1g_2)^{-1})) 
\in \ni(\lrang{g_1,g_2}, (\C_1,\C_2,\C_3))$ is
obstructed in going to $\ni(U_k',(\C_1,\C_2,\C_3))$. We continue this topic in 
\S\ref{genHMreps}. 

\subsubsection{$\sh$-incidence, spin separation and the dominance of H-M reps.}
\S\ref{shinc1} found the data for computing the
$\sh$-incidence matrix for orbit
$O_1^+$ in two stages. Going to higher levels gives further reason to carefully 
organize information
from
$\sH_{k+1}^\rd \to \sH_k^\rd$. It makes sense to combine $\sh$-incidence data 
and the use of
Riemann-Hurwitz as in Lem.~\ref{RHOrbitLem}. 

Prop.~\ref{spinSepLst} finds all degree 120 spin separating representations of 
$G_1$. That gave
subgroups $H_0\le A_5=G_0$ and $H_1\le G_1$ with the natural map
$G_1\to A_5$ inducing $H_1\twoheadrightarrow H_0$. Suppose we have a sequence
$\{H_k\}_{k=0}^\infty$, with each $H_k$ inducing a spin separating 
representation of $G_k$,
$k=0,1,\dots$. Regard the projective limit
$\lim_{\leftarrow k}H_k=\tilde H$ as giving the compatible system of spin 
separating
representations of $\tG p2$. The proof of Prop.~\ref{spinSepLst} gives hints of 
the
following sort. We can expect to find 
pro-2 $\tilde H$ containing a pro-cyclic subgroup $C$ with generator $\alpha\in 
\tG 2$
where 
\begin{triv} $\alpha^{2^{k+1}} \bmod \ker_k\in G_k$ is conjugate (in $G_k$) to 
no other
element of $H_k$. \end{triv}  

Lem.~\ref{tough2-4} and Lem.~\ref{sixCent}
has difficult arguments, respectively about shortening the respective type 
(6,12) and type
(10,20) $q_2$ orbits. These used the precise form of the level 1
characteristic module $M(A_5)$ from Cor.~\ref{A5Morbs}. 
We have yet to capture the essence of $Q''$ shortening
from $p$-group representations. Lem.~\ref{4HM4NHM} says an $H_4$ orbit 
containing
an H-M rep., at level $k\ge 1$ when $p=2$,  has
faithful $\sQ''$ action. 

Note, too, the relation between near H-M and H-M reps.~exploited in such
places as Prop.~\ref{nearHMreps-count}. Further, the notion of complements also 
goes with 
H-M reps. 
Consider any Nielsen classes 
$\{\ni_k=\ni(G_k,\bfC)\}_{k=0}^\infty$ attached to a Modular Tower (for the 
prime $p$). Let
$\{\bg_k\}_{k=0}^\infty$ be a projective system of Nielsen class
representatives. Suppose the conclusion of
Prop.~\ref{nearHMreps-count} holds for this system: 
$\mpr(\bg_{k+1})=p\mpr(\bg_k)$ for $k\ge 0$. If
$p=2$, Def.~\ref{compDef} would give $(\bg_{k+1})q_2^{\mpr(\bg_{k+1})/2}$ as the 
complement of 
$\bg_{k+1}$ in $\ni_{k+1}$. Both elements lie over the same element of $\ni_k$. 
A similar
notion with a general prime
$p$ replacing 2 would give
$p-1$ such elements as 
$(\bg_{k+1})q_2^{\mpr(\bg_{k+1})/p}$, $(\bg_{k+1})q_2^{2\mpr(\bg_{k+1})/p}$, 
\dots,
$(\bg_{k+1})q_2^{(p\nm1)\mpr(\bg_{k+1})/p}$. These are examples of H-M
reps.~indicating a moduli (boundary) behavior that informs about the components 
in which they
appear. 

\subsection{Generalizing the notion of H-M reps.} \label{genHMreps} The 
following problem appears in \S\ref{jumps} in a special
case. To start with assume $\tG p$ is the univeral $p$-Frattini cover of $G_0$, 
with notation as usual for the characteristic
quotients. Let $\C_1,\C_2$ be $p'$ conjugacy classes in $G_0$. Refer to 
$(g_1,g_2)
\in G_0^2$ as {\sl $p$-divisible at level $k$\/} if all lifts
$(\hat g_1,
\hat g_2)$ to $G_k^2\cap (\C_1,\C_2)$  have $p$ dividing the order of $\hat g_1 
\hat
g_2$. Refer to $(g_1,g_2)$ as $p$-divisible if it is $p$-divisible at some level 
$k$. If $\lrang{g_1,g_2}$ is a $p'$ group, then there is a lift of the whole 
group to $G_k$ for each $k$ (special case of
Schur-Zassenhaus). That gives $p\,|\,|\lrang{g_1,g_2}|$ as necessary condition 
for $p$-divisibility. Is it sufficient?  

\begin{prob} \label{p-divisprob} Is $p$ dividing $|\lrang{g_1,g_2}|$ sufficient 
for $(g_1,g_2)$ being $p$-divisible?
\end{prob}

The following says the $p$-divisible property holds for lifting to the universal 
$p$-Frattini cover of $G_0$ if and only if it
holds  for lifting to the universal $p$-Frattini cover of $\lrang{g_1,g_2}$. 
Recall the universal $p$-Frattini cover of a $p'$
group is just the group itself, so the lemma applies even if 
$\lrang{g_1,g_2}=G^{1,2}$ is a $p'$ group. 

\begin{lem}[Local $p'$-lift] That $(g_1,g_2)$ is not $p$-divisible is equivalent 
to their being no lift $(\hat g_1,
\hat g_2)$ of $(g_1,g_2)$ to ${\tG p}^2\cap (\C_1,\C_2)$ with $\hat g_1\hat g_2$ 
having the same order as $g_1g_2$. Suppose
$g_1,g_2\in G_0$ are
$p'$ elements. Let
$G^{1,2}_k$ be the 
$k$th characteristic quotient of the universal $p$-Frattini cover of $G^{1,2}$. 
Then, $(g_1,g_2)$ is $p$-divisible at level $k$ for
lifts to $G_k$ if and only if 
it is $p$-divisible at level $k$ for lifts to $G^{1,2}_k$. \end{lem}

\begin{proof} Suppose at each level $k$ there are lifts $(\hat g_{1,k},\hat 
g_{2,k})$ to $G_k$ as above,   where 
$\hat g_{1,k}, \hat g_{2,k}$ has the same order as $g_1  g_2$. The sets $L_k$ of 
such lifts in $G_k^2$ form a projective system
of nonempty closed subsets. So, there is an element in the projective limit of 
the $A_k\,$s giving the desired lift $(\hat
g_1,\hat g_2)\in {\tG p}^2$. Now consider the relation between lifts to $G_k^2$ 
to that of lifts to $(G^{1,2}_k)^2$.  

Apply Prop.~\ref{pprojchar} to the pullback $\tilde G^*$ of $\lrang{g_1,g_2}$ in 
${}_p\tilde G_0$. Then, $\tilde G^*$ is a
$p$-projective cover of
$\lrang{g_1,g_2}$. So, it maps surjectively to the $p$-Frattini cover of 
$\lrang{g_1,g_2}$.
A lifting $(\hat g_{1,k},\hat g_{2,k})$ to $G_k$, as in the statement where 
$\hat g_{1,k}, \hat g_{2,k}$ has the
same order as $g_1  g_2$, will map to such a lifting in the image. Conversely, 
the universal $p$-Frattini cover of
$G^{1,2}$ has a map back through $\tilde G^*$ over 
$ G_{1,2}$.  Again, a lifting of the type in the lemma to $G^{1,2}_k$ produces 
one in $G_k$.
\end{proof} 

\begin{exmp} Consider $(A_5,\bfC_{5_+^{2}5_-^{2}},p=2)$ at level 0. The two 
components separate 
by lifting invariant values. In the argument for Princ.~\ref{5-5-5princ},  
$\bfC_{5_+5_-3}$ has lifting invariant +1. This was a nontrivial case of 
nonobstruction for
$r=3$. We still, however, don't know if in the $(A_5,\bfC_{5_+5_-3},p=2)$ 
Modular Tower there is 
obstruction at some level. In the language here, is the following pair 2-
divisible: $((1\,2\,3\,4\,5), (1\,2\,5\,4\,3))$? 
\end{exmp}

Reminder of the positive effect of obstructed components: 
They have nothing above them at higher levels, so that eliminates the
component from consideration in the Main Conjecture for any value of $r$. In the 
problem here the
negative effect of the appearance of $\one\,$s is toward providing sufficient 
ramification over cusps in the Main Conjecture. Here
is a projective version of it. 

\begin{prob}[Questions related to Prop.~\ref{serLift}] Let $\tilde \bg$ be a 
projective system of
Nielsen class representatives as in Lem.~\ref{cuspWidth}. Characterize when 
$\jmp(\tilde \bg)=
\infty$.  What generalization of Serre's Prop. could deduce for a projective
system above $\bg_0$ what  might be the values of $\jmp(\tilde \bg)$.  Even if 
the genus
0 condition does not hold in Prop.~\ref{serLift}, is there a general procedure 
for checking the
possible values of the lifting invariant by coalescing branch cycles, as does 
\S\ref{onlytwo},  to the
case of  genus 0?  \end{prob} 

The main effect of Harbater-Mumford representatives on our Main Conjecture has 
been to give lifts of pairs from level 0 to level
$k$ where the order of the product does not change. This has the effect of 
giving projective systems of cusps on a Modular
Tower with relative ramification of degree 1. For this property (depending on 
the outcome of Prob.~\ref{p-divisprob}), we suggest
the natural generalization of H-M reps. is to consider elements
$\bg\in
\ni(G,\bfC)$ with having the following form: 
$$ \bg=(g_1,\dots,g_{s_1},g_{s_1\np 1},\dots, g_{s_1\np s_2},\dots,g_{s_1\np 
s_2\np \cdots \np
s_{t-1}\np 1},\dots, g_{s_1\np s_2\np \cdots \np s_t})$$
\begin{triv} \label{p'property} where, $\lrang{g_{s_1\np s_2\np \cdots \np
s_{i}\np 1},\dots, g_{s_1\np s_2\np \cdots \np s_{i+1}}}$ is a $p'$ group, 
$i=1,\dots,t-1$, and for all $i$, $s_i\ge 2$.  
\end{triv}

\begin{exmp} Consider ${}_7\bg$, ${}_8\bg$ and ${}_9\bg$ from Table 
\ref{lista5}. None are H-M reps. Yet, for $p=5$,
each satisfies the $p'$ property \eqref{p'property}; $(g_1,g_2,g_3,g_4)$ has 
$\lrang{g_1,g_2}$ and $\lrang{g_3,g_4}$ isomorphic to $A_4$ (5 not dividing the 
order). 
\end{exmp} 

H-M reps.~were simple enough to allude to them by just the shape of their 
Nielsen class representives. Some geometry applications
of H-M reps.~require not only that 
\eqref{p'property} holds. Suppose also for each $i$, $(g_{s_1\np s_2\np \cdots 
\np
s_{i}\np 1},\dots, g_{s_1\np s_2\np \cdots \np s_{i+1}})$ has entries with 
product 1, satisfying the genus 0 condition (on each
orbit). Such situations aid in analyzing  
geometry behind components of reduced Hurwitz spaces. 

As an example need for treating Prob.~\ref{p-divisprob}, note the difficulties 
in Lem.~\ref{4HM4NHM} in
establishing that all $H_4$ orbits at level 1 of $(G_1,\bfC_{3^4})$ have 
faithful $\sQ''$ action. This
used details from Cor.~\ref{A5Morbs}. 

\begin{appendix} 
\renewcommand{\labelenumi}{{{\rm (\teql \alph{enumi})}}} 

\section{Coset description of a fundamental domain for
$\sH_k^\rd$} \label{HkH1} \S\ref{fundDomainCusp} suggests finding good cosets to 
go from
fundamental a domain for level 0 of the $(A_5,\bfC_{3^4},p=2)$ to a fundamental 
domain for level 1. 
This includes inspecting the monodromy groups $H_{k,0}$ of the Galois closure of
$\bar\sH(G_{k},\bfC_{3^4})\to
\bar\sH(G_{0},\bfC_{3^4})$ for $k$ large. This will mostly be a
2-group (for a general Modular Tower, mostly a $p$-group).  When $k=1$ the map
has degree 16. The monodromy group, $H_{1,0}$ has these properties.
\begin{edesc} \item  $|H_{1,0}| = 192 = 2^6 \cdot 3$ and its center has order
2.
\item  It has a subnormal series with abelian quotients: 
$T_0\norm T_1\norm T_2\norm H_{1,0}$ with $|T_2|=2^5\cdot 3$, $|T_1|=2^5$ and 
$|T_0|=2$. 
 \end{edesc}  

Note: $H_{1,0}$ is not entirely a 2-group, having a copy of $S_3$ 
at the top. This is from a $q\in \bar M_4$
inducing an automorphism on a 2-Sylow of $G_1$, though $q$ does not act
as an automorphism of $G_1$ through its Nielsen class action. \cite{AGR}
documents this very Galois-like (as in \cite{GalLife}) and subtle 
phenomenon.

\section{Serre's half-canonical invariant} \label{hcan}
For any a projective curve $X$, denote
the divisor class of a divisor $D$ on $X$ by $[D]$. Our concern here is for 
divisors $D$ that are {\sl
half-canonical\/}: $[2D]$ is the canonical class (of divisors of meromorphic 
differentials) on $X$. 
Consider the lifting invariant $s(\bg)$ for a cover $\phi:X\to \prP^1_z$ with a 
branch cycle description $\bg$. 
Prop.~\ref{serLift} discusses the special case of the lifting invariant that 
\S\ref{NielSep} calls
$s_{H,\bfC_H}(\bg)$. For that case, $H\to G$ comes from a spin separating 
representation and $\bfC_H$ is
a lift of $2'$ conjugacy classes.  Particular cases are in Prop.~\ref{expInv} 
and Prop.~\ref{expInv2}. Much of our 
technical work is to establish the former proposition is a case of spin 
separation.  
For such a spin separating representation, under special hypotheses,  
\cite[p.~479]{SeLiftAn} gives a formula
for computing 
$s(\bg)$, and \S\ref{onlytwo} gives examples of loosening the genus 0 hypothesis 
to make this computation.
The formula of  \cite[p.~548 and p.~550]{SeTheta} ties the lifting invariant to 
an invariant of Riemann 
for describing many delicacies about Riemann surfaces. We explain its potential 
to create
automorphic functions, produced by 
$\theta$-nulls, on certain Hurwitz space components. Such functions produced by 
natural moduli properties will reveal 
properties of the moduli space, especially about the nature of the cusps and 
divisors supported in the cusps. 

\subsection{Half-canonical classes from odd ramification} 
Assume the cover 
$\phi: X\to
\prP^1_z$,
$\deg(\phi)=n$ has odd order ramification. At each point $x'\in X$, lying over 
$z'\in \prP^1_z$, express $\phi$
locally as $t_{z'}^e+z'$ with $(z-z')^{1/e}=t_{z'}$ (if $z'=\infty$ use $z^{-
1/e}=t'$) with $e$ odd. As the multiplicity
of zeros and poles of $d\phi$ is independent of the local uniformizing 
paramater, the differential
$d\phi$ has all multiplicities of its zeros and poles even.  Serre treats the 
case with a general curve $Y$ replacing
$\prP^1_z$. Then, the monodromy group of the cover may not be in
$A_n$. That requires extra care
describing the precise group cover
$\hat S_n$ that one uses for computing the lifting invariant. Back to 
$Y=\prP^1_z$. 

Let $D_\phi$ be the divisor
$\frac{(d\phi)}{2}$, the divisor of a half-canonical class. Replacing $\phi$ by 
$\alpha\circ\phi$ with
$\alpha\in \PSL_2(\bC)$ replaces $D$ be a divisor linearly equivalent to $D$. 
The formula is 
$s(\bg)=(-1)^{\ell(D_\phi)}(-1)^{\omega(\phi)}$.  Notations are as follows: 
\begin{edesc} \item $\ell(D)$ is the
dimension of the linear system of the divisor; and  
\item $\omega(\phi)$ is the sum over ramification indices $e_x$ of $\frac{e_x^2-
1}{8}$ of
ramified $x\in X$. \end{edesc} Formulas using ramification indices depend on the 
permutation representation
giving a cover. For either of the Nielsen classes
$\ni(A_5,\bfC_{3^4})^\abs $ or
$\ni(A_4,\bfC_{3_+^23_-^2})^\abs$ (absolute classes here mean for the standard 
representations of $A_n$,
$n=4$ or 5),
$\omega(\phi)$ is even. So, 
$s(\bg)=(-1)^{\ell(D_\phi)}$ (see \S\ref{modtowPullback}). This also holds for 
inner classes in this case. More
significantly, it holds in replacing $G_0=A_5$ (resp.~$A_4$) by $G_k$, the $k$th 
characteristic quotient of the
universal 2-Frattini cover of $G_0$, and the standard representation of $G_0$ by 
some spin separating
representation
$T_k$ of $G_k$. 

For $\bp\in \sH^\inn$ and $(G,T)$ a group with a faithful representation,
$G(T,1)$ the stabilizer of 1 (\S\ref{setupNC}), we have always to keep in mind 
the relation of the Galois cover
$\phi_\bp: \hat X_\bp\to \prP^1_z$ and the cover
$\phi_{\bp,T}:\hat X_{\bp}/G(T,1)\to \prP^1_z$ on the absolute space $\sH^\abs$. 
It is the latter we are attaching
a half-canonical divisor class (and
$\theta$ function below). Still, we suppress $T$ in the continuing discussion. 
Regard $\bp$ as a point of
$\sH(G,\bfC,T)^{\abs,\rd}$. To simplify notation, refer to a cover associated 
with $\bp$ as $X_\bp$. 

\begin{prop}[Even $\theta$ nulls] \label{evenThetas} Each representation 
$T:G_1\to A_{\deg(T)}$ from 
Prop.~\ref{spinSepLst} and Prop.~\ref{spinSepLst2} provides an even $\theta$ 
function
$\theta_{\bp,1}$ on  the cover of degree $\deg(T)$ corresponding to a point 
$\bp\in
\sH_1^+\subset
\sH_1^{\inn,\rd}$. Similarly, there is an even $\theta$ function 
$\theta_{\bp,0}$ corresponding to
the Galois cover $\phi_\bp$ with group $A_5$ attached to $\bp\in 
\sH_0^{\inn,\rd}$ (see
\S\ref{applyA5}). 
\end{prop} 

\subsection{Coordinates from $\theta$ functions} 
The final version of \cite{FrLInv} gives normalizations of these theta functions 
(in their
dependency on $\bp$). That will refer to this exposition for its setup. 
Typically the normalization is to integrate the
$\theta$ function with respect to a volume form, and set it equal to 1. 
Prop.~\ref{evenThetas} implies the geometric
realization of the unramified spin cover of $G$  comes from a $G$ invariant 
nontrivial 2-torsion point $u_\bp$ on
$\Pic(\hat X_\bp)^{(0)}$;  with this description $u_\bp$ is unique. In some 
cases, the ratio 
$|\theta_{\bp,0}(0)/\theta_{\bp,0}(u_\bp)|$ (or its log) appears in a nontrivial 
measure of the distance from an H-M
cusp, so long as that ratio makes sense, the technical point of the duration of 
this section.

Let $\sM_g$ be the moduli space of
curves of genus
$g$. Suppose
$\sH$ is a Hurwitz space of covers $\phi:X\to \prP^1_z$ having genus $g$. Then, 
let $[X]\in \sM_g$
denote the isomorphism class of the curve. A corollary of one of the main 
theorems on moduli of
curves is that the map $\bp\in \sH \mapsto [X_\bp]\in \sM_g$ is a morphism of 
quasi-projective varieties. 

For our applications, consider $\theta_{\bp,0}$: An {\sl even\/} function on the 
universal
covering space $\widetilde \Pic(X_\bp)^{(0)}$ of $\Pic(X_\bp)^{(0)}$, the 
Jacobian (divisor classes of degree 0) of
$X_\bp$. Denote by
$X_\bp^{(t)}$ the positive divisors on $X_\bp$ of degree $t\ge 1$ (different 
than in \S\ref{sequiv} where it was a fiber
product). 
\cite[Lecture III]{CurveJac} guided this, though detecting which bits of this 
construction depend only on $\bp$ is
less obvious.  

We need a basis $\pmb \omega=(\row \omega g)$ for global holomorphic 
differentials on $X_\bp$ to {\sl see\/}
$\theta_{\bp,0}$ in its classical coordinates $\pmb \bw=(\row w g)$. Choose a
positive degree $g$ divisor $D'=\row {x'} g$. Then, consider any
$g$-tuple of paths 
$(\row P g)$ on $X_\bp$ where $P_i$ starts at $x_i'$ (denote its endpoint 
$x_i$),
$i=1,\dots, g$, so $D_{\pmb P}=\sum_{i=1}^g x_i\in X_\bp^{(g)}$ is a positive 
divisor of degree $g$. The Jacobi
inversion theorem says the following about the map
$\mu_{D',\pmb \omega}: \pmb P=(\row P g)\mapsto
[D_{\pmb P}-D']\in \Pic(X_\bp)^{(0)}$.  
\begin{edesc} \item $\mu_{D',\pmb \omega}$ is surjective. 
\item It factors through $\pmb P \to (\sum_{i=1}^g \int_{P_i}\omega_1,\dots,
\sum_{i=1}^g \omega_g)\in \bC_\bw$. 
\item The induced map $\bC_\bw\to \Pic(X_\bp)^{(0)}$ identifies $\bC_{\bw}$ with 
the universal covering space
$\widetilde \Pic(X_\bp)^{(0)}$. 
\end{edesc} This is the classical version of {\sl Jacobi Inversion\/} 
\cite[p.~281]{Springer}. If we
know $\Pic(X_\bp)^{(0)}$ is a projective variety, then this complex analytic map 
$X_\bp^{(g)}\to
\Pic(X_\bp)^{(0)}$ is algebraic (Chow's Lemma \cite[p.~115]{MumRB} or 
\cite[Exer.~4.10, p. 107]{Hart}).
It is a birational morphism, identifying the two function fields. The 
differentials $(\row \omega
g)$ identify with the differentials $(\row {dw} g)=d\bw$, which live on  
$\Pic(X_\bp)^{(0)}$.  

Further, for any integer $t$ and any divisor $D''$ of degree $t$, a similar map 
$$\mu_{D'',\pmb
\omega}: \pmb P=(\row P t)\mapsto
[D_{\pmb P}-D'']\in \Pic(X_\bp)^{(0)}$$ factors through $X_\bp^{(t)}$.  Taking 
$t=1$ and $D''=x_0\in X_\bp$
gives an embedding 
$\mu_{x_0,\pmb
\omega}: X_\bp\to \Pic(X_\bp)^{(0)}$ (assuming $g\ge 1$). Again by Chow's Lemma, 
this gives coordinates on
$X_\bp$ compatible with those of $\Pic(X_\bp)^{(0)}$. Recover the basis $\pmb 
\omega$ as 
restriction of $d\bw$ to  $\mu_{x_0,\pmb\omega}(X_\bp)$. Finally, let $D_0$ be a 
positive divisor of degree
$g-1$. Then, $\Theta_{D_0,\pmb \omega}=\mu_{D_0,\pmb
\omega}(X_\bp^{(g-1)})$ is a divisor on $\Pic(X_\bp)^{(0)}$.  

Forming this data, however, seemed 
to require choices in $\pmb\omega$ and the positive divisors. Under what 
circumstances can one use this
to construct objects analytically varying in the coordinates of $\bp$ (alone 
\wsp not dependent on our choices)? We are
looking for a global construction tightly tying the curves $X_\bp$ in the family 
to coordinates for
$\Pic(X_\bp)^{(0)}$. That is the topic for the rest of this section.

\subsection{$\theta$ functions along the moduli space}
Even the most relaxed conditions won't give explicit $\theta$ functions in 
classical coordinates from a uniform
construction.  Here are reasons. 
\begin{edesc} \label{nonanalReas} \item  \label{nonanalReasa} A nontrivial 
family of curves will have
nontrivial monodromy on the 1st homology (or cohomology) of its fibers. 
\item  \label{nonanalReasb} Analytic coordinates in $\bp$ for an embedding
$\mu_{\bp,t}: X_\bp^{(t)} \to
\Pic(X_\bp)^{(0)}$ come from an analytic choice $[D_{\bp,t}]$ of divisor class 
of degree $t$.
\end{edesc}    From \eql{nonanalReas}{nonanalReasa}, in nontrivial algebraic 
families it is impossible to get a uniform
basis $\pmb\omega_\bp$ of holomorphic differentials. Since $X_\bp^{(t)}$ 
naturally embeds in
$\Pic(X_\bp)^{(t)}$, use $D\in X_\bp^{(t)}\mapsto [D-D_{\bp,t}]\in
\Pic(X_\bp)^{(0)}$ to give analytic coordinates in 
\eql{nonanalReas}{nonanalReasb}. It is especially
significant for
$t=1$ and $t=g-1$. A $[D_{\bp,t}]$ analytic in $\bp\in
\sH$ need not have a representative divisor (positive or not) with coordinates 
in $\bp$.
The divisor class is sufficient. Finding, however,  such an analytically varying
divisor usually requires it appear from the map $\phi_\bp$ through ramification. 
That isn't a formal
definition, though several illustrative situations occur in \cite{DFrVarFam}.  
An example of it for 
\eql{nonanalReas}{nonanalReasb} when $t=g-1$ is having an analytically varying 
half-canonical
divisor class. In turn that produces a precise copy in $\Pic(X_\bp)^{(0)}$ of 
the divisor
$\Theta_{D_{\bp,g-1}}$. 

Let $h$
be any meromorphic function on $X_\bp$  of degree $u$. Then,
$h:X_\bp\to \prP^1_v$  
has zeros $\row {x^0} u$ and poles
$\row {x^\infty} u$. 

\subsubsection{Abel's Theorem motivates finding nondegenerate odd $\theta$-
nulls}
Suppose we have an embedding of $X_\bp$ in $\Pic(X_\bp)^{(0)}=J(X_\bp)$. 
So, each zero $x_i^0$ and pole of $h$ produces a point in $J(X_\bp)$. List
these as $\bw_{x^0_i},\bw_{x^\infty_i}$, $i=1,\dots,u$.
Yet, finding these $\bw\,$s doesn't require
giving $h$. It only needs points $\row {x^0} u$ and $\row {x^\infty}
u$ on $X_\bp$ viewed as inside $J(X_\bp)$. 
Define $$[D_{\bx}]=[D(\bp_{x^0_i},\bw_{x^\infty_i},i=1,\dots,u)]$$ 
to be the sum of all the $\bw_{x^0_i}\,$s minus  
all the $\bw_{x^\infty_i}\,$s on $J(X_\bp)$. To say $[D_\bx]$ is
zero means it is the origin of 
$J(X_\bp)$. 
\begin{thm}[Abel's Theorem (generalized)] Existence of $h$
with these zeros and poles characterizes exactly when
$[D]$ is zero.\end{thm}

If $h$ exists, consider the logarithmic derivative
$dh/h$. This is a meromorphic differential of 3rd kind with pure imaginary
periods. Even if $h$ doesn't exist, given the divisor $D_{\bx}$ above, the
following always holds. 
\begin{trivl}\label{3rdKindDiff} There
is a unique differential $\omega_{\bx}$ with residue divisor $D_{\bx}$
having pure imaginary periods. \end{trivl}
 
Suppose $[D]$ is not zero, but $m[D]$ is zero on $J(X_\bp)$
for some integer $m$. Then, repeating all the zeros and poles
$m$ times produces a function $h$ on $X_\bp$. The $m$th root
of $h$ defines a cyclic unramified degree $m$ cover $Y\to X_\bp$. 

\subsubsection{Abel's Theorem and nondegenerate odd $\theta\,$s}
Riemann produced $\theta=\theta_{X_\bp}$
functions to provide uniformizing coordinates for this
construction. {\sl Torelli space\/} $\sT_g$ is a cover of $\sM_g$. The points of 
$\sT_g$ lying over $[X_\bp]$
consists of a canonical basis for $H_1(X_\bp,\bZ)$ for which we borrow the 
notation of 
$\ba=(\row a g)$ and ${\pmb b}=(\row b g)$ generators from \eqref{poinGens}. See 
\cite[p.~151]{ShAbVar} for the full
abelian variety context of the next statement. 

\begin{lem}[The Symplectic group] \label{Symp} The fiber of $\sT_g$ over 
$[X_\bp]$ is a homogeneous space for
$\Sp_{2g}(\bZ)$ the group $$\Bigl\{T=\smatrix A B C D \mid A,B,C,D\in 
\bM_g(\bZ),\ {}^\tr T\smatrix {0_g}{-I_g}{I_g}{0_g}
T=\smatrix {0_g}{-I_g}{I_g}{0_g}\Bigr\}.$$  \end{lem}

The standard normalization is to choose $\row \omega g=\pmb 
\omega=\pmb\omega_\ba$ so the matrix
$A_{\pmb
\omega_{\ba}}$ with $(i, j)$ entry
$\int_{a_j}\omega_i$ is the identity matrix. The matrix $B_{\pmb \omega_\ba}$  
with $(i, j)$ entry
$\int_{a_j}\omega_i$ is the {\sl period matrix}. Then $B_{\pmb \omega_\ba, \pmb 
b}$ is symmetric and has
positive definite imaginary part (Riemann's bilinear relations: 
\cite[p.~254]{Springer}).  Given this, there is a
unique
$\theta$ function $\theta=\theta_{\bp,\ba,\pmb b}$ attached to $(\ba,\pmb b)$. 

Many mathematical items on $X_\bp$ appear constructively from this.
This includes functions and meromorphic differentials 
(with particular zeros
and poles). This was a central goal in generalizing Abel's
Theorem: To provide Abel(-Jacobi) constructions for a
general Riemann surface. For the function $h$ it has
this look: 
\begin{equation} \label{constFuncts} h(x)=\prod_{i=1}^u
\theta(\int_{x_i^0}^x \pmb \omega)/\prod_{i=1}^u 
\theta(\int_{x_i^\infty}^x \pmb \omega).\end{equation} 
In $\theta$ you see $g$ coordinates; the $i$th entry is $\int_{x_i^0}^x 
\omega_i$, 
an integral over a path on $X_\bp$. 
The integrals make sense up to integration around closed paths. 
So, they define a point in $J(X_\bp)$. 

Even if $h$ doesn't exist, 
the logarithmic differential of \eqref{constFuncts} 
does. It gives the third kind differential from \eqref{3rdKindDiff}. 
This gives the differential equation defining $\theta$ functions. 
In expressing $h$, replace 
$\int_{x_i^\infty}^x \pmb \omega$ by a vector $\bw$ in the universal
covering space of the Jacobian. Form the logarithmic differential of it: 
$d\theta(\bw)/\theta(\bw)$. To construct his $\theta$ functions, Riemann 
following Abel's
case for cubic equations, imposed the condition that translations by periods 
would change the logarithmic
differential only by addition of a constant. With
$\nabla$ the gradient in
$\bw$,   
\begin{trivl} \label{diffinbw} $\nabla(\nabla\theta(\bw)/\theta(\bw))$ is 
invariant under the lattice of
periods. \end{trivl}

The main point about $\theta$ functions (of 1st order for dimension $g$) is that 
there is essentially just one:
$\theta(\bw,B)$ where
$B$ runs over symmetric $g\times g$ matrices with imaginary part positive 
definite. The $\Theta$
divisor determines
$\theta(\bw,B)$ for fixed $B$ {\sl up to a constant multiple}, though we can 
always translate the $\Theta$ divisor on
$\Pic(X_\bp)^{(0)}$. The space of translated divisors is a homogeneous space for 
$\Pic(X_\bp)^{(g-1)}$. 
In varying the $\theta$ function with
$\bp$, our concerns include that choice of constant, and also the analytic 
choice in $\bp$ of $[D_\bp]$ in 
$\Pic(X_\bp)^{(g-1)}$. It depends on our goals (one element of depth in 
Riemann's results) what type of
choice we want for $[D_\bp]$. We illustrate with Riemann's original goal. 

To represent all functions
$h$ appropriately in
\eqref{constFuncts}, the following is necessary: 
\begin{trivl} $\theta$ must be an {\sl odd\/} function so as to give 0 at 
$\bw=\pmb 0$ with  
multiplicity 1. \end{trivl}

\newcommand{\tpi}[5]{{\theta\Bigl[\!\!\begin{array}{c}\raise-
2pt\hbox{{$#1$}}\\\raise2pt\hbox{{$#2$}}
\end{array}\!\!\Bigr](#3,\Pi_{#4,#5})}}

\subsubsection{Fay's normalizations and Riemann's Theorem}
\cite[p.~1]{Fay} makes a choice to replace the matrices $A=I_g$ and $B$ 
respectively by $2\pi iI_g$ and
$\Pi_{\ba,\pmb b}=2\pi iB$ is symmetric and its real part is negative definite. 
This scaling doesn't change the
isomorphism class of
$\bC^g_\bw/L_{\ba,\pmb b}$ and it has numerous advantages. With this notation, 
and all these choices, rewrite Riemann's
classical function 
$\theta(\bw,B)$ as
\begin{equation} \theta(\bw,\Pi_{\ba,\pmb b})= \tpi{\pmb 0}{\pmb 
0}{\bw}{\ba}{\pmb b}=\sum_{\bm\in
\bZ^g}\exp({\scriptscriptstyle\frac 1 2}
\bm \Pi_{\ba,\pmb b} {}^\tr\bm + \bm {}^\tr \bw).\end{equation}

One can see the difference in the formulas by contrasting with the notation of 
\cite[Chap.~2]{FKraTheta}. The
$2\pi i$ factor appears in their expression for a $\theta$ function, and their 
version of the {\sl heat
equation\/}. \cite{FKraTheta} only considers the case
$g=1$.  From our viewpoint, that is because they only consider moduli spaces of 
genus 1 curves where
the $j$-invariant for a point $\bp$ is the $j$-invariant of the curve it 
represents: These are
modular curves. Their $\theta$-nulls are close to  
automorphic functions on these curves, and the relations among the $\theta$-
nulls with various characteristics
are algebraic equations describing these curves. 
\cite[Chap.~2,\S8]{FKraTheta}: 
\begin{quote} In fact, partial motivation for the results discussed so far is to 
better understand the
multivariable case.\end{quote}
Even for $r=4$, Modular Towers forces  relating moduli spaces of curves of
infinitely many genuses parametrized by reduced 1-dimensional Hurwitz spaces.  
When
$r=4$, the moduli spaces are upper half plane quotients. So comparing with 
$\theta$-nulls from 
modular curves is inevitable. Our period matrices
$\Pi_{\ba,\pmb b}=\Pi_{\ba_\bp,{\pmb b}_\bp}$ for the curve $X_\bp$ are 
therefore (matrix) functions of $\tau\in \bH$,
rather than just being
$\tau$. 

Notice: Because a lattice is invariant under multiplication by -1, 
$\theta(\bw,\Pi_{\ba,\pmb b})$ is an even
function of $\bw$. This isn't the function that works in \eqref{constFuncts}. 
Since the divisors for the $\theta$
functions differ only by a translation, giving an appropriate function for 
\eqref{constFuncts} comes by translating
$\bw$. Further, the even $\theta$ from Prop.~\ref{evenThetas}, for $g>1$, is not 
likely
to be this one. As in \S\ref{shincidence} use ${}^t A$ for the transpose of a 
matrix. 

\newcommand{\pdelta}{{\pmb \delta}}  \newcommand{\pepsilon}{{\pmb \epsilon}}
Let $\be=2\pi i {\pepsilon}+{\pmb \pdelta}\Pi_{\ba,\pmb b}$ with 
${\pepsilon},{\pdelta}\in
\bR^g$, represent an arbitrary point of $\bC^g$. A special notation \wsp given 
by characteristics \wsp
expresses the $\theta$ translated by $\be$ (then multiplied by an exponential 
factor):
$\tpi{\pdelta}{\pepsilon}{\bw}{\ba}{\pmb b}=$
\begin{equation} \label{thetaChar} \begin{array}{rl} 
&\exp({\scriptscriptstyle\frac 1 2}
\pdelta \Pi_{\ba,\pmb b} {}^\tr\pdelta + (\bw+2\pi i \pepsilon) {}^\tr 
\pdelta)\tpi{\pmb 0}{\pmb
0}{\bw+\be}{\ba}{\pmb b}\\&=
\sum_{\bm\in
\bZ^g}\exp({\scriptscriptstyle\frac 1 2}
(\bm+\pdelta) \Pi_{\ba,\pmb b} {}^\tr(\bm+\pdelta) + (\bw+2\pi i \pepsilon) 
{}^\tr
(\bm+\pdelta)).\end{array}\end{equation}

There are $2^{2g-1}+2^{g-1}$ choices of $(\pepsilon,\pdelta)\in \frac 1 2\bZ 
\mod \bZ$ giving even $\theta\,$s
($4\pdelta\cdot {}^\tr\pdelta\equiv 0 \mod 2$ \wsp {\sl even characteristics}) 
and $2^{2g-1}-2^{g-1}$ choices of
$(\pepsilon,\pdelta)$ giving odd
$\theta\,$s.  Riemann's Theorem identifies the
$\theta$ from Prop.~\ref{evenThetas} with an even characteristic 
\cite[Thm.~1.1]{Fay}. 

\subsection{Monodromy,  characteristics and Pryms} \label{monhalfcan} 
Again consider $\bp\in \sH^\rd$, and let $P$ be a closed path in $\sH^\rd$ based 
at $\bp$. Suppose 
$\ba_\bp,{\pmb b}_\bp$ is a choice of canonical homology basis  for 
$H_1(X_\bp,\bZ)$. It is classical
that analytic continuation around $P$ induces a linear transformation $\psi_P$ 
on $H_1(X_\bp,\bZ)$:
$(\ba_\bp,{\pmb b}_\bp)\mapsto \psi_P(\ba_\bp,{\pmb b}_\bp)=(\tilde 
\ba_\bp,\tilde {\pmb b}_\bp)$. Using Lem.~\ref{Symp},
this represents
$\psi_P$ as an element $T_P\in \Sp_{2g}(\bZ)$.  For the modular curve case, the 
identification of 
$H_1(X_\bp,\bZ)$ with $\lrang{2\pi i,\tau}$ (Fay's notation) makes this 
transparent because each path identifies with an
element of $\SL_2(\bZ)$. 

\subsubsection{Monodromy acting on a $\theta$} \cite[Chap.~2]{FKraTheta}, 
\cite[p.~7]{Fay} and \cite[p.~179]{ShAbVar} use
a transformation formula that explains how $\theta$ functions change with an 
application 
of $\psi_P$ a homology basis. Here are is a qualitative description. 
\begin{edesc} \item Change $\bw$ to $\tilde \bw$ so analytic continuation around 
$P$ gives the new normalized period
matrix
$\tilde \Pi= \Pi_{\tilde \ba_\bp,\tilde {\pmb b}_\bp}$.
\item Express a translate of $\tilde \bw$ that gives the new $\theta$, a 
function of $(\tilde \bw, \Pi_{\tilde
\ba,\tilde{\pmb b}})$, as an explicit multiple of $ 
\tpi{\pdelta}{\pepsilon}{\bw}{\ba}{\pmb b}$. 
\end{edesc} 
 
Express $\Pi_{\tilde \ba,\tilde{\pmb b}}$ as $(A_P\Pi_{\ba,{\pmb b}}+ 2\pi 
iB_P)( C_P\Pi_{\ba,{\pmb b}}+ 2\pi
iD_P)^{-1}$. Denote  $C_P\Pi_{\ba,{\pmb b}}+ 2\pi
iD_P$ by $M_P=M_{T_P}(\Pi_{\ba,{\pmb b}})=M$.  To simplify, refer to $T_P\in 
\Sp_{2g}(\bZ)$ as $T$. Then,
$^\tr(\tilde
\ba_\bp,\tilde {\pmb b}_\bp)=\smatrix D C B A {}^\tr(\ba,\pmb b)$ and
$2\pi i\bw =\tilde
\bw M$. For a $g\times g$ matrix $U$, use bracket notation $\{U\}$ for 
the vector of diagonal elements (as in \cite{ShAbVar}). The result:  
\begin{equation} \label{thetaAction} \begin{array}{rl} \text{with}
&\Bigl[\begin{array}{c}\raise-2pt\hbox{{$\tilde \pdelta$}}\\ 
\raise2pt\hbox{{$\tilde \pepsilon$}}
\end{array}\Bigr] = \smatrix D {-C} {-B} A
\Bigl[\begin{array}{c}\raise-2pt\hbox{{$\pdelta$}}\\ 
\raise2pt\hbox{{$\pepsilon$}}
\end{array}\Bigr]+\frac 12 \Bigl[\begin{array}{c}\raise-2pt\hbox{{$\{C{}^\tr 
D\}$}}\\ \raise2pt\hbox{{$\{A{}^\tr B\}$}}
\end{array}\Bigr],\ \tpi{\tilde \pdelta}{\tilde \pepsilon}{\tilde \bw}{\tilde 
\ba}{\tilde {\pmb b}}=\\
&K_{\pdelta,\pepsilon,T}\det(M)^{\frac 1 2} \exp(\frac 1 2 (\bw \nabla_\Pi \ln 
(\det(M)) {}^\tr \bw)\tpi{
\pdelta}{\pepsilon}{\bw}{ \ba}{ \pmb b}. \\\end{array}\end{equation} Here
$K_{\pdelta,\pepsilon,T}$, having absolute value 1, depends on which branch of 
$\det(M)^{\frac 1 2}$ we have chosen. When 
$g=1$,
\cite[p.~101]{FKraTheta} has $K_{\pdelta,\pepsilon,T}$ explicitly, recognizing 
$K_{\pmb 0,\pmb 0,T}$ as an 8th root 
of 1. View
$\nabla_\Pi$ as producing the matrix whose
$(i,j)$ entry is the partial with respect to the variable in the $(i,j)$ 
position of $\Pi$. 

\subsubsection{Monodromy on a Hurwitz fiber} Suppose $\sH$ is a Hurwitz family 
(any $r$) for some 
absolute Nielsen class from $(G\le S_n,\bfC)$ (or inner cases). 

\begin{defn} Call $(G\le S_n,\bfC)$ {\sl g(alois)c(losure)-ramified\/} if for 
any cover $\phi:X\to \prP^1_z$ in the
Nielsen class, there is no unramified cover $Y\to X$, of degree exceeding 1, 
fitting in a diagram $\hat X\to Y\to X$
with the Galois closure $\hat X$ of
$\phi$.
\end{defn}

The notion gc-ramified (without the name) is from \cite[Lem.~3.3]{FrCCMD} which 
characterizes it using a
branch cycle description $\bg$ for $\phi$. It is usual that $\phi$ is gc-
ramified (true in our cases), though
nonpathological examples abound where it is not \cite[Ex.~3.4]{FrCCMD}. Let
$H$ be any subgroup of
$G(1)$. Consider the permutation representation $T_H:G\to S_{(G:H)}$ on $H$ 
cosets. There is a map from disjoint cycles
of
$T_H(g)$, $g\in G$, to disjoint cycles of $g$: Any $g$ orbit of $H$ cosets goes 
to the corresponding orbit of $G(1)$
cosets.  Call
$g$ unramified in
$H$ if the lengths of the disjoint cycle of $T_H(g)$ and the corresponding 
disjoint cycle of $g$ are equal. 

\begin{lem}  Let
$H_t(\bg)$ be minimal among $H\le G(1)$ with $g_i$  unramified in $H$, 
$i=1,\dots, r$. Then,
$(G,\bfC)$ is gc-ramified if and only if $H_t(\bg)=G(1)$. \end{lem} In the 
notation above, let
$\bz$ be the branch points of
$\phi$, $\row P r$ a set of classical generators of $\pi_1(U_\bz,z_0)$ 
(\S\ref{classGens}), and  $g$ the genus of
$X$.  By a $\theta$-null we mean a function of the form 
$\tpi{\pdelta}{\pepsilon}{\pmb 0}{ \ba}{ \pmb b}\ne 0$, a
$\theta$ with characteristic evaluated at $\pmb 0$. Let $S^1$ be the 
multiplicative group of complex numbers having
absolute value 1. The expression $f(\bp) \mod S^1$ refers to a function at 
$\bp\in \sH^\rd$ up to multiplication by
elements in
$S^1$. 

\begin{prop} \label{prodThetNull} Assume  $(G\le S_n,\bfC)$ is gc-ramified. Then 
there is an effective procedure for
computing
$H_r$ on
$\pi_1(X)$. Assume $r=4$. This action extends to $\bar M_4$ acting on the 
reduced Nielsen classes. For 
$\bg\in \ni(G,\bfC)^\rd/N$, let $\bar M_\bg$ be the stabilizer in $\bar M_4$ of 
the reduced class of $\bg$. Then, there is
an effective procedure for computing $\bar M_\bg$ acting on
$\pi_1(X)$ (extending the action of \S\ref{j-covers} and Prop.~\ref{j-Line}). 
This induces a map $\phi_*: \bar M_\bg\to
\Sp_{2g}(\bZ)$ through which $\bar M_\bg$ acts as in
\eqref{thetaAction}. 

Suppose elements in the classes $\bfC$ all have odd order and there is an even 
half-canonical class on $X_\bp$
given by $d\phi_\bp$ according to the conclusion of Prop.~\ref{evenThetas}. Let 
$\Bigl[\begin{array}{c}\raise-2pt\hbox{{$
\pdelta_\bp$}}\\ \raise2pt\hbox{{$\pepsilon_\bp$}}
\end{array}\Bigr]$ be the characteristic
corresponding to it (relative to $\ba_\bp,{\pmb b}_\bp$) at the
point $\bp\in
\sH^\rd$.  If
$\tpi{\pdelta_\bp}{\pepsilon_\bp}{\pmb 0}{
\ba_\bp}{
\pmb b_\bp}$ is not 0 at all points
$\bp$, then the action of $\phi_*$ gives a $\theta$-null mod $S^1$, an 
automorphic function on $\sH^\rd$. 
\end{prop} 

\begin{proof}[Comments] The effective procedure for computing the monodromy is 
\cite[\S3.2]{FrCCMD}. For $\gamma,
\gamma'\in \bar M_4$, consider
$T_{\gamma},T_{\gamma'}$ the corresponding elements in $\Sp_{2g}(\bZ)$, so  and 
$M_{\gamma}(\Pi_{\ba,{\pmb b}})$ and
$M_{\gamma'}(\Pi_{\ba,{\pmb b}})$ are the corresponding $M$ factors. Then, 
$$M_{\gamma'\gamma}(\Pi_{\ba,{\pmb
b}})=M_{\gamma'}(\gamma(\Pi_{\ba,{\pmb b}})M_{\gamma}(\Pi_{\ba,{\pmb b}}) \mod 
S^1 \text{\cite[p.~177]{ShAbVar}}.$$ So,
mod $S^1$, the $M\,$s are factors of automorphy. 

From \cite[Thm.~27.2]{ShAbVar}, however, finding factors of automorphy without 
modding out by
$S^1$ improves if we know $T_{\gamma}$ is in the level 2 congruence subgroup of 
$\Sp_{2g}(\bZ)$. Then, the
4th power of the $\theta$-null is an automorphic function.   
\end{proof}

\subsubsection{Automorphic functions and Pryms along $\sH^\rd$} 
\label{autFunctHrd} Again, suppose elements in the
classes
$\bfC$ all have odd order and $d\phi_\bp$ gives an even half-canonical class on 
$X_\bp$ as in the conclusion of
Prop.~\ref{evenThetas}. Let $\hat G\to G$ be the spin cover. Let $\hat 
\bw_\bp\in \bC^{\hat g}$ represent the unique
$G$ invariant nontrivial 2-division point on
$\Pic^{(0)}(\hat X_\bp)$ corresponding to the unramified spin cover $\psi_\bp: 
Y_\bp\to \hat X_\bp$ corresponding to
spin separation. Denote the genus of $\hat X_\bp$ by $\hat g$, global 
holomorphic differentials on $Y_\bp$ by
$\Gamma(\Omega^1(Y_\bp))$. As in
\cite[p.~91]{CurveJac},
$\Pic(Y_\bp)^{(0)}$ is isogenous to
$\Pic(\hat X_\bp)^{(0)}\times \sA_\bp$ with $\sA_\bp$ (the {\sl Prym\/}) an 
abelian variety of dimension $\hat g-1$.
Properties of this construction: 
\begin{edesc} \item  The tangent space at ${\pmb 0}\in \sA_\bp$ identifies with 
the $-1$ eigenspace,
of the involution of $Y_\bp$ commuting with $\psi_\bp$, on 
$\Gamma(\Omega^1(Y_\bp))$. 
\item The isogenies between $\Pic(\hat X_\bp)^{(0)}\times \sA_\bp$ and 
$\Pic(Y_\bp)^{(0)}$ have degree $2^{2g-1}$
in each direction,  composing to multiplication by 2, with $\hat \bw_\bp$ 
generating $\ker(\Pic(\hat
X_\bp)^{(0)}\to
\Pic(Y_\bp)^{(0)})\eqdef\bZ/2(\bp)$. 
\end{edesc} More than just the involution of
$Y_\bp$, all of $\hat G$ acts on this construction. Analysis of this {\sl 
prym\/} varying
canonically with $\bp$ will be in the final version of \cite{FrLInv} (see 
\S\ref{notModular} for the easiest appearance
of a Prym). 

Consider again $\bp\in \sH(A_5,\bfC_{3^4})^{\inn,\rd}$. Is there a nondegenerate  
$\theta$-null
(nonzero and varying with $\bp$) through
$X_\bp$ represented as a point on the moduli space of curves of genus 21? This 
is what Prop.~\ref{prodThetNull} has as
an hypothesis.  This  requires at least that the image $[X_\bp]$ of the curve 
for
$\bp\in
\sH_0^{\inn,\rd}$ vary  nontrivially in the moduli space $\sM_{21}$ as $\bp$ 
varies.
\cite[Main Thm]{FrKK} shows the moduli of all families of inner 
$(A_n,\bfC_{3^r})$
covers, $r\ge n$ as in Ex.~\ref{full3cycleList}, varies nontrivially. 

Our case, however,
is inner covers with $r=n-1=4$, so this does not apply directly. It is more 
elementary.
Let $\sH(G,\bfC)^{\inn,\rd}$ be an inner reduced Hurwitz space with $r\ge 4$, 
where the
genus
$g$ of the covers exceeds 1. 

\begin{lem} The moduli map $\bp\in \sH(G,\bfC)^{\inn,\rd}\mapsto [X_\bp]\in
\sM_g$ is nonconstant.   
\end{lem}

\begin{proof} Restrict the argument to any geometrically connected component of
$\sH(G,\bfC)^{\inn,\rd}$. Assume the result is false and $X=X_\bp$ is the 
constant
image of the moduli map. So, with $\bp$ running over $\sH(G,\bfC)^{\inn,\rd}$, 
there are an infinite number of
inequivalent Galois covers
$\phi_\bp: X
\to \prP^1_z$  with group $G_\bp$ (isomorphic to $G$, though not necessarily 
equal
to it). Consider the group $G^*$ these $G_\bp\,$s generate in the automorphism 
group of $X$. Since
$g(X)\ge 2$, this is a finite group. So there are only finitely many covers up 
to
equivalence in this class. The family is connected: There is only one cover up 
to
equivalence, contradicting the above.  
\end{proof}

\cite[Cor.~3.2]{Fay} provides coordinates for period matrix degeneration as 
points on $\sM_g$ approach its boundary.
The simplest corollary considers being near that part of the boundary with 
period matrices close to
products of period matrices of elliptic curves. All the even $\theta$-nulls are 
nonzero and all the odd
$\theta$-nulls are nondegenerate on an elliptic curve. (\S\ref{A434casejaware} 
therefore shows 
$\sH(A_4,\bfC_{3^4})^{\rd,+}$ has a nonzero even $\theta$-null.) So,  these 
coordinates show the same holds for
jacobians close to such points. For H-M components of Modular Tower levels, 
although the degeneration isn't quite so
simple, this type of analysis also works. This is the Mumford half of the reason 
for calling them Harbater-Mumford
representatives (continuing in
\cite{FrLInv}).  
  
\section{Covers and lifting invariants} \label{NielSep} 
Computations toward Main Conjecture Prob.~\ref{MPMT} \wsp  
showing the genus of components of levels of a reduced Modular Tower go up \wsp 
require 
recognizing the organization of $\gamma_\infty$ orbits within the $\bar M_4$
orbits on the Nielsen class. Every curve cover in a Nielsen class deforms to
a cover associated with collections of cusps. Equivalencing covers according to 
which
cusps are in their deformation range (call this {\sl d(eformation)-
equivalence\/}) is equivalent
to finding $\bar M_4$ orbits. 

Acting with a braid generator on Nielsen classes is quite elementary. Yet, as 
the main examples of this paper
show, finding $\bar M_4$ orbits is usually difficult. Further, even when we
find from $\bar M_4$ action that there are several components, we need to know 
if these
components are conjugate under $G_\bQ$. In significant cases we expect intrinsic 
geometry (of covers in our
case) to inform $G_\bQ$  about d-equivalence. Our examples
show that augmenting bare braid group computations with lifting
invariants often reveals such intrinsic geometry. Two covers certainly have 
differences in their d-equivalence
classes if one has unramified covers with a given group the other does not have. 
While
there are other ways to say they have different geometries, this and its 
generalizations are d-equivalence invariants. The remainder of this section
formulates this purely group theoretically.   

\subsection{The lifting invariant and Nielsen separation} Up to
d-equivalence, assume  two covers (in the same Nielsen class) whose $H_r$ orbits 
we wish to
compare have the same branch points $\bz$ so can compare them as covers by their 
branch cycles using
the same classical generators of $\prP^1_z\setminus \{\bz\}$. If the monodromy 
group of the cover is $G$, let
these branch cycles be $\bg=(\row g r)\in G^r$ and $\bg'=(\row {g'} r)\in G^r$ 
(both in $\ni(G,\bfC)$; 
\S\ref{setupNC}). 

Consider $\psi_H: H\to G$, and a choice of 
$\bfC_H$: $r$ conjugacy classes lifting to $H$ those of
$\bfC$. 

For $\bg\in \ni(G,\bfC)$, consider the $(H,\bfC_H)$-lifting
invariant: $$s_{H,\bfC_H}(\bg)= \{
\Pi(\bg^*)\mid \bg^*=\bg \bmod
\ker(\psi_H),\ \bg^*\in
\bfC_H,\  
\lrang{\bg^*}=H\}.$$ 

Use $\sY_\bg$ for the collection of $s_{H,\bfC_H}(\bg)$
running over all $(H,\bfC_H)$ to separate d-inequivalent elements in 
$\ni(G,\bfC)$. Two cases appear often in this paper.
\begin{edesc}\label{exLiftInv} \item \label{exLiftInva} $G\le A_n$ for some 
integer $n$ and
$H=\hat G$ is the pullback of
$G$ to $\Spin_n$ in Prop.~\ref{serLift}.  \item \label{exLiftInvb} $G=G_k$ is 
the $k$
characteristic quotient of the universal
$p$-Frattini cover of the group $G_0$, and $H=R_k'$ in Prop.~\ref{RkGk}. 
\end{edesc}  

When $p=2$, and $G_0=A_n$, we may compare \eql{exLiftInv}{exLiftInva} with
\eql{exLiftInv}{exLiftInvb}, and ask if the latter is a case of the former 
through some
embedding of $G_k$ in an alternating group. That is exactly what happens in our 
main
example at level 1 (Prop.~\ref{expInv}). The following easy lemma slightly 
generalizes
\cite[Lem.~3.12]{FrMT}

\begin{lem} \label{HrInv} $\sY_\bg$ is an
$H_r$  invariant:
$$s_{H,\bfC_H}(\bg)=s_{H,\bfC_H}((\bg)Q) \text{ for }Q\in H_r,$$ any
$(H,\bfC_H)$. It is a d-equivalence invariant. \end{lem} 

\begin{defn} Call $\bg$ and $\bg'$ in $\ni(G,\bfC)$ {\sl Nielsen
separated\/} if the collection $\sY_\bg$ differs from
$\sY_{\bg'}$. \end{defn} 

\begin{quest} If $\bg$ and $\bg'$ are in separate d-equivalence classes, are 
$\sY_\bg$
and $\sY_\bg'$ provably (computationally, significantly)  different so as to 
detect this?
\end{quest} 

\subsection{Frattini central extensions} 
The case of Nielsen separation from  this paper is where 
$\psi: H\to G$ is a Frattini {\sl central\/} extension with $\ker(\psi)$
of order prime to the orders of elements in $\bfC$. Construct $\bfC_H$ as the
conjugacy classes $\bfC$ of $G$, each lifted (uniquely) to $H$ to have the same
order. The next lemma follows easily from the technique of \cite[App]{FrVMS}. 
 
\begin{lem} \label{CenFratLift} For $\bg\in \ni(G,\bfC)$, 
$s_{H,\bfC_H}(\bg)=s(\bg)$ is a single element. 
Assume $\ni(G,\bfC)\ne \emptyset $ and 
\begin{trivl} conjugacy classes
in $\bfC$ appear
many times.\end{trivl} Then
$\{s(\bg)\mid \bg\in \ni(G,\bfC)\}=\ker(\psi)$. 
So, there are at least $|\ker(\psi)|$ orbits
for $H_r$ on $\ni(G,\bfC)$. \end{lem}

General Nielsen separation uses complicated 
witnessing pairs $(H,\bfC_H)$. In the last example, the $(H,\bfC_H)$
witnessing separation is a small cover of $(G,\bfC)$. 
Lem.~\ref{CenFratLift} shows that embeddings in alternating groups are only a
tiny part of the story of general Nielsen separation of orbits of $H_r$. 

For instance, consider the universal $p$-Frattini cover ${}_p\tilde G$ of any 
simple group $G$ ($p\,|\,|G|$;
Thm.~\ref{densityOnes}) or
$p$-split,
$p$-perfect group with non-cyclic $p$-Sylow (Prop.~\ref{psplitOne}). Then, 
${}_p\tilde G$ has infinitely many
quotients $G'$ with
$H'\to G'$ a central (Frattini) extension with
$\ker(H'\to G')=\bZ/p$. In any of these circumstances, consider such a group 
$G'$ and 
corresponding $p'$ conjugacy classes $\bfC'$ satisfying the hypotheses of 
Lem.~\ref{CenFratLift}. Then, the Hurwitz
space $\sH(G',\bfC')$ has at least $p$ components. For,
however, a fixed $G$, $p$, and $r$ (say, $r=4$) consider for  which $r$-tuples 
of
$p'$  conjugacy classes
$\bfC$, does the collection of covers
$H'\to G'$ Nielsen separate $H_r$ (equivalently $\bar M_4$) orbits on 
$\ni(G',\bfC)$ (with $G'$ running over the
mentioned quotients of $\tG p$). 

\section{The Open Image Theorem} \label{openIm} Use the notation of
Prop.~\ref{GKaction}.  Figuring  the image of 
$G_K$ (up to subgroups of finite index) in $\Hom^*(\pi_1^{(p)},\ker_0)/\ker_0$ 
or in 
$\Hom^*(\bT_p,\ker_0)/(\ker_0,\ker_0))$ generalizes 
the goal of Serre's Theorem in \cite{SeAbell-adic}. One 
goal would be to decide (find situations) when this action has an open orbit. 

\subsection{Interpreting and extending the Open Image Theorem} From the
Hurwitz  viewpoint, 
\cite{SeAbell-adic} is the case 
$G_0=D_p$, where $X_{\bp_0}$ is an elliptic curve. It 
says, the image  of  
$G_K$ is an open subgroup of the general linear group acting on the Tate module 
of $X_{\bp_0}$, if 
$X_{\bp_0}$ has no 
complex multiplications. Otherwise, its image is an open subgroup of the 
Tate module automorphisms  commuting with the ring of endomorphisms. 

Much of \cite{SeAbell-adic} describes the 
complex multiplication case. Still, the guiding theorem is 
\cite[IV-20]{SeAbell-adic}: 
If $j(E)$ is not an algebraic integer (at $p$), then the image of $G_K$ in the 
Tate
module of
$E$ is the full general linear group. Further, in its closing pages, 
\cite{SeAbell-adic}
outlines the case when $j(E)$ is integral and $E$ has no complex multiplication.
Prop.~\ref{tGppts} has a framework for considering these results applied to a 
Modular
Tower where direct analogs  would start with one prime $\ell$ (with $\ell$ not
dividing $G$). We'll use 
$\bp_0\in \sH(A_5,\bfC_{3^4})^\inn(\bar \bQ)$ for explicit comments. 

Since $\sH(A_5,\bfC_{3^4})^{\inn,\rd}$ covers the $j$-line, interpret $\bp_0\in
\sH(A_5,\bfC_{3^4})^\inn(\bar \bQ)$ as being $\ell$-adically close to
$\infty$ if its $j$-invariant is sufficiently negative for the $\ell$-adic 
valuation. A first goal
is to find the largest (up to equivalence according to the Lie algebra of the
image) possible action of
$G_K$ on projective systems of points lying over $\bp_0\in 
\sH(A_5,\bfC_{3^4})^\inn(K)$. Finding how to estimate
how
$\ell$-adically close to $\infty$ would be part of this.  

We make simplifying assumptions to
propose a program. 
We haven't yet found a good analog on a {\sl general\/}
($p$-perfect) Modular Tower. There are two parts.
\S\ref{outMon} reminds of the place where a $p$-Frattini covering property 
enters in Serre's Open Image Theorem.
\S\ref{findFrob} gives a precise test case for starting from Serre's situation 
using $j$ values $\ell$-adically close
to $\infty$ for some prime $\ell$. 

\newcommand{\HMP}{\text{\rm HMP}}

\subsection{$p$-part of the outer monodromy group} \label{outMon} Let 
$\HMP(A_5,\bfC_{3^4})$ be the subset of projective systems 
$\{\bp_k\}_{k=0}^\infty$ (points) on
the
$(A_5,\bfC_{3^4},2)$ Modular Tower lying over $\bp_0\in 
\sH(A_5,\bfC_{3^4})^{\inn,\rd}$ (so
they correspond to elements of $\ni(\tG p,\bfC)^{\inn}$ (Prop.~\ref{tGppts}) 
with this additional
property:  
\begin{trivl} Each $\bp_k$ lies in a component at level $k$ containing H-M 
reps.\end{trivl} An
absolute Galois group preserves the collection of H-M components \cite[proof 
Lem.~3.22]{FrMT}. So, it is
a special case of the problem to describe the H-M components at every level (see 
Rem.~\ref{countHMcomp}).
 
Let $K$ be a number field. Consider an arbitrary reduced Modular Tower (possibly 
inner, or absolute)
$\{\sH^\rd_k\}_{k=0}^\infty$. Apply the notation of \S\ref{HkH1} to a projective 
system of $K$ components
$\{\sH'_k\}_{k=0}^\infty$. Specifically, let $H_{k,j}'$ be the (geometric) 
monodromy group of the map
$\sH'_k\to
\sH'_j$ with $H_k'$ the (geometric) monodromy group from $\sH'_k\to 
\prP^1_j\setminus\{\infty\}$. 
The \cite{SeAbell-adic} case is  $\cdots \to Y_0(p^{k+1})\to Y_0(p^k)\to \cdots 
\to
\prP^1_j\setminus\{\infty\}$, the arithmetic monodromy groups are 
$\GL_2(p^{k+1})$ while the geometric
monodromy groups are $\SL_2(p^{k+1})$ \cite[p.~47]{Se-GT}. 
Lem.~\ref{basicA5cover} computes the geometric and
arithmetic monodromy groups for level 0 of the $(A_5,\bfC_{3^4},p=2)$ Modular 
Tower. The modular curve
tower has $p$-groups for $H_{k,0}$.  \S\ref{HkH1} notes that for the 
$(A_5,\bfC_{3^4},p=2)$ Modular
Tower, $H_{1,0}'$ (for either component) is not a 2-group as its order is 
divisible by 1 power of 3.

A basic lemma is that
$\PSL_2(\bZ_p)\to \PSL_2(\bZ/p)$ is a Frattini cover if $p\ge 5$ \cite[IV-
23]{SeAbell-adic}. For $p=3$: $\PSL_2(3)$
is $A_4$ \wsp acting on the 4 points of the projective plane. So, it is not 3-
perfect and its universal 3-Frattini cover
is
$K_4\xs \bZ_3$. 

For $p\ge 5$, the result is a very strong Hilbert's Irreducibility Theorem for 
this situation hidden in
\cite{SeAbell-adic}. It starts by assuming
$\bp_0\in Y_0(p)$, and the Galois closure of $K(\bp)/K$ contains $\SL_2(\bZ/p)$.  
The conclusion: For
any $k\ge 0$ and $\bp_k\in Y_0(p^{k+1})$ over $\bp_0$, the Galois closure of 
$K(\bp_{k+1})/K$ equals the full
arithmetic monodromy group of  $Y_0(p^{k+1})\to \prP^1_j\setminus\{\infty\}$. 
This is the strongest conclusion
supporting the open image theorem. Considering  Modular Towers 
requires addressing the next problem, especially in the $p$-split case (the 
topic of \cite{AGR}). 

\begin{prob} \label{pFrattiniTower} Assume $G_0$ is $p$-perfect. Characterize 
when for some $j\ge 1$, $\lim_{\infty
\leftarrow k} H_k'$ is a $p$-Frattini cover of $H_j'$?  \end{prob}

To understand a general Modular Tower requires
considering geometric components other than H-M types  that could provide
projective (unobstructed) systems of components. Here is an example. Suppose  
a collection of conjugacy classes  $\bfC$ satisfies the following
properties. 
\begin{edesc} \item $\ni(G_0,\bfC)$ contains no H-M reps (say, if $r$ is odd).
\item $\ni(G_k,\bfC)$ is nonempty for all $k\ge 0$.
\end{edesc} 
\begin{exmp}[Geometry of other than H-M components] 
Example: Suppose $r=3t$ with $t$ odd and with all conjugacy classes equal to 
$\C$, a rational $p'$ conjugacy
class of an element of order
$t$. Assume also: $G_0=\lrang{g_1,g_2,g_3}$ with $g_i\in \C$, $i=1,2,3$. If 
$\tilde g_i$ is a lift of $g_i$
to $\C$ in $G_k$, then $$(\underbrace{\tilde g_1,\dots, \tilde g_1}_{t \rm\; 
times},\underbrace{\tilde
g_2,\dots,
\tilde g_2}_{t  \rm\; times},
\underbrace{\tilde g_3,\dots,
\tilde g_3}_{t  \rm\; times})\in \ni(G_k,\bfC).$$
\end{exmp} 

\begin{rem}[How many H-M components at each level?] \label{countHMcomp} It would 
help if there
were just one H-M component at each level of a Modular Tower. For a general 
Modular Tower a 
simple assumption guarantees this \cite[Th.~3.21]{FrMT}. This assumption, 
however, holds in none
of the cases when
$r=4$ in this paper. The conclusion can even be false: Ex.~\ref{A4C34-wstory} 
has two H-M components at level 1.
\cite{FrLInv} gives an explicit $r_0=r_0(n)$
for the Modular Towers associated to
$(A_n,\bfC_{3^r},2)$ so that each level of the Modular Tower has exactly one H-M 
component (though
possibly many components if the level $k$ is large) if
$r\ge r_0$. If $r=n-1$ or $n$, however, the number of H-M components is still a 
mystery. \end{rem}

\subsection{Finding the Frobenius} \label{findFrob} 

Consider $\prP^4(\overline{\bQ((t))})\setminus D)\times \prP^1_z$ over the field
$\bQ((t))=K$ of Laurent series. Regard $\PSL_2(\overline{K})$ as acting as in 
\S\ref{SL2-act} (except extending
to the
$\prP^1_z$ component). Let $d$ be a squarefree integer. Start with a tangential 
base point
convenient for such an $\ell$-adic investigation. A convenient example is that 
from
\cite[\S2]{WeFM}: 
$$\Spec(\bQ((t))) \to 
\prP^4(\overline{\bQ((t))})\setminus D)\times 
\prP^1_z/\PSL_2(\overline{\bQ(t)})$$ by using the unordered four
points $S=\{\pm \sqrt{d}\pm i \sqrt { t} \}=\{\row z 4\}$ in the first position  
and $z_0=0$ (a fixed base point) in the
last position. To extend the map to $\Spec(\bQ[[t]])=\Spec(R)$ requires forming 
a special fiber in place of
$\prP^1_z$ over
$j=\infty$. Extending the {\sl Stable Compactification Theorem\/} appropriately 
to
covers replaces $\prP^1_z$ by three copies $\prP^1_z$, $\prP^1_{w_5}$ and 
$\prP^1_{w_6}$ tied
together as on the left side of Table 5.

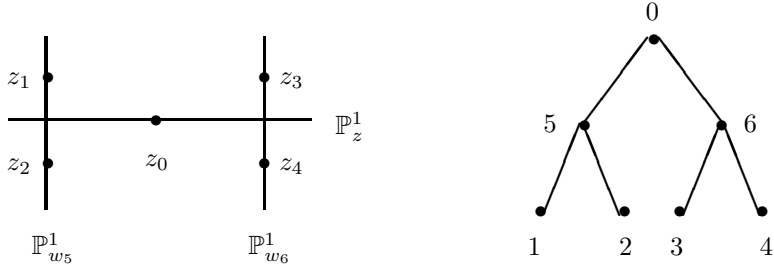
\begin{table}[h]  \label{stabComp}
\caption{The configuration of $\prP^1\,$s at a special fiber} 
\setlength{\unitlength}{0.0005in}
\thicklines
\begin{picture}(8068,2520)(0,-10)
\put(1670,1377){$\bullet$}
\put(545,1827){$\bullet$}
\put(545,927){$\bullet$}
\put(2795,1827){$\bullet$}
\put(2795,927){$\bullet$}
\put(6845,2227){$\bullet$}
\put(6125,1327){$\bullet$}
\put(7550,1327){$\bullet$}
\put(5670,427){$\bullet$}
\put(6540,427){$\bullet$}
\put(7115,427){$\bullet$}
\put(7975,427){$\bullet$}
\put(200,1460){\line(1,0){3150}}
\put(590,520){\line(0,1){1800}} 
\put(2865,520){\line(0,1){1800}}
\put(5765,475){\line(2,5){375}}
\put(6555,465){\line(-2,5){360}}
\put(7220,465){\line(2,5){360}}
\put(8015,465){\line(-2,5){360}}
\put(6165,1400){\line(3,4){680}}
\put(7645,1400){\line(-3,4){680}}

\put(185,1827){$z_1$}
\put(185,927){$z_2$}
\put(3020,1827){$z_3$}
\put(3020,927){$z_4$}
\put(1620,972){$z_0$}
\put(445,57){$\prP^1_{w_5}$}
\put(2695,57){$\prP^1_{w_6}$}
\put(3605,1307){$\prP^1_{z}$}
\put(5765,1327){$5$}
\put(7850,1327){$6$}
\put(5605,50){$1$}
\put(6555,50){$2$}
\put(7085,50){$3$}
\put(8015,50){$4$}
\put(6830,2497){$0$}
\end{picture}
\end{table}

Each of the two pairs of points coming together (complex
conjugate pairs) is wandering in $t$ space. Add a line of directions for the 
approach 
of each pair. Each pair degenerates to the common limit of the two points with 
the addition of the
direction of their coming together in the copy of $\prP^1_z$. The downward 
pointing tree on the right in Table 5
has root 0 labeled for the component $\prP^1_z$, one edge each for 
$\prP_{w_5}^1$ and 
$\prP_{w_6}^1$ (with corresponding vertices 5 and 6) and four {\sl leaf\/} edges 
terminating in vertices corresponding
to the points $\row z 4$.  

Computing with this amounts to finding explicit action of
$G_{\bQ((t))}$ on the elements of the Nielsen classes 
$\ni(G_k,\bfC_{3^4})^{\inn,\rd}$ representing germs of covers
approaching the special fiber.   For example, in the complex conjugate pairs 
business above, identify each
component with a choice of coordinate
$w$. That identifies the function field of the copy of $\prP^1_{K}$ (modulo the 
action of $\PGL_2(R)$). For the root
component use $z$. For $\prP^1_{w_5}$, the component from  $z_1$ and $z_2$ 
coming together, choose
$w_5=(z+\sqrt{d})/\sqrt t$. For $\prP^1_{w_6}$ choose $w_6=(z-\sqrt{d})/\sqrt 
t$. Then, $$z_1=-\sqrt d+i\sqrt t,\
z_2=-\sqrt d-i\sqrt t,\ z_3=\sqrt d+i\sqrt t,\ z_4=\sqrt d-i\sqrt t.$$ As $t 
\mapsto 0$, $z_1$ and $z_2$ have
limit $-\sqrt d$ on $\prP^1_z$. On, however,  $\prP^1_{w_5}$, the $w_5$ values 
representing the limit of $z_1$ and
$z_2$ are
$i$ and
$-i$.  Similarly for
$z_3$ and
$z_4$ coming together, given by
$\sqrt d$ on $\prP^1_z$ and $w_6$ values $i$ and $-i$. 

Let $q_0$ be the generator of $G(\overline{\bQ((t))}/\bar \bQ((t)))$. The 
coordinates on each component give the $q_0$
action on the tree.  Extend it to the fundamental group of the tree. The 
assumptions for this explicit
action are in \cite{WeFM}. First: Identify $t^{1/n} \mapsto \zeta_nt^{1/n}$ as a  
braid by
recognizing its action on the tree. As in our example, this is from writing 
$\row z r \in \bC(t^{1/n})$ explicitly as
$z_i(t^{1/n})$.  You get the image of the braid in $S_r$ just from its action on 
the power series for $\row z r$.
Vertices $v$ of the tree consist of $S\cup\{z_0\}$ together with the components 
of the special fiber. Edges are
from a point of $S$ lying on a component, or from the meeting of two components. 
Suppose an edge
$e=(v_1,v_2)$ represents the component meeting at a singular point
$x_e\in X'$ (like
$\prP^1_{w_5}$ and
$\prP^1_z$ meeting). Then, the local ring around that point is
$\sO_{R,x_e} = R[[u_1,u_2]]/(u_{e,1}u_{e,2}=t^{n_e})$ for some well-determined 
$n_e$ where $u_{e,i }=u_i= 0$ defines the
two components, $i=1,2$. Fix $t_0<\epsilon$, and form a neighborhood $V_e$ using 
the coordinates  $(u_{e,1},u_{e,2},t)$
to include the arc $$c_e: s\in [0,1] \mapsto (t_0^{n_e/2n}e^{2\pi is},\ 
 t_0^{n_e/2n}e^{-2\pi is},\ t_0^{1/n}).$$ Let $U_{t_0}=\prP^1_z\setminus 
\{{z_1(t_0^{1/n})},\dots,
{z_r(t_0^{1/n})}\}$.

For a leaf edge $e$  take the clockwise path around $z_i$. 
For $v\in V(T)\setminus \{v_0\}$, set $e=(\Re(v),v)\in E(T)$ and let $U_v$ be 
the connected component of
$U_{t_0}\setminus (\cup_{e'} c_{e'})$ containing the annulus 
$$\{(u_{e,1},u_{e,2},t)\in V_e \mid |u_{1,e}|>
|u_{2,e}|\}.$$ This is the topological realization of the ordered tree. 

Reminder of those details: \cite[\S2.3]{WeFM} produces the presentation 
compatible with a given graph. 
For $v\in V(T)\setminus\{v_0\}$, let the unique path $p_v$ connecting $v_0$ to 
$v$ be denoted 
$(v_0,v_1,\dots,v_{k\nm1},v)$. Let
$e=(v_{k\nm1},v)$,
$z_v=c_e(0)$ and
$z_{v_0}=z_0$. Choose a simple path $a_v:[0,1] \to \bar U_v$ leading from 
$z_{v_{k\nm1}}$ to $z_v$. Let
$b_v=a_{v_1}\cdots a_v$, and 
$\gamma_v = b_vc_e b_v^{-1}$ (an element of $\pi_1(U,z_0)$). 

As usual, use classical generators $\row \gamma 4$, clockwise around the $\row z 
4$ from
$z_0$. Add these paths: 
$\gamma_5=\gamma_1\gamma_2, \gamma_6=\gamma_3\gamma_4$. We act on the opposite 
side of Wewers, following the notation
of \S\ref{groth}. So,
$q_0$ acts as 
$$ \gamma_1^{q_0}=\gamma_1^{-1}\gamma_2\gamma_1, \gamma_2^{q_0}=\gamma_1, 
\gamma_3^{q_0}=\gamma_3^{-1}\gamma_4\gamma_3, \gamma_4=\gamma_3^{q_0}.$$ Then, 
with $\Pi_5=\lrang{\gamma_1,\gamma_2}$
and $\Pi_6=\lrang{\gamma_3,\gamma_4}$, and $\sigma \in G_\bQ$, 
$$ \gamma_i{}^\sigma=\beta_{(i)\sigma}^{-
1}\gamma_{(i)\sigma}^{\chi(\sigma)}\beta_{(i)\sigma} \text{ with } 
{\beta_i} \in \Pi_5, i=1,2, \text{ and } {\beta_i} \in {\Pi_6}, {i=3,4}.$$ The 
action of $\sigma$ on the
subscripts of the $\beta_i\,$s is the obvious one from our the action of 
$\sigma$ on the Puiseux expression of
the paths in $S$. Further, with
$\epsilon: G_\bQ\to
\{\pm 1\}$,
$\gamma_i^\sigma=\gamma_i^{\epsilon(\sigma)\chi(\sigma)}$, $i=5,6$. 

\cite[Lem.~2.19]{WeFM} shows that if $\ell$ is prime to $|G|$ and the branch 
locus doesn't degenerate at $\ell$, we
may take $t_0$ in the maximal ideal of the valuation ring for $\bQ_\ell$ and 
specialize all the actions. 
There would be two steps in the
program corresponding to two parts of Serre's result. Still, we restrict to the 
case of H-M components for the
$(A_5,\bfC_{3^4},p=2)$ Modular Tower. 
\begin{edesc} \label{concOpenImage} \item  \label{concOpenImagea} For $a=t_0\in 
\bar\bQ$, giving $j(t_0)$ suitably
$\ell$-adically close to
$\infty$ draw deductions about the action of $\bQ_\ell$ on $\bp\in 
\sH(G_k,\bfC_{3^4})^{\inn,\rd}$ for $\bp$ lying over
$j(a)$ in an H-M component. \item  \label{concOpenImageb} Estimate from this the 
$\ell$-adic domain in the $j$-line to
which the one may analytically continue from $j(a)$ the Lie algebra action of 
$G_{\bQ_\ell}$ on the $\ell$-adic
Grassmanian of
\S\ref{grassman}. 
\end{edesc}   Toward \eql{concOpenImage}{concOpenImagea}, an accomplishment 
would be to differentiate between the action
of $G_{\bQ_\ell}$ on H-M reps.~ and on near H-M reps., exactly as we did with 
$G_\bR$ at all levels of the Modular
Tower in Prop.~\ref{HMnearHMLevel}, and Wewers did $\ell$-adically at level 0 of 
the same Modular Tower.

\section{Simple branching and $j$-awareness} \label{jawareMT} Two situations 
about a Modular Tower require 
extending considerations of this paper. The most classical situation doesn't 
satisfy 
$p$-perfectness. \S\ref{non-p-perfect} discusses a reasonable remedy.  Further, 
given a reduced Modular
Tower $\{\sH_k^\rd\}_{k=0}^\infty$, some significant part of the tower may come 
from a tower of
modular curves. While monodromy groups of the levels over the $j$-line often 
preclude this,
it is difficult to compute the monodromy groups even should that be
sufficient to set the matter straight. \S\ref{modtowPullback} illustrates by 
example the problems. 

\subsection{More complicated Hurwitz spaces for non-$p$-perfect $G_0$} 
\label{non-p-perfect} The
pairing of dihedral groups and modular curves generalizes neatly into Modular 
Towers to include, for example,
simple groups and any prime that divides the order of such a group. Still, the 
most classical case of moduli of
covers does not: {\sl simple branching\/} 
$(S_n,\bfC_{2^r}, p=2)$ with  
$\bfC_{2^r}$ being $r\ge 2(n-1)$ involutions in $S_n$ with $n\ge 4$. Since
$\phi: S_n\twoheadrightarrow
\bZ/2$ has kernel $A_n$, $S_n$ is not 2-perfect, and 2 divides the orders of the 
conjugacy
classes. \cite{AGR} explains this situation in more generality (as in 
Ex.~\ref{SnAnbcycles}).  

The remedy, however, is clear. Regard  each level 0 curve in the reduced Hurwitz 
space as
coming with a map to a hyperelliptic curve. Then view the (ramified) $A_n$ cover 
of hyperelliptic curves
(of genus $g=\frac{r-2(n-1)}2$) as the starting point for unramified projective 
systems of
covers coming from the universal 2-Frattini cover of $A_n$. 
Prop.~\ref{serLift} applies to this situation extending its application to 
Modular Towers starting with covers of
$\prP^1_z$. \cite{AGR} illustrates this with $n=5$ using just the group theory 
from this paper.

\subsection{Modular curve quotients from a Modular Tower} \label{modtowPullback} 
Continue 
with $r=4$ as at the beginning of this section.  Ex.~\ref{SnAnbcycles} gives 
situations automatically implying
the jacobian of
$\hat X_\bp$,  
$\bp\in
\sH^\rd$ maps through the elliptic curve with $j$-invariant $j(\bp)$. (This is  
equivalent to $\hat  X_\bp$ itself
mapping onto this elliptic curve.) This is the situation we call $j$-aware. Even 
if $\hat X_\bp$ maps through an
elliptic curve, it may not be the elliptic curve with $j$-invariant $j(\bp)$. 
For a general Modular Tower the jacobians
of the corresponding curves at level $k$ go up in rank very quickly. Given the 
rank
$\rk_0$ of
$\ker_0$, inductively compute the rank of $\ker_{k+1}$ as $\rk_{k+1}=1+(\rk_k-
1)p^{\rk_k}$
from the Shreier formula. Given the genus $g_0$ of $\hat X_0$ (a curve at level 
0), computed from the branch
cycles as a cover of $\prP^1_z$, then $g_{k+1}-1=p^{\rk_k}(g_k-1)$. For 
$(A_5,\bfC_{3^4},p=2)$, $g_0=21$ and
$\rk_0=5$, so $g_1=1+5\cdot 2^7, \rk_1=1+2^7$, $g_2=1+5\cdot 
2^{14},\rk_2=1+2^{1+7+2^7}$,
\dots. 

\subsubsection{Extending Ex.~\ref{A4C34-wstory}} \label{extA4C34} 
\S\ref{hcan} continues the Nielsen class 
$\ni(A_4,\bfC_{3_+^23_-^2})$ from Ex.~\ref{A4C34-wstory}. 
We show level 0 is not a modular curve.  

The natural map from modding out by the Klein 4-group in $A_4$ gives  $A_4 \to  
A_3$. Regard
$C_4=(\pm1)^4/\lrang{(-1,-1,-1,-1)}$ as a multiplicative group of order 8. Any 
element of
$\ni(A_4,\bfC_{3_+^23_-^2})$ (the level 0 Nielsen class) gives a symbol $(\row 
\epsilon 4)\in
C_4$ where two entries are +1 and two entries are -1. Let
$\sH_0^{\rd,+}$ be the H-M component of this reduced Hurwitz space at level 0 
(as in 
Ex.~\ref{full3cycleList}). The absolute and inner (nonreduced) spaces, 
$\sH_0^{\abs,+}$ and $\sH_0^{\inn,+}$, at level 0, have a natural degree 2 map 
$\psi^{\inn,\abs}: \sH_0^{\inn,+}\to
\sH_0^{\abs,+}$. The element $\sh^2q_1q_3^{-1}\in \sQ''$ takes the 4-tuple $( 
(1\, 2\, 3), (1\, 3\, 2), (1\, 3\,
4), (1\, 4\, 3) )$ to the 4-tuple $( (1\, 4\, 3), (1\, 3\, 4), (1\, 3\, 2), (1\, 
2\, 3) )$. The latter is
conjugation of the former by $(2\,4)$. This is exactly why $\psi^{\inn,\abs}$ 
has degree 2 (illustrating
Thm.~\ref{FrVMS}).  Then, 
$\sH_0^{\abs,+}$ is a family of genus 1 curves, while the curves of
$\sH_0^{\inn,+}$ have genus $g_{\inn,0}$ with $2(12+g_{\inn,0}-1)=4\cdot 4\cdot 
2$: $g_{\inn,0}=5$. 
Since, however, we mod out by $\sQ''$ to get reduced Nielsen classes, this 
induces an
isomorphism of
$\sH_0^{\inn,\rd,+}$ with $\sH_0^{\abs,\rd,+}$. As in \S\ref{fineModDil} and
Prop.~\ref{HMnearHMLevel}, the spaces still represent different moduli problems. 
Neither reduced space is
a fine (or even b-fine) moduli space according to Prop.~\ref{redHurFM}. Still, 
to simplify, refer to both as
$\sH_0^{\rd,+}$. 

Choose a basepoint $j_0\in (1,+\infty)$ from $\prP^1_j(\bR)$ (as in 
Lem.~\ref{locj}). Then, 
$\bar M_4$ acts on allowable symbols from $C_4$, giving a cover $\sH(A_3)\to 
\prP^1_j\setminus \{\infty\}$ with
this branch cycle description:   $\gamma_0' = (1\, 2\, 3)$, $\gamma_1' = (1\, 
3)$ and $\gamma_\infty' =  (1\, 2)$. 
One H-M rep.~$\bg^{+-+-}\in\ni(A_4,\bfC_{3_+^23_-^2})$ lies over $(+1,-1,+1,-
1)\in C_4$ and 
and one over $(+1,-1,-1,+1)$, $\bg^{+--+}\in\ni(A_4,\bfC_{3_+^23_-^2})$: 
$\mpr(\bg^{+-+-})=2$ while
$\mpr(\bg^{+--+})=3$ (as in \S\ref{HMrepRubric}).  

With $(+1,-1,+1,-1)$ over $j_0$ as a basepoint for
the cover $\sH_0^{\rd,+}\to \sH(A_3)$, 
the fiber over this point has representatives
$$\begin{array}{rl}\bg_1=& ( (1\, 2\, 3), (1\, 3\, 2), (1\, 3\, 4), (1\, 4\, 3) 
), \bg_2= ( (1\, 2\, 3), (1\,
2\, 4), (1\, 3\, 4), (2\, 3\, 4) ),\\ \bg_3=& ( (1\, 2\, 3), (1\, 2\, 4), (1\, 
4\, 2), (1\, 3\, 2)
).\end{array}$$

It is easy from this to form a branch cycle description for $\bar 
\sH_0^{\rd,+}\to \prP^1_j$:
$$\begin{array}{rl} \gamma_0''=&(1,6,9)(2,5,8)(3,4,7), \\
\gamma_1''=&(1,8)(2,7)(3,9)(4,6),\\
\gamma_\infty''=&(1,5,2,4)(3,6)(7,8,9).\end{array}$$ 
The group generated by these has size 648 and is centerless. Denote
$\lrang{\gamma_0'',\gamma_1'',\gamma_\infty''}$ by $H_0$. As the geometric 
monodromy group of a chain of two
degree 3 extensions, it would naturally be a subgroup of the wreath product 
$S_3\wr S_3$, as in the argument of
Lem.~\ref{basicA5cover}. The full wreath product has order 1296, while the 
kernel for the projection $H_0
\to S_3$ is the subgroup of $(a,b,c)\in S_3^3$ for which $abc\in A_3$. 

\subsubsection{$\sH(A_4,\bfC_{3_+^23_-^2})^{\abs,\rd,+}$ is not a modular curve} 
\label{notModular} 
Our problem was to decide, for a general cover $\phi_\bp: X_\bp \to \prP^1_z$ in 
this family, lying over $j(\bp)$,  if
the Galois closure $\hat X_\bp$ of the cover could factor through the elliptic 
curve with $j$ invariant $j(\bp)$.  This
is equivalent to $\Pic^{(0)}(\hat X_\bp)$ having that elliptic curve as an 
isogeny factor. The generality of the next
argument in $r$ will appear in later applications. Let $K_4$ be the copy of the 
Klein 4-group in $A_4$. 

\begin{lem}  \label{isogenyFactors} For $\bp$ generic, $\Pic^{(0)}(\hat X_\bp)$ 
is isogenous to $$\Pic^{(0)}(\hat
X_\bp/K_4)\times
\Pic_0(X_\bp)\times E_{\bp,1}\times E_{\bp,2}$$ with $\Pic^{(0)}(\hat 
X_\bp/K_4)$ a simple abelian variety of dimension
2, and $E_{\bp,1}$ and $E_{\bp,2}$ Prym (elliptic) varieties for the 
hyperelliptic curve $\hat X_\bp/K_4$. \end{lem}

\newcommand{\GU}{\text{{\rm GU}}} \newcommand{\SU}{\text{{\rm SU}}}
\begin{proof}[Terms and outline of proof] Riemann-Hurwitz shows $\hat X_\bp/K_4$ 
has genus 2.
\cite[p.~167]{VCycCov} shows for any
$r\ge 4$, the general cyclic degree 3 triple cover with $r$ branch points has a 
simple Jacobian. This result applies to
available for
$A_4$ covers no matter the number of branch points. The  proof considers $K$, 
the function field $\bQ(U_r)$ of the branch point locus, with the field of 
definition for the generic cover $C\to
\prP^1_z$ and its group of automorphisms adjoined. Then,  for $\ell\ge 5$ a 
prime, $G_K$ acts on the $\ell$-division
points $\Pic(C)^{(0)}_\ell=J_\ell$ according to this recipe:
\begin{edesc} \label{simpleFactors} \item If $\ell\equiv 1 \mod 3$ then 
$J_\ell=V\oplus V'$ with $V$ and $V'$ a standard
isotropic decomposition with respect to the Weil pairing, $G_K$ fixing $V$ and 
$V'$, inducing a group between
$\SL_{r-2}(\bZ/\ell)$ and
$\GL_{r-2}(\bZ/\ell)$; and 
\item if $\ell\equiv -1 \mod 3$, then $\Aut(C/\prP^1_z)$ makes $J_\ell$ an 
$\bF_{\ell^2}$ vector space, with $G_K$
inducing a group between $\SU_{r-2}(\bZ/\ell)$ and
$\GU_{r-2}(\bZ/\ell)$.
\end{edesc}
Suppose $J(C)=J$ is not simple. Let $K_1$ be a finite extension over which $J$ 
is isogenous to a product of abelian
varieties, so $G_{K_1}$ acts through a reducible normal subgroup of one of the 
simple groups occurring in
\eqref{simpleFactors}. Then, the finite group $G_K/G_{K_1}$ would contain 
infinitely many nonisomorphic simple groups, 
according to the composition factors appearing in the cases of 
\eqref{simpleFactors}. 

Finally, the cover $\hat X_\bp\to \hat X_\bp/K_4$ is unramified, and for each 
copy $A$ of a $\bZ/2$ in $K_4$, $\hat
X_\bp/A\to \hat X_\bp/K_4$ is a Prym cover. So, $\Pic^{(0)}(\hat X_\bp/A)$ has a 
degree 2 isogeny to $\Pic^{(0)}(\hat
X_\bp/K_4)\times E_A$ with $E_A$ an elliptic curve (\S\ref{autFunctHrd} and 
\cite[Lecture IV, p.~90-94]{CurveJac}). If
$A_1$ and $A_2$ generate
$K_4$, then their associated elliptic curves $E_i$, $i=1,2$, give the curves 
appearing in the statement of the lemma.
These isogeny factors account for 
$\Pic^{(0)}(\hat X_\bp)$ up to isogeny.  
\end{proof} 

Many different maps from
$\PSL_2(\bZ)$ to
$\PSL_2(\bZ/n)$ induce projective systems of covers of the $j$-line having the 
same geometric monodromy as the
projective system of curve covers called modular curves. Even, however, if some 
of these covers did have significant
moduli problems describing them, they would not be the same as those describing 
modular curves. For example, having
monodromy group
$S_3=\PSL_2(\bZ/2)$ does not imply an upper half plane quotient is the 
$\lambda$-line cover of the
$j$-line. That is the case with our space $\sH^{\rd,+}$. 

\begin{lem} The cover $\bar\psi^\rd: \bar
\sH^{\rd,+}\to \prP^1_j$ does not factor through $\prP^1_\lambda$.  So, 
$\sH^{\rd,+}$ is
not a modular curve. 
\end{lem} 

\begin{proof} Since the moduli
space  $\sH^{\rd,+}$ is a family of $r=4$ branch point covers,  
Ex.~\ref{conjClasspairscont} shows it does not satisfy the
moduli interpretation for factoring through $\prP^1_\lambda$. As the monodromy 
group has $S_3$ as a quotient, this means
$\sH^{\rd,+}\to \prP^1\setminus\{\infty\}$ is not a quotient of a modular curve. 
\end{proof} 

\subsubsection{Inspecting if $\sH^{\rd,+}$ is $j$-aware} \label{A434casejaware} 
Suppose $g\in A_4$ defines the degree 2
cover $\hat X_\bp/\lrang{g}\to \hat X_\bp/K_4$. The other degree two covers (for 
example, $E_{\bp,2}$) in
Lem.~\ref{isogenyFactors} have the  form $\hat X_\bp/\lrang{hgh^{-1}}\to \hat 
X_\bp/K_4$ with $h$ having order 3. So, if
for a generic $j$ value $\hat X_\bp/\lrang{g}$ is isogenous to $E_{\psi(\bp)}$, 
then so will be $\hat
X_\bp/\lrang{hgh^{-1}}$. According to Lem.~\ref{isogenyFactors}, to decide if
$\sH^{\rd,+}$ is $j$-aware requires answering the following question.   
\begin{prob} For general $\bp\in
\sH^{\rd,+}$, is either $X_\bp$ or $E_{\bp,1}$ isogenous to the elliptic curve 
with $j$-invariant $\psi(\bp)$. 
\end{prob} 

We leave this question after the following observations. Since the question is 
for general $\bp$, and we know from
\cite[Main Thm]{FrKK} we may assume $X_\bp$ has no complex multiplication. The
Riemann-Roch Theorem for genus 1 curves gives the dimension of linear systems of 
2-division points as the dimension of
$H^0(\kappa_\bp-D_{\phi_\bp})$ where
$\kappa_\bp$ is the canonical class, which on a genus 1 curve is trivial. This 
dimension is 1 (odd) if $D_{\phi_\bp}$
is linearly equivalent to 0, and 0 (even) if not. So, according to 
Prop.~\ref{evenThetas}, the linear
equivalence class of $D_{\phi_\bp}$ defines a non-trivial $\bQ(\bp)$ 2-division 
point on 
$\Pic^0(X_\bp)$.  

\newcommand{\so}{\text{\sl o}}
\newcommand{\sgo}{\text{\sl g}}
\newcommand{\sro}{\text{\sl r}}
\newcommand{\sho}{\text{\sl h}}
\newcommand{\spo}{\text{\sl po}}

\section{A representation difference between $\bar M_4$ on $\ni_1^+$ and
$\ni_1^-$} \label{GKS} The second author saw an S.~Sternberg talk in Jerusalem 
based on \cite{GKS}. 
In that talk, three semi-simple lie group representations had appeared as 
associated with the
existence of elementary particles, and the point was to see them as a natural 
multiplet. The shape
of the final result was this. 

They have $\sro \subset \sgo$ two complex Lie algebras with $\sgo$ semi-simple, 
$\sro$ reductive
and both having the same Carton subalgebra $\sho$. So the following properties 
hold. 
\begin{itemize} \item $\sro$ and $\sgo$ have the same rank. 
\item The Weyl group of $\sro$ is a subgroup of the Weyl group of $\sgo$.
\item With a choice of positive roots for each assume the positive Weyl chamber 
of $\sgo$ is in the
positive Weyl chamber of $\sro$.
\item The killing form on $\sgo$ gives an orthogonal splitting as $\sho\oplus 
\spo$. 
\end{itemize} The last item 
allows embedding $\sro$ in the orthogonal Lie algebra acting on
$\spo$. Further, this embedding has two representations $S^+$ and
$S^-$ associated with the plus and minus representations of the orthogal Lie 
algebra on the even
and odd part of the Clifford algebra for the orthogonal product (as in 
Prop.~\ref{liftEven}). 
Consider the elements
$C_{\sgo/\sro}$ of the Weyl group of $\sgo$ that map the positive Weyl chamber 
of $\sgo$ into the
positive Weyl chamber of $\sro$. The {\sl highest weights\/} $\lambda$ in the 
positive Weil chamber of $\sgo$
correspond to irreducible representations $V_\lambda$ of $\sgo$ 
\cite[p.~49]{Fuchs}. 

In the Grothendieck group of
$\sro$ representations, form
$V_\lambda\otimes S^+- V_\lambda\otimes S^-$ as a representation of $\sro$. 
Linear algebra magic
occurs here: For each $c\in C$, \cite{GKS} constructs a highest weight for 
$\sro$ from the pair
$(c,\lambda)$. Denote the corresponding representation by $U_{c,\lambda}$. Their 
formula is 
\begin{equation} V_\lambda\otimes S^+- V_\lambda\otimes S^-=\oplus_{c\in C} \pm
U_{c,\lambda},\end{equation} with $\pm$ indicating some choice of sign. 

The analogy for us starts with a group $\bar M_4=\SL_2(\bZ)$, and  $\bg_0\in
\ni(A_5,\bfC_{3^4})$ representing an H-M rep.~at level 0. Let $M_{\bg_0}$ be the 
subgroup of $\bar M_4$
stabilizing the reduced class of $\bg_0$. Then, Prop.~\ref{expInv} gives two 
representations $S^+$ and $S^-$ of
$H_{\bg_0}$, respectively corresponding to the actions on the sets of 
\eqref{conjSep} and \eqref{conjSep2}. The analogy
would go further if there was a representation $V_\lambda$ of $\bar M_4$ that 
fit in the formula $V_\lambda\otimes S^+-
V_\lambda\otimes S^-$ to give the difference between the two representations on 
\eqref{conjSep} and \eqref{conjSep2}.
The  analogy misses a lot, though we can easily construct $p$-adic Lie algebras 
from the
projective situation.  Since there are infinitely many levels of a
Modular Tower, our goal would be to understand if there are such paired 
representations at 
higher levels of the Modular Tower. This would give a direct translation of the 
Clifford (Spin)
invariant of Prop.~\ref{expInv}, and its  relation to
$\theta$-nulls for each component of $\sH_1=\sH({}_2^1\tilde
A_5,\bfC_{3^4})^{\inn,\rd}$ (as in \S\ref{hcan}). 

\section{Where are those regular realizations?} \label{galSummary}
Let $\{G_k\}_{k=0}^\infty$ the characteristic quotients of the universal $p$-Frattini cover $\tG
p$ of a $p$-perfect group $G=G_0$. Our contributions to Main Conjecture (\cite[Main 
Conj.~0.1]{FrKMTIG}, Prob.~\ref{MPMT},
\S\ref{diophMT})  with $r\le 4$ branch points suggest that no bounded set of 
branch points will allow regular
realizations over a fixed number field of the complete set 
$\{G_k\}_{k=0}^\infty$. 

Further, suppose for some $r_0\ge 1$, there are realizations of  
$\{G_k\}_{k=0}^\infty$, each requiring no more than
$r_0$ branch points. Then, Thm.~\ref{thm-rbound} says there must exist a Modular 
Tower for $(G_0,\bfC,p)$ with $\bfC$ a
set of $p'$ conjugacy classes having cardinality $r\le r_0$, and rational points 
at every level. The tower levels
are manifolds. So, having  $\bQ$ points at every level implies the Modular Tower 
levels  have a projective system of
$\bQ$ components. \cite[Thm.~3.21]{FrMT} shows that if there is just one 
absolutely irreducible H-M component at every
level of a $\bQ$ Modular Tower, then this will be defined over $\bQ$. There are 
variants on this, though this is
the clearest case establishing this necessary condition. Further, assume 
$\ni(G_0,\bfC)$ contains H-M reps. For any integers $1\le i< j\le r$, call 
$(i,j)$ an H-M pair if $\C_i$ is
inverse to the class $\C_j$. In this situation let 
$\bfC_{i,j}$ be $\bfC$ with $\C_i$ and $\C_j$ removed. The only successful 
hypotheses guaranteeing  one
H-M component at every level are mildly variant on the other half of 
\cite[Thm.~3.21]{FrMT}. There will be 
just one H-M component if these (H-Mgcomplete) hypotheses hold:
\begin{trivl}  For each H-M pair $(i,j)$,  and $\bh\in \bfC_{i,j}$, $\lrang{\row 
h {r-2}}=G_0$. \end{trivl} 
\noindent For example, this would hold if for each H-M pair $(i,j)$, 
$\bfC_{i,j}$ contains every nontrivial conjugacy
class of $G_0$.  

There are results that encourage using H-Mgcomplete reps.: For any prime $\ell$, 
\cite[Thm.~4.2]{DebDes} finds a
$\bQ_\ell$ regular realization of $\tG p$  using the H-Mgcomplete condition. The 
suggested generalization of the Open
Image Theorem  in  \S\ref{outMon} would apply especially well to a sequence of 
H-Mgcomplete components.    Reminder:
Shafarevich's method with solvable groups does not give regular realizations. 
Only some version of the braid rigidity
method has systematically produced regular realizations, and going beyond 
solvable groups the story has been regular
realizations all the way. That even includes Shih's results using modular 
curves: These came from the first version of
the method the second author told Shimura when he first worked it out at the 
Institute for
Advanced Study at Princeton in 1967-69. We consider
Modular Towers a general invention, extending through the
analogy with modular curves, the potential for many
applications. What, however, of someone who wants such regular realizations? 
Does this say anything of 
the next reasonable step in the search for them? 

\subsection{Reminder of the setup} For a large integer $k$, where would one find 
a regular
realization of $G_k$. From Thm.~\ref{thm-rbound} combined with the conclusion of 
the Main Conjecture, you
will need many branch points, and so you need a model for where to get an 
indication of how many, and
what kind of conjugacy classes (with elements of order divisible by $p$) you 
will use.  

When $G_0$ is $p$-split, with cyclic
$p$-Sylow, it is tempting to imitate realization of dihedral groups. No one 
knows how to get regular
{\sl involution\/} realizations of dihedral groups (when $r=4$, the moduli of 
these relates 
special Hurwitz spaces and modular curves \S\ref{compMC}). Yet, any elementary 
algebra
book shows how to get regular realizations of dihedral groups. You must use $p$-cycles as branch
cycles, and you require many branch points (see the implications of this 
elementary analysis in
\cite{Fr-Se2} or \cite{Fr-Se1}).  Yet, for a general
finite group $G$ and prime $p$, no one has produced similar branch cycles. 
Results from \S\ref{genFrattini} show
divide the case where $M_0$ is cyclic from those where it is not. We refer to 
the former case as {\sl dihedral-like\/}. 
It is the latter case that is serious. We consider finding such branch
cycles when the situation is not dihedral-like, our basic assumption from here.  
So, we are seeking 
dihedral-like realizations in non-dihedral-like situations.   

Recorded in \cite{FrVMS} is the following possibility: Given $(G_0,\bfC,p)$, 
there might 
be infinitely many $\bfC'$ formed by   suitably increasing the multiplicity of 
appearance of
conjugacy classes in $\bfC$, so that  $\sH(G_k,\bfC)^{\inn}$ will be a uni-
rational variety. 
(That result would make a mockery of the inverse Galois problem.) Our 
replacement would consider using special
realizations at level 0, formed, say, from understanding simple groups.  We 
acknowledge the Main Conjecture might be false; 
it has phenomenal implications for $r\ge 5$ (see the comments about
\cite{FKVo} in \S\ref{usebcl}).  

\subsection{From Chevalley realizations to the whole $p$-Frattini tower} An 
approach of Thompson-V\"olklein
(\cite{TVo} has an example, or the Thompson tuples of Ex.~\ref{ThomTuples}) uses 
intricate knowledge of
Chevalley simple groups to locate conjugacy classes that show some Chevalley 
series over each finite field $\bF_q$ has
regular realizations excluding finitely many
$q$. This has worked for many  series, extending V\"olklein's production of high 
rank Chevalley group realizations.
(Prior to his approach, a less abstract use of the main idea of 
\cite{FrVAnnals}, there had been essentially no higher
rank realizations.) We consider as a starting point some successful realization 
cases where the Hurwitz space 
$\sH(G_0,\bfC)^{\inn}$ is uni-rational (it might even be an abelian cover of  
$\prP^r\setminus D_r$ as in the
Thompson-V\"olklein examples). Take a prime
$p$ not dividing the orders of elements in
$\bfC$. 

Form, for each $k$, a sequence $\bfC_k$ of conjugacy classes, using the
following principles. Reminder: Def.~\ref{ratunion} explains rationalization of 
any conjugacy classes
$\bfC'$. 
\begin{edesc} \label{pchoices} \item $\bfC_k=(\bfC,\bfC_k')$ is a rational union 
of conjugacy classes in
$G_k$ with $\bfC_k'$ (possibly empty) conjugacy classes whose elements have $p$ 
power order, $\bfC$ are
$p'$ classes, $\bfC_k'$ are conjugacy classes of $p$-power order, and 
$\sH(G_0,(\bfC,\bfC_0'))^\inn$ is
a rational (or uni-rational) variety over $\bQ$.
\item $\bfC_k'$ consists of $(\bfC_k^*,\bfC_k^{**})$ where $\bfC_k^*$ consists 
of the rationalization
of $\bfC_{k-1}^\dagger$ whose terms consist of a lift to $G_k$ of the conjugacy 
class for $G_{k-1}$ for
each member of $\bfC_k'$, and $\bfC_k^{**}$ is a rational union of conjugacy 
classes (for $G_{k}$) in 
$\ker(G_{k+1}\to G_{k})=M_k$. 
\end{edesc} 

Since $G_k\to G_0$ is a Frattini cover, any lifts of generators of $G_0$ to 
$G_k$ generate $G_k$. So,
as with the universal $p$-Frattini cover, it is only the product 1 condition 
that governs the
possibility of braid orbits. There will be no obvious minimal choice. 

Given (nontrivial) 
$h\in \bfC_{k-1}'$, any lift $\hat h$ to $G_{k}$ will have $p$ times the order 
of $h$ (Lem.~\ref{FrKMTIG}). Further,
since the situation is not dihedral-like,  there will be many
$G_{k-1}$ orbits on the collection of lifts of $h$. We use the notation 
established in the rest of the paper for
Nielsen classes at levels of a Modular Tower. 

\begin{exmp}[Start from $(A_n,\bfC_{3^{n-1}},p=2)$] \label{AnchooseConj} 
Consider
$\bfC_1=(\bfC_{3^{n-1}},\C_1')$ with
$\C_1'$ the conjugacy class of an element $c$ in $M_0=M_0(A_n)$. For 
$\ni(G_1,\bfC_1)$ to be nonempty, the image
of $c$ in
$\ker(\Spin_n\to A_n)$ must be $(-1)^{n-1}$ (see Ex.~\ref{full3cycleList}). When 
$n=5$ we know this is
{\sl sufficient\/} to guarantee $\ni(G_1,\bfC_1)$ is nonempty; the argument 
right after Princ.~\ref{obstPrinc}. (For
$n=5$, there are two conjugacy class choices for such $c$, and several reasons 
(like Prop.~\ref{expInv}) to see them as
significantly different, though not for this particular obstruction.) For, 
however,
$n\ge 6$, appearance of $\one_{A_n}$ beyond the first Loewy layer of $M_0$ might 
obstruct
$\ni(G_1,\bfC_1)$ from being nonempty for some choices of conjugacy class 
$\C_1$. That is, we don't know for certain if
$\{\Pi(\bg)\mid \bg \in \ni(G_1,\bfC_1)\}$ fills out all conjugacy classes lying 
over the conjugacy class of
$(-1)^{n-1}$ in $\Spin_n$.  
\end{exmp}

Assume the non-dihedral-like case, with the usual $p$-perfect assumption,  and 
consider  
the list \eqref{pchoices}. The next problem asks if it is possible to remove 
choices after level 1
and still be certain higher level Nielsen classes are nonempty. 

\begin{prob} Is it always possible to find some $\bfC_1$ so the following holds. 
For $k\ge 2$, with
$\bfC_k^{**}$ empty and any choice of a single lift of each class in $\bfC_{k-
1}'$ to $G_k$ to give
$\bfC_{k-1}^\dagger$ and its rationalization 
$\bfC_k^{*}$, then $\ni(G_k,\bfC_k)$ is nonempty? \end{prob} 
\end{appendix}


\begin{thebibliography}{Ser90b}

\bibitem[Ah79]{Ahlfors} L.~Ahlfors, \emph{Introduction to the {T}heory of 
{A}nalytic 
{F}unctions of {O}ne}, 3rd edition, Inter. Series in Pure and Applied Math., 
McGraw-Hill
{C}omplex {V}ariable, 1979. 

\bibitem[AW67]{AW}
M.F. Atiyah and C.T.C. Wall, \emph{Cohomology of groups}, vol. Algebraic Number
  Theory, pp.~94--115, Thompson Book Co. and Academic Press, 1967, Editor
  J.~W.~S.~Cassels.

\bibitem[Ben91]{Ben1}
D.J. Benson, \emph{I: Basic representation theory of finite groups and
  associative algebras}, Cambridge Studies in advanced math., vol.~30,
  Cambridge U. Press, Cambridge, 1991.

\bibitem[BF98]{BeFr}
G.~Berger and M.~Fried, \emph{Rational cusps on noncongruence towers of the
  $j$-line}, preprint (Sept.~1998), 1--28.

\bibitem[BF82]{BFr}
R.~Biggers and M.~Fried, \emph{Moduli spaces of covers and the {H}urwitz
  monodromy group}, Crelles Journal \textbf{335} (1982), 87--121.

\bibitem[Bi75]{Birman}
J.~Birman, \emph{Braids, links and mapping class groups}, 
Based on notes of J.~Cannon: Princeton U. Press
\textbf{82}, 1975. 

\bibitem[Boh47]{Bohnen}
F.~Bohnenblust, \emph{The algebraical braid group}, Ann.~of Math. \textbf{(2)
  48} (1947), 127--136.

\bibitem[Ch44]{Chat1} 
F.~Ch\^atelet, 
\emph{Variations sur un th\`eme de Poincar\'e}, 
Annales ENX 3${}^e$ s\'erie {\bf 61} (1944), 
249--300.

\bibitem[CoH85]{CoombesHarb} K.~Coombes  and D.~Harbater, \emph{Hurwitz
families and  arithmetic Galois groups}, Duke Math.~J. \textbf{52} 
(1985),  821--839.

\bibitem[CT88]{CT} 
J.~Colliot-Th\'elene,  
\emph{Les Grands th\`emes de Fran\c cois Ch\^ atelet},
L'Ens. Math.   {\bf 34} (1988),  387--405. 

\bibitem[DF90]{DFrRRCF}
P.~D\`ebes and M.~Fried, \emph{Rigidity and real 
residue class fields}, Acta Arith.  \textbf{56} (1990), 13--45.

\bibitem[DF90b]{DFrVarFam} 
P.~D\`ebes and M.~Fried, 
\emph{Arithmetic variation of fibers in 
families: Hurwitz  monodromy criteria for rational points \dots}, 
Crelles J. \textbf{409} 
(1990), 106--137. 


\bibitem[Deb95]{DebP1Qp}
P.~D\`ebes, \emph{Covers of ${\prP}^1$ over the
${\bQ}_p$-adics}, Cont.~Math. AMS, Recent Developments in
the Inverse Galois 
Problem, Seattle 1995, Editors
M.~Fried, et.~al., \textbf{186}, 217--238. 


\bibitem[DF99]{DFrIntSpec}
P.~D\`ebes and M.~Fried, \emph{Integral specialization 
of families of rational functions},  PJM {\bf 190} (1999) 45--85
(typos corrected at www.math.uci/\~{}mfried/\#math).

\bibitem[Deb01]{DDescTh} 
P.~D\`ebes, 
\emph{Descent Theory for Algebraic Covers}, this volume.

\bibitem[DDes01]{DebDes} P.~D\`ebes et B.~Deschamps,
\emph{Corps $\psi$-libres et th\'eorie inverse de Galois
infinie}, preprint as of October, 2001. 


\bibitem[DDE00]{DDEm} 
P.~D\`ebes, J.-C.~Douai and M.~Emsalem, 
\emph{Familles de Hurwitz et Cohomlogies non Ab\'elienne}, 
Annales de l'institut Fourier, 
{\bf 50} (2000), 113--149.  

\bibitem[Fal83]{FaMordCon}
G.~Faltings, \emph{Endlichkeitss\"atze \"uber 
Zahlk\"orpern}, Invent.~Math.
  \textbf{73} (1983), 349--366.

\bibitem[FaK01]{FKraTheta}
H.~Farkas and I.~Kra, \emph{Theta Constants, Riemann
Surfaces and the Modular Group}, AMS graduate text
series \textbf{37}, 2001.


\bibitem[Fay73]{Fay} 
J.~Fay, 
\emph{Theta {F}unctions on {R}iemann {S}urfaces}, 
Lecture notes in Mathematics 
\textbf{352}, Springer Verlag, Heidelberg, 1973. 

\bibitem[FKVo99]{FKVo}
G.~Frey, E.~Kani and H.~V\"olklein,
\emph{Curves with infinite $K$-rational geometric fundamental group}, 
   in:  Aspects of Galois Theory, edts. H. Voelklein et. al.,
   London Math. Soc. Lecture Notes   {\bf 256},   
   Cambridge Univ. Press (1999), 85--118.     


\bibitem[FJ86]{FrJ}
M.~Fried and M.~Jarden, \emph{Field arithmetic}, Ergebnisse
der Mathematik
  {III}, vol.~11, Springer Verlag, Heidelberg, 1986.

\bibitem[FK97]{FrKMTIG}
M.~Fried and Y.~Kopeliovic, \emph{Applying Modular Towers
to the inverse
  {G}alois problem}, Geometric {G}alois {A}ctions {II} {D}essins d'{E}nfants,
  Mapping Class Groups and Moduli, vol. 243, Cambridge U.~Press, 1997, 
London  Math.~Soc. Lecture Notes, pp.~172--197.

\bibitem[FKK01]{FrKK} 
M.~Fried, E.~Klassen, Y.~Kopeliovic, 
\emph{Realizing alternating groups
as monodromy groups of  genus one covers}, 
PAMS \textbf{129} (2000), 111--119.

\bibitem[Fri70]{FrSchur}
M.~Fried, \emph{On a conjecture of {S}chur}, Mich.~Math.~J.
\textbf{17}
  (1970), 41--55.

\bibitem[Fri75]{FrBGJac} M.~Fried, \emph{Brauer groups and
Jacobians}, Preprint from a talk at Oberwolfach Alg. No.
Theory session, July 1975. 

\bibitem[Fri77]{FrHFGG}
M.~Fried, \emph{Fields of definition of function fields and
{H}urwitz
  families and groups as {G}alois groups}, Communications in Algebra \textbf{5}
  (1977), 17--82.

\bibitem[Fri78]{FrGGCM}
M.~Fried, \emph{{G}alois groups and {C}omplex {M}ultiplication}, 
Trans.A.M.S. {\bf 235} (1978) 141--162. 

\bibitem[Fri89]{FrCCMD}
M.~Fried, \emph{Combinatorial computation of moduli
dimension of {N}ielsen
  classes of covers}, Cont.~Math., vol.~89, pp.~61--79, AMS, Providence, RI,
  1989.

\bibitem[Fri90]{Fr3-4BP}
M.~Fried, \emph{Arithmetic of 3 and 4 branch point covers:
{A} bridge
  provided by 
noncongruence subgroups of ${\SL}_2({\bZ})$}, Progress in Math.~Birkhauser 
\textbf{81} (1990), 77--117.

\bibitem[Fri94]{Fr-Se2} 
M.~Fried, \emph{review---{T}opics in {G}alois Theory, {J.-P}.
{S}erre},  BAMS \textbf{30 \#1}  (1994), 124--135.  

\bibitem[Fri95a]{FrMT}
M.~Fried, \emph{Introduction to Modular Towers: {\sl
Generalizing the relation between  dihedral groups and
modular curves\/}}, Proceedings AMS-NSF  Summer Conference,
vol.  186, 1995, Cont. Math series, Recent Developments in
  the Inverse {G}alois Problem, pp.~111--171.

\bibitem[Fri95b]{FrExtConst}
M.~Fried, \emph{Extension of constants, rigidity, and the
  {C}howla-{Z}assenhaus conjecture}, Finite Fields and their applications;
  Carlitz Volume \textbf{1} (1995), 326--359.

\bibitem[Fri95c]{Fr-Se1}
M.~Fried, \emph{Enhanced review: Serre's {\rm {T}opics in
{G}alois Theory}}, Proceedings of the Recent
developments in  the Inverse {G}alois Problem 
conferenc, AMS Cont. Math.~ \textbf{186} (1995), 15--32.  


\bibitem[Fri96]{FrLInv}
M.~Fried, \emph{Alternating groups and lifting invariants},
Preprint 
  07/01/96 (1996), 1--34.

\bibitem[Fri99]{Fr-Schconf}
M.~Fried, \emph{Separated variables polynomials and moduli
spaces}, Number
  Theory in Progress (Berlin-New York) (J.~Urbanowicz K.~Gyory, H.~Iwaniec,
  ed.), Walter de  Gruyter, Berlin-New York (Feb. 1999), Proceedings of the
Schinzel Festschrift, Summer 1997, pp.~169--228,
http://www.math.uci.edu/{$\tilde{\phantom u}$mfried/\#math}.

\bibitem[Fr01]{AGR}
M.~Fried and D.~Semmen, \emph{Relatively nilpotent algebraic extensions 
and work of  Abel, Galois and Riemann}, in preparation. 

\bibitem[Fr02]{FrRIMS}  M.~D.~Fried, \emph{Moduli of
relatively nilpotent extensions\/}, Institute of
Mathematical Science Analysis 1267, June 2002, Comm. in
Arithmetic Fundamental Groups, 70--94. 

\bibitem[FV91]{FrVMS}
M.~Fried and H.~V{\"o}lklein, \emph{The inverse {G}alois
problem and rational
  points on moduli spaces} , Math.~Annalen
\textbf{290} (1991), 771--800.

\bibitem[FV92]{FrVAnnals}
M.~Fried and H.~V{\"o}lklein, \emph{The embedding
problem  over an  {H}ilbertian-{PAC} field}, Annals of Math
\textbf{135} (1992), 469--481.

\bibitem[Fu92]{Fuchs}
J.~Fuchs, \emph{Affine Lie algebras and quantum groups},
Cambridge Univ. Press, Monographs on Mathematics physics, 
1992.

\bibitem[GB68]{GilBus}
P.~Gillette and van Buskirk, \emph{The word problem and its consequences for the 
braid groups
and mapping class groups of the 2-sphere}, TAMS \textbf{131
No.~2} (1968), 277--296.

\bibitem[GS78]{GriessFrat} R.~Griess and P.~Schmid, \emph{The
Frattini module}, Archiv.~Math. \textbf{30} (1978), 256Ñ266. 

\bibitem[GKS99]{GKS}
B.~Gross, B.~Kostant and S.~Sternberg, \emph{Gainesville Lectures on Kostant's 
Dirac 
Operator}, by Shlomo Sternberg, Feb.~23--26, 1999. 

\bibitem[Gu67]{gunInd}
R.C.~Gunning, \emph{Special Coordinate Covering of Riemann Surfaces}
Math.~Ann. \textbf{170} (1967), 67--86. 

\bibitem[Ha63]{HaThGps}
M.~Hall, \emph{The theory of groups}, Macmillan, 1963.

\bibitem[Har77]{Hart}
R.~Hartshorne, \emph{Algebraic {G}eometry},
Graduate Texts in Math. \textbf{52}, Springer-Verlag,
1977.

\bibitem[Ha80]{HaTF} D.~Harbater, \emph{Deformation theory
and the tame fundamental group},   TAMS,
{\bf 262} (1980), 399-415.

\bibitem[Ha84]{HaMC}  D.~Harbater,  \emph{Mock covers and
Galois extensions},  J.~Algebra {\bf 91} (1984), 281-293.


\bibitem[Iha86]{IharProBraids}
Y.~Ihara, \emph{Profinite braid groups}, Annals of Math. \textbf{123} (1986),
  43--106.

\bibitem[Iha91]{IharIntCong}
Y.~Ihara, \emph{Braids, {G}alois groups and some arithmetic functions},
  Proceedings of the International Congress, vol. Kyoto 1990, pp.~99--120,
  Springer-Verlag, Tokyo, 1991.

\bibitem[Iha00]{IharSC} 
Y.~Ihara, \emph{Shimura curves over finite fields and their rational
points}, Cont.~Math., proceedings of AMS-IMS-SIAM  Summer
Conference 1997 on Applications of Curves over finite fields,
Cont.~Math. {\bf 245} (1999), 15--23.


\bibitem[IM95]{IharMat}
Y.~Ihara and M.~Matsumoto, \emph{On Galois actions on profinite completions of
  braid groups}, Proceedings AMS-NSF Summer Conference, vol. 186, 1995, Cont.
  Math series, Recent Developments in the Inverse {G}alois Problem,
173--200.

\bibitem[Ja99]{JarProj} M.~Jarden, \emph{The projectivity of
the fundamental group of an affine line}, Turk.~J.~Math.
\textbf{23} (1999), 531--547. 

\bibitem[Je41]{Je} 
S.A.~Jennings, \emph{The structure of the group ring of a
$p$-group over a modular field}, TAMS \textbf{50} (1941),
175--185.

\bibitem[KMS66]{KMS}
A.~Karrass, W.~Magnus and D.~Solitar,
\emph{Combinatorial group theory}, Interscience Publishers,
J. Wiley and Sons, 1966. 

\bibitem[Ma34]{Mag} 
W.~Magnus,
\emph{{\"U}ber {A}utomorphismen von {F}undamentalgruppen berandeter
{F}l\"achen}, Math.~Ann. 109, 1934, 617--646.

\bibitem[Mar45]{Markov}
A.~Markoff, \emph{Foundations of the algebraic theory of tresses}, Trav. Inst.
  Math. Steklov, vol.~16, Steklov, Moscow, 1945.

\bibitem[MM99]{MM}  
G.~Malle and B.H.~Matzat,
\emph{Inverse {G}alois Theory},
ISBN 3-540-62890-8, Monographs in Mathematics, 
Springer,1999. 

\bibitem[Mes90]{Me}
J.-F. Mestre, \emph{Extensions r\'eguli\`eres de ${\bQ}(t)$ de groupe de
  {G}alois $\tilde {A}_n$}, J.~of Alg. \textbf{131} (1990), 483--495.

\bibitem[Mes94]{Me2}
J.-F. Mestre, \emph{Annulation, par changement de variable,
d'\'el\'ements de $\text{\rm Br}_2(k(x))$ ayant quatre
p\v oles}, CR.~Acad.~Sci., Paris \textbf{319} (1994),
781--782.

\bibitem[Moc96]{MochpUnif}
S.~Mochizuki, \emph{A theory of ordinary $p$-adic curves}, Publ.~RIMS
  \textbf{32 \#6} (1996), 957--1151.

\bibitem[Mu66]{MumRB}
D.~Mumford, \emph{Introduction to Algebraic Geometry; {T}he
{R}ed {B}ook}, Harvard Lecture Notes, 1966. 

\bibitem[Mu72]{MuDC} D.~Mumford, \emph{An analytic
construction of degenerating curves over complete local
rings}, Comp.~Math. {\bf 24} (1972), 129--174.

\bibitem[Mu76]{CurveJac}
D.~Mumford, \emph{Curves and their {J}acobians},
{Progress in Mathematics}, Univ.~of Mich.~Press, Ann Arbor,
1976.


\bibitem[Nor62]{NCHA}
D.~G. Northcott, \emph{An introduction to homological algebra}, Camb.~U.~press,
  1962.


\bibitem[Na98]{tanBasePts} 
H.~Nakamura, \emph{Tangential base points and Eisenstein series}, 
preprint, 1998, 1-16. 

\bibitem[Qu68]{Qu} D.G.~Quillen, \emph{On the associated
graded ring of a group ring}, J.~Alg. \textbf{10} (1968),
411-418. 

\bibitem[Rig96]{GalLife} L.T.~Rigatelli, \emph{Evariste Galois:
1811-1832}, Vol. 11,
translated from the Italian by John Denton, Vita Mathematica,
Birkh\"auser, 1996. 

\bibitem[Rob73]{RoEC}
A.~Robert, \emph{Elliptic curves}, Lecture Notes, vol. 326, Springer Verlag,
  Heidelberg, 1973; second edition 1986.

\newcommand{\etalchar}[1]{$^{#1}$}
\bibitem[S{\etalchar{+}}95]{GAP}
 Martin Sch{\accent127 o}nert et~al.
\newblock {\em {GAP} -- {Groups}, {Algorithms}, and {Programming}}.
\newblock Lehrstuhl D f{\accent127 u}r Mathematik,
Rheinisch Westf{\accent127 a}lische Technische Hoch\-schule,
Aachen, Germany, fifth edition, 1995.

\bibitem[Se04b]{Darren} D. Semmen, \emph{Jennings' theorem
for
$p$-split groups}, J.~Alg. \textbf{285} (2005), 730--742. 


\bibitem[Ser68]{SeAbell-adic}
J.-P. Serre, \emph{Abelian $\ell$-adic representations and elliptic curves},
  1st ed., McGill University Lecture Notes, Benjamin, New York $\bullet$
  Amsterdam, 1968, in collaboration with Willem Kuyk and John Labute.

\bibitem[Ser72]{SeElliptic}
J.-P. Serre, \emph{Propri\'et\'es galoisiennes des points d'ordre fini des
  courbes elliptiques}, Inv.~Math. \textbf{15} (1972), 259--331.

\bibitem[Ser77]{SeArbres}
J.-P. Serre, \emph{Arbres, amalgames, $\SL_2$}, vol.~46,
Ast\'erisque, New
  York, Amsterdam, 1977.

\bibitem[Ser90a]{SeLiftAn}
J.-P. Serre, \emph{Rel\^evements dans $\tilde A_n$}, C.~R.~Acad. Sci. Paris
  \textbf{311} (1990), 477--482.

\bibitem[Ser90b]{SeTheta}
J.-P. Serre, \emph{Rev\^etements a ramification impaire et
  th\^eta-caract\'eristiques}, C.~R.~Acad. Sci. Paris \textbf{311} (1990),
  547--552.

\bibitem[Ser92]{Se-GT}
J.-P. Serre, \emph{Topics in {G}alois theory}, 
no. ISBN \#0-86720-210-6,
  Bartlett and Jones Publishers, notes taken by 
H.~Darmon,
1992.

\bibitem[Ser96]{SeGalCoh}
J.-P.~Serre, \emph{Galois Cohomology},
4th Edition, Springer: Translated from the French by Patrick Ion;
original from Springer LN5 (1964), 1996. 


\bibitem[Sh01]{Sh} A.~Shalev, \emph{Asymptotic Group Theory}, Notices of the AMS 
{\bf April 2001},
383--389. 

\bibitem[Sh94]{ShAutFuncts}
G.~Shimura, \emph{Introduction 
to the arithmetic of
automorphic functions}, Publications of Math.~Soc.~Japan,
Princeton 
U.~Press, 2nd Edition, First printing 1971, 
1994.

\bibitem[Sh98]{ShAbVar}
G.~Shimura, \emph{Abelian Varieties with Complex
Multiplication and Modular Functions}, Princeton  U.~Press,
1998.

\bibitem[Sp57]{Springer} 
G.~Springer,
\emph{Introduction to {R}iemann {S}urfaces}, Addison-Wesley
Mathematics Series \textbf {11}, 1957.

\bibitem[St94]{StAlg} M.~Steinberger, \emph{Algebra}, PWS
publishing, 1994. 

\bibitem[Ta58]{TateWC} 
J.~Tate, \emph{WC-groups over $p$-adic fields},
S\'em.~Bourbaki, Expos\'e {\bf 156}, 1957--1958. 

\bibitem[TVo99]{TVo}
J.G.~Thompson, and H.~V\"olklein,
\emph{Braid abelian tuples in $\Sp_n(K)$}, 
   in:  Aspects of Galois Theory, edts. H. Voelklein et. al.,
   London Math. Soc. Lecture Notes   {\bf 256},   
   Cambridge Univ. Press (1999), 85--118.     

\bibitem[{V\"o}95a]{VCycCov}
H.~V{\"o}lklein, \emph{Cyclic covers of $\prP^1$ and Galois
action on their Division Points\/},
Proceedings AMS-NSF  Summer Conference, vol.  186, 1995,
Cont. Math series, Recent Developments in
  the Inverse {G}alois Problem, 91--107.


\bibitem[V{\"o}96]{VB}
H.~V{\"o}lklein, \emph{Groups as {G}alois 
{G}roups} {\bf 53},
Cambridge Studies in Advanced Mathematics,
Camb.~U.~Press, Camb.~England, 1996.


\bibitem[V{\"o}01]{VoBC}
H.~V{\"o}lklein, \emph{Galois realizations of profinite
projective linear groups}, this volume.

\bibitem[We56]{WeilField} 
A.~Weil, \emph{The field of definition of a variety},
Amer.~J.~Math. {\bf 78} (1956), 509--524.  

\bibitem[We98]{WeTh}
S.~Wewers, \emph{Construction of {H}urwitz spaces}, {Thesis},
Institut f\"ur Experimentelle Mathematik {\bf 21}  
(1998), 1--79.


\bibitem[We01]{WeFM}
S.~Wewers, \emph{Field of moduli and field of definition of
Galois covers}, this volume.
\end{thebibliography}
\providecommand{\bysame}{\leavevmode\hbox to3em{\hrulefill}\thinspace}

\end{document}